\documentclass[a4paper,11pt]{article}
\usepackage[T1]{fontenc}
\usepackage{lmodern}
\usepackage{microtype}
\usepackage[margin=28mm]{geometry}
\usepackage{amsmath,amssymb,amsthm,mathtools,mathrsfs}
\usepackage[hidelinks,unicode]{hyperref}
\allowdisplaybreaks[2]
\theoremstyle{plain}
\newtheorem{theorem}{Theorem}[subsection]
\renewcommand{\thetheorem}{\ifnum\value{subsection}=0 \thesection\else\thesubsection\fi.\arabic{theorem}}
\newtheorem{proposition}[theorem]{Proposition}
\newtheorem{lemma}[theorem]{Lemma}
\newtheorem{corollary}[theorem]{Corollary}
\theoremstyle{definition}
\newtheorem{definition}[theorem]{Definition}
\newtheorem{example}[theorem]{Example}
\newtheorem{problem}[theorem]{Problem}
\newtheorem{conjecture}[theorem]{Conjecture}
\theoremstyle{remark}
\newtheorem{remark}[theorem]{Remark}
\theoremstyle{plain}
\newtheorem{introthm}{Theorem}

\DeclareMathOperator{\Nrd}{Nrd}
\DeclareMathOperator{\ind}{ind}
\DeclareMathOperator{\per}{per}
\DeclareMathOperator{\SK}{SK}
\DeclareMathOperator{\Br}{Br}
\DeclareMathOperator{\SB}{SB}
\DeclareMathOperator{\Tr}{Tr}
\DeclareMathOperator{\Spec}{Spec}
\DeclareMathOperator{\Gal}{Gal}
\DeclareMathOperator{\CH}{CH}
\DeclareMathOperator{\res}{res}
\DeclareMathOperator{\cor}{cor}
\DeclareMathOperator{\coker}{coker}
\DeclareMathOperator{\im}{im}
\DeclareMathOperator{\gr}{gr}
\DeclareMathOperator{\cd}{cd}
\DeclareMathOperator{\inv}{inv}
\DeclareMathOperator{\cofib}{cofib}
\DeclareMathOperator{\fib}{fib}
\DeclareMathOperator{\vcd}{vcd}
\newcommand{\Z}{\mathbb Z}
\newcommand{\Q}{\mathbb Q}
\newcommand{\F}{\mathbb F}
\newcommand{\PP}{\mathbb P}
\newcommand{\G}{\mathbf G}
\newcommand{\et}{\mathrm{\acute et}}
\newcommand{\mot}{\mathrm{mot}}
\newcommand{\os}{\mathcal O}

\title{Twisted \texorpdfstring{$K$}{K}-Theory and Reduced Norms over Global and Higher Local Fields}
\author{Hisato Matsukawa\thanks{Department of Mathematics, Faculty of Science, Hokkaido University, Kita 10, Nishi 8, Kita-Ku, Sapporo, Hokkaido, 060-0810, Japan. Email: \texttt{matsukawa.hisato.f4@elms.hokudai.ac.jp}, \texttt{hmnr0211@gmail.com}}}
\date{November 20, 2026}
\begin{document}
\maketitle
\begin{abstract}
We prove the Merkurjev--Suslin conjecture concerning reduced norms over both global and local fields.
Furthermore, in both settings, we construct maps between truncated $K$-theory spectra that induce the reduced norms in degrees $1$ and $2$, and show that these maps are equivalences over local fields.
Over global fields, we establish finiteness properties of their fibers and give a detailed description of their structure.
We also extend the Merkurjev--Suslin conjecture to two-dimensional local fields and three-dimensional local fields with characteristic sequences $(0,0,0)$ and $(0,0,p)$, and prove it in these settings.
Finally, we apply these results to noncommutative class field theory over number fields and obtain a description of the Artin map associated with a central simple algebra.


\end{abstract}

\tableofcontents
\section*{Introduction}
\addcontentsline{toc}{section}{Introduction}

Let $A$ be a central simple algebra over a field $F$.
In degrees $i=1,2$, the groups $K_i(A)$ are compared with $K_i(F)$ by the reduced norms $\Nrd_i\colon K_i(A)\to K_i(F)$ \cite{MerkurjevSuslin1982,Suslin1987}.
We study these maps over nonarchimedean local fields, global fields, and two- and three-dimensional local fields.

\subsection*{Background}

Let $F$ be a local or global field.
Merkurjev and Suslin conjectured that the sequence
\begin{equation}\label{eq:intro-MS}
 0\longrightarrow K_2(A)\xrightarrow{\ \Nrd_2\ }K_2(F)
 \xrightarrow{\ (s_v)_{v}\ }(\Z/2)^{S_{\mathbb R}}\longrightarrow0
\end{equation}
is exact, where $S_{\mathbb R}$ is the set of real places of $F$ at which $A$ is ramified and $s_v$ is the Hilbert symbol at $v$ \cite[Remark 17.5]{MerkurjevSuslin1982}.
They proved~\eqref{eq:intro-MS} when the index of $A$ is squarefree and prime to the characteristic of $F$ \cite[Theorem 17.4]{MerkurjevSuslin1982}.
Merkurjev and Suslin and, independently, Kahn and Levine proved the injectivity of $\Nrd_2$ for $A$ of squarefree index over every field of characteristic zero \cite[Theorem 2.3]{MerkurjevSuslin2010}, \cite[Corollary 2]{KahnLevine2010}.
For global fields, Bak and Rehmann identified the image of $\Nrd_2$ with the kernel of $(s_v)_v$ \cite[\S2, Theorem 3]{BakRehmann1984}, so the main point of~\eqref{eq:intro-MS} is the injectivity of $\Nrd_2$.
This injectivity over global fields has sometimes been attributed to \cite{BakRehmann1984}; as pointed out in \cite[Remark 10.4(3)]{BHRAdeles}, however, \cite{BakRehmann1984} does not claim it.

The corresponding statements in degree one are classical \cite{NakayamaMatsushima1943,Wang1950}, \cite[Corollary 4.10(iii)]{HazratWadsworth}, \cite[Theorem 33.15]{ReinerMaximalOrders}, \cite[Corollary 2]{Yanchevskii1975}.

Let $F$ be a nonarchimedean local field with residue characteristic $p$, and let $D$ be the division algebra in the Brauer class of $A$.
Suslin and Yufryakov showed that, if $\deg D$ is prime to $p$, then for every prime $\ell\ne p$ and every $j\ge1$ there is an isomorphism $\pi_jK(D)^\wedge_\ell\cong\pi_jK(F)^\wedge_\ell$ whose product with $\deg D$ is induced by the norm $N_{D/F}$ \cite[Theorem 3]{SuslinYufryakov}, and Hesselholt, Larsen and Lindenstrauss proved the analogous statement for $\ell=p$ when $p$ is odd \cite[Theorem A]{HesselholtLarsenLindenstrauss2020}.

Clausen described the Artin map of class field theory through the $K$-theory of locally compact $F$-modules \cite{ClausenArtin}.
Braunling, Henrard and van Roosmalen developed a noncommutative analogue, which relates the id\`eles of $A$ to $K_1$ of the category $\operatorname{LCA}_A$ of locally compact $A$-modules \cite{BHRAdeles}.
Reduced norms are closely related to noncommutative $K$-theoretic class field theory; for example, the Merkurjev--Suslin conjecture implies the isomorphism $K_2(\operatorname{LCA}_A)\cong K_2(\operatorname{LCA}_F)$ of \cite[Conjecture 10.2]{BHRAdeles}; see \cite[Theorem 10.5]{BHRAdeles}.

\subsection*{Main results}

\begin{introthm}[Merkurjev--Suslin conjecture, Theorems~\ref{thm:paper-local-norm} and~\ref{thm:paper-global-norm}]\label{thm:intro-K2}
Let $F$ be a nonarchimedean local field or a global field and $A$ a central simple algebra over $F$.
Then $\Nrd_2$ is injective.
In particular, the sequence
\[
 0\longrightarrow K_2(A)\xrightarrow{\ \Nrd_2\ }K_2(F)
 \xrightarrow{\ (s_v)_{v}\ }(\Z/2)^{S_{\mathbb R}}\longrightarrow0
\]
is exact.
\end{introthm}

Hesselholt, Larsen and Lindenstrauss remark that it is not known whether their reduced norm isomorphism is induced by a map of $K$-theory spectra, and that such a map may well not exist \cite[Introduction]{HesselholtLarsenLindenstrauss2020}; compare \cite[Proposition 4]{Merkurjev1995}.
We construct a map of truncated $K$-theory spectra that induces the reduced norms in degrees $1$ and $2$.
Moreover, over local fields we obtain an equivalence of spectra.

\begin{introthm}[Theorems~\ref{thm:paper-local-norm} and~\ref{thm:paper-local-spectrum}]\label{thm:intro-local}
Let $F$ be a nonarchimedean local field, and let $A$ be a central simple algebra over $F$.
There is an equivalence of spectra
\[
 \nu\colon\tau_{\ge1}K(A)\xrightarrow{\ \sim\ }\tau_{\ge1}K(F).
\]
On $\pi_j$ it induces $\Nrd_j$ for $j=1,2$.
\end{introthm}

\begin{introthm}[Theorems~\ref{thm:paper-global-comparison}, \ref{thm:paper-global-split}, \ref{thm:paper-global-function-split}, \ref{thm:paper-global-fiber-bounds} and~\ref{thm:paper-global-fiber-gem}, Corollaries~\ref{cor:global-fibre-connectivity} and~\ref{cor:paper-global-fiber-prime-four}]\label{thm:intro-global}
Let $F$ be a global field and let $A$ be a central simple algebra over $F$.
If $A$ splits at every real place of $F$, all statements are integral; if $A$ is ramified at some real place, all spectra and groups are taken with $\Z[1/2]$-coefficients, and the prime $\ell=2$ is excluded.
\begin{enumerate}
\item There is a map of spectra
\[
 g\colon\tau_{\ge1}K(A)\longrightarrow\tau_{\ge1}K(F)
\]
that agrees with $\Nrd_1$ and $\Nrd_2$ in degrees one and two.
\item Let $Z=\fib(g)$.
There is an equivalence
\[
 \tau_{\ge2}K(A)\simeq\tau_{\ge2}K(F)\oplus Z
\]
whose first component is $\tau_{\ge2}g$, and $Z\in\mathrm{Sp}_{\ge4}$.
Thus, for $j\ge2$, the map $g_j=\pi_jg$ is split surjective with kernel $\pi_jZ$, and $g_2$ and $g_3$ are isomorphisms.
\item For a prime $\ell\mid\ind A$ with $\ell\ne\operatorname{char}F$, put $e_\ell=\ell^{v_\ell(\ind A)}$, and let $B_\ell(A)$ be the fiber of the localization $K(A)_{(\ell)}\to L_{E(1),\ell}K(A)$.
Then
\[
 Z\simeq\prod_{\ell}\tau_{\ge2}B_\ell(A),\qquad
 e_\ell\,\pi_jB_\ell(A)=0\ \ (j\ge1),\qquad
 B_\ell(A)\in\mathrm{Sp}_{\le e_\ell-3},
\]
and $\tau_{\ge2}B_\ell(A)=0$ if $e_\ell=\ell$.
\item
If $F$ is a global function field, then for $j\ge2$
\[
 \pi_jZ\cong K_j(A)\Big/\sum_{L}\theta_LK_j(L),
\]
where $L$ runs over the finite separable extensions of $F$ that split $A$ and $\theta_L\colon K(L)\to K(A)$ is the splitting transfer.
\item If $F$ is a number field, there is a noncanonical equivalence
\[
 Z\simeq\bigoplus_{\ell\mid\ind A}\ \bigoplus_{j=4}^{e_\ell-3}\Sigma^jH\bigl(\pi_jZ_{(\ell)}\bigr),
\]
where each $\pi_jZ_{(\ell)}$ is a finite $\ell$-group, namely the kernel of $(g_j)_{(\ell)}\colon K_j(A)_{(\ell)}\to K_j(F)_{(\ell)}$.
\end{enumerate}
\end{introthm}

Consequently, the $\ell$-primary part of $\pi_jZ$ can be nonzero only if $4\le j\le e_\ell-3$ and $e_\ell\ne\ell$.
In particular, $g$ is an equivalence if $e_\ell\in\{\ell,4\}$ for every prime $\ell\mid\ind A$ with $\ell\ne\operatorname{char}F$.

Inverting $2$ is necessary: an explicit counterexample is given in Example~\ref{ex:paper-global-k3-real-ramification}.

In degree three the comparison is explicit, and the following descriptions are integral, also in the presence of real ramification (Theorem~\ref{thm:global-k3} and Corollary~\ref{cor:global-k3-structure}).
If $F$ is a global function field with constant field $\F_q$, then $K_3(A)\cong K_3(F)\cong\Z/(q^2-1)$.
If $F$ is a number field, then $K_3(A)\cong K_3(F)/J_F^{\mathrm{ram}}$, where $J_F^{\mathrm{ram}}$ is the part of the image of $K_3^M(F)\to K_3(F)$ supported at the real places at which $A$ is ramified.
Explicitly, let $r_2$ be the number of complex places of $F$, let $r_{\mathrm{sp}}$ be the number of real places at which $A$ splits, and put $w_2=|H^0(F,\Q/\Z(2))|$; then $K_3(A)$ has free rank $r_2$, and its torsion subgroup is $\Z/w_2$ if $r_{\mathrm{sp}}=0$ and $\Z/2w_2\oplus(\Z/2)^{r_{\mathrm{sp}}-1}$ if $r_{\mathrm{sp}}>0$.

Up to equivalence, the fiber $Z$ depends only on the cyclic subgroup of $\Br(F)$ generated by $[A]$.

\begin{introthm}[Proposition~\ref{prop:paper-brauer-power-fiber}]\label{thm:intro-fiber}
Let $F$ be a global field and $A$ a central simple algebra over $F$, use the coefficient convention of Theorem~\ref{thm:intro-global}, and write $Z_A$ for the fiber of the comparison $g$ of Theorem~\ref{thm:intro-global} for $A$.
If $a\ge1$ and $\gcd(a,\per A)=1$, then there is a noncanonical equivalence $Z_A\simeq Z_{A'}$ for every central simple algebra $A'$ over $F$ with $[A']=a[A]$ in $\Br(F)$.
\end{introthm}

\begin{introthm}[Merkurjev--Suslin conjecture for higher local fields, Theorems~\ref{thm:second-two-surj}, \ref{thm:second-two-vanishing} and~\ref{thm:second-three-coker}, Corollaries~\ref{cor:review-two-local-sk1-zero} and~\ref{cor:second-three-charzero}]\label{thm:intro-higher}
Let $F$ be a two- or three-dimensional local field, let $A$ be a central simple algebra over $F$, put $\alpha=[A]\in\Br(F)$, and let $N_2(\alpha)\subseteq K_2(F)$ be the subgroup generated by the norms $N_{E/F}K_2(E)$ from all finite extensions $E/F$ that split $\alpha$.
\begin{enumerate}
\item If $F$ is two-dimensional, then $\Nrd_1$ is injective, $N_2(\alpha)=K_2(F)$, and $\Nrd_2\colon K_2(A)\to K_2(F)$ is an isomorphism.
\item If $F$ is three-dimensional, then
\[
 N_2(\alpha)=\Nrd_2K_2(A)=\ker\bigl(\alpha\cup-\colon K_2(F)\to H^4(F,\Q/\Z(3))\bigr),
\]
and $\inv_F(\alpha\cup-)$ induces an isomorphism $\coker(\Nrd_2)\cong\frac1{\per\alpha}\Z/\Z$.
If moreover the first residue field of $F$ has characteristic zero, then $\Nrd_2$ is injective.
\end{enumerate}
\end{introthm}

Together with Theorem~\ref{thm:intro-K2}, Theorem~\ref{thm:intro-higher} shows that $\SK_2(A)=\ker(\Nrd_2)$ vanishes over every local field of dimension at most two (in dimension zero, $A$ is split by Wedderburn's theorem).
In dimension four and higher the vanishing fails (Remark~\ref{rem:paper-sk2-dimension-four}).
This parallels the behavior of $\SK_1(A)=\ker(\Nrd_1)$, which vanishes over local fields of dimension at most two (\cite{NakayamaMatsushima1943,HazratWadsworth} and Theorem~\ref{thm:intro-higher}(1)) but not in dimension three \cite{Platonov1976}.
We conjecture that $\SK_2(A)=0$ over every three-dimensional local field (Conjecture~\ref{conj:paper-three-local-sk2}) and two- and three-dimensional global fields (Conjecture~\ref{conj:paper-global-sk2}).

We apply these results to noncommutative class field theory.

\begin{introthm}[Theorems~\ref{thm:paper-lca-comparison}, \ref{thm:paper-lca-fiber-splitting} and~\ref{thm:paper-lca-artin}, Corollary~\ref{cor:paper-lca-gem}]\label{thm:intro-LCA}
Let $F$ be a number field, let $A$ be a central simple algebra over $F$, and use the coefficient convention of Theorem~\ref{thm:intro-global}.
\begin{enumerate}
\item The comparisons of Theorems~\ref{thm:intro-local} and~\ref{thm:intro-global} induce a map $c^{\mathrm{lc}}_A\colon\tau_{\ge1}K(\operatorname{LCA}_A)\to\tau_{\ge1}K(\operatorname{LCA}_F)$, compatible with the localization sequences of \cite{BHRAdeles}, such that $\fib(c^{\mathrm{lc}}_A)\simeq\Sigma Z$.
There is a noncanonical equivalence
\[
 \tau_{\ge2}K(\operatorname{LCA}_A)\simeq\tau_{\ge2}K(\operatorname{LCA}_F)\oplus\Sigma Z.
\]
In particular, $c^{\mathrm{lc}}_A$ is an isomorphism in degrees $1\le n\le4$.
\item For every prime $p$, odd if $A$ is ramified at a real place, the comparison $c^{\mathrm{lc}}_A$ is compatible with Clausen's pro-$p$ Artin maps \cite{ClausenArtin}, and the Artin map of $A$ vanishes on $\fib(c^{\mathrm{lc}}_A)$.
On id\`eles of $A$, the resulting Artin map sends $(a_v)_v$ to $\operatorname{rec}_{F,p}\bigl((\Nrd_{A_v/F_v}(a_v))_v\bigr)$, where $\operatorname{rec}_{F,p}$ is the global Artin map followed by the projection to the maximal pro-$p$ quotient.
\end{enumerate}
\end{introthm}

Note that the Merkurjev--Suslin conjecture implies the isomorphism $K_2(\operatorname{LCA}_A)\cong K_2(\operatorname{LCA}_F)$ \cite[Theorem 10.5]{BHRAdeles} (Theorem~\ref{thm:paper-lca-k2}).

\subsection*{Organization}

Section~\ref{sec:local} treats nonarchimedean local fields, comparing $K(D)$ with its \'etale hypersheafification by the methods of Clausen and Mathew \cite{ClausenMathew}, with input from \cite{SuslinYufryakov}, \cite{HesselholtLarsenLindenstrauss2020}, and a Frobenius comparison of twisted cycle complexes (Appendix~\ref{app:second-frobenius}).
Section~\ref{sec:splitting-transfers} proves Theorem~\ref{thm:intro-K2} over global fields from the norm theorem of Bak and Rehmann \cite{BakRehmann1984}, cyclic Hilbert~90 for $K_2$ \cite[Theorem 14.1]{MerkurjevSuslin1982}, and a local--global principle for $K_2(D)$ modulo splitting transfers, proved with Picard norms on curves.
Section~\ref{sec:global} constructs the global comparison with the twisted motivic spectral sequence of Kahn and Levine \cite{KahnLevine2010}, arithmetic squares, transfers from homotopy orbits, Selmer $K$-theory \cite{ClausenMathew}, and Quillen's decomposition of the $K$-theory of Severi--Brauer varieties.
Sections~\ref{sec:two-local} and~\ref{sec:three-local} use Kato's higher local class field theory \cite{KatoExistence,KatoLocalII,Kurihara2000} and, for wild ramification, relative topological cyclic homology of graded covers (Appendix~\ref{app:second-tc}).
Section~\ref{sec:lca} treats locally compact modules with the localization sequences of \cite{BHRAdeles} and Clausen's Selmer $K$-homology \cite{ClausenArtin}.
Appendix~\ref{app:second-frobenius} treats fields $F$ of characteristic $p$ with $[F:F^p]=p$: it compares twisted and ordinary cycle complexes when the period is a power of $p$ and shows that $K_i(D;\Z/p^a)=0$ for $i\ge2$, which is the input at the characteristic prime in Sections~\ref{sec:local}--\ref{sec:global}.
Appendix~\ref{app:second-tc} compares relative $K$-theory of maximal orders with relative $\mathrm{TC}$ at finite truncations and represents the resulting classes by integral classes, as needed for wild ramification in Section~\ref{sec:two-local} and for deep principal units in Section~\ref{sec:three-local}.
Appendix~\ref{app:review-cycle-transfer} collects facts on zero-cycles with coefficients.

\section*{Acknowledgements}
\addcontentsline{toc}{section}{Acknowledgements}
I am deeply grateful to Seidai Yasuda for his invaluable support during the preparation of this paper.
I also thank Ryunosuke Nakano for his help with proofreading.
This work is supported by JSPS KAKENHI Grant Number JP26KJ0476.

\section*{Use of AI}
\addcontentsline{toc}{section}{Use of AI}
All proofs in this paper were given by OpenAI's ChatGPT (GPT-6.1 Sol).
This manuscript was written with the assistance of Anthropic's Claude (Claude Opus 5.5) and OpenAI's ChatGPT (GPT-6.1 Sol and GPT-6.2 Sol).
The author takes full responsibility for the entire content of this manuscript.

\section*{Conventions}
\addcontentsline{toc}{section}{Conventions}

\paragraph{Twisted $K$-theory.}
We use nonconnective $K$-theory.
For a central simple algebra $A$ over $F$, we write $K(A)$ for the $K$-theory spectrum of finitely generated projective right $A$-modules.
We write $\per A$ and $\ind A$ for its period and index, that is, the order of $[A]$ in $\Br(F)$ and the degree of the division algebra \(D\) in its class; thus $\ind A=\deg D$.
Morita equivalence identifies $K(A)$ with $K(D)$, where $D$ is the division algebra in the Brauer class of $A$, and all our statements are invariant under this identification.
From now on, we therefore identify $A$ with $D$ in this way and state all results for the division algebra $D$.

The reduced rank \(\Nrd_0\) identifies $K_0(D)$ with the subgroup $\deg D\cdot\Z$ of $K_0(F)=\Z$, the reduced norm $\Nrd_1$ sends the class of $x\in D^\times$ to the determinant of $x\otimes1$ under $D\otimes_FL\simeq M_{\deg D}(L)$ for a splitting field $L$, and $\Nrd_2$ is defined through the $K$-cohomology of the Severi--Brauer variety of $D$ \cite{MerkurjevSuslin1982,Suslin1987}.
For a real place $v$ of $F$, we write the Hilbert symbol $s_v\colon K_2(F)\to\Z/2$ additively: $s_v\{a,b\}=1$ if $a$ and $b$ are both negative at $v$, and $s_v\{a,b\}=0$ otherwise.

\paragraph{Local fields.}
An \emph{$h$-dimensional local field} is a field $F=F_h$ for which there are fields $F_{h-1},\ldots,F_0$ such that $F_0$ is finite and, for $1\le i\le h$, the field $F_i$ is complete discretely valued with residue field $F_{i-1}$.
We call $F_{h-1}$ the \emph{first residue field} and $F_0$ the \emph{last residue field} of $F$, and $(\operatorname{char}F_h,\ldots,\operatorname{char}F_1)$ its \emph{characteristic sequence}.
For a three-dimensional local field $F$, we write $\inv_F\colon H^4(F,\Q/\Z(3))\xrightarrow{\sim}\Q/\Z$ for Kato's isomorphism \cite[Theorem 2]{KatoExistence}.
For a field $F$ of characteristic $p$, we set
\[
 H^{a+1}(F,\Z/p^r(a))=H^1(F,W_r\Omega^a_{\log}),
\]
where $W_r\Omega^a_{\log}$ is the logarithmic de Rham--Witt sheaf.

\paragraph{Spectra.}
We write $\mathrm{Sp}$ for the category of spectra, $\mathbb S$ for the sphere spectrum, $\otimes=\wedge$ for the smash product, and $\Sigma$ for suspension.
The full subcategories $\mathrm{Sp}_{\ge a}$ and $\mathrm{Sp}_{\le b}$ consist of the spectra whose homotopy groups vanish below $a$ and above $b$, respectively, and $\tau_{\ge a}$ and $\tau_{\le b}$ denote the $(a-1)$-connected cover and the $b$th Postnikov truncation.
For an abelian group $M$, we write $HM$ for its Eilenberg--Mac Lane spectrum.
For a spectrum $T$, a prime $\ell$, and an integer $n\ge1$, we write
\(
 T_{\Q}=T\otimes H\Q\), \(T_{(\ell)}=T\otimes\mathbb S_{(\ell)}\), \(T[1/n]\)
for the rationalization, the localization at $\ell$, and the inversion of $n$, and
\(
 T^\wedge_\ell=\operatorname*{holim}_rT/\ell^r
\)
for the derived $\ell$-completion, where $T/\ell^r$ is the cofiber of multiplication by $\ell^r$.
We distinguish $\pi_j(T^\wedge_\ell)$ from the $\ell$-adic completion $\varprojlim_r\pi_j(T)/\ell^r$ of the group $\pi_j(T)$.
The functor $L_{K(1),\ell}$ is Bousfield localization with respect to height-one Morava $K$-theory at $\ell$, and $L_{E(1),\ell}$ is the corresponding localization with respect to height-one Johnson--Wilson theory.
For a finite group $G$ acting on $T$, we write
\[
 T_{hG}=\operatorname*{colim}_{BG}T,\qquad T^{hG}=\operatorname*{lim}_{BG}T
\]
for the homotopy orbits and homotopy fixed points, and for a $G$-module $M$ we write $M_G$ for the ordinary coinvariants.

\section{Twisted \texorpdfstring{$K$}{K}-theory over local fields}\label{sec:local}
This section focuses on nonarchimedean local fields.
We first identify the positive-degree $K$-groups of a central division algebra with those of its center (Section~\ref{sec:paper-local-norm}).
We then lift these isomorphisms to an equivalence of positive-degree spectra (Section~\ref{sec:paper-local-spectrum}).
\subsection{\'Etale descent and reduced norms}\label{sec:paper-local-norm}
The main result of this subsection is Theorem~\ref{thm:paper-local-norm}, which identifies the positive-degree $K$-groups of a central division algebra over a nonarchimedean local field with those of the field, compatibly with splitting transfers and field norms.
We first fix the conventions on norms and algebras that are used in the rest of the paper.
Let $R\to S$ be a homomorphism of rings such that $S$ is perfect as a right $R$-module.
Restriction of scalars $\operatorname{Perf}(S)\to\operatorname{Perf}(R)$ induces the \emph{norm}
\[
 N_{S/R}\colon K(S)\longrightarrow K(R).
\]
The norms commute with flat base change: if $R\to S$ is a homomorphism of algebras over a commutative ring $k$ and $k'$ is a flat commutative $k$-algebra, then $N_{S\otimes_kk'/R\otimes_kk'}\circ(-\otimes_kk')=(-\otimes_kk')\circ N_{S/R}$.
For composable maps $R\to S\to T$ with $S$ perfect as a right $R$-module and $T$ perfect as a right $S$-module, we have $N_{T/R}=N_{S/R}\circ N_{T/S}$.

Let $F$ be a field, let $A$ be a central simple algebra over $F$, and let $D$ be the division algebra in the Brauer class of $A$.
Since Morita equivalence identifies $K(A)$ with $K(D)$, we work with $D$ instead of $A$ from now on.
We write $i_D\colon K(F)\to K(D)$ for extension of scalars.
We write $\Nrd_i\colon K_i(D)\to K_i(F)$ for the reduced norms ($i=1,2$).
We write $F^{\mathrm{sep}}$ for a separable closure of $F$.
For a field extension $E/F$, put $D_E=D\otimes_FE$, a central simple $E$-algebra that need not be a division algebra, and let $D'$ be the division algebra in its Brauer class.
When we apply a result stated for central division algebras to $D_E$, we identify $K(D_E)$ with $K(D')$ by the Morita equivalence given by a simple $(D',D_E)$-bimodule with compatible $E$-actions.
Norms are still taken for $D_E$ itself: if $D_E\simeq M_r(D')$, then $N_{D_E/E}=r\,N_{D'/E}$ under this identification.
If $E/F$ is finite and splits $D$, we identify $K(E)$ with $K(D_E)$ by the Morita equivalence given by a simple $(E,D_E)$-bimodule on which the two actions of $E$ agree; this identification does not depend on the choice of the simple bimodule.
The splitting transfer is the composite
\[
 \theta_E\colon K(E)\simeq K(D_E)\xrightarrow{\ N_{D_E/D}\ }K(D).
\]
A simple $(E,D_E)$-bimodule has dimension $\deg D$ over $E$, so $N_{D_E/E}=\deg D$ under this identification, and transitivity of norms gives
\[
 N_{D/F}\circ\theta_E=\deg D\cdot N_{E/F}.
\]

Throughout this subsection, $F$ is a nonarchimedean local field with residue field $k$ of characteristic $p$.

\begin{theorem}\label{thm:paper-local-norm}
Let $F$ be a nonarchimedean local field, and let $D$ be a central division algebra over $F$.
For each $j\ge1$ there is an isomorphism
\[
 \nu_{j,D}\colon K_j(D)\xrightarrow{\ \sim\ }K_j(F)
\]
satisfying $\nu_{j,D}\theta_E=N_{E/F}$ for every finite splitting field $E/F$ of $D$.
For a finite separable extension $F'/F$ it satisfies $\nu_{j,D_{F'}}\circ\res_{F'/F}=\res_{F'/F}\circ\nu_{j,D}$ and $\nu_{j,D}\circ N_{D_{F'}/D}=N_{F'/F}\circ\nu_{j,D_{F'}}$, with the convention above for $D_{F'}$, and $\deg D\cdot\nu_{j,D}=\pi_jN_{D/F}$.
In degrees $1$ and $2$ it is $\Nrd_1$ and $\Nrd_2$, respectively.
\end{theorem}
The proof of Theorem~\ref{thm:paper-local-norm} is given at the end of this subsection.

Let $K^{\et}(D)$ be the spectrum of global sections of the hypersheafification of $L\mapsto K(D_L)$ on the small \'etale site of $F$, and let
\[
 c_D\colon K(D)\longrightarrow K^{\et}(D)
\]
be the \'etale comparison map.
For a finite separable extension $E/F$, the commutative \'etale comparison $c_E\colon K(E)\to K^{\et}(E)$ induces isomorphisms on $\pi_i$ for $i\ge0$ \cite[Theorems 7.12(2) and 7.13]{ClausenMathew}.
At each prime $\ell$ the dimension bound in that theorem is at most two: it is $\vcd_\ell(E)$ for $\ell\ne\operatorname{char}E$, and it is $1+\log_p[E:E^p]=2$ for $\ell=p=\operatorname{char}E$.
Write $G_F=\Gal(F^{\mathrm{sep}}/F)$ and $G_k=\Gal(k^{\mathrm{sep}}/k)$.
The proof of Theorem~\ref{thm:paper-local-norm} compares $K(D)$ with its \'etale sheafification, whose geometric stalks agree with those of $F$ by Morita equivalence.
Lemma~\ref{lem:paper-local-complete-high} shows that the comparison is an isomorphism after completion in degrees at least two.
Lemma~\ref{lem:paper-local-galois-vanishing} gives the vanishing range of the Galois cohomology that enters integral descent, and Lemma~\ref{lem:paper-local-low-descent} uses it to show that $\pi_0$ of the \'etale theory is torsion-free and that the comparison on $\pi_1$ is $\Nrd_1$.
Lemma~\ref{lem:paper-local-compare} combines Lemmas~\ref{lem:paper-local-complete-high} and~\ref{lem:paper-local-low-descent} into an integral comparison, and Lemma~\ref{lem:paper-local-even-coinvariants} describes the even degrees by coinvariants.

\begin{lemma}\label{lem:paper-local-complete-high}
For every prime $\ell$, the map
\[
 (c_D)^\wedge_\ell\colon K(D)^\wedge_\ell\longrightarrow K^{\et}(D)^\wedge_\ell
\]
induces an isomorphism on $\pi_j$ for every $j\ge2$.
If $\ell\ne\operatorname{char}F$, the localization unit
\[
 K(D)^\wedge_\ell\longrightarrow L_{K(1),\ell}K(D)
\]
induces an isomorphism on $\pi_j$ for every $j\ge2$.
\end{lemma}
\begin{proof}
We first compute the targets after derived $\ell$-completion.
Let $\ell\ne\operatorname{char}F$ and $a\ge1$.
By Morita equivalence, the geometric stalks are
\[
 \pi_t\bigl(K(D_{F^{\mathrm{sep}}})/\ell^a\bigr)\simeq\pi_t\bigl(K(F^{\mathrm{sep}})/\ell^a\bigr)\simeq
 \begin{cases}
  \Z/\ell^a(b), & t=2b\ge0,\\
  0, & \text{otherwise},
 \end{cases}
\]
by Gabber--Suslin rigidity and the computation of $K(F^{\mathrm{sep}})/\ell^a$ \cite{SuslinClosed}, and for $L_{K(1),\ell}K(D_{F^{\mathrm{sep}}})/\ell^a$ the same formula holds with $\Z/\ell^a(b)$ in every even degree $t=2b\in\Z$.
The height-one localization satisfies \'etale hyperdescent, since it is the value of a weakly localizing invariant with tensor coefficients $\operatorname{Perf}(D)$ \cite[Example 5.12 and Theorem 7.14]{ClausenMathew}.
The descent spectral sequences therefore have
\[
 E_2^{s,2b}=H^s(F,\Z/\ell^a(b)),\qquad E_2^{s,t}=0\quad(t\text{ odd}).
\]
Since the geometric degrees are even and $\cd_\ell(F)\le2$, the differentials vanish, and since these groups are finite, their inverse systems are Mittag--Leffler.
By hyperdescent of the height-one localization, the localization unit factors as $K(D)^\wedge_\ell\xrightarrow{(c_D)^\wedge_\ell}K^{\et}(D)^\wedge_\ell\to L_{K(1),\ell}K(D)$.
The descent calculations for $K^{\et}(D)/\ell^a$ and $L_{K(1),\ell}K(D)/\ell^a$ agree in positive degrees, so, after passage to the inverse limit, the second map induces isomorphisms on $\pi_j$ for $j\ge1$.
Thus the two assertions of the lemma are equivalent for $\ell\ne\operatorname{char}F$.
Passing to the inverse limit over $a$ gives, for $m\ge1$,
\[
 \pi_{2m-1}\bigl(K^{\et}(D)^\wedge_\ell\bigr)\simeq H^1(F,\Z_\ell(m)),\qquad
 \pi_{2m}\bigl(K^{\et}(D)^\wedge_\ell\bigr)\simeq H^2(F,\Z_\ell(m+1));
\]
the other term in degree $2m$, namely $H^0(F,\Z_\ell(m))$, vanishes, since the $\ell$-adic cyclotomic character of $F$ has infinite image.
Under geometric Morita equivalence, $N_{D/F}$ acts on these groups as multiplication by $\deg D$, and $\theta_E$ acts as $\cor_{E/F}$.

\emph{Case $\ell\ne p$.}
Let $E/F$ be the unramified extension of degree $\deg D$, which splits $D$.
We first show that
\[
 \cor_{E/F}\colon\pi_j\bigl(K^{\et}(D_E)^\wedge_\ell\bigr)\longrightarrow\pi_j\bigl(K^{\et}(D)^\wedge_\ell\bigr)
\]
is surjective for $j\ge2$.
Let $I_F\subset G_F$ be the inertia subgroup, so that $G_F/I_F\simeq G_k$ and $\cd(G_k)=1$.
Since $\ell\ne p$, we have $H^1(I_F,\Z_\ell(m))\simeq\Z_\ell(m-1)$ as $G_k$-modules.
For $m\ge2$ we have $H^0(k,H^1(I_F,\Z_\ell(m)))=H^0(k,\Z_\ell(m-1))=0$, and hence
\[
 H^1(F,\Z_\ell(m))\simeq H^1(k,\Z_\ell(m))\quad(m\ge2),\qquad
 H^2(F,\Z_\ell(m+1))\simeq H^1(k,\Z_\ell(m))\quad(m\ge1).
\]
On $H^1(k,\Z_\ell(m))$, corestriction along the unramified extension $E/F$ is the quotient map onto Frobenius coinvariants, so it is surjective.
To deduce surjectivity of $\pi_j((c_D)^\wedge_\ell)$ for $j\ge2$, we lift a class along the isomorphism $\pi_j((c_E)^\wedge_\ell)$, where $K(D_E)\simeq K(E)$ by Morita equivalence, and use
\[
 c_D\circ\theta_E=\cor_{E/F}\circ c_E.
\]

For injectivity, suppose first that $p\nmid\deg D$.
Suslin and Yufryakov compute $K_*(D)$ by reduction of the maximal order and localization \cite[Lemma 2 and Theorem 3]{SuslinYufryakov}.
It follows that, for $j\ge2$,
\[
 \#\varprojlim_aK_j(D)/\ell^a=\#\pi_j\bigl(K^{\et}(D)^\wedge_\ell\bigr)<\infty.
\]
Their decomposition and the description of $K_1(D)$ show that $K_{j-1}(D)[\ell^\infty]$ is finite, so its Tate module $\varprojlim_aK_{j-1}(D)[\ell^a]$ in the Bockstein sequence vanishes.
Since the finite-coefficient groups are finite, the Moore towers satisfy the Mittag--Leffler condition and the Milnor $\varprojlim^1$ terms vanish; the surjection of the preceding paragraph is therefore a surjection between finite groups of the same order, hence an isomorphism.
For general $D$, let $F'/F$ be the unramified extension of degree $p^{v_p(\deg D)}$, so that the index of $D_{F'}$ is prime to $p$.
Since $\cor_{F'/F}\circ\res_{F'/F}=[F':F]\in\Z_\ell^\times$, restriction is injective on $\pi_j(K(D)^\wedge_\ell)$, so the case over $F'$ shows that the kernel of the map induced by $(c_D)^\wedge_\ell$ vanishes over $F$.

\emph{Case $\operatorname{char}F=0$ and $\ell=p$.}
In this case we use the trace comparison with closed support of Hesselholt--Larsen--Lindenstrauss \cite[Theorem B(1) and Theorem 2.5]{HesselholtLarsenLindenstrauss2020}, which gives
\[
 \tau_{\ge2}K(D)^\wedge_p\simeq\tau_{\ge2}K(F)^\wedge_p,
\]
also in residue characteristic two by the paragraph following \cite[Theorem B]{HesselholtLarsenLindenstrauss2020}.
Since $L_{K(1),p}HM\simeq0$ for every abelian group $M$, a connective spectrum $T$ satisfies $L_{K(1),p}\tau_{\ge2}T\simeq L_{K(1),p}T$.
Hence, for connective $T$, the localization unit induces an isomorphism on $\pi_j$ for $j\ge2$ if and only if it does so for $\tau_{\ge2}T$.
By the Quillen--Lichtenbaum comparison for $F$ \cite[Theorems 6.13(1) and 7.13--7.14]{ClausenMathew}, the map $K(F)^\wedge_p\to L_{K(1),p}K(F)$ induces an isomorphism on $\pi_j$ for $j\ge2$.
The result for $D$ follows from naturality of the localization unit and the displayed equivalence.
By the factorization of the localization unit above, which holds also for $\ell=p$ in characteristic zero \cite[Example 6.3]{ClausenMathew}, $(c_D)^\wedge_p$ induces an isomorphism on $\pi_j$ for $j\ge2$.

\emph{Case $\operatorname{char}F=\ell=p$.}
Then $[F:F^p]=p$, and Proposition~\ref{prop:app-frobenius}(2) gives
\[
 \pi_j\bigl(K(D)/p^a\bigr)=0\qquad(j\ge2).
\]
Since $[F^{\mathrm{sep}}:(F^{\mathrm{sep}})^p]=p$ as well, the same vanishing holds for $K(D_{F^{\mathrm{sep}}})/p^a$ after split Morita equivalence.
This upper bound is preserved by global sections of the \'etale hypersheaf and by the inverse limit.
\end{proof}

\begin{lemma}\label{lem:paper-local-galois-vanishing}
For $t\ge2$ and $s\ge1$, we have
\[
 H^s\bigl(F,K_t(F^{\mathrm{sep}})\bigr)=0
\]
unless $s=1$ and $t$ is odd.
\end{lemma}
\begin{proof}
For $t\ge2$ the groups $K_t(F^{\mathrm{sep}})$ are divisible, and for $b\ge1$
\[
 K_{2b}(F^{\mathrm{sep}})\otimes\Q\simeq K_{2b}(F^{\mathrm{sep}}),\qquad
 K_{2b+1}(F^{\mathrm{sep}})_{\mathrm{tors}}\simeq\bigoplus_{\ell\ne\operatorname{char}F}\Q_\ell/\Z_\ell(b+1),
\]
where in mixed characteristic the sum includes $\ell=p$, and $K_{2b+1}(F^{\mathrm{sep}})/K_{2b+1}(F^{\mathrm{sep}})_{\mathrm{tors}}$ is a rational vector space.
For an algebraic closure of $F$ this is Suslin's theorem \cite{SuslinClosed}, and in characteristic zero $F^{\mathrm{sep}}$ is algebraically closed.
In characteristic $p$, restriction from $F^{\mathrm{sep}}$ to its algebraic closure is an isomorphism after inverting $p$, since for every finite purely inseparable extension $F'/F^{\mathrm{sep}}$
\[
 N_{F'/F^{\mathrm{sep}}}\res_{F'/F^{\mathrm{sep}}}=\res_{F'/F^{\mathrm{sep}}}N_{F'/F^{\mathrm{sep}}}=[F':F^{\mathrm{sep}}].
\]
For the $p$-part, $F^{\mathrm{sep}}$ and its finite purely inseparable extensions $F'$ satisfy $[F':F'^p]=p$, so $K_t^M(F')/p=0$ for $t\ge2$ by the logarithmic differentials, and the calculation of Geisser and Levine together with the Bockstein sequences in degrees $t$ and $t+1$ shows that $K_t(F^{\mathrm{sep}})$ is uniquely $p$-divisible for $t\ge2$ \cite{GeisserLevine}.

Positive-degree cohomology of discrete rational coefficients vanishes by averaging, which proves the assertion for $t$ even.
For $t=2b+1$ odd and $s\ge2$, the long exact sequence for the torsion subgroup and the rational quotient reduces the assertion to the torsion subgroup.
There $H^s$ vanishes for $s\ge3$ since $\cd_\ell(F)\le2$, and for $s=2$ local Tate duality \cite{TateDuality} and $H^0(F,\Z_\ell(-b))=0$ give
\[
 H^2\bigl(F,\Q_\ell/\Z_\ell(b+1)\bigr)=0.
\]
\end{proof}

For a central simple algebra $B$ over a field $F$, the \emph{geometric rank} of a class in $\pi_0K^{\et}(B)$ is its image in $\pi_0K^{\et}(B_{F^{\mathrm{sep}}})\simeq\pi_0K(F^{\mathrm{sep}})=\Z$ under Morita equivalence; similarly, for $\ell\ne\operatorname{char}F$, the geometric rank of a class in $\pi_0L_{K(1),\ell}K(B)$ is its image in $\pi_0L_{K(1),\ell}K(B_{F^{\mathrm{sep}}})\simeq\pi_0L_{K(1),\ell}K(F^{\mathrm{sep}})=\Z_\ell$.
Classes of geometric rank $1$ are the basic input for the comparisons of Section~\ref{sec:paper-local-spectrum} and Section~\ref{sec:global-spectrum}.

\begin{lemma}\label{lem:paper-local-low-descent}
The group $\pi_0K^{\et}(D)$ is torsion-free, and under the identification $\pi_1K^{\et}(D)\simeq F^\times$ given by the edge map, the map $\pi_1(c_D)$ is the reduced norm isomorphism
\[
 \Nrd_1\colon K_1(D)\xrightarrow{\ \sim\ }F^\times.
\]
Moreover $(c_D)^\wedge_\ell$ induces an isomorphism on $\pi_1$ for every prime $\ell$.
If $\ell\ne\operatorname{char}F$, the localization unit $K(D)^\wedge_\ell\to L_{K(1),\ell}K(D)$ also induces an isomorphism on $\pi_1$.
\end{lemma}
\begin{proof}
We use the integral descent spectral sequence
\[
 E_2^{s,t}=H^s\bigl(F,K_t(F^{\mathrm{sep}})\bigr)\Longrightarrow\pi_{t-s}K^{\et}(D),
\]
whose coefficients are identified with those of $F$ by Morita equivalence.
For any discrete module $W$, we first separate its torsion subgroup and then rationalize the torsion-free quotient; since $\cd_\ell(F)\le2$ for torsion coefficients, we get $H^s(F,W)=0$ for $s>3$, so the spectral sequence converges.
By Lemma~\ref{lem:paper-local-galois-vanishing}, in the rows $t\ge2$ the only possible nonzero terms with $s\ge1$ are $H^1(F,K_t(F^{\mathrm{sep}}))$ with $t$ odd.
In total degree one, $H^1(F,K_2(F^{\mathrm{sep}}))$ and $H^2(F,K_3(F^{\mathrm{sep}}))$ therefore vanish, as do the targets $E_r^{r,r}$ of the differentials leaving $E_r^{0,1}$.
Thus the edge map gives $\pi_1K^{\et}(D)\simeq F^\times$.
The composite of $\pi_1(c_D)$ with this edge map extends scalars to $D_{F^{\mathrm{sep}}}\simeq M_{\deg D}(F^{\mathrm{sep}})$ and then takes the determinant, so it is $\Nrd_1$.
The map $\Nrd_1\colon K_1(D)\to F^\times$ is surjective: $D$ has ramification index $\deg D$ over $F$ and contains an unramified maximal subfield, so $\Nrd_1$ takes values of every valuation, and the norm from the unramified maximal subfield is surjective on units \cite[V, \S2, Proposition 3(b) and its Corollary; XII, \S2, Propositions 1--2 and Exercises 2--3]{SerreLocal}.
It is injective since $\SK_1(D)=0$ \cite[Corollary 4.10(iii)]{HazratWadsworth}; the hypotheses there hold because $D$ is tame and semiramified, with residue field cyclic of degree $\deg D$ over $k$, and the Tate group $\widehat H^{-1}$ of the cyclic residue Galois group appearing there vanishes by Hilbert 90.
In total degree zero,
\[
 H^1\bigl(F,(F^{\mathrm{sep}})^\times\bigr)=0,\qquad H^2\bigl(F,K_2(F^{\mathrm{sep}})\bigr)=0,
\]
by Hilbert 90 and by Lemma~\ref{lem:paper-local-galois-vanishing}, respectively.
Hence the geometric rank map $\pi_0K^{\et}(D)\to H^0(F,\Z)=\Z$ is injective, and $\pi_0K^{\et}(D)$ is torsion-free.

Let $T$ be $K(D)$ or $K^{\et}(D)$.
Since $\pi_0T$ is torsion-free, the Bockstein sequence gives $\pi_1(T/\ell^a)\simeq\pi_1(T)/\ell^a$.
The Milnor sequence for the derived completion is
\[
 0\to\varprojlim\nolimits^1_a\pi_2(T/\ell^a)
 \to\pi_1(T^\wedge_\ell)\to\varprojlim_a\pi_1(T/\ell^a)\to0.
\]
The transition maps of $\{\pi_2(T)/\ell^a\}_a$ are surjective, so the derived inverse limit of the degree-two Bockstein sequence gives
\[
 \varprojlim\nolimits^1_a\pi_2(T/\ell^a)\simeq\varprojlim\nolimits^1_a\pi_1(T)[\ell^a].
\]
The isomorphism $\pi_1(c_D)$ induces isomorphisms on both outer terms, which proves the claim in degree one.

Let $\ell\ne\operatorname{char}F$.
Then $\pi_1L_{K(1),\ell}K(D)\simeq H^1(F,\Z_\ell(1))$, and Kummer theory and Hilbert 90 give $H^1(F,\mu_{\ell^a})\simeq F^\times/(F^\times)^{\ell^a}$.
Since $F^\times[\ell^\infty]$ is finite, the inverse system $\{F^\times[\ell^a]\}_a$ of $\ell$-power roots of unity is Mittag--Leffler, and $\varprojlim^1_aF^\times[\ell^a]=0$.
Thus
\[
 H^1(F,\Z_\ell(1))\simeq\varprojlim_a F^\times/(F^\times)^{\ell^a}.
\]
Together with $\pi_1(c_D)=\Nrd_1$, these isomorphisms show that $K(D)^\wedge_\ell\to L_{K(1),\ell}K(D)$ induces an isomorphism on $\pi_1$ as well.
\end{proof}

\begin{lemma}\label{lem:paper-local-compare}
The \'etale comparison map $c_D\colon K(D)\to K^{\et}(D)$ induces an isomorphism on $\pi_j$ for every $j\ge1$.
For $\ell\ne\operatorname{char}F$, the localization unit $K(D)^\wedge_\ell\to L_{K(1),\ell}K(D)$ also induces an isomorphism on $\pi_j$ for every $j\ge1$.
\end{lemma}
\begin{proof}
Set $Y=\cofib(c_D)$.
Let $E/F$ be a finite separable splitting field and $Y_E=\cofib(c_{D_E})$.
By Morita equivalence and the comparison for fields above, $\pi_iY_E=0$ for $i\ge1$.
Since $\cor_{E/F}\circ\res_{E/F}=[E:F]$ on $\pi_iY$, we obtain
\[
 [E:F]\cdot\pi_iY=0\qquad(i\ge1).
\]
By Lemmas~\ref{lem:paper-local-complete-high} and~\ref{lem:paper-local-low-descent}, $\pi_i(Y^\wedge_\ell)=0$ for $i\ge2$ and $\pi_1(c_D)=\Nrd_1$ is an isomorphism.
The Bockstein sequence for $Y/\ell\simeq Y^\wedge_\ell/\ell$ then gives $\pi_i(Y)[\ell]=0$ for $i\ge2$, and the bounded exponent implies $\pi_iY=0$ for $i\ge2$.
For $\ell\ne\operatorname{char}F$, the same two lemmas give the height-one comparison on $\pi_j$ for $j\ge2$ and on $\pi_1$, respectively.
\end{proof}

To define the even-degree isomorphisms of Theorem~\ref{thm:paper-local-norm}, we describe splitting transfers from Galois extensions in terms of coinvariants.

\begin{lemma}\label{lem:paper-local-even-coinvariants}
Let $L/F$ be a finite Galois extension with group $G$, and let $m\ge1$.
Then the norm $N_{D_L/D}$ induces an isomorphism
\[
 K_{2m}(D_L)_G\xrightarrow{\ \sim\ }K_{2m}(D).
\]
Moreover $K_{2m}(D)/K_{2m}(D)_{\mathrm{tors}}$ is a rational vector space, and $N_{D_L/D}$ maps $(K_{2m}(D_L)_{\mathrm{tors}})_G$ onto $K_{2m}(D)_{\mathrm{tors}}$.
All these assertions hold with $D$ replaced by $F$.
\end{lemma}
\begin{proof}
Set
\[
 M=K_{2m+1}(F^{\mathrm{sep}}),\qquad P=K_{2m}(F^{\mathrm{sep}}).
\]
By integral descent, Lemma~\ref{lem:paper-local-compare} over $L$, and Lemma~\ref{lem:paper-local-galois-vanishing}, we have a natural exact sequence
\[
 0\longrightarrow H^1(L,M)\longrightarrow K_{2m}(D_L)
 \longrightarrow H^0(L,P)\longrightarrow0.
\]
The right-hand term is a rational vector space and the left-hand term is torsion, so $H^1(L,M)$ is the torsion subgroup of $K_{2m}(D_L)$.
For $L=F$ this shows that $K_{2m}(D)/K_{2m}(D)_{\mathrm{tors}}$ is a rational vector space.
Restriction and transfer are compatible with these exact sequences and act on the outer terms by Galois restriction and corestriction, since at geometric points they are the diagonal map and the sum over the factors of the separable base change, respectively.
Under local Tate duality, the dual of the corestriction
\[
 H^1(L,\Q_\ell/\Z_\ell(m+1))_G\longrightarrow H^1(F,\Q_\ell/\Z_\ell(m+1))
\]
is restriction \cite{TateDuality}.
Since the dual coefficients are $\Z_\ell(-m)$ and $H^0(L,\Z_\ell(-m))=0$, inflation--restriction identifies this restriction map with an isomorphism onto the invariants.
For $M/M_{\mathrm{tors}}$ and $P$ we use averaging over rational coefficients: on rational vector spaces with $G$-action, the norm induces an isomorphism from coinvariants to invariants.
We have
\[
 H^1(L,M)\simeq\coker\bigl(H^0(L,M/M_{\mathrm{tors}})\to H^1(L,M_{\mathrm{tors}})\bigr),\qquad H_1(G,H^0(L,P))=0,
\]
so corestriction induces isomorphisms on the coinvariants of the two outer terms, and hence on $K_{2m}(D_L)_G$.
Since the left-hand terms are the torsion subgroups, the isomorphism on their coinvariants, $(K_{2m}(D_L)_{\mathrm{tors}})_G\xrightarrow[\ \sim\ ]{N_{D_L/D}}K_{2m}(D)_{\mathrm{tors}}$, is the second assertion.
The same calculation applies to $F$.
\end{proof}

\begin{proof}[Proof of Theorem~\ref{thm:paper-local-norm}]
We first treat odd degrees $j=2m-1\ge3$.
By Lemma~\ref{lem:paper-local-galois-vanishing}, the only contribution to total degree $j$ in descent is $H^0(F,K_j(F^{\mathrm{sep}}))$, and the differentials leaving it land in vanishing groups.
By Lemma~\ref{lem:paper-local-compare} applied to $D$ and to $F$, the \'etale comparison maps give $K_j(D)\simeq\pi_jK^{\et}(D)$ and $K_j(F)\simeq\pi_jK^{\et}(F)$.
The edge maps
\[
 e_D\colon\pi_jK^{\et}(D)\longrightarrow H^0(F,K_j(F^{\mathrm{sep}})),\qquad
 e_F\colon\pi_jK^{\et}(F)\longrightarrow H^0(F,K_j(F^{\mathrm{sep}}))
\]
identify the groups and are compatible with restriction and transfer.
We define $\nu_{j,D}$ as the composite
\[
 K_j(D)\xrightarrow{\ \pi_j(c_D)\ }\pi_jK^{\et}(D)\xrightarrow{\ e_D\ }H^0(F,K_j(F^{\mathrm{sep}}))\xrightarrow{\ e_F^{-1}\ }\pi_jK^{\et}(F)\xrightarrow{\ \pi_j(c_F)^{-1}\ }K_j(F).
\]
For the normalization $\nu_{j,D}\theta_E=N_{E/F}$, the splitting field $E$ need not be separable.
Let $E_s$ be its maximal separable subextension.
The algebra $E\otimes_FF^{\mathrm{sep}}$ is a product, indexed by the embeddings of $E_s$, of purely inseparable field extensions of $F^{\mathrm{sep}}$, over each of which $D$ splits.
Under the common geometric Morita identification, the component of $\theta_E$ at each factor is the corresponding field norm.
Flat base change and injectivity of the edge map therefore give $\nu_{j,D}\theta_E=N_{E/F}$.
Under the common geometric Morita identification for $e_D$ and $e_F$, restriction of scalars at a geometric point is multiplication by $\deg D$, so injectivity of $e_F\circ\pi_j(c_F)$ gives $\deg D\cdot\nu_{j,D}=\pi_jN_{D/F}$ in odd degrees.

Next let $j=2m\ge2$ be even, and let $L/F$ be a finite Galois splitting field with group $G$.
By Lemma~\ref{lem:paper-local-even-coinvariants}, the maps $\theta_L$ and $N_{L/F}$ induce isomorphisms
\[
 K_{2m}(D)\xleftarrow[\ \sim\ ]{\ \theta_L\ }K_{2m}(L)_G\xrightarrow[\ \sim\ ]{\ N_{L/F}\ }K_{2m}(F),
\]
and we define $\nu_{2m,D}=N_{L/F}\circ\theta_L^{-1}$, where $\theta_L$ and $N_{L/F}$ denote the induced maps on coinvariants.
The induced identification $K_{2m}(L)\simeq K_{2m}(D_L)$ is $G$-equivariant, since conjugation takes the defining simple bimodule to an isomorphic one.
This definition does not depend on $L$: two Galois splitting fields can be compared over their compositum, using the surjectivity of the field norm from the compositum and the transitivity of transfer.
For any finite splitting field $E/F$, the extension $LE/E$ is Galois, so its even-degree field norm is surjective by Lemma~\ref{lem:paper-local-even-coinvariants} applied over $E$, and transitivity of transfer gives $\nu_{2m,D}\theta_E=N_{E/F}$.
For a finite separable extension $F'/F$, choose a Galois splitting field $L/F$ containing $F'$; then $\nu_{2m,D}$ is compatible with restriction and transfer along $F'/F$.
Indeed, the Mackey formula for restriction of tensor products gives the restriction identity on the image of $\theta_L$, and transitivity of transfer gives the transfer identity on the image of the splitting transfer $\theta_{L/F'}\colon K(L)\to K(D_{F'})$ for $D_{F'}$.
The identity $N_{D/F}\theta_L=\deg D\cdot N_{L/F}$ gives $\deg D\cdot\nu_{2m,D}=\pi_{2m}N_{D/F}$.
In degree two, the transfer formula $\Nrd_2\theta_L=N_{L/F}$ \cite[Proposition 6.6.1(1)]{KahnLevine2010} and the surjectivity of $\theta_L$ give $\nu_{2,D}=\Nrd_2$.

It remains to treat degree one.
We set $\nu_{1,D}=\Nrd_1$, which is an isomorphism by Lemma~\ref{lem:paper-local-low-descent}.
Compatibility with restriction follows from base change for reduced norms.
The transfer formula for reduced norms \cite[Proposition 6.6.1(1)]{KahnLevine2010} gives $\Nrd_1\theta_E=N_{E/F}$ and $\Nrd_1\circ N_{D_{F'}/D}=N_{F'/F}\circ\Nrd_1$.
For $x\in D^\times$, the determinant of multiplication by $x$ on the $F$-vector space $D$ is $\Nrd_1(x)^{\deg D}$, so
\[
 \pi_1N_{D/F}=\deg D\cdot\Nrd_1\quad\text{on }K_1(D).
\]
\end{proof}

\subsection{Normalized equivalences of positive-degree spectra}\label{sec:paper-local-spectrum}
In this subsection, we lift the isomorphisms of Theorem~\ref{thm:paper-local-norm} to an equivalence of positive-degree spectra.
We then show that such an equivalence cannot in general be obtained by dividing the norm by $\deg D$, and we record the rank-$1$ lemma and the naturality over finite Galois extensions that are used later.

\begin{theorem}\label{thm:paper-local-spectrum}
Let $F$ be a nonarchimedean local field, and let $D$ be a central division algebra over $F$.
There is an equivalence of spectra
\[
 g\colon\tau_{\ge1}K(D)\xrightarrow{\ \sim\ }\tau_{\ge1}K(F)
\]
inducing all the isomorphisms of Theorem~\ref{thm:paper-local-norm} simultaneously.
\end{theorem}
The proof of Theorem~\ref{thm:paper-local-spectrum} is given after Lemma~\ref{lem:paper-rational-correction} below.

We glue $N_{D/F}/\deg D$ on the rational part to the inverse rank-$1$ equivalences of Lemma~\ref{lem:paper-etale-rank-one} at primes $\ell\ne\operatorname{char}F$.
At the characteristic prime we use $\Nrd_1$ on the $p$-completed factor, which is concentrated in degree one.
We then use the boundary of the arithmetic square to correct the even-degree maps.

The rank-$1$ input is the following lemma, which we state also over global fields, since it is used again in Section~\ref{sec:global-spectrum}.
We record coherence over finite Galois extensions with spectral Mackey functors.
For a finite group $G$, write $\mathscr M_G$ for the category of spectral Mackey functors on the effective Burnside category of finite $G$-sets; the abelian category of ordinary Mackey functors is denoted by $\mathcal M_G$.
For a finite Galois extension $L/F$ with group $G$ and a central simple algebra $B$ over $F$, the functor $G/H\mapsto K(L^H)$, with tensor products and norms along finite \'etale algebras, is a commutative spectral Green functor, and $G/H\mapsto K(B_{L^H})$ is a module over it.
After height-one localization, we denote by $\mathbf R_\ell$ and $\mathbf M_\ell$ the corresponding functors with values in $L_{K(1),\ell}$-local spectra, and we use the same notation $N_{B_{F'}/B}$ for the norm after height-one localization.

\begin{lemma}\label{lem:paper-etale-rank-one}
Let $F$ be a nonarchimedean local field or a global field, let $B$ be a central simple algebra over $F$, and suppose $\ell\ne\operatorname{char}F$.
If $F$ is a global field, we assume that $B$ splits at each real place of $F$.
Let $x_\ell\in\pi_0L_{K(1),\ell}K(B)$ have geometric rank $1$, and for a finite separable extension $F'/F$ put $x_{\ell,F'}=\res_{F'/F}x_\ell$.
\begin{enumerate}
\item Multiplication by $x_\ell$ is an equivalence of $L_{K(1),\ell}K(F)$-modules
\[
 m_{x_\ell}\colon L_{K(1),\ell}K(F)\xrightarrow{\ \sim\ }L_{K(1),\ell}K(B).
\]
\item For a finite Galois extension $L/F$ with group $G$, the maps $m_{x_{\ell,L^H}}$ are the components at the orbits $G/H$ of an equivalence $\mathbf R_\ell\to\mathbf M_\ell$ of $\mathbf R_\ell$-modules in $\mathscr M_G$.
In particular, they respect restriction, conjugation, and the projection formula for transfer.
\end{enumerate}
\end{lemma}
\begin{proof}
(1) Under Morita equivalence over $F^{\mathrm{sep}}$, multiplication by a geometric rank-$1$ class acts as the identity on the \'etale $KU_\ell$-module, and hence it is an equivalence on every geometric stalk with finite coefficients.
\'Etale hyperdescent and derived inverse limits give the equivalence $m_{x_\ell}$ \cite{ClausenMathew}.

(2) Since $\mathbf R_\ell$ is the free $\mathbf R_\ell$-module on one generator, we have
\[
 \operatorname{Map}_{\operatorname{Mod}_{\mathbf R_\ell}}
       (\mathbf R_\ell,\mathbf M_\ell)
 \simeq\Omega^\infty\mathbf M_\ell(G/G).
\]
A representative of $x_\ell$ thus defines a map $m_{x_\ell}\colon\mathbf R_\ell\to\mathbf M_\ell$ of $\mathbf R_\ell$-modules whose component at $G/H$ is multiplication by $x_{\ell,L^H}=\res_{L^H/F}x_\ell$.
In particular, for every intermediate field $F'$,
\[
 N_{B_{F'}/B}\circ m_{x_{\ell,F'}}\simeq m_{x_\ell}\circ N_{F'/F}.
\]
The homotopy is part of the structure of the module map.
Since the map is an equivalence at every orbit, its inverse in the same category retains compatibility with restriction, transfer, conjugation, and their composites even for a virtual rank-$1$ class \cite{BarwickSpectralMackey}, \cite[Sections 3--5, Construction 9.7]{BGS2017}.
\end{proof}

\begin{lemma}\label{lem:paper-rank-zero-product}
Let $F$ be a nonarchimedean local field or a global field, let $\ell\ne\operatorname{char}F$ be a prime, and let $a\in\pi_0L_{K(1),\ell}K(F)$ have geometric rank $0$.
If $\ell=2$ and $F$ has real places, assume in addition that $a$ restricts to zero at the algebraic real closure $F_v^{\mathrm{rc}}$ of every real place $v$.
Then $a\cdot y=0$ for every $y\in\pi_jL_{K(1),\ell}K(F)$ with $j\ge1$.
\end{lemma}
\begin{proof}
We call the case where $\ell=2$ and $F$ has real places the real-relative case; there we write $H^i_{\mathrm{rel}}$ for the cohomology of the fiber of restriction to the algebraic real closures at all real places.
Let $\star$ be empty in the absolute case and $\star=\mathrm{rel}$ in the real-relative case.
In the absolute case, the cohomology with finite $\ell$-primary coefficients vanishes above degree two because $\cd_\ell(F)\le2$.
In the real-relative case, for finite coefficients, restriction to the real places is surjective on $H^2$ and an isomorphism on $H^{\ge3}$ \cite[I, Theorem 4.10(c) and Corollary 4.16]{MilneADT}, so the long exact sequence of the fiber gives $H^i_{\mathrm{rel}}=0$ for $i>2$.

Let $\mathcal R_r$ be the periodic \'etale $K$-theory hypersheaf modulo $\ell^r$.
Since $F(\mu_{\ell^\infty})/F$ is an infinite extension, the $\ell$-adic cyclotomic character of $F$ has infinite image (for $\ell$ different from the residue or constant characteristic this is visible on a Frobenius; for a mixed-characteristic local field with $\ell$ equal to the residue characteristic, $F/\Q_\ell$ is finite and $[\Q_\ell(\mu_{\ell^r}):\Q_\ell]$ is unbounded), so
\[
 H^0(F,\Z_\ell(n))=0\qquad(n>0).
\]
Hence each $y\in\pi_jL_{K(1),\ell}K(F)$ with $j>0$ lifts along
\[
 \pi_j\varprojlim_r R\Gamma(F,\tau_{\ge j+1}\mathcal R_r)
       \longrightarrow \pi_j L_{K(1),\ell}K(F).
\]
Indeed, the global sections of the truncation cofiber $\tau_{\le j}\mathcal R_r$ have zero $\pi_{j+1}$, and their $\pi_j$ is $H^0(F,\Z/\ell^r(j/2))$ for $j$ even and zero for $j$ odd; the inverse limit of these groups over $r$ is zero.
Likewise, $a$ lifts to sections of $\tau_{\ge2}\mathcal R_r$, by the description of the geometric stalks and the vanishing of its geometric rank.
In the real-relative case, $a$ lifts to relative sections, and the relative sections of $\tau_{\le1}\mathcal R_r$ have zero $\pi_0$, since $H^0_{\mathrm{rel}}(F,\Z/\ell^r)=0$ and there are no odd geometric layers.

To multiply the two lifts, we use the direct limit of the duals of the Moore tower,
\[
 \varinjlim_r (\mathbb S/\ell^r)^\vee
        \simeq\fib(\mathbb S\to\mathbb S[1/\ell]).
\]
It has a coalgebra structure from the counit of the smashing colocalization onto $\ell$-primary torsion.
By compactness of finite Moore spectra, on each finite stage the comultiplication factors through a tensor product of two sufficiently late stages.
Dualizing gives products on the inverse system, and additivity of connectivity sends the product of the $\tau_{\ge2}$ and $\tau_{\ge j+1}$ stages to $\tau_{\ge j+3}$.
This multiplication is defined on the tower, without a ring structure on each finite Moore spectrum.
In the real-relative case, the fiber of real restriction is a module over absolute sections, so the same product uses a relative and an absolute stage.

For the $\star$-sections of this stage, the geometric degrees satisfy $t\ge j+3$ and the cohomological degrees satisfy $s\le2$, so the total degrees satisfy $t-s\ge j+1$.
Hence $\pi_j$ modulo $\ell^r$ vanishes, and the only possible term in $\pi_{j+1}$, occurring when $t-s=j+1$, is $H^2_\star(F,\Z/\ell^r(n))$ with $2n=j+3$.
The transition maps of these groups are surjective because $H^3_\star=0$, so $\varprojlim^1$ vanishes in the Milnor sequence, and $a\cdot y=0$ in the completed group.
\end{proof}

The correction step uses the following general lemma, which is also used in Proposition~\ref{prop:paper-finite-galois-naturality} and in Section~\ref{subsec:paper-lca-comparison}.

\begin{lemma}\label{lem:paper-rational-correction}
Let $Y\simeq Y_1\times_WY_2$ be a pullback of spectra in which $W$ is rational, and let $\partial\colon\Omega W\to Y$ be its boundary.
For $j\ge1$ put $V_j=\pi_{j+1}W$ and $P_j=\im\bigl(\pi_{j+1}Y_1\oplus\pi_{j+1}Y_2\to V_j\bigr)$, so that $\partial$ induces an injection $V_j/P_j\hookrightarrow\pi_jY$.
Let $X$ be a spectrum with rationalization $q_X\colon X\to X_\Q$, let $g^0\colon X\to Y$, and let $\nu_j\colon\pi_jX\to\pi_jY$ $(j\ge1)$ be homomorphisms such that $\nu_j=\pi_jg^0$ for odd $j$ and, for even $j$, the difference $\nu_j-\pi_jg^0$ factors as
\[
 \pi_jX\longrightarrow\pi_jX\otimes\Q\longrightarrow V_j/P_j\longrightarrow\pi_jY.
\]
\begin{enumerate}
\item If $\operatorname{Ext}^1_\Z(\Q,P_j)=0$ for every even $j\ge2$, then there is $b\colon X_\Q\to\Omega W$ such that $g=g^0+\partial\,b\,q_X$ satisfies $\pi_jg=\nu_j$ for all $j\ge1$, and $g$ and $g^0$ have the same composites with $Y\to Y_1$ and $Y\to Y_2$.
\item The same holds for spectral Mackey functors for a finite group, with homotopy groups replaced by homotopy Mackey functors, provided the factorizations are maps of Mackey functors and $\operatorname{Ext}^1_\Z(\Q,P_j)=0$ at every orbit.
\item Let $M$ be a finitely generated $\Z_p$-module, put $V=M\otimes_{\Z_p}\Q_p$, let $\Lambda$ be the image of $M\to V$, and let $\Gamma\subseteq V$ be a $p$-divisible subgroup.
Then $\operatorname{Ext}^1_\Z(\Q,\Gamma+\Lambda)=0$.
\end{enumerate}
\end{lemma}
\begin{proof}
(1) Put $U_j=\pi_jX\otimes\Q$.
Since $\operatorname{Ext}^1_\Z(-,P_j)$ turns direct sums into products, $\operatorname{Ext}^1_\Z(U_j,P_j)=0$, so $U_j\to V_j/P_j$ lifts to a homomorphism $U_j\to V_j$, which is $\Q$-linear because both groups are rational.
Rational spectra are determined by their homotopy groups, so these lifts in even degrees, together with zero in all other degrees, are realized by a map $b\colon X_\Q\to\Omega W$.
Then $\pi_j(\partial\,b\,q_X)=\nu_j-\pi_jg^0$ for every $j\ge1$.
The boundary $\partial$ becomes null after composition with $Y\to Y_1$ and $Y\to Y_2$, so the correction does not change these composites.

(2) Let $\mathcal M_{G,\Q}$ be the rational counterpart of $\mathcal M_G$.
The exact forgetful functor $\mathcal M_{G,\Q}\to\mathcal M_G$ has a right adjoint, given orbitwise by $\operatorname{Hom}_\Z(\Q,-)$.
Evaluation at an orbit preserves injectives because its left adjoint tensors with the representable Mackey functor, whose Burnside groups are finite free.
Hence the derived right adjoint is computed orbitwise by $R\operatorname{Hom}_\Z(\Q,-)$.
Derived adjunction and semisimplicity of $\mathcal M_{G,\Q}$ give, for rational $U$,
\[
 \operatorname{Ext}^1_{\mathcal M_G}(U,P_j)
 \simeq \operatorname{Hom}_{\mathcal M_{G,\Q}}
 \bigl(U,\operatorname{Ext}^1_\Z(\Q,P_j)\bigr)=0,
\]
where $\operatorname{Ext}^1_\Z$ on the right is taken orbitwise.
Thus $U_j\to V_j/P_j$ lifts to a map of Mackey functors $U_j\to V_j$.
Semisimplicity in the rational spectral category realizes all these lifts, together with zero in the other degrees, by a single map $b$ \cite{BarnesRational2009}, and the rest of the argument is as in~(1).

(3) The group $\Lambda$ is a free $\Z_p$-module of finite rank, hence derived $p$-complete, and $p$ is invertible on $\Q$, so $\operatorname{Ext}^1_\Z(\Q,\Lambda)=0$ \cite[Tag 091N]{StacksProject}.
The quotient $(\Gamma+\Lambda)/\Lambda$ is a $p$-divisible subgroup of the $p$-primary torsion group $V/\Lambda$; since primes other than $p$ act invertibly on it, it is divisible and hence injective.
Since $\operatorname{Ext}^2_\Z=0$, the long exact sequence for $0\to\Lambda\to\Gamma+\Lambda\to(\Gamma+\Lambda)/\Lambda\to0$ gives the assertion.
\end{proof}

\begin{proof}[Proof of Theorem~\ref{thm:paper-local-spectrum}]
We first construct the completed comparisons.
Set
\[
 X=\tau_{\ge1}K(D),\qquad Y=\tau_{\ge1}K(F).
\]
For $\ell\ne\operatorname{char}F$, set
\[
 R_F=L_{K(1),\ell}K(F),\qquad R_D=L_{K(1),\ell}K(D).
\]
Descent modulo $\ell^a$ in degree zero gives
\[
 0\longrightarrow H^2(F,\Z/\ell^a(1))\longrightarrow
 \pi_0(R_D/\ell^a)\xrightarrow{\mathrm{rk}_{\mathrm{geom}}}\Z/\ell^a\longrightarrow0.
\]
Since coefficient reduction is surjective on $H^2$, we can choose compatible finite-coefficient classes of rank $1$, and the surjection $\pi_0R_D\to\varprojlim_a\pi_0(R_D/\ell^a)$ of the Milnor sequence lifts them to a class $x\in\pi_0R_D$.
By Lemma~\ref{lem:paper-etale-rank-one}, multiplication by $x$ is an equivalence
\[
 m_x\colon R_F\xrightarrow{\ \sim\ }R_D.
\]
The projection formula gives $N_{D/F}\circ m_x\simeq m_{N_{D/F}(x)}$, where $m_{N_{D/F}(x)}$ is multiplication by $N_{D/F}(x)\in\pi_0R_F$.
The element $N_{D/F}(x)-\deg D$ has geometric rank $0$, so Lemma~\ref{lem:paper-rank-zero-product} gives
\[
 \bigl(N_{D/F}(x)-\deg D\bigr)\cdot y=0\qquad\text{for }y\in\pi_iR_F,\ i\ge1.
\]
Let $u_\ell\colon X^\wedge_\ell\to Y^\wedge_\ell$ be the completed equivalence obtained from $m_x^{-1}$ and Lemma~\ref{lem:paper-local-compare}.
Then
\[
 \pi_i\bigl((N_{D/F})^\wedge_\ell\bigr)=\deg D\cdot\pi_iu_\ell\qquad(i\ge1).
\]
For $\ell=\operatorname{char}F$, the completed positive-degree spectrum is concentrated in degree one by Proposition~\ref{prop:app-frobenius}(2), and there we take the reduced norm $\Nrd_1$.

Next we glue the completed comparisons.
For the extension of scalars $i_D\colon K(F)\to K(D)$, we have $N_{D/F}\circ i_D=(\deg D)^2$ because $D$ has $F$-dimension $(\deg D)^2$.
The composite $i_D\circ N_{D/F}$ is tensoring with the bimodule $D\otimes_FD$, which is a sum of $(\deg D)^2$ copies of $D$ since $D\otimes_FD^{\mathrm{op}}\simeq\operatorname{End}_F(D)$, so additivity gives $i_D\circ N_{D/F}=(\deg D)^2$ as well.
Hence we can take $N_{D/F}/\deg D$ after inverting $\deg D$.
For $S=X,Y$, put $S^\wedge=\prod_\ell S^\wedge_\ell$ and $Z_S=(S^\wedge)_{\Q}$.
The arithmetic square for $S$ gives a fiber sequence
\[
 \Omega Z_S\xrightarrow{\ \partial_S\ }S\longrightarrow S_{\Q}\oplus S^\wedge
 \longrightarrow Z_S.
\]
That this square is a pullback can be checked by rationalization and comparison modulo every prime.
On homotopy groups, $N_{D/F}/\deg D$ agrees with the rationalization of the product of the completed equivalences.
Two maps between rational spectra are homotopic if they have the same effect on homotopy groups, so gluing gives an equivalence $g^0\colon X\to Y$.

In odd degrees $j=2m-1$, we have $\pi_jg^0=\nu_{j,D}$.
Indeed, agreement on the completed factors follows from injectivity of restriction of $H^1(F,\Z_\ell(m))$ to a splitting field, by inflation--restriction and $H^0(E,\Z_\ell(m))=0$.
The boundary of the arithmetic square introduces no ambiguity because $\pi_{j+1}Z_Y=0$: the positive even completed groups are finite, and only finitely many primes contribute in any fixed degree.

For even $j=2m\ge2$, we first compare on integral torsion.
By the second assertion of Lemma~\ref{lem:paper-local-even-coinvariants}, $\theta_E$ from a Galois splitting field $E/F$ is surjective on torsion.
With $x_E=\res_{E/F}x$, the projection formula gives
\[
 m_x^{-1}\theta_E=N_{E/F}\,m_{x_E}^{-1}
\]
at height one.
Descent modulo $\ell^a$ in degree $j+1$ has the single term $H^1(E,\Z/\ell^a(m+1))$, on which multiplication by the rank-$1$ class $x_E$ acts as the identity.
By the naturality of the finite-coefficient comparison and of the integral Bockstein sequence, $\pi_jg^0\circ\theta_E=N_{E/F}$ on $K_j(E)[\ell^a]$.
Torsion at the characteristic prime vanishes by Proposition~\ref{prop:app-frobenius}(2).
Thus
\[
 \pi_jg^0=\nu_{j,D}\quad\text{on }K_j(D)_{\mathrm{tors}}.
\]
The difference $\nu_{j,D}-\pi_jg^0$ factors through the rational vector space given by the torsion-free quotient of $K_j(D)$.
Since each even completed group is finite, the difference has zero image in the completed factors.
The rational components agree by construction.

We correct this difference by Lemma~\ref{lem:paper-rational-correction}, applied to the arithmetic square $Y\simeq Y_\Q\times_{Z_Y}Y^\wedge$.
Write
\[
 V_j=\pi_{j+1}Z_Y,\qquad
 \Gamma_j=\im\bigl(\pi_{j+1}Y_{\Q}\to V_j\bigr),\qquad
 \Lambda_j=\im\bigl(\pi_{j+1}Y^\wedge\to V_j\bigr).
\]
The difference gives a map $K_j(D)\otimes\Q\to V_j/(\Gamma_j+\Lambda_j)$.
Here $\Gamma_j$ is a rational subspace.
In equal characteristic $V_j=0$, so $\Gamma_j+\Lambda_j=0$; in mixed characteristic $V_j$ and $\Lambda_j$ arise from the finitely generated $\Z_p$-module $\pi_{j+1}Y^\wedge_p$ as in Lemma~\ref{lem:paper-rational-correction}(3).
In both cases $\operatorname{Ext}^1_{\Z}(\Q,\Gamma_j+\Lambda_j)=0$.
Lemma~\ref{lem:paper-rational-correction} therefore gives $b\colon X_{\Q}\to\Omega Z_Y$, and with the rationalization $q_X\colon X\to X_{\Q}$ we put
\[
 g=g^0+\partial_Y\,b\,q_X.
\]
\end{proof}

In general, the norm $N_{D/F}$ cannot be divided by $\deg D$ on the level of spectra; this is why the proof above glues completed comparisons.

\begin{proposition}\label{prop:paper-local-no-division}
Let $F$ and $D$ be as in Theorem~\ref{thm:paper-local-spectrum}.
\begin{enumerate}
\item If some prime $\ell\mid\deg D$ satisfies $\ell\ne\operatorname{char}F$, then for every $J\ge1$ the map $\tau_{\ge J}N_{D/F}$ is not $\deg D$ times any map of spectra $\tau_{\ge J}K(D)\to\tau_{\ge J}K(F)$.
\item In equal characteristic $p$ with $\deg D$ a power of $p$, the map $\tau_{\ge1}N_{D/F}$ is $\deg D$ times an equivalence $\tau_{\ge1}K(D)\to\tau_{\ge1}K(F)$.
\end{enumerate}
\end{proposition}
\begin{proof}
For~(1), let $\ell\mid\deg D$ with $\ell\ne\operatorname{char}F$, put
\[
 R_F=L_{K(1),\ell}K(F),\qquad R_D=L_{K(1),\ell}K(D),
\]
and choose a class $x\in\pi_0R_D$ of geometric rank $1$ as in the proof of Theorem~\ref{thm:paper-local-spectrum}.
By Lemma~\ref{lem:paper-etale-rank-one}, multiplication by $x$ is an equivalence $m_x\colon R_F\to R_D$, and by the projection formula $N_{D/F}\circ m_x\simeq m_{N_{D/F}(x)}$.
In integral \'etale descent the boundary of geometric rank $1$ is the scalar $2$-cocycle $[D]$ obtained by gluing the split Morita modules (compare the proof of \cite[Proposition 6.9.1]{KahnLevine2010}); since $[D]$ has order $\deg D$ in $\Br(F)$ by the local invariant, we obtain
\[
 \pi_0K^{\et}(D)=\deg D\cdot\Z,\qquad \pi_{-1}K^{\et}(D)=\Br(F)/\langle[D]\rangle.
\]
The map $\Br(F)/\langle[D]\rangle\to\Br(F)$ induced by $N_{D/F}$ is multiplication by $\deg D$, and it is an isomorphism by the local invariant.
The localization unit induces an isomorphism $\pi_0(K^{\et}(D)/\ell^a)\simeq\pi_0(R_D/\ell^a)$: in the finite-coefficient descent spectral sequences only the stalk degrees $0$ and $2$ contribute in total degree zero, the localization unit is an isomorphism on both, and since $\cd_\ell F\le2$ and the odd stalk degrees vanish, no differentials enter or leave.
Comparison modulo $\ell^a$ in degrees zero and one therefore gives the exact sequence
\[
 0\longrightarrow \deg D\cdot\Z_\ell\longrightarrow\pi_0R_D
 \longrightarrow T_\ell\bigl(\Br(F)/\langle[D]\rangle\bigr)\longrightarrow0.
\]
If $\ell\mid\deg D$, the image of the rank-$1$ class $x$ in the right-hand group is not divisible by $\ell$, since a lift of that divisibility would make the rank of $x$ a multiple of $\ell$.
The same sequence for $F$ reads $0\to\Z_\ell\to\pi_0R_F\to T_\ell\Br(F)\to0$, whose kernel of geometric rank is $H^2(F,\Z_\ell(1))$, and $N_{D/F}$ maps the two sequences to each other, inducing on the right $T_\ell$ of the isomorphism $\Br(F)/\langle[D]\rangle\xrightarrow{\deg D}\Br(F)$.
Thus
\[
 \beta=N_{D/F}(x)-\deg D\in H^2(F,\Z_\ell(1))
\]
is not divisible by $\ell$.

We take $n=\ell$ for odd $\ell$ and $n=4$ for $\ell=2$, so $n\not\equiv2\bmod4$ and the identity of the Moore spectrum $\mathbb S/n$ has additive order $n$ \cite[p.~257]{Oka1984}.
Choose a sufficiently large multiple $m$ of $\varphi(n)$ so that $\Z/n(m)$ has trivial action.
In descent modulo $n$ in degree $2m$, the outer terms are
\[
 H^2(F,\Z/n(m+1))\qquad\text{and}\qquad H^0(F,\Z/n(m)),
\]
both cyclic of order $n$.
Since the middle term is killed by $n$, the sequence splits as a sequence of abelian groups, so an element of the left-hand term that is not divisible by $\ell$ there is not divisible by $\ell$ in the middle term.
If $y$ lifts a generator of the right-hand term, then $\beta y$ is not divisible by $\ell$ in the left-hand term, since it is the cup product of the coefficient reduction of $\beta$ under the local invariant with that generator.
Thus
\[
 N_{D/F}\,m_x(y)=\deg D\cdot y+\beta y
\]
is not divisible by $\ell$ in the middle term.
By Lemma~\ref{lem:paper-local-compare}, for $2m-1\ge J$ the groups $\pi_{2m}(\tau_{\ge J}K(D)/n)$ and $\pi_{2m}(R_D/n)$ agree, and similarly for $F$; if $\deg D\cdot u=\tau_{\ge J}N_{D/F}$ for some map $u$, then $N_{D/F}(m_xy)$ would be divisible by $\ell$, contradicting the displayed formula.

For~(2), let $p=\operatorname{char}F$, assume that $\deg D$ is a power of $p$, and put $X=\tau_{\ge1}K(D)$.
The completed positive-degree spectrum at $p$ is concentrated in degree one.
We use the reduced norm $\Nrd_1$ at $p$ and $N_{D/F}/\deg D$ after inverting $p$.
Their overlap is also concentrated in degree one, and the two maps agree there because $\pi_1N_{D/F}=\deg D\cdot\Nrd_1$ on $K_1(D)$.
Every map from $X$ to the loop spectrum of the overlap vanishes, so the two pieces glue to a map $u$, and the identity $\deg D\cdot u=\tau_{\ge1}N_{D/F}$, which holds on the two pieces, holds in the integral homotopy category as well.
The map $u$ is an equivalence: on the $p$-completed part it is the isomorphism $\Nrd_1$ in degree one, and after inverting $p$ the map $N_{D/F}/\deg D$ has inverse $i_D/\deg D$, since $N_{D/F}i_D$ and $i_DN_{D/F}$ are multiplication by $(\deg D)^2$.
\end{proof}

The comparisons of Theorem~\ref{thm:paper-local-spectrum} can moreover be chosen coherently over the intermediate fields of a finite Galois extension.

\begin{proposition}\label{prop:paper-finite-galois-naturality}
Let $F$ be a nonarchimedean local field, let $D$ be a central division algebra over $F$, and let $L/F$ be a finite Galois extension.
The comparisons
\[
 g_{F'}\colon\tau_{\ge1}K(D_{F'})\longrightarrow\tau_{\ge1}K(F')
\]
for the intermediate fields $F'$ of $L/F$ can be chosen simultaneously and coherently with respect to restriction, transfer, conjugation, and their composites.
Each $g_{F'}$ is an equivalence inducing the isomorphisms of Theorem~\ref{thm:paper-local-norm}.
\end{proposition}
\begin{proof}
Let $G=\Gal(L/F)$.
We apply positive connective covers to the Mackey functors introduced before Lemma~\ref{lem:paper-etale-rank-one} and obtain objects of $\mathscr M_G$,
\[
 X(G/H)=\tau_{\ge1}K(D_{L^H}),\qquad
 Y(G/H)=\tau_{\ge1}K(L^H).
\]
Connective covers preserve finite direct sums, so these remain spectral Mackey functors.
We use the norms $N_{D_{F'}/F'}$ of the algebras $D_{F'}$, which have degree $\deg D$ even when the index drops over $F'$.

For $S\in\mathscr M_G$, put
\[
 S^\wedge_{\deg D}=\prod_{\ell\mid\deg D}S^\wedge_\ell.
\]
Rationalization, derived completion, products, and pullbacks are computed orbitwise.
The arithmetic square for the primes dividing $\deg D$ gives
\[
 S\simeq S[1/\deg D]\mathop{\times}_{S^\wedge_{\deg D}[1/\deg D]}S^\wedge_{\deg D}.
\]
This can be checked orbitwise by rationalization and reduction modulo every prime.
Primes other than $\ell$ are invertible in an $\ell$-complete object, so $S^\wedge_{\deg D}[1/\deg D]$ is rational.

For each $\ell\mid\deg D$ distinct from $\operatorname{char}F$, choose over $F$ a class $x_\ell$ of geometric rank $1$ as in the proof of Theorem~\ref{thm:paper-local-spectrum}, and restrict it to every intermediate field.
From the inverse $m_{x_\ell}^{-1}$ of Lemma~\ref{lem:paper-etale-rank-one} and the natural \'etale comparison, we obtain a map
\[
 u_\ell\colon X^\wedge_\ell\longrightarrow Y^\wedge_\ell\quad\text{in }\mathscr M_G.
\]
In characteristic $p$ the $p$-completions of $X$ and $Y$ are concentrated in degree one, so we use the map of Mackey functors given by the reduced norm $\Nrd_1$.
After inverting $\deg D$, we take
\[
 u_{\deg D}=N/\deg D\colon X[1/\deg D]\longrightarrow Y[1/\deg D],
\]
where $N$ is the map of Mackey functors given by the norms $N_{D_{F'}/F'}$ with their restriction and transfer structure.

We glue these maps over the rational overlap of the arithmetic square in $\mathscr M_G$.
The product calculation in the proof of Theorem~\ref{thm:paper-local-spectrum} identifies the rationalized $u_\ell$ with $N/\deg D$ on positive homotopy groups over every intermediate field.
Spectral Mackey functors are equivalent to genuine $G$-spectra \cite{BarwickSpectralMackey}, and rational spectral Mackey functors are described by the product of categories of complexes over the rational Weyl group algebras $\Q[W_GH]$, indexed by conjugacy classes of subgroups $H\le G$ \cite{BarnesRational2009}.
These algebras are semisimple, so rational maps inducing the same map of graded Mackey homotopy groups are homotopic.
Choosing a homotopy on the overlap, we obtain $g^0\colon X\to Y$ from the pullback; it is an equivalence because both arithmetic components are.

We correct $g^0$ to obtain the prescribed maps $\nu_{j,D_{F'}}$, denoted collectively by $\nu_j$, while retaining the arithmetic components.
The projection-formula calculation in the proof of Theorem~\ref{thm:paper-local-spectrum} gives agreement in odd degrees and on torsion, and it applies to the set of primes dividing $\deg D$ as well.
The components after inverting $\deg D$ agree since $\deg D\cdot\nu_j=\pi_jN$.
Thus the difference $\nu_j-\pi_jg^0$ is confined to even degrees $j\ge2$ and maps to zero in both arithmetic components.

In the abelian category of ordinary Mackey functors, images, kernels, and quotients are computed orbitwise.
Define
\[
 \begin{aligned}
 V_j&=\pi_{j+1}\bigl(Y^\wedge_{\deg D}[1/\deg D]\bigr),\\
 \Gamma_j&=\im\bigl(\pi_{j+1}Y[1/\deg D]\to V_j\bigr),\\
 \Lambda_j&=\im\bigl(\pi_{j+1}Y^\wedge_{\deg D}\to V_j\bigr),
 \qquad P_j=\Gamma_j+\Lambda_j.
 \end{aligned}
\]
The arithmetic fiber sequence yields the natural exact sequence
\[
 0\longrightarrow V_j/P_j
  \xrightarrow{\ \bar\partial\ }\pi_jY
  \longrightarrow\pi_jY[1/\deg D]\oplus\pi_jY^\wedge_{\deg D}.
\]
Consequently $\nu_j-\pi_jg^0$ defines a unique map of Mackey functors into $V_j/P_j$.
It factors through the rational Mackey functor $U_j=\pi_jX\otimes\Q$, since it vanishes on torsion and, for a local field, the torsion-free quotient of a $K$-group in positive even degree is a rational vector space.

We now compute the obstruction to a lift at each intermediate field.
Let $p$ be the common residue characteristic.
In equal characteristic, or in mixed characteristic with $p\nmid\deg D$, the completed groups in degree $j+1\ge3$ give $V_j=0$.
In mixed characteristic with $p\mid\deg D$, let
\[
 M_j(F')=\pi_{j+1}\bigl(Y(F')^\wedge_p\bigr).
\]
This is a finitely generated $\Z_p$-module, and the corresponding completions for $\ell\ne p$ are finite.
Thus
\[
 V_j(F')=M_j(F')\otimes_{\Z_p}\Q_p,\qquad
 \Lambda_j(F')=M_j(F')/M_j(F')_{\mathrm{tors}}.
\]
Since $\Gamma_j(F')$, the image of a $\Z[1/\deg D]$-module, is $p$-divisible, Lemma~\ref{lem:paper-rational-correction}(3) gives $\operatorname{Ext}^1_\Z(\Q,P_j(F'))=0$.

Lemma~\ref{lem:paper-rational-correction}(2), applied to the arithmetic square in $\mathscr M_G$, now gives a map
\[
 b\colon X_\Q\longrightarrow\Omega\bigl(Y^\wedge_{\deg D}[1/\deg D]\bigr)
\]
such that, with the arithmetic boundary $\partial$ and the rationalization $q_X\colon X\to X_\Q$, the map
\[
 g=g^0+\partial\,b\,q_X
\]
satisfies $\pi_jg=\nu_j$ at every orbit.
By Lemma~\ref{lem:paper-rational-correction}(2), the correction preserves both arithmetic components, and it is coherent because it is constructed in $\mathscr M_G$.
\end{proof}

\begin{remark}\label{rem:paper-galois-splitting-transfer}
Let $E$ be an intermediate field of the Galois extension $L/F$ of Proposition~\ref{prop:paper-finite-galois-naturality} (and its global analogue in Section~\ref{sec:global-spectrum}) that splits $D$, and let $\ell$ be a prime at which a class of geometric rank $1$ was chosen, so that $\ell\mid\deg D$, $\ell\ne\operatorname{char}F$, and $\ell\ne2$ under real ramification.
At the $\ell$-completed factor, write $x_{\ell,E}=1+\alpha_{\ell,E}$ under the standard Morita equivalence, and write $g_{F,\ell}^{\wedge}$ for the $\ell$-adic completion of the component $g_F$ over $F$ of the map $g$ of that proposition.
On the $\ell$-completed positive-degree parts, the splitting transfer then satisfies
\[
 g_{F,\ell}^{\wedge}\circ\theta_E\simeq N_{E/F}\circ m_{x_{\ell,E}}^{-1}.
\]
The element $\alpha_{\ell,E}$ has geometric rank $0$; if $\ell=2$ and $E$ has real places, it also restricts to zero at the real closures, since there $\pi_0(KO^\wedge_2)=\Z_2$ is detected by rank.
Hence by Lemma~\ref{lem:paper-rank-zero-product} it acts trivially on the positive-degree homotopy groups of the $\ell$-completed factor, which gives the transfer identity on groups.
With finite coefficients, the action of $\alpha_{\ell,E}$ must be retained in the displayed formula.
\end{remark}

These choices are made separately for each finite Galois extension.

\begin{problem}\label{prob:paper-coherent-choices}
Can the choices of Proposition~\ref{prop:paper-finite-galois-naturality} (and its global analogue in Section~\ref{sec:global-spectrum}) be made coherently over all finite extensions at once, and compatibly with tensor products of central simple algebras?
\end{problem}

\section{Splitting transfers and local--global principles}\label{sec:splitting-transfers}
This section studies the subgroup of $K_2(D)$ generated by transfers from splitting fields.
We introduce the corresponding quotient $Q_F(D)$ and show that, over fields complete for a nontrivial rank-$1$ nonarchimedean valuation, evaluations of $K_2$-classes are locally constant modulo the transfers from a fixed finite separable splitting field (Section~\ref{sec:paper-splitting}).
Over global fields, we prove that $Q_F(D)$ is detected at the places where $D$ does not split (Section~\ref{sec:paper-Q-local-global}).
We then deduce that $\SK_2(D)=0$ over global fields and determine the image of the reduced norm (Section~\ref{sec:paper-global-weight-one}).
\subsection{Splitting transfers and rigidified Picard norms}\label{sec:paper-splitting}
In this subsection, we introduce the splitting transfer subgroup $I_F(D)$ and the quotient $Q_F(D)=K_2(D)/I_F(D)$, which measures how far $K_2(D)$ is from being generated by transfers from splitting fields.
The main result is Proposition~\ref{prop:paper-picard-eval}: over a field complete for a nontrivial rank-$1$ nonarchimedean valuation, the evaluations at rational points of a $K_2$-class on a smooth affine scheme are locally constant modulo the transfers from a fixed finite separable splitting field.
We prove it using rigidified Picard norms on curves.
It is used, through Corollary~\ref{cor:paper-congruence-splitting}, in Appendix~\ref{app:second-tc}, and Remark~\ref{rem:paper-local-completion} records a further application; Section~\ref{sec:paper-Q-local-global} uses the reciprocity on curves and a global analogue of the Picard argument.
Throughout this subsection, $F$ is a field and $D$ is a central division algebra over $F$.

\begin{definition}\label{def:paper-transfer}
For a central division algebra $D$ over $F$, we define the splitting transfer subgroup $I_F(D)$, the splitting transfer quotient $Q_F(D)$, and the norm subgroup $N_2([D])$ by
\[
 \begin{split}
 I_F(D)&=\sum_{D_E\,\mathrm{split}}\theta_E K_2(E),\\
 Q_F(D)&=K_2(D)/I_F(D),\\
 N_2([D])&=\sum_{D_E\,\mathrm{split}}N_{E/F}K_2(E)\subseteq K_2(F).
 \end{split}
\]
For any $\alpha\in\Br(F)$, set $N_2(\alpha)=\sum_{\alpha_E=0}N_{E/F}K_2(E)$.
In all these sums, $E/F$ ranges over finite extensions, including inseparable ones.
\end{definition}

For a field extension $L/F$ we write $\res_{L/F}\colon K(D)\to K(D_L)$ for extension of scalars and, if $L/F$ is finite, $\cor_{L/F}=N_{D_L/D}\colon K(D_L)\to K(D)$.
If $E/L$ is a finite splitting field of $D_L$, we write $\theta_{E/L}\colon K(E)\to K(D_L)$ for the splitting transfer of $D_L$, normalized as in Section~\ref{sec:paper-local-norm}; thus $\theta_E=\theta_{E/F}$ and $\cor_{L/F}\circ\theta_{E/L}=\theta_{E/F}$.
With these conventions, $I_F(A)$ and $Q_F(A)$ are defined for every central simple algebra $A$ over $F$.

The transfer formula $\Nrd_2\theta_E=N_{E/F}$ gives $N_2([D])=\Nrd_2(I_F(D))$.
With $C_2(D/F)=\coker(\Nrd_2)$, the snake lemma for the reduced norm and its restriction $I_F(D)\to N_2([D])$ gives the exact sequence
\[
 0\longrightarrow\SK_2(D)\cap I_F(D)\longrightarrow\SK_2(D)
 \longrightarrow Q_F(D)\longrightarrow K_2(F)/N_2([D])
 \longrightarrow C_2(D/F)\longrightarrow0.
\]

\begin{lemma}\label{lem:paper-Q-functorial}
Let $F'/F$ be a field extension.
Then restriction maps $I_F(D)$ into $I_{F'}(D_{F'})$, where the splitting transfers for the central simple algebra $D_{F'}$ are normalized by Morita equivalence as in Section~\ref{sec:paper-local-norm}; hence it induces a map $Q_F(D)\to Q_{F'}(D_{F'})$.
Moreover,
\[
 \deg D\cdot Q_F(D)=0.
\]
\end{lemma}
\begin{proof}
Let $E/F$ be a finite splitting field, and write
\[
 E\otimes_F F'=\prod_i R_i,
\]
where $\kappa_i$ and $\mathfrak m_i$ are the residue field and the radical of the Artinian local ring $R_i$.
Filtering restriction of scalars after base change by the powers $\mathfrak m_i^j$ on each factor, we obtain as $j$th graded piece the composite of restriction to $\kappa_i$ and splitting transfer from $\kappa_i$ to $F'$, with multiplicity $\dim_{\kappa_i}(\mathfrak m_i^j/\mathfrak m_i^{j+1})$.
The simple module of the split algebra remains simple after base change to $\kappa_i$, so the graded pieces respect the Morita normalization.
Additivity \cite{QuillenHigher} therefore gives
\[
 \res_{F'/F}\theta_E
  =\sum_i \operatorname{length}_{R_i}(R_i)\,
          \theta_{\kappa_i/F'}\,\res_{\kappa_i/E},
\]
which proves the first assertion; the lengths may exceed one when $E/F$ is inseparable.

For the second assertion, let $L\subset D$ be a separable maximal subfield, so that $[L:F]=\deg D$ and $L$ splits $D$.
For $x\in K_2(D)$ we have
\[
 \deg D\cdot x=N_{D_L/D}\bigl(\res_{L/F}x\bigr)\in\theta_LK_2(L)\subseteq I_F(D),
\]
since $K(D_L)$ is identified with $K(L)$ and $\theta_L=N_{D_L/D}$ under this identification \cite{KahnLevine2010}.
\end{proof}

The lemma on transfers of mixed products is used in part~(1) of the reciprocity law on curves.
Part~(2), for functions with residue one along the boundary, is used in Proposition~\ref{prop:paper-picard-eval}.

\begin{lemma}\label{lem:paper-mixed-products}
If $F$ is infinite, then for every finite extension $L/F$, the transfer to $F$ maps the mixed products $L^\times\cdot K_1(D_L)\subset K_2(D_L)$ into $I_F(D)$.
\end{lemma}

\begin{proof}
Let $D'$ be the division algebra in the Brauer class of $D_L$.
By Morita equivalence, we represent an element of $K_1(D_L)$ by $b\in D'^\times$.
Since $L$ is infinite, we can choose $c$ such that both $c$ and $cb$ lie in the nonempty Zariski open set of units generating separable maximal subfields.
Each factor of $b=c^{-1}\cdot(cb)$ then has the form $\theta_{E/L}(u)$ in $K_1(D_L)$ for a separable maximal subfield $E\subset D'$ and some $u\in E^\times$.
For $a\in L^\times$, the projection formula and transitivity of transfer give
\[
 \cor_{L/F}\bigl(\{a\}\cdot\theta_{E/L}(u)\bigr)
 =\cor_{L/F}\theta_{E/L}(\{a|_E,u\})=\theta_E(\{a|_E,u\}),
\]
which lies in $I_F(D)$ since $E$ splits $D$.
By bilinearity, the transfer of $\{a\}\cdot b$ to $F$ lies in $I_F(D)$.
\end{proof}

\begin{lemma}\label{lem:paper-curve-reciprocity}
Let $Y/F$ be a proper geometrically integral curve, let $V\subset Y$ be a smooth open subset, and let $\alpha\in K_2(V;D)$, where $K(V;D)$ is the $K$-theory of the Azumaya algebra $\mathcal O_V\otimes_FD$ over $V$.
Let $f$ be a rational function on the normalization of $Y$ that is a unit at every branch of the boundary $Y\setminus V$, and put
\[
 \Sigma(f)=\sum_{z\in V\ \mathrm{closed}}v_z(f)\,\cor_{\kappa(z)/F}\bigl(\alpha(z)\bigr)\in K_2(D).
\]
\begin{enumerate}
\item If $F$ is infinite, then $\Sigma(f)\in I_F(D)$.
\item If $f$ has residue $1$ at every branch of the boundary, then $\Sigma(f)=0$.
\item Suppose that $F$ is infinite, and let $P,P'\in V(F)$ satisfy $P-P'=\operatorname{div}(f)+Z$ on the normalization, where $f$ is as above and $Z$ is supported on closed points of $V$ whose residue fields split $D$.
Then $[\alpha(P)]=[\alpha(P')]$ in $Q_F(D)$.
\end{enumerate}
\end{lemma}
\begin{proof}
On the normalization of $Y$, which is unchanged over $V$, we form $f\cdot\alpha$ at the generic point.
The composite of the localization boundary and proper pushforward is zero \cite{QuillenHigher}.
At a closed point $z$ of $V$, the contribution is $v_z(f)\cor_{\kappa(z)/F}(\alpha(z))$.
At the boundary, the contributions are mixed products of the residue $\overline f$ with the $K_1$ boundary of $\alpha$; they lie in $I_F(D)$ by Lemma~\ref{lem:paper-mixed-products} if $F$ is infinite, which gives~(1), and they vanish if $\overline f=1$ at every branch, which gives~(2).
Assertion~(3) follows from~(1), since the contributions of split closed points lie in $I_F(D)$.
\end{proof}

\begin{proposition}\label{prop:paper-picard-eval}
Let $F$ be complete for a nontrivial rank-$1$ nonarchimedean valuation, and let $E/F$ be a finite separable splitting field of $D$.
For a smooth affine $F$-scheme $X$ and $\alpha\in K_2(X;D)$, the map
\[
 X(F)\longrightarrow K_2(D)/\theta_E K_2(E),\qquad P\longmapsto[P^*\alpha],
\]
is locally constant.
\end{proposition}

We change one coordinate at a time and compare the values on a curve using the rigidified Picard norm and Lemma~\ref{lem:paper-curve-reciprocity}(2).
The rigidification along the boundary makes the boundary contributions vanish, which gives constancy modulo the single subgroup $\theta_EK_2(E)$ rather than modulo $I_F(D)$.

\begin{proof}
We first set up a family of curves with rigidified Picard spaces.
In the projective closure of an affine embedding of $X$, let $H_\infty$ be the hyperplane at infinity.
We take a projective flat family $Y\to B$ of geometrically integral curves, parametrizing the complete-intersection curves constructed in the last paragraph, meeting $X$ in every fiber.
The boundary $M=Y\cap H_\infty$ is a finite flat relative Cartier divisor over $B$.
The forgetful map of the relative Picard space $\mathcal P_M$ with trivializations along $M$ has kernel
\[
 \operatorname{Res}_{M/B}\G_m/\G_m,
\]
which is smooth because it is the open locus of units in a relative projective bundle.
The Picard space is smooth since $H^2(\mathcal O)$ vanishes for curves.
On each proper geometrically integral fiber, regular functions come from the base, so a rigidified line bundle has trivial automorphism group, since $M$ is a nonempty finite flat boundary; hence $F$-points of $\mathcal P_M$ are represented by rigidified line bundles on the fibers.

Next we solve the norm equation near the diagonal.
The Picard norm along $E$ is smooth because it becomes addition after finite \'etale splitting.
For two points $P$ and $P'$ on the smooth locus, the line bundle $\mathcal O(P-P')$ has a canonical trivialization along $M$.
Starting from the trivial solution along the diagonal, we obtain nearby solutions by the \'etale inverse function theorem, which follows here from Hensel's lemma for a complete rank-$1$ valued field.
Let $\mathcal M$ be such a solution: a line bundle on $Y_E$, rigidified along $M$, whose norm is $\mathcal O(P-P')$ with its canonical rigidification.
Since line bundles over semilocal rings are free and units lift from the finite boundary, we can choose a rational section of $\mathcal M$ that takes the value $1$ along the boundary under the rigidification.
If $Z$ is the divisor of the section, then on the normalization
\[
 P-P'-\operatorname{Nm}(Z)=\operatorname{div}(f),\qquad f|_M=1.
\]
We also choose the section to be invertible at the points where the affine curve is not smooth over $F$, so that the support of $\operatorname{Nm}(Z)$ lies in the smooth open subset.
At the branches over the non-smooth points of the affine curve, $f$ is a unit, and $\alpha$ extends along the pullback to the normalization, so its $K_1$ boundary vanishes; at the branches at infinity, the rigidification makes the residue of $f$ equal to $1$.
Thus the boundary terms vanish by the localization calculation in the proof of Lemma~\ref{lem:paper-curve-reciprocity}, and the remaining terms are transfers through $E$.
Indeed, for a closed point $w$ in the support of $Z$ with image $z$, the pushforward of $[w]$ contributes $[\kappa(w):\kappa(z)]\cdot\cor_{\kappa(z)/F}\alpha(z)=\cor_{\kappa(w)/F}\alpha(w)$, which factors through $\theta_E$ since $E\subset\kappa(w)$.

Finally, we apply this calculation near an arbitrary smooth point.
The connected component of $X$ containing that point is smooth with an $F$-point, hence geometrically integral, and we replace $X$ by it; if it has dimension zero, local constancy is immediate.
By Bertini's theorem we choose geometrically integral complete-intersection curves through that point tangent to the coordinate axes; sections of sufficiently large degree remain free away from the point when the point and the tangent direction are prescribed.
For two nearby points differing only in one \'etale coordinate $t_i=r_i/s_0$, we correct each defining section by multiples of $s_0$ and $r_i$ so that it vanishes at both points.
Since the central curve is tangent to the $i$th axis, the correction is small, so flatness, geometric integrality, and the open norm condition persist.
In a product coordinate neighborhood, we move from one point to another by changing one coordinate at a time and apply the preceding calculation at each step.
\end{proof}

To apply local constancy to Steinberg words, we use the universal $K_2$-class on the smooth locus of the identity fiber of the word map defined below.
Denote the standard generators of the stable Steinberg group by $x_{ij}(a)$.
For a stable Steinberg word over $D$, we use an independent affine variable for each occurrence of a coefficient and take the matrix product in $\operatorname{SL}_n(D)$, the algebraic group of elements of reduced norm one in $M_n(D)$, for sufficiently large $n$; we call the resulting morphism the word map.

\begin{lemma}\label{lem:paper-steinberg-fiber}
Append to a stable Steinberg word over $D$ the letters
\[
 x_{ij}(t)\qquad\text{and}\qquad x_{ij}(u)\,x_{ji}(t)\,x_{ij}(-u),
\]
for all $i\ne j$ and all $u$ in a basis of $D$ containing $1$, with an independent $D$-valued variable $t$ for each of these letters.
Then the identity fiber of the word map is smooth at every point where the added variables are zero, and on the smooth open locus of the identity fiber the universal Steinberg word defines a $K_2$-class.
Moreover, the word map has an \'etale local section near the identity, sending the identity to the point at which the original variable coefficients and the added variables $t$ vanish.
\end{lemma}
\begin{proof}
At $t=0$ the word is unchanged.
The differentials of the added letters generate the off-diagonal entries and the elements $ute_{ii}-tue_{jj}$.
The diagonal differences for $u=1$, together with the linear span of the commutators $ut-tu$, give the entire subspace of reduced trace zero, since by nondegeneracy of the reduced trace pairing the linear span of the commutators in $D$ is $\ker\operatorname{Trd}$, also when the characteristic divides $\deg D$.
Thus the differential of the matrix product is surjective where the added variables vanish, and the identity fiber is smooth there.
On its smooth open locus, the universal Steinberg word has matrix product one and defines a $K_2$-class.
The same surjectivity gives an \'etale local section of the word map near the identity through the point where the original variable coefficients and the added variables $t$ vanish.
There the fixed letters $x_{ij}(u)$ and $x_{ij}(-u)$ cancel in the Steinberg group, so the word is trivial.
\end{proof}

\begin{corollary}\label{cor:paper-congruence-splitting}
Let $F$ be complete for a nontrivial rank-$1$ nonarchimedean valuation, and let $E/F$ be a finite separable splitting field of $D$.
For $n$ sufficiently large, let $\widetilde G\to E_n(D)$ be the pullback of the stable universal central extension by $K_2(D)$ to the elementary subgroup $E_n(D)$ of $\operatorname{SL}_n(D)$, and push it out to $K_2(D)/\theta_E K_2(E)$.
For $\gamma>0$ put $J_\gamma=\{x\in D:v_D(x)>\gamma\}$, where $v_D$ is the valuation of $D$.
Then this central extension splits over $E_n(D)\cap(1+\operatorname{Mat}_n(J_\gamma))$ for all sufficiently large $\gamma$.
In the discrete case, this gives a splitting over $E_n(D)\cap(1+\pi_F^m\operatorname{Mat}_n(\mathcal O_D))$ for all sufficiently large $m$, where $\pi_F$ is a uniformizer of $F$ and $\mathcal O_D$ is the valuation ring of $D$.
\end{corollary}
\begin{proof}
We choose an affine \'etale neighborhood $\rho\colon V\to\operatorname{SL}_n(D)$ of the identity on which the \'etale local section of Lemma~\ref{lem:paper-steinberg-fiber} (applied to the empty word) is defined and sends the identity to the trivial word, and we lift points near the identity by Hensel's lemma.
The multiplicative error of this section is a universal $K_2$-class on the closed subscheme of $V^3$ where $\rho(c)=\rho(a)\rho(b)$, which is smooth and affine near the identity triple, and we apply Proposition~\ref{prop:paper-picard-eval} to it.
The error vanishes at the identity triple, and hence, by local constancy, on a neighborhood of it.
Choosing $\gamma$ so that the product of the subgroup $E_n(D)\cap(1+\operatorname{Mat}_n(J_\gamma))$ with itself lies in this neighborhood, we obtain the splitting over this subgroup, and hence over its subgroups for larger $\gamma$.
\end{proof}

\begin{remark}\label{rem:paper-local-completion}
We have also proved the following extension of Theorem~\ref{thm:paper-local-norm} in degree two; it is not used in this paper, and we omit the proof.
Let $H$ be an algebraic extension of a nonarchimedean local field, and let $\widehat H$ be its completion.
Then, for every central division algebra $D$ over $\widehat H$, the reduced norm $\Nrd_2\colon K_2(D)\to K_2(\widehat H)$ is an isomorphism.
The proof descends $D$ to a finite subextension and applies Proposition~\ref{prop:paper-picard-eval} over $\widehat H$.
\end{remark}

\subsection{Local--global injectivity for the splitting transfer quotient}\label{sec:paper-Q-local-global}
In this subsection we prove the local--global principle for the splitting transfer quotient $Q_F(D)$ over global fields (Theorem~\ref{thm:paper-Q-local-global}) by comparing two global evaluations of a single universal $K_2$-class.
Proposition~\ref{prop:paper-global-picard} and Lemma~\ref{lem:paper-adelic-eval} are the global counterparts of the curve argument of Section~\ref{sec:paper-splitting}.
Here the norm equation must be solved simultaneously at all places $v\in S$.
Since only classes modulo $I_F(D)$ are needed, Lemma~\ref{lem:paper-curve-reciprocity}(1) handles the boundary terms, and no rigidification is required.

Throughout this subsection, $F$ is a global field and $D$ is a central division algebra over $F$.
Let $S$ be the finite set of places where $D$ does not split, including ramified real places, and write $D_v=D\otimes_FF_v$.

\begin{proposition}\label{prop:paper-global-picard}
Let $V$ be a nonempty smooth open subset of a projective geometrically integral curve $C/F$.
Let $P_D(C,V)\subset\operatorname{Pic}(C)$ be the subgroup generated by the divisor classes of the closed points of $V$ whose residue fields are separable over $F$ and split $D$.
Restriction induces an injection
\[
 \operatorname{Pic}(C)/P_D(C,V)\hookrightarrow
 \prod_{v\in S}\operatorname{Pic}(C_{F_v})/P_{D_v}(C_{F_v},V_{F_v}).
\]
\end{proposition}
\begin{proof}
Let $\mathcal L$ be a line bundle on $C$ whose class vanishes in every local quotient.
At each $v\in S$, we lift the separable split closed points to rational points over their residue fields and write
\[
 \mathcal L_{F_v}\simeq\bigotimes_iN_{E_{v,i}/F_v}(\mathcal M_{v,i})
\]
for line bundles $\mathcal M_{v,i}$ on $C_{E_{v,i}}$, where the $E_{v,i}/F_v$ are splitting fields of $D_v$; if this decomposition is empty, we add one trivial factor from a splitting field.
Put $n_v=\sum_i[E_{v,i}:F_v]$, choose a very ample line bundle $\mathcal H$ on $C$, and choose a sufficiently large common multiple $m$ of the $n_v$.
Twisting each $\mathcal M_{v,i}$ and each trivial bundle on $C_{E_{v,i}}$ by $\mathcal H^{\otimes m/n_v}$, we obtain
\[
 \bigotimes_iN_{E_{v,i}/F_v}\bigl(\mathcal M_{v,i}\otimes\mathcal H^{\otimes m/n_v}\bigr)\simeq\mathcal H^{\otimes m}\otimes\mathcal L,
 \qquad
 \bigotimes_iN_{E_{v,i}/F_v}\bigl(\mathcal H^{\otimes m/n_v}\bigr)\simeq\mathcal H^{\otimes m}.
\]

We next choose split divisors locally.
We take general sections of the twisted line bundles on the extensions $C_{E_{v,i}}$, and require their zero divisors to be geometrically reduced, to avoid the boundary $C\setminus V$, and to be disjoint from one another and from all their distinct conjugates, over all indices $i$ and all embeddings of the fields $E_{v,i}$.
Over $F_v^{\mathrm{sep}}$ the Weil restriction of the section space gives independent section spaces for the conjugates.
For two sections, vanishing at the same point imposes two independent linear conditions, so even when the point varies on the curve this is a proper closed condition.
Disjointness ensures that every point of a pushed-forward divisor has a residue field containing the corresponding splitting field, so we obtain split divisors in both localized complete linear systems $|\mathcal H^{\otimes m}\otimes\mathcal L|$ and $|\mathcal H^{\otimes m}|$ at $v$.

Finally we globalize.
Over the locus of divisors with simple zeros, the family of divisors is finite \'etale.
The inverse function theorem for the finite \'etale scheme of isomorphisms with a chosen local divisor then shows that nearby sections have the same local splitting type.
Weak approximation in the global section spaces gives two global divisors with the required conditions at every $v\in S$.
A closed point of either global divisor has a finite separable residue field $\kappa$ over $F$, and $D$ splits over every completion of $\kappa$; by the Brauer--Hasse--Noether theorem, $D_\kappa$ is split.
The difference of the two divisors represents $\mathcal L$.
When $S$ is empty, $D$ is split and general sections suffice.
\end{proof}

\begin{lemma}\label{lem:paper-adelic-eval}
Let $X/F$ be smooth quasi-projective and let $\alpha\in K_2(X;D)$.
For any $(b_v)\in\prod_{v\in S}X(F_v)$, the evaluation $[y^*\alpha]\in Q_F(D)$ is the same for all global points $y$ in a sufficiently small product neighborhood of $(b_v)$.
\end{lemma}
\begin{proof}
If the neighborhood contains a global point, we shrink it so that all its global points lie in a single connected component of $X$; this component is smooth and has an $F$-point, hence is geometrically integral, and we work there.
We fix a very ample line bundle on a projective closure of $X$; for a curve we use the fixed projective closure itself.
Dimension zero and the case $S=\varnothing$ are trivial.

We first construct curves through two nearby points.
In dimension at least two, we take \'etale coordinates
\[
 t_i=r_i/s_0
\]
at each $b_v$, where $r_i$ and $s_0$ are sections of the same sufficiently high power of the fixed very ample line bundle.
For every tangent direction $u$, we choose a geometrically integral complete-intersection curve through $b_v$ tangent to $u$; in sufficiently high degree, the sections satisfying the prescribed point and first-order conditions remain free away from the point, so Bertini's theorem applies.
We take the local curves to be cut out by sections of one sufficiently high degree, so that they form a projective flat family of geometrically integral curves over a common open set of parameters: the geometrically non-Cohen--Macaulay locus has codimension at least two, Bertini's theorem lets the curves avoid it, and the defining sections then form a regular sequence along them.
For the local norm we fix a finite separable splitting field.
The relative Picard norm is smooth, since it becomes addition after finite \'etale base change and $H^2(\mathcal O)=0$.
Since the curves have local rational points, Picard points are represented by line bundles.
By smoothness at the trivial solution on the diagonal, the condition that $\mathcal O(y-y')$ is a norm defines an open neighborhood of the diagonal in the family of curves with two smooth marked points $y,y'$.
We use this open neighborhood when the defining sections are moved by weak approximation.

Suppose that the secant direction of nearby points $y,y'$ is close to $u$.
Let $s_j$, $j\ge1$, be the defining sections of the central curve, and choose $i$ such that, for a fixed constant $c>0$,
\[
 |t_i(y')-t_i(y)|\ge c\,\|t(y')-t(y)\|.
\]
Sections vanishing at both points are given by
\[
 \beta_j=\frac{(s_j/s_0)(y')-(s_j/s_0)(y)}{t_i(y')-t_i(y)},\qquad
 s'_j=s_j-\bigl((s_j/s_0)(y)-\beta_j t_i(y)\bigr)s_0-\beta_jr_i.
\]
Since $d(s_j/s_0)(b_v)$ annihilates $u$, the correction $s'_j-s_j$ is small when the direction and the endpoints are close to $u$ and $b_v$, respectively; this follows by analytic expansion in \'etale coordinates over both real and nonarchimedean completions.
Compactness of the projective space of directions over $F_v$ gives a neighborhood of points that is uniform in the direction.

Now let $y,y'$ be global points in this neighborhood.
Weak approximation on the space of global sections vanishing at both points gives a single global curve, geometrically integral (an open condition preserved by approximation) and smooth at $y$ and $y'$, on which $\mathcal O(y-y')$ is a local norm from a splitting field at every $v\in S$.
As in the proof of Proposition~\ref{prop:paper-global-picard}, twisting by a very ample line bundle and pushing forward general sections shows that the local class of $\mathcal O(y-y')$ lies in the subgroup generated by split closed points.
Applied to the intersection $V$ of the smooth locus of the global curve with $X$, the proposition gives a divisor $Z$ supported on split closed points of $V$ with $[y-y']=[Z]$ in the Picard group.
Pulling back to the normalization gives $y-y'=\operatorname{div}(f)+Z$ for a rational function $f$.
Thus $\operatorname{div}(f)$ is supported in $V$, and $f$ is a unit at every branch of the boundary.
By Lemma~\ref{lem:paper-curve-reciprocity}(3),
\[
 y^*\alpha-y'^*\alpha\in I_F(D).
\]
\end{proof}

\begin{lemma}\label{lem:paper-local-I-single-field}
Let $k$ be a nonarchimedean local field and $B$ a central simple algebra over $k$, or let $k=\mathbb R$ and $B$ a nonsplit central simple algebra over $\mathbb R$.
For every finite separable splitting field $L/k$ of $B$,
\[
 I_k(B)=\theta_LK_2(L).
\]
\end{lemma}
\begin{proof}
Let $E/k$ be any finite splitting field of $B$, and put $M=EL$; the extension $M/E$ is separable.
For $k$ nonarchimedean, the field norm $N_{M/E}\colon K_2(M)\to K_2(E)$ is surjective \cite{BassTate}, \cite[p.~4]{BakRehmann1984}, so transitivity of transfer gives
\[
 \theta_E K_2(E)=\theta_M K_2(M)\subseteq\theta_LK_2(L).
\]
For $k=\mathbb R$, every splitting field is $\mathbb C$, and the same identity holds.
\end{proof}

For a unit $u$, we use the standard generators $x_{ij}(a)$ of the stable Steinberg group from Section~\ref{sec:paper-splitting} to define
\[
 w_{ij}(u)=x_{ij}(u)\,x_{ji}(-u^{-1})\,x_{ij}(u),\qquad h_{ij}(u)=w_{ij}(u)\,w_{ij}(-1).
\]
The matrix image of $h_{ij}(u)$ is diagonal with entries $u$ and $u^{-1}$ in positions $i$ and $j$, respectively, and $h_{ij}(1)=1$.

\begin{theorem}\label{thm:paper-Q-local-global}
For a global field $F$ and a central division algebra $D$ over $F$, the diagonal restriction
\[
 Q_F(D)\longrightarrow\prod_{v\in S}Q_{F_v}(D_v)
\]
is injective, where $S$ is the finite set of places where $D$ does not split.
\end{theorem}

For a class in the kernel, we lift its local images through the completions of one separable maximal subfield and compare two global evaluations of a $K_2$-class built from Steinberg words of a common shape.

\begin{proof}
Let $z\in K_2(D)$ have zero class in every local quotient.
If $S=\varnothing$, then $D$ is split, and there is nothing to prove.

We first show that the image of $z$ in the product of the local algebras comes from a separable maximal subfield.
Let $L\subset D$ be a separable maximal subfield, and put
\[
 D_S=\prod_{v\in S}D_v,\qquad L_S=L\otimes_F\prod_{v\in S}F_v=\prod_{v,\,w\mid v}L_w.
\]
For the base change $L_S\to D_S$ of $L\subset D$, the $(v,w)$-component of extension of scalars $K_2(L_S)\to K_2(D_S)$ is tensoring with the $(L_w,D_v)$-bimodule $e_wD_v$, where $e_w$ is the idempotent of the factor $L_w$.
The bimodule $e_wD_v$ is the base change to $L_w$ of the simple module $D$ over $L\otimes_F D^{\mathrm{op}}\simeq M_{\deg D}(L)$, so this component is $\theta_{L_w}$.
By Lemma~\ref{lem:paper-local-I-single-field}, $I_{F_v}(D_v)=\theta_{L_w}K_2(L_w)$ for every $v\in S$ and $w\mid v$.
Thus the image of $z$ in $K_2(D_S)$ comes from $K_2(L_S)$.

Let $\mathsf Z$ be a freely reduced Steinberg word representing $z$, and let $\mathsf T$ be a word in Matsumoto symbols representing a lift to $K_2(L_S)$ of the image of $z$ in $K_2(D_S)$.
We represent a symbol $\{a,b\}$ by $h_{12}(a)h_{12}(b)h_{12}(ab)^{-1}$ \cite{Matsumoto}, and a symbol in one factor of $L_S$ by taking the units equal to $1$ in the other factors.
Since $\mathsf Z=\mathsf T$ in the stable Steinberg group over $D_S$, we can write
\[
 \mathsf Z=\mathsf R\mathsf T
\]
in the free group on the elementary generators, where $\mathsf R$ is a finite product of conjugates of defining relations and their inverses.
We work at one matrix rank containing all the indices used, fix the indices and signs in $\mathsf R\mathsf T$, and use an independent variable for the coefficient of each occurrence of a letter.
Since $D_S$ is a product ring, this gives a single word shape at all $v\in S$, and the original coefficients give local parameter points $b_v$.

We construct two global points for this word shape.
In the first, we insert common approximating coefficients in the pairs of letters that cancel when $\mathsf R\mathsf T$ freely reduces to $\mathsf Z$, and we use the original coefficients of $\mathsf Z$ in the remaining letters; the value at this point is $z$.
In the second, we approximate the independent parameters of the defining relations, compute the sum and product coefficients after approximation, and approximate the units in $\mathsf T$ simultaneously at all $L_w$ by weak approximation on $L$; the value at this point lies in $\theta_L K_2(L)$.
At both points the matrix word is exactly the identity.

Finally, the additional letters of Lemma~\ref{lem:paper-steinberg-fiber} make the identity fiber smooth at all $b_v$, and the universal word defines a $K_2$-class on its smooth open subset.
If the two global points are close enough to the $b_v$, their evaluations $[z]$ and $0$ agree in $Q_F(D)$ by Lemma~\ref{lem:paper-adelic-eval}; hence $[z]=0$.
This proof does not use the vanishing of $\SK_2$ over local fields.
\end{proof}

\subsection{Twisted higher Chow groups and the image of the reduced norm}\label{sec:paper-global-weight-one}
In this subsection, we combine the local--global principle of Section~\ref{sec:paper-Q-local-global} with a cyclic splitting field to show that $\SK_2(D)=0$ over global fields, and we recover the formula of Bak and Rehmann \cite[\S2, Theorem 3]{BakRehmann1984} for the image of the reduced norm (Theorem~\ref{thm:paper-global-norm}).
We also introduce twisted higher Chow groups and prove the vanishing in weight one that is used in Section~\ref{sec:global-spectrum} (Theorem~\ref{thm:paper-global-weight-one}).

\begin{lemma}\label{lem:paper-cyclic-real-split}
Let $F$ be a global field and $D$ a central division algebra over $F$.
There is a cyclic maximal subfield $L_0\subset D$ that splits completely at every real place of $F$ at which $D$ splits.
\end{lemma}
\begin{proof}
If $F$ is a global function field or $\deg D$ is odd, a cyclic maximal subfield of $D$ has the required property, since a Galois extension of odd degree splits completely at all real places.
Let $F$ be a number field and $\deg D$ even.
Starting from a cyclic maximal subfield, we correct its class field theory character by a quadratic character that is trivial at every place where $D$ does not split and cancels the signs at the real places where $D$ splits; weak approximation of square classes gives these local conditions.
The corrected character has order dividing $\deg D$ because the original and quadratic characters do.
At the places where $D$ does not split, the local degrees are unchanged and are multiples of the local indices, whose least common multiple is $\deg D$.
Hence the corrected cyclic extension $L_0/F$ has degree $\deg D$ and splits $D$, so it embeds in $D$ as a maximal subfield.
\end{proof}

\begin{theorem}\label{thm:paper-global-norm}
For a global field $F$ and any central division algebra $D$ over $F$, we have $\SK_2(D)=0$.
If $S_{\mathbb R}$ is the set of ramified real places of $D$, then the image of the reduced norm is
\[
 R_D=\ker\left(K_2(F)\xrightarrow{(s_v)_{v\in S_{\mathbb R}}}
                        (\Z/2)^{S_{\mathbb R}}\right).
\]
Here $s_v\{a,b\}=1$ if both $a$ and $b$ are negative at $v$, and is zero otherwise.
Consequently $C_2(D/F)\simeq(\Z/2)^{S_{\mathbb R}}$.
\end{theorem}

We use cyclic Hilbert~90 with the field $L_0$ of Lemma~\ref{lem:paper-cyclic-real-split} to prove injectivity of $\Nrd_2$ on $I_F(D)$.
Local generation and Theorem~\ref{thm:paper-Q-local-global} then give $K_2(D)=I_F(D)$, whose reduced norm image is $N_{L_0/F}K_2(L_0)$.

\begin{proof}
Let $L_0/F$ be as in Lemma~\ref{lem:paper-cyclic-real-split}.
We first show that $\Nrd_2$ is injective on $I_F(D)$.
Every real place of a finite splitting field $E/F$ lies over a real place of $F$ at which $D$ splits, and $L_0$ splits completely there.
Hence no real place of $E$ becomes complex in $EL_0$, so the field norm $K_2(EL_0)\to K_2(E)$ is surjective \cite[\S2, Theorem 2]{BakRehmann1984}.
Transitivity of transfer gives
\[
 I_F(D)=\theta_{L_0}K_2(L_0).
\]
Let $\sigma$ be a generator of $G=\Gal(L_0/F)$.
By Skolem--Noether, $\theta_{L_0}\bigl((1-\sigma)K_2(L_0)\bigr)=0$, so $\theta_{L_0}$ factors through $K_2(L_0)_G$ and we obtain
\[
 K_2(L_0)_G\xrightarrow{\ \theta_{L_0}\ }I_F(D)\xrightarrow{\ \Nrd_2\ }K_2(F),
 \qquad \Nrd_2\circ\theta_{L_0}=N_{L_0/F}.
\]
If $\operatorname{char}F\nmid\deg D$, cyclic Hilbert 90 \cite[Theorem 14.1]{MerkurjevSuslin1982} gives $\ker N_{L_0/F}=(1-\sigma)K_2(L_0)$, so the composite is injective on $K_2(L_0)_G$.
If $\operatorname{char}F=p$ divides $\deg D$, let $P$ be a Sylow $p$-subgroup of $G$ and $M=L_0^P$.
Since $[L_0:L_0^p]=[M:M^p]=p$, the groups $K_2(L_0)$ and $K_2(M)$ are uniquely $p$-divisible \cite[Theorem 8.1]{GeisserLevine}, and the identities $\res_{L_0/M}\circ N_{L_0/M}=\sum_{g\in P}g$ and $N_{L_0/M}\circ\res_{L_0/M}=|P|$ give $K_2(L_0)_P\simeq K_2(M)$.
Cyclic Hilbert 90 for the cyclic extension $M/F$ of degree prime to $p$ then gives the same injectivity.
Since the first map is surjective, $\Nrd_2$ is injective on $I_F(D)$.

Finally we show that $K_2(D)=I_F(D)$.
At every finite place $v$, the group $K_2(D_v)$ is generated by the images of splitting transfers: for a Galois splitting field $L/F_v$ of $D_v$, the norm $N_{(D_v)_L/D_v}$ is $\theta_L$ under the split Morita identification, so Lemma~\ref{lem:paper-local-even-coinvariants} gives
\[
 \theta_L\colon K_2(L)_{\Gal(L/F_v)}\xrightarrow{\ \sim\ }K_2(D_v).
\]
At ramified real places, the splitting transfer $\theta_{\mathbb C}$ surjects onto $K_2$ of the real quaternion algebra as well \cite[Proposition 6.6.1(2)]{KahnLevine2010}.
Hence the target of Theorem~\ref{thm:paper-Q-local-global} vanishes, so $K_2(D)=I_F(D)$, and the injectivity of $\Nrd_2$ on $I_F(D)$ gives $\SK_2(D)=0$.
The image of $\Nrd_2$ is $N_{L_0/F}K_2(L_0)$, which equals $R_D$ by \cite[\S2, Theorem 2]{BakRehmann1984}.
Weak approximation realizes any real sign vector on a symbol $\{a,-1\}$, so the displayed sign map is surjective.
\end{proof}

\begin{definition}\label{def:paper-twisted-chow}
Let $X/F$ be a smooth equidimensional scheme, and put
\[
 D_X=D\otimes_F\mathcal O_X,\qquad \Delta_F^m=\Spec F[t_0,\ldots,t_m]\Big/\Bigl(\sum_i t_i-1\Bigr).
\]
An integral closed subscheme $W$ of codimension $r$ in $X\times\Delta_F^m$ is an admissible support if it meets every face $X\times\{t_{i_1}=\cdots=t_{i_s}=0\}$ properly.
We define
\[
 z^r(X,m;D)=\bigoplus_{W\ \text{admissible}}
                   K_0(D_{\kappa(W)}).
\]
Write $\partial_i^*$ for pullback along the $i$th face and $\partial=\sum_{i=0}^m(-1)^i\partial_i^*$ for the differential.
The pullbacks are constructed by $K$-theory with support and satisfy the face and degeneracy identities, so $z^r(X,\bullet;D)$ is a chain complex.
We set
\[
 \CH^r(X,m;D)=H_m\bigl(z^r(X,\bullet;D)\bigr)
\]
and abbreviate this to $\CH^r(F,m;D)$ when $X=\Spec F$.
For $D=F$ these are Bloch's higher Chow groups \cite{LevineHigherChow}, \cite[Definitions 5.6.1 and 5.6.3]{KahnLevine2010}.

The reduced rank gives an injection
\[
 K_0(D_{\kappa(W)})\hookrightarrow\Z,
 \qquad\im=\ind(D_{\kappa(W)})\Z.
\]
It induces the \emph{reduced cycle map} $\rho$, which realizes the twisted cycle complex as a subcomplex of the ordinary cycle complex.
Under this identification the coefficient differential is face restriction with the usual intersection multiplicities.
For finite coefficients we use the derived tensor product of the complex, and we distinguish the homology of this derived tensor product from the quotient of integral homology by the coefficient integer \cite[Sections 5.6 and 5.9, Theorem 6.5.3]{KahnLevine2010}.
\end{definition}

The same methods give the following vanishing for twisted higher Chow groups in weight one: in degree two we use the local--global principle, and in degree three the field $L_0$ and the reduced cycle map.

\begin{theorem}\label{thm:paper-global-weight-one}
Let $F$ be a global field and let $D$ be a central division algebra over $F$.
Then $\CH^1(F,2;D)=\CH^1(F,3;D)=0$.
\end{theorem}

\begin{proof}
\emph{Step 1: identification with $Q_k(D)$.}
Over a perfect field $k$, we use the edge exact sequence of \cite[Theorem 6.7.1(2)]{KahnLevine2010},
\[
 0\to\CH^1(k,3;D)\to\CH^2(k,2;D)\xrightarrow{e_2}K_2(D)
 \to\CH^1(k,2;D)\to0.
\]
Transfers from splitting fields generate $\CH^2(k,2;D)$ \cite[Lemmas 6.6.2--6.6.4]{KahnLevine2010}: at each admissible closed point $x$, we choose a maximal subfield of the division algebra in the Brauer class of $D_{\kappa(x)}$, and pushforward from the resulting point gives the generator of the $K_0$ coefficient.
Over a splitting field $E$, the map $e_2$ is the isomorphism $\CH^2(E,2)\simeq K_2(E)$ \cite{NesterenkoSuslin}, compatible with transfer.
Thus
\[
 \im e_2=I_k(D),\qquad \CH^1(k,2;D)=Q_k(D),
\]
naturally in field extensions.

\emph{Step 2: imperfect fields.}
For a field $k$ of characteristic $p$, we first invert $p$ and pass to the perfect closure of $k$.
For a finite purely inseparable extension $M/k$, we have
\[
 \cor_{M/k}\circ\res_{M/k}=\res_{M/k}\circ\cor_{M/k}=[M:k].
\]
For cycles, the factor $[M:k]$ is the product of the support multiplicity and the residue degree.
For $K$-theory, we filter the local Artinian ring of the self-base-change by the powers of its radical and use additivity \cite{QuillenHigher}.
Twisted cycle complexes, $K$-theory, and the quotient by splitting transfers commute with filtered colimits of fields \cite{KahnLevine2010}.
Hence restriction to the perfect closure is an isomorphism after inverting $p$, and
\[
 \CH^1(k,2;D)[1/p]\simeq Q_k(D)[1/p].
\]

\emph{Step 3: the $p$-primary part.}
For the $p$-primary part, restriction and corestriction along an extension of degree prime to $p$ that splits the prime-to-$p$ component of the Brauer class reduce the problem to a $p$-primary class.
Since $[F:F^p]=p$ for local fields $F$ of characteristic $p$ and for global function fields $F$, the rank-$1$ Frobenius comparison of Proposition~\ref{prop:app-frobenius} makes the reduced cycle map in positive weights a quasi-isomorphism.
Hence the $p$-primary parts in weight one and degrees two and three vanish by the commutative calculation \cite[Proposition 6.8.1]{KahnLevine2010}.

\emph{Step 4: degree two over $F$.}
The proof of Theorem~\ref{thm:paper-global-norm} gives $Q_F(D)=0$.
In characteristic zero, Step~1 gives $\CH^1(F,2;D)\simeq Q_F(D)=0$; in positive characteristic, Step~2 gives the same after inverting $p$, and Step~3 shows that the $p$-primary torsion vanishes.

\emph{Step 5: degree three over $F$.}
The generation of the diagonal group $\CH^2(F,2;D)$ by transfers from splitting fields holds over any field, so this step takes place over $F$; only the edge exact sequence below requires a perfect field.
For splitting fields $E$, the isomorphism $\CH^2(E,2)\simeq K_2(E)$ is compatible with transfer \cite{NesterenkoSuslin}.
For the field $L_0$ of Lemma~\ref{lem:paper-cyclic-real-split}, the field norms $K_2(EL_0)\to K_2(E)$ for splitting fields $E$ are surjective, as in the proof of Theorem~\ref{thm:paper-global-norm}, so the pushforward
\[
 \psi\colon\CH^2(L_0,2)\longrightarrow\CH^2(F,2;D)
\]
is surjective.
Then $\psi$ factors through the coinvariants of $G=\Gal(L_0/F)$ because the structure map $\Spec L_0\to\Spec F$ is $G$-invariant, and $\rho\circ\psi=N_{L_0/F}$ is the compatibility of $\rho$ with pushforward \cite[Lemma 6.6.4]{KahnLevine2010}.
Hence cyclic Hilbert 90, combined in characteristic $p$ with the averaging argument in the proof of Theorem~\ref{thm:paper-global-norm}, makes $\rho\circ\psi$ injective on coinvariants, and
\[
 \rho\colon\CH^2(F,2;D)\to\CH^2(F,2)
\]
is injective.
When $F$ is perfect, under $\CH^2(L_0,2)\simeq K_2(L_0)$ and $\CH^2(F,2)\simeq K_2(F)$ we have
\[
 \Nrd_2\circ e_2\circ\psi=\Nrd_2\circ\theta_{L_0}=N_{L_0/F}=\rho\circ\psi,
\]
and the surjectivity of $\psi$ gives $\Nrd_2\circ e_2=\rho$ \cite[Proposition 6.6.5]{KahnLevine2010}.
Hence the injectivity of $\rho$ and the edge exact sequence give
\[
 \CH^1(F,3;D)=\ker e_2\subseteq\ker\rho=0.
\]
In positive characteristic, restriction to the perfect closure commutes with $\rho$ and is an isomorphism on $\CH^2(-,2;D)[1/p]$ and $\CH^2(-,2)[1/p]$, so the injectivity of $\rho$ over $F$ shown above gives its injectivity over the perfect closure after inverting $p$.
There the edge exact sequence and the identity $\Nrd_2\circ e_2=\rho$ give $\CH^1(F,3;D)[1/p]=0$, while the $p$-primary part vanishes by the comparison above.
\end{proof}

\section{Twisted \texorpdfstring{$K$}{K}-theory over global fields}\label{sec:global}
This section compares $K(D)$ with $K(F)$ over a global field $F$ in positive degrees.
We first prepare the comparisons of twisted and ordinary higher Chow groups used in low degrees (Section~\ref{sec:global-chow}).
We then construct the comparison map, split it in degrees at least two, and study its fiber (Section~\ref{sec:global-spectrum}); when $D$ is ramified at a real place, the statements about the comparison, its splitting, and its fiber hold after inverting $2$.
\subsection{Structural pushforward and twisted higher Chow groups}\label{sec:global-chow}
This subsection prepares the cycle-theoretic input for Section~\ref{sec:global-spectrum}.
We relate twisted higher Chow groups to ordinary ones through the pushforward along the structure map of a Severi--Brauer variety: this pushforward factors through the twisted higher Chow groups, the factored map is an isomorphism on the diagonal and surjective on the next row, and composing it with the reduced cycle map gives the usual pushforward.
After $\ell$-localization, the reduced cycle map is an isomorphism on both rows in weights $r\ge2$, except for $\ell=2$ over number fields with real places (Proposition~\ref{prop:global-two-chow-rows}); for $\CH^2(F,3;D)$ this exception is removed in Proposition~\ref{prop:global-real-weight-two}.
These isomorphisms are used for $K_3(D)$ and for the group $\pi_4$ of the comparison fiber.

Let $F$ be a global field and $D$ a central division algebra over $F$.
At a closed support $T$, the reduced rank identifies $K_0(D_{\kappa(T)})$ with $\ind(D_{\kappa(T)})\Z$, inducing the reduced cycle map $\rho$ of Definition~\ref{def:paper-twisted-chow}.
We compare $\rho$ with the \emph{structural trace}, the pushforward along the structure map of a Severi--Brauer variety.
Write $\mathcal F^rK(D)$ for the homotopy coniveau tower.
The homotopy group $\pi_m$ of the weight-$r$ layer is $\CH^r(F,m;D)$, and $\pi_i\mathcal F^rK(D)=0$ for $i<r$.
The weight-zero layer has only the rank group in degree zero \cite[Corollary 6.1.4 and Section 6.6]{KahnLevine2010}.
Over a global function field, which is imperfect, the identification of the layers is available by Lemma~\ref{lem:app-twisted-tower}.

\subsubsection{Structural pushforward and reduced rank}
Our first goal is Lemma~\ref{lem:integral-structural-push}, which expresses the twisted higher Chow groups on the diagonal and the next row through higher Chow groups of a Severi--Brauer variety, compatibly with the reduced cycle map.
Its proof uses the following boundary correction.

\begin{lemma}\label{lem:paper-structural-boundary-correction}
Let $F$ be an infinite field, let $B$ be a central simple algebra over $F$, put $X=\SB(B)$ and $n=\dim X$, let $r\ge1$, and write $f_*\colon z^{n+r}(X,\bullet)\to z^r(F,\bullet)$ for the proper pushforward along the structure map $f\colon X\to\Spec F$.
If a chain $c$ of degree $r$ satisfies $f_*c=0$, then there is a chain $W$ of degree $r+1$ with
\[
 f_*W=0,\qquad\partial W=c.
\]
\end{lemma}
\begin{proof}
The chain $c$ has codimension $n+r=\dim(X\times\Delta^r)$, so its components are closed points, and proper intersection with the faces forces their images in $\Delta^r$ to avoid all faces.
Over each interior closed point $q\in\Delta^r$, the component of $c$ is a zero-cycle of degree zero on $X_{\kappa(q)}$.
The degree map on $\CH_0$ of a Severi--Brauer variety is injective \cite[Corollary 7.3]{ChernousovMerkurjev2006}, so the component is a sum of pushforwards of divisors of rational functions on proper normal integral curves $C\to X_{\kappa(q)}$.
By normality of $C$, a rational function $\psi$ with $\operatorname{div}(\psi)\ne0$ defines a nonconstant morphism $\psi\colon C\to\PP^1_{\kappa(q)}$.
This morphism of proper curves is finite, and it is flat because $\psi_*\mathcal O_C$ is torsion-free over the regular curve $\PP^1_{\kappa(q)}$.

Set $L=\kappa(q)$ and place the corresponding $L$-point $a=(a_0,\ldots,a_{r+1})$ on the face $t_0=0$ of $\Delta^{r+1}_L$.
Thus $a_0=0$ and $a_j\ne0$ for $j>0$.
We choose $b_j\in L$ with $\sum_j b_j=0$, $b_j\ne0$, and $a_i b_j-a_j b_i\ne0$ for $i\ne j$; this avoids finitely many proper linear subspaces of the hyperplane $\sum_jb_j=0$ over the infinite field $L$.
The line $t_j=a_j+b_j u$, parametrized by $u\in\mathbb A^1_L$, meets the codimension-one faces transversely at distinct points $u_j=-a_j/b_j$ and avoids all faces of higher codimension.
We interpolate a polynomial $h$ with $h(u_0)=1$ and $h(u_j)=0$ for $j>0$ and set $\varphi(u)=[1-h(u):h(u)]$.
Pullback along $\varphi$ and along the constant map with value $\infty$ gives
\[
 C\times_{\PP^1_L,\varphi}\mathbb A^1_L,\qquad
 C\times_{\PP^1_L,\infty}\mathbb A^1_L,
\]
both finite flat of the same degree over the line, where $0=[0:1]$ and $\infty=[1:0]$ in $\PP^1_L$.
The difference of their fundamental cycles gives an admissible chain $W_\psi$ by proper pushforward to $X_L\times\Delta^{r+1}_L$ followed by finite pushforward to $F$.
Equal degrees give $f_*W_\psi=0$.
By transversality and finite flatness, restriction to $t_j=0$ gives the difference of the fiber cycles at $u_j$, namely $\psi^{-1}(0)-\psi^{-1}(\infty)$ for $j=0$ and zero for $j>0$.
The pushforward has no extra degree factor since $L=\kappa(q)$, and the boundary sign of $t_0=0$ is $+1$.
Summing the chains $W_\psi$ over the divisors, we obtain a chain $W$ with
\[
 f_*W=0,\qquad\partial W=c.
\]
\end{proof}

\begin{lemma}\label{lem:integral-structural-push}
Let $F$ be an infinite field, let $D$ be a central division algebra over $F$, let $B$ be a central simple algebra over $F$ with $[B]=[D]$, and put $X=\SB(B)$ and $n=\dim X$.
For every $r\ge1$ the structure map $f\colon X\to\Spec F$ induces
\[
 \CH^{n+r}(X,r+1)\twoheadrightarrow\CH^r(F,r+1;D),\qquad
 \CH^{n+r}(X,r)\xrightarrow{\sim}\CH^r(F,r;D).
\]
The composite with $\rho$ is the usual proper pushforward along the structure map.
These maps commute with restriction and transfer for finite extensions.
\end{lemma}

\begin{proof}
We first show that $f_*$ factors through the twisted cycle complex.
In each simplicial degree, let $W\subset X\times\Delta^m$ be an integral cycle with image $T\subset\Delta^m$.
Then
\[
 f_*[W]=
 \begin{cases}
  [\kappa(W):\kappa(T)]\,[T], & \dim T=\dim W,\\
  0, & \dim T<\dim W.
 \end{cases}
\]
In the first case $X$ has a $\kappa(W)$-rational point, so $[\kappa(W):\kappa(T)]$ is divisible by $\ind(D_{\kappa(T)})$.
Consequently $f_*$ factors through the twisted cycle complex $z^r(F,\bullet;D)\subset z^r(F,\bullet)$ embedded by reduced rank, compatibly with face restrictions, flat base change, and proper transfer.

We next show surjectivity on chains.
For $m=r,r+1$, an admissible support $T$ is a closed point or an integral curve, respectively.
To lift the coefficient generator $\ind(D_{\kappa(T)})[T]$, we take the closure of a closed point of $X_{\kappa(T)}$ of degree $\ind(D_{\kappa(T)})$.
When $T$ is a curve, the closure is proper over $T$, and its closed fibers are proper closed subsets of the closure, an integral curve; hence they are zero-dimensional, and the closure is finite over $T$.
It inherits the proper intersection conditions with faces from $T$.
This proves surjectivity on chains in both degrees.

We lift a cycle in degree $r+1$ of the twisted cycle complex to a chain $Z$; then $f_*\partial Z=0$.
For the chain $W$ of Lemma~\ref{lem:paper-structural-boundary-correction} with $c=\partial Z$, the difference $Z-W$ is a cycle lifting the given cycle.
In degree $r$ the chain group one degree lower vanishes, giving surjectivity.
Suppose that a degree-$r$ chain $Y$ and a chain $V$ satisfy $f_*Y=\partial V$.
We choose a lift $\widetilde V$ of $V$ and subtract $\partial\widetilde V$ from $Y$.
Then $f_*(Y-\partial\widetilde V)=0$, so $Y-\partial\widetilde V$ is also a boundary by Lemma~\ref{lem:paper-structural-boundary-correction} with $c=Y-\partial\widetilde V$.
This proves injectivity in degree $r$.
\end{proof}

\subsubsection{Reduced norms on the diagonal and the next row}
Next we use Lemma~\ref{lem:integral-structural-push} to study the reduced cycle map $\rho$ on $\CH^r(F,m;D)$ for $m=r,r+1$.

\begin{proposition}\label{prop:global-two-chow-rows}
Let $F$ be a global field and $\ell$ a prime.
If $F$ is a number field, we assume that $\ell$ is odd or that $F$ has no real places.
For every $r\ge2$,
\[
 \rho\colon\CH^r(F,m;D)_{(\ell)}
    \xrightarrow{\sim}\CH^r(F,m)_{(\ell)}
    \qquad(m=r,r+1).
\]
The isomorphisms commute with restriction and transfer.
\end{proposition}
\begin{proof}
First let $\ell=\operatorname{char}F$, so that $[F:F^\ell]=\ell$.
We split the prime-to-$\ell$ component of the Brauer class over an extension $F'/F$ of degree prime to $\ell$.
Since $\cor_{F'/F}\res_{F'/F}=[F':F]$, the kernel and cokernel of $\rho_{(\ell)}$ are detected after extension.
There we apply the rank-$1$ Frobenius comparison of Proposition~\ref{prop:app-frobenius}.

Now let $\ell\ne\operatorname{char}F$, and put $X=\SB(D)$, $n=\dim X$, and $R=\Z_{(\ell)}$.
We consider the fiber of the structural trace
\[
 \mathcal C_r=\fib\!\left[
 R\Gamma_{\mot}(X,R(n+r))[2n]
   \xrightarrow{f_*}R\Gamma_{\mot}(F,R(r))\right].
\]
Since $\CH^r(F,m)=H^{2r-m}_{\mot}(F,\Z(r))$ and $\CH^{n+r}(X,m)=H^{2r-m}\bigl(R\Gamma_{\mot}(X,\Z(n+r))[2n]\bigr)$, the rows $m=r$ and $m=r+1$ correspond to degrees $r$ and $r-1$ of the trace.
By the hypothesis on $\ell$, $\cd_\ell F\le2$.
At a point of dimension $a$ and codimension $c=n-a$ in $X$, purity gives weight $a+r$, and the residue field has $\ell$-cohomological dimension at most $a+2\le a+r$ (as $r\ge2$).
The layers of the support filtration therefore lie in the range of the Beilinson--Lichtenbaum comparison \cite{Voevodsky}, \cite[Theorem 2.3(2)]{Yamazaki2013}.
Since the cycle class maps commute with proper trace, this identifies the finite-coefficient trace with the \'etale trace in the required degrees.

Over $F^{\mathrm{sep}}$, $X$ is $\PP^n$.
The trace shifted by $[2n]$ is an isomorphism on the top geometric group $\Z/\ell(r)$.
The remaining groups lie in geometric cohomological degrees $-2,-4,\ldots,-2n$.
Modulo $\ell$, Hochschild--Serre gives
\[
 H^q(\mathcal C_r/\ell)=0\qquad(q\ge1).
\]
With rational coefficients, we descend the projective bundle decomposition from a finite splitting field by restriction and transfer, using the hyperplane class $c_1(\omega_X^{-1})/(n+1)$.
The trace is projection onto the top component, so
\[
 (\mathcal C_r)_\Q\simeq
 \bigoplus_{k=1}^{n}R\Gamma_{\mot}(F,\Q(r+k))[2k].
\]
For number fields the weight decomposition of Borel's calculation \cite{Borel} via the Chern character gives $H^i_{\mot}(F,\Q(b))=0$ for $b\ge2$ and $i\ne1$ (see the proof of \cite[Theorem 3.2]{Scholbach2011}).
For global function fields the groups in this range vanish by Harder's theorem \cite{Harder}; this is consistent with localization, since the positive $K$-groups of the finite residue fields are torsion \cite{QuillenFinite}.
Thus $H^q((\mathcal C_r)_\Q)=0$ for $q\ge0$.

Combining the two computations, the coefficient sequence shows that multiplication by $\ell$ is invertible on $H^q(\mathcal C_r)$ for $q\ge2$.
These groups are $\ell$-primary torsion and hence zero.
The group $H^1(\mathcal C_r)$ is $\ell$-divisible.
Consequently the structural trace is an isomorphism in degrees $q\ge2$ and is surjective with $\ell$-divisible kernel in degree one.
Write $\widetilde f_*$ for the maps $\CH^{n+r}(X,m)\to\CH^r(F,m;D)$ of Lemma~\ref{lem:integral-structural-push}, so that $f_*=\rho\circ\widetilde f_*$.
On the diagonal, $\widetilde f_*$ is an isomorphism and identifies $\rho$ with the trace in degree $r$.
On the next row, $f_*$ and $\widetilde f_*$ are both surjective, so $\rho$ is surjective and $\ker\rho=\widetilde f_*(\ker f_*)$.
This kernel vanishes for $r\ge3$, since $f_*$ is injective in degree $r-1\ge2$, and it is $\ell$-divisible for $(r,m)=(2,3)$.
In the case $(r,m)=(2,3)$, restriction and transfer show that $\ker\rho$ is killed by the degree of a finite splitting field.
The prime-to-$\ell$ part of this degree is a unit, so this $\ell$-divisible group is killed by a power of $\ell$ and is therefore zero.
\end{proof}

\subsubsection{The trace relative to the real places}
Proposition~\ref{prop:global-two-chow-rows} excludes the prime two over number fields with real places.
For $\CH^2(F,3;D)$ we remove this restriction in Proposition~\ref{prop:global-real-weight-two}, using the structural trace relative to the real places.

\begin{proposition}\label{prop:global-real-weight-two}
For a number field $F$ with a real place and any central division algebra $D$ over $F$,
\[
 \rho\colon\CH^2(F,3;D)_{(2)}\xrightarrow{\sim}\CH^2(F,3)_{(2)}.
\]
\end{proposition}
\begin{proof}
Set $X=\SB(D)$, $n=\dim X$, and $R=\Z_{(2)}$.
The structural trace in degree one is
\[
 f_*\colon H^{2n+1}_{\mot}(X,R(n+2))\longrightarrow H^1_{\mot}(F,R(2)),
\]
and the target of $\rho$ is $\CH^2(F,3)_{(2)}=K_3^{\mathrm{ind}}(F)_{(2)}$.

We write $R\Gamma_{\mot,\mathrm{rel}}$ and $H^i_{\mathrm{rel}}$ for the fiber of restriction to the algebraic real closures $F_v^{\mathrm{rc}}$ at all real places $v$, and its cohomology.
For $M=R,\Q,\Q/R$, set
\[
 \mathcal C_{\mathrm{rel}}(F,M)=\fib\left(
 R\Gamma_{\mot,\mathrm{rel}}(X,M(n+2))[2n]
 \xrightarrow{f_*} R\Gamma_{\mot,\mathrm{rel}}(F,M(2))
 \right).
\]
We first show that $H^0\mathcal C_{\mathrm{rel}}(F,\Q/R)$ and $H^1\mathcal C_{\mathrm{rel}}(F,\Q/R)$ vanish.

We apply finite-coefficient motivic--\'etale comparison to the support filtration.
At a point of dimension $a$ on $X$, purity shifts the weight to $a+2$.
Without real relativization the virtual $2$-cohomological dimension is at most $a+2$.
In the relative calculation, the constant field in the residue field is a finite number field, and its real and complex branches account for the restriction targets.
By Jannsen \cite[Proposition 1.2 and Lemma 1.3]{Jannsen}, restriction to the real branches is surjective in degree $a+2$ and an isomorphism above that degree.
The complex branches have cohomological dimension at most $a$.
Thus the comparison cofiber above the Beilinson--Lichtenbaum range \cite{Voevodsky} vanishes after taking the relative fiber.
For $x\in X$, let $k_x$ be the algebraic closure of $F$ in $\kappa(x)$; the residue fields on the real-closure side are the factors of $\kappa(x)\otimes_FF_v^{\mathrm{rc}}$.
Closed supports over a real closure descend to finite intermediate fields and lie in the base change of their finite images; since codimension is preserved, base-changed supports are cofinal, so the relative layers can be computed residue field by residue field \cite[Proposition 1.2, Lemma 1.3(b), Remark 1.4(a)]{Jannsen}.
Taking the direct limit over finite coefficients, we compute $\mathcal C_{\mathrm{rel}}(F,\Q/R)$ by \'etale cohomology.

Over $F^{\mathrm{sep}}$ the trace is an isomorphism on $(\Q/R)(2)$ in top geometric degree $2n$, so after shifting by $[2n]$ the geometric fiber has cohomology only in even degrees at most $-2$, with $(\Q/R)(3)$ in degree $-2$ for $n\ge1$.
The relative cohomological dimension of the base is at most two, so positive degrees vanish and the only possible contribution in degree zero is $H^2_{\mathrm{rel}}(F,(\Q/R)(3))$.
The motivic coefficient triangle gives
\[
 H^2(F,(\Q/R)(3))\simeq K_3^M(F)\otimes R.
\]
Indeed, the rational groups $H^2_{\mot}(F,\Q(3))$ and $H^3_{\mot}(F,\Q(3))$ on either side of the $\Q/R$-coefficient term vanish by $K_4(F)_\Q=0$ \cite{Borel} and $K_3^M(F)_\Q=0$, respectively.
The same holds at algebraic real closures by taking the direct limit over finite number fields.
Using the description of $K_3^M(F)$ by signs at the real places \cite{BassTate}, we have
\[
 H^2(F,(\Q/R)(3))\simeq\prod_{v\,\mathrm{real}}H^2(F_v^{\mathrm{rc}},(\Q/R)(3)).
\]
Complex conjugation acts on $(\Q/R)(3)$ by $-1$, so
\[
 H^1\bigl(F_v^{\mathrm{rc}},(\Q/R)(3)\bigr)=(\Q/R)(3)/2(\Q/R)(3)=0.
\]
The relative exact sequence gives $H^2_{\mathrm{rel}}(F,(\Q/R)(3))=0$.
Hence $H^i\mathcal C_{\mathrm{rel}}(F,\Q/R)=0$ for $i\ge0$.

Next we compute with rational coefficients.
Using the hyperplane class above, we descend the rational projective bundle decomposition.
Since the structural trace projects onto the top component,
\[
 \mathcal C_{\mathrm{rel}}(F,\Q)\simeq
 \bigoplus_{j=3}^{n+2}
 R\Gamma_{\mot,\mathrm{rel}}(F,\Q(j))[2j-4].
\]
For $i\ge1$ the absolute group $H^{2j+i-4}_{\mot}(F,\Q(j))$ vanishes since the degree exceeds the weight except for $j=3$ and $i=1$, where it is $K_3^M(F)_\Q=0$.
The preceding real-closure group $H^{2j+i-5}_{\mot}(F_v^{\mathrm{rc}},\Q(j))$ is handled in the same way.
The exceptions within the weight bound are $H^2_{\mot}(-,\Q(3))$, $K_3^M(-)_\Q$, and $K_4^M(-)_\Q$, all zero by the number field calculation and filtered colimits.
Therefore $H^i\mathcal C_{\mathrm{rel}}(F,\Q)=0$ for $i\ge1$, and the coefficient triangle gives $H^i\mathcal C_{\mathrm{rel}}(F,R)=0$ for $i\ge1$.

We now return to the absolute trace.
For a field extension $k/F$ and $M$ as above, write $X_k=X\times_Fk$ and define
\[
 \mathcal C_{\mathrm{abs}}(k,M)=\fib\left(
 R\Gamma_{\mot}(X_k,M(n+2))[2n]
 \xrightarrow{f_*} R\Gamma_{\mot}(k,M(2))
 \right).
\]
The exact sequence with the absolute trace fiber gives
\[
 H^1\mathcal C_{\mathrm{abs}}(F,R)
 \xrightarrow{\sim}\prod_{v\,\mathrm{real}}H^1\mathcal C_{\mathrm{abs}}(F_v^{\mathrm{rc}},R),
 \qquad
 H^2\mathcal C_{\mathrm{abs}}(F,R)
 \hookrightarrow\prod_{v\,\mathrm{real}}H^2\mathcal C_{\mathrm{abs}}(F_v^{\mathrm{rc}},R).
\]
By Lemma~\ref{lem:integral-structural-push} with $r=2$, the structure map induces a surjection $\CH^{n+2}(X,3)\twoheadrightarrow\CH^2(F,3;D)$ whose composite with $\rho$ is $f_*$.
The same holds over $F_v^{\mathrm{rc}}$ and over finite extensions of $F$, and these surjections commute with restriction and transfer.
Over a real closure, the index is $1$ or $2$; after passing by Morita equivalence to degree $2$ in the second case, $\rho$ is an isomorphism in weight two \cite[Theorems 6.5.3 and 8.1.4(1)]{KahnLevine2010}.
Hence $f_*$ is surjective onto $H^1_{\mot}(F_v^{\mathrm{rc}},R(2))$, so the connecting map $H^1_{\mot}(F_v^{\mathrm{rc}},R(2))\to H^2\mathcal C_{\mathrm{abs}}(F_v^{\mathrm{rc}},R)$ vanishes.
By naturality and the injectivity in $H^2$ displayed above, the connecting map $H^1_{\mot}(F,R(2))\to H^2\mathcal C_{\mathrm{abs}}(F,R)$ also vanishes.
Thus $f_*$ is surjective over $F$, and since $f_*$ factors through $\rho$, so is $\rho$ onto $K_3^{\mathrm{ind}}(F)_{(2)}$.
For injectivity, lift a kernel element in weight two along the surjection of Lemma~\ref{lem:integral-structural-push} to $H^{2n+1}_{\mot}(X,R(n+2))$ and then to an element $y$ of $H^1\mathcal C_{\mathrm{abs}}(F,R)$.
Choose $k=F(\alpha)$ where the minimal polynomial of $\alpha$ is a monic polynomial of degree $\deg D$ chosen by weak approximation: at each finite place where $D$ is ramified, close to an irreducible polynomial defining the unramified extension of degree $\deg D$ (Krasner's lemma), and at each real place, close to a polynomial with distinct real roots (simple roots form an open condition); then all real places split completely in $k$, and the local degrees at the finite ramified places are divisible by the index.
If no finite place is ramified, then $\ind D\le2$ and $k=F$ suffices.
Brauer--Hasse--Noether gives $\ind(D_k)\le2$.
Place the component of $y$ at each real place on one real branch of $k$ and use zero on the others.
Since the relative computation applies over $k$ to $X_k=\SB(D_k)$, the isomorphism in $H^1$ onto the real components, applied over $k$, lifts the resulting tuple to $y_k\in H^1\mathcal C_{\mathrm{abs}}(k,R)$.
At each real place $v$ of $F$, corestriction adds the components over the real places of $k$ above $v$, so $\cor_{k/F}y_k$ and $y$ have the same real components, and the injectivity in $H^1$ gives $\cor_{k/F}y_k=y$.
Since $\rho$ is an isomorphism at index at most two, after passing by Morita equivalence to degree $2$ at index two \cite[Theorems 6.5.3 and 8.1.4(1)]{KahnLevine2010}, $y_k$ has zero image in the twisted weight-two group.
Transferring back shows that the original kernel element is zero.
\end{proof}

\subsection{The positive-degree comparison and its fiber over global fields}\label{sec:global-spectrum}
In this subsection, we construct the comparison $g\colon\tau_{\ge1}K(D)\to\tau_{\ge1}K(F)$ over a global field (Theorem~\ref{thm:paper-global-comparison}) and study its fiber $Z$.
We split $g$ in degrees at least two, first over number fields and then over global function fields, and compute $K_3(D)$.
We then show that $Z$ has finitely many primary factors with bounded homotopy (Theorem~\ref{thm:paper-global-fiber-bounds}), that, up to equivalence, it depends only on the cyclic subgroup of the Brauer group generated by $[D]$ (Proposition~\ref{prop:paper-brauer-power-fiber}), and that over number fields it is noncanonically equivalent to a finite sum of shifted Eilenberg--Mac Lane spectra (Theorem~\ref{thm:paper-global-fiber-gem}).
When $D$ is ramified at a real place, the statements about the comparison, its splitting, and its fiber are made after inverting $2$; the computation of $K_3(D)$ is integral.

\subsubsection{Construction of the positive-degree comparison}
The comparison is glued in the arithmetic square from the inverse of multiplication by a class of geometric rank $1$ at each prime $\ell\ne\operatorname{char}F$ dividing $\deg D$ (Proposition~\ref{prop:global-rank}), from $\Nrd_1$ at $p=\operatorname{char}F$ when $p\mid\deg D$, and from $N_{D/F}/\deg D$ after inverting $\deg D$.

Let $F$ be a global field.
In the construction below we assume that $D$ splits at all real places; the real-ramified case is treated after inverting $2$ in Theorem~\ref{thm:paper-global-comparison}.
As in Proposition~\ref{prop:global-real-weight-two}, $H^i_{\mathrm{rel}}$ denotes the cohomology of the fiber of restriction to the algebraic real closures at all real places.
Fix a prime $\ell\ne\operatorname{char}F$.
In \'etale descent, the geometric stalks with coefficients modulo $\ell^r$ are
\[
 \Z/\ell^r(n)\ \text{in degree }2n,\qquad 0\ \text{in odd degrees}.
\]
Given a compatible geometric rank-$1$ class $x\in\pi_0L_{K(1),\ell}K(D)$, multiplication by $x$ is an equivalence
\[
 m_x\colon L_{K(1),\ell}K(F)\xrightarrow{\ \sim\ }L_{K(1),\ell}K(D)
\]
by Lemma~\ref{lem:paper-etale-rank-one} (see \cite[Theorems 6.18 and 7.14]{ClausenMathew}).

\begin{proposition}\label{prop:global-rank}
Let $F$ be a global field, let $\ell\ne\operatorname{char}F$ be a prime, and assume that $D$ splits at all real places of $F$.
Then one can choose a geometric rank-$1$ class $x\in\pi_0L_{K(1),\ell}K(D)$.
For a number field with real places and for $\ell=2$, one can prescribe the rank-$1$ Morita class at every real closure.
\end{proposition}
\begin{proof}
We have $\cd_\ell F\le2$ unless $\ell=2$ and $F$ has a real place; this exceptional case is the real-relative case treated below.

Suppose first that $\cd_\ell F\le2$.
Then the rank lifts along the edge map in total degree zero modulo $\ell^r$, and the difference of two lifts lies in $H^2(F,\Z/\ell^r(1))$.
The coefficient transition maps
\[
 H^2(F,\Z/\ell^{r+1}(1))\longrightarrow H^2(F,\Z/\ell^r(1))
\]
are surjective, so we can choose compatible lifts.

In the real-relative case, the obstruction to lifting the prescribed real elements lies in total degree $-1$ of the relative groups, namely in $H^1_{\mathrm{rel}}$ in weight zero.
It is the relative boundary of the tuple
\[
 (1)_v\in\prod_{v\ \mathrm{real}}H^0(F_v^{\mathrm{rc}},\Z/2^r),
\]
which is the restriction of the diagonal $1$, so this boundary vanishes.
The differences of lifts lie in relative $H^2$ in weight one, whose transition maps are surjective because relative $H^3$ vanishes.
The Milnor exact sequence for the homotopy inverse limit therefore supplies $x$.
Since the corresponding groups on the real side modulo $2^r$ are finite, the restriction of $x$ is still the prescribed one in the limit.
\end{proof}

We also write $N_{D/F}\colon L_{K(1),\ell}K(D)\to L_{K(1),\ell}K(F)$ for the map induced by the norm, and $m_{N_{D/F}(x)}$ for multiplication by $N_{D/F}(x)\in\pi_0L_{K(1),\ell}K(F)$.

\begin{lemma}\label{lem:paper-global-product-vanishing}
Let $\ell\ne\operatorname{char}F$, and let $x\in\pi_0L_{K(1),\ell}K(D)$ be a class of geometric rank $1$ as in Proposition~\ref{prop:global-rank}; if $\ell=2$ and $F$ has real places, we assume that $x$ is the rank-$1$ Morita class at every real closure.
Then
\[
 \bigl(N_{D/F}(x)-\deg D\bigr)\,y=0\qquad\text{for }y\in\pi_iL_{K(1),\ell}K(F),\ i\ge1,
\]
so $N_{D/F}\circ m_x$ acts as multiplication by $\deg D$ on $\pi_i$ for $i\ge1$.
Moreover, for every finite separable splitting field $E/F$, multiplication by $x_E-1$ vanishes on $\pi_iL_{K(1),\ell}K(E)$ for $i\ge1$, where $x_E=\res_{E/F}x$ and $K(D_E)$ is identified with $K(E)$ by Morita equivalence.
\end{lemma}
\begin{proof}
By the projection formula,
\[
 N_{D/F}\circ m_x\simeq m_{N_{D/F}(x)}.
\]
The classes $N_{D/F}(x)-\deg D\in\pi_0L_{K(1),\ell}K(F)$ and $x_E-1\in\pi_0L_{K(1),\ell}K(E)$ have geometric rank $0$.
In the real-relative case they restrict to zero at every real closure, because $x$, and hence $x_E$, is the rank-$1$ Morita class there.
Lemma~\ref{lem:paper-rank-zero-product}, applied over $F$ and over $E$, gives both assertions.
\end{proof}

\begin{lemma}\label{lem:paper-global-no-divisible}
Let $F$ be a global field and $j\ge1$.
Then $K_j(F)$ has no nonzero divisible subgroup.
The same holds for $K_j(F)_{(\ell)}$ for every prime $\ell$ and for $K_j(F)[1/2]$.
\end{lemma}
\begin{proof}
In degree one, we use the map to the direct sum of the valuations at the finite places or closed points; its kernel is the finitely generated unit group.
In higher degrees, we use localization for the ring of integers or the smooth proper curve: the direct sum of the groups $K_{j-1}$ of the residue fields has no nonzero divisible subgroup, so any divisible subgroup lies in the image of the finitely generated group $K_j$ of the model \cite{QuillenHigher,QuillenCurve,QuillenFinite}, \cite[Theorem 1]{QuillenFiniteGeneration}.
After localization at $\ell$ or after inverting $2$, the direct sum of the value groups or of the residue $K$-groups still has no nonzero divisible subgroup, and the rest lies in a finitely generated $\Z_{(\ell)}$-module or $\Z[1/2]$-module, so the same argument applies.
\end{proof}

\begin{theorem}\label{thm:paper-global-comparison}
If $D$ splits at all real places of $F$, there is an integral map of positive-degree spectra
\[
 g\colon\tau_{\ge1}K(D)\longrightarrow\tau_{\ge1}K(F)
\]
that carries the splitting transfer $\theta_E$ from each finite separable splitting field $E/F$ to $N_{E/F}$ on every $\pi_j$.
In degrees $1$ and $2$, the map $g$ agrees with the reduced norms $\Nrd_1$ and $\Nrd_2$.
With real ramification the same construction works after inverting $2$.
\end{theorem}

For $j\ge1$, we write $g_j=\pi_jg$.

\begin{proof}
We first define the components of $g$ in the arithmetic square.
For $\ell\mid\deg D$ with $\ell\ne\operatorname{char}F$, we use the composite
\[
 K(D)\longrightarrow L_{K(1),\ell}K(D)\xrightarrow{\ m_x^{-1}\ }L_{K(1),\ell}K(F),
\]
together with the Quillen--Lichtenbaum comparison $\tau_{\ge1}L_{K(1),\ell}K(F)\simeq(\tau_{\ge1}K(F))^\wedge_\ell$ for $F$.
If $p=\operatorname{char}F$ divides $\deg D$, then $[F:F^p]=p$, and the rank-$1$ Frobenius comparison of Proposition~\ref{prop:app-frobenius} gives
\[
 K_i(D;\Z/p)=K_i(F;\Z/p)=0\qquad(i\ge2).
\]
Thus the completed positive-degree spectrum at $p$ is concentrated in degree one, where we use the global reduced norm $\Nrd_1$.
After inverting $\deg D$, we use $N_{D/F}/\deg D$.

These components agree on the overlap.
For $\ell\ne\operatorname{char}F$, Lemma~\ref{lem:paper-global-product-vanishing} gives $N_{D/F}m_x=\deg D$ in positive degrees, so
\[
 m_x^{-1}=N_{D/F}/\deg D
\]
on the rationalized positive homotopy groups of the $\ell$-completed factor.
If $p=\operatorname{char}F$ divides $\deg D$, then for every finite separable splitting field $E/F$
\[
 N_{D/F}\theta_E=\deg D\cdot N_{E/F}=\deg D\cdot\Nrd_1\theta_E,
\]
and the splitting transfers generate $K_1(D)$, as shown below, so $N_{D/F}=\deg D\cdot\Nrd_1$ on $K_1(D)$.
Since the overlap of the arithmetic square is rational, maps into it that agree on homotopy groups are homotopic, so we can glue the completed maps and $N_{D/F}/\deg D$ to obtain $g$.
With real ramification, after inverting $2$ only odd primes $\ell$ occur, and then $\cd_\ell F\le2$.
For odd $\ell$, the rank-$1$ lift, the multiplicative equivalence and the vanishing of the products in Lemma~\ref{lem:paper-global-product-vanishing} follow by the same argument from absolute \'etale descent; the geometric rank argument uses no prescription at the real places.

Next we check the transfer identity on the components.
For a finite separable splitting field $E$, we identify $K(D_E)$ with $K(E)$ by the Morita equivalence in the definition of $\theta_E$, so that $\theta_E=N_{D_E/D}$.
On the completed factors, the projection formula of Lemma~\ref{lem:paper-etale-rank-one} gives $N_{D_E/D}m_{x_E}\simeq m_xN_{E/F}$, and hence
\[
 m_x^{-1}\theta_E=N_{E/F}m_{x_E}^{-1}.
\]
By Lemma~\ref{lem:paper-global-product-vanishing}, multiplication by $x_E-1$ vanishes in positive degrees over $E$, so the right-hand side induces the field norm.
For the only positive degree of the completed characteristic-$p$ factor, the transfer identity for $\Nrd_1$ gives the same conclusion.
After inverting $\deg D$, we use $N_{D/F}\theta_E=\deg D\cdot N_{E/F}$.
Consequently, on integral homotopy groups $g_j\theta_E-N_{E/F}$ takes values in the divisible image of the boundary from the rational overlap of the arithmetic square.
This image vanishes by Lemma~\ref{lem:paper-global-no-divisible}, applied to $K_j(F)$ or, with real ramification, to $K_j(F)[1/2]$, so the transfer identity holds with the coefficients of the theorem.

Finally we identify $g_1$ and $g_2$.
In degree one, the units generating separable maximal subfields form an open subset, and every unit is the ratio of two elements of this subset; consequently the splitting transfers generate $K_1(D)$.
In degree two, the proof of Theorem~\ref{thm:paper-global-norm} shows that $K_2(D)$ is generated by transfer from one separable splitting field.
Since the reduced norms satisfy the same transfer identities $\Nrd_j\theta_E=N_{E/F}$ for $j=1,2$, we obtain $g_j=\Nrd_j$.
\end{proof}

\begin{lemma}\label{lem:paper-global-low-degrees}
Let $F$ be a global field, and let $g$ be the comparison of Theorem~\ref{thm:paper-global-comparison}.
If $D$ splits at all real places of $F$, then $g_1$ and $g_2$ are isomorphisms; otherwise they are isomorphisms after inverting $2$.
\end{lemma}
\begin{proof}
By Theorem~\ref{thm:paper-global-comparison}, $g_j=\Nrd_j$ for $j=1,2$.
Over a global function field, $\SK_1(D)=0$ by \cite[Corollary 2]{Yanchevskii1975}, whose proof also gives surjectivity of $\Nrd_1$.
Over a number field, $\Nrd_1$ is injective because $\SK_1(D)=0$ \cite{Wang1950}.
By the Hasse--Schilling--Maass norm theorem \cite[Theorem 33.15]{ReinerMaximalOrders}, its image consists of the elements of $F^\times$ that are positive at every real place where $D$ is ramified.
Hence $\Nrd_1$ is surjective if $D$ splits at all real places.
If $D$ is ramified at a real place, the cokernel of $\Nrd_1$ is killed by $2$, since squares are positive at every real place, so $\Nrd_1$ is an isomorphism only after inverting $2$.
By Theorem~\ref{thm:paper-global-norm}, $\Nrd_2$ is injective and its cokernel is $(\Z/2)^{S_{\mathbb R}}$, where $S_{\mathbb R}$ is the set of real places where $D$ is ramified.
Thus $\Nrd_2$ is an isomorphism if $D$ splits at all real places, and after inverting $2$ in general.
\end{proof}

Over number fields, the comparisons of Theorem~\ref{thm:paper-global-comparison} can be chosen coherently over the intermediate fields of a finite Galois extension, as in the local case (Proposition~\ref{prop:paper-finite-galois-naturality}).

\begin{proposition}\label{prop:paper-finite-galois-naturality-global}
Let $F$ be a number field, let $D$ be a central division algebra over $F$, and let $L/F$ be a finite Galois extension.
If $D$ is ramified at some real place, we invert $2$ throughout.
Then the comparisons
\[
 g_{F'}\colon\tau_{\ge1}K(D_{F'})\longrightarrow\tau_{\ge1}K(F')
\]
of Theorem~\ref{thm:paper-global-comparison} for the intermediate fields $F'$ of $L/F$ can be chosen simultaneously and coherently with respect to restriction, transfer, conjugation, and their composites.
\end{proposition}
\begin{proof}
Let $G=\Gal(L/F)$, and define $X,Y\in\mathscr M_G$ by $X(G/H)=\tau_{\ge1}K(D_{L^H})$ and $Y(G/H)=\tau_{\ge1}K(L^H)$, as in the proof of Proposition~\ref{prop:paper-finite-galois-naturality}.
As there, we use the norms $N_{D_{F'}/F'}$ of the algebras $D_{F'}$, which have degree $\deg D$ even when the index drops over $F'$, and the arithmetic square in $\mathscr M_G$ for the primes dividing $\deg D$.
Only if $D$ is ramified at some real place do we replace $X$ and $Y$ by their localizations away from two at the outset; then the two-completed factor is zero, and we omit it in choosing rank-$1$ classes.

For each remaining $\ell\mid\deg D$, choose a class $x_\ell$ of geometric rank $1$ over $F$ and restrict it to every intermediate field.
If $D$ splits at all real places, the class is supplied by Proposition~\ref{prop:global-rank}, and Lemma~\ref{lem:paper-etale-rank-one} applies; if moreover $\ell=2$, we choose $x_2$ with the rank-$1$ Morita class prescribed at every real closure, as Proposition~\ref{prop:global-rank} allows.
This prescription is inherited by the real closures of the intermediate fields and of the splitting fields.
Under real ramification only odd $\ell$ remain; there we do not apply Proposition~\ref{prop:global-rank} or Lemma~\ref{lem:paper-etale-rank-one} directly, which assume that $D$ splits at the real places, but lift the rank-$1$ class by absolute \'etale descent as in the proof of Theorem~\ref{thm:paper-global-comparison}, and construct the equivalence by the geometric-point and free-module arguments in the proof of Lemma~\ref{lem:paper-etale-rank-one}.
On each orbit, $u_\ell$ (as in the proof of Theorem~\ref{thm:paper-local-spectrum}) is the composite of the natural localization unit on the twisted side, $m_{x_\ell}^{-1}$, and the inverse of the Quillen--Lichtenbaum comparison on the field side; the latter is natural and an equivalence on each orbit, so it can be inverted in the Mackey category.
The gluing then follows the argument of Proposition~\ref{prop:paper-finite-galois-naturality}.
Here the product calculation in the proof of Theorem~\ref{thm:paper-global-comparison} identifies the rationalized $u_\ell$ with $N/\deg D$ on positive homotopy groups over every intermediate field.
We obtain $g\colon X\to Y$ from the pullback.

For number fields the correction of the local case is unnecessary.
Once the arithmetic components are fixed, the maps on homotopy groups in degrees $j\ge1$ are determined up to the boundary image of the rational group $\pi_{j+1}\bigl(Y(F')^\wedge_{\deg D}[1/\deg D]\bigr)$, and this image is divisible and hence zero by Lemma~\ref{lem:paper-global-no-divisible}, applied to $K_j(F')$ or, under real ramification, to $K_j(F')[1/2]$.
For each intermediate field $F'$ and each finite separable splitting field $E/F'$, the projection formula and the vanishing of products with classes of geometric rank $0$, as in the proof of Theorem~\ref{thm:paper-global-comparison}, give the transfer identity on the completed components.
After inverting $\deg D$, the identity follows from $\pi_j(N_{D_{F'}/F'}\circ\theta_E)=\deg D\cdot\pi_jN_{E/F'}$.
The vanishing of the boundary image therefore gives the transfer identity on homotopy groups with the stated coefficients.
Applying the generation by splitting transfers, as in the proof of Theorem~\ref{thm:paper-global-comparison}, to the division algebra representing $D_{F'}$ over each intermediate field gives the degree-one and degree-two normalizations.
Hence $g_{F'}$ satisfies the transfer identities and the degree-one and degree-two normalizations of Theorem~\ref{thm:paper-global-comparison}.
\end{proof}

\subsubsection{Derived norms and splitting over number fields}
Over number fields, the splitting of Theorem~\ref{thm:paper-global-split} uses the derived norm equivalence of Proposition~\ref{prop:paper-derived-norm} for a finite Galois extension $E/F$, with every real place of $F$ split completely in $E$ unless $2$ is inverted.

For a finite Galois extension $E/F$ with group $G$, the field norm defines a map from homotopy orbits using the canonical isomorphism between restriction of scalars and its composite with Galois conjugation.
We include the coherence of these isomorphisms in the action defining the orbits.

We use \'etale cohomology of rings of $S$-integers, as follows.
Let $E/F$ be a finite Galois extension of number fields with group $G$.
Fix $\ell$, take a finite set $S$ of places containing the places over $\ell$ and those ramified in $E/F$, and put
\[
 V_L=\Spec\mathcal O_{L,S}\qquad(L\in\{E,F\}).
\]
Then $V_E\to V_F$ is a finite \'etale $G$-cover.
For a finite $\ell$-primary local system $M$, we use ordinary $R\Gamma$ if $\ell$ is odd or $F$ is totally imaginary.
If $\ell=2$ and $F$ has real places, we use
\[
 R\Gamma_{\mathrm{rel}}(V_L,M)=
 \fib\!\left(R\Gamma(V_L,M)\longrightarrow
              \bigoplus_{v\ \mathrm{real}}R\Gamma(L_v,M)\right).
\]
The real terms are ordinary group cohomology of $\Z/2$.
We write $R\Gamma_\star(V_L,-)$ for the chosen functor.
The cohomology of $R\Gamma_\star(V_L,M)$ lies in $[0,2]$ for finite $\ell$-primary local systems because $S$ contains the places above $\ell$.
For $R\Gamma_{\mathrm{rel}}(V_L,M)$, this bound uses the isomorphism of real restriction in degrees at least three and surjectivity in degree two.
Surjectivity in degree two follows from Poitou--Tate duality: choose a place $v\mid\ell$ of $L$; local duality makes $H^2(L_v,M)\to H^0(V_L,M^D)^\vee$ surjective, where $M^D=\operatorname{Hom}(M,\mu_{\ell^\infty})$, so the obstruction for any tuple at the real places can be cancelled at that finite place.
For continuous $\Z_\ell$ coefficients, finiteness at each finite stage gives the Mittag--Leffler condition and preserves this bound \cite{KahnPositive,ChinburgEtal}.

\begin{lemma}\label{lem:paper-derived-tate-norm}
Let $M$ be a finite local system of $\Z/\ell^a$-modules on $V_F$, and let $\omega\colon V_E\to V_F$ be the cover.
In the real-relative case, assume that every real place of $F$ splits completely in $E$.
Then the trace induces an equivalence
\begin{equation}\label{eq:paper-derived-tate-norm}
  \Z/\ell^a\otimes_{\Z/\ell^a[G]}^{\mathbf L} R\Gamma_\star(V_E,\omega^*M)
      \xrightarrow{\ \sim\ } R\Gamma_\star(V_F,M).
\end{equation}
\end{lemma}
\begin{proof}
For the cover $\omega$, the stalks of $\omega_*\omega^*M$ are induced $G$-modules with derived coinvariants $M$ via trace.
This is the tensor product calculation with the free group-ring resolution for a regular permutation action.
We have $R\Gamma(V_F,\omega_*\omega^*M)=R\Gamma(V_E,\omega^*M)$, compatibly with restriction to real places in the real-relative case because they split completely.
A free resolution then gives~\eqref{eq:paper-derived-tate-norm}: finite truncations of the bar resolution have finite free terms, so sections commute with the finite sums and cofibers.
The omitted tail contributes nothing to a fixed degree past a sufficiently late stage by the cohomological upper bound two.
The same reasoning applies to relative sections.
\end{proof}

\begin{proposition}\label{prop:paper-derived-norm}
For a finite Galois extension $E/F$ of number fields in which every real place splits completely, the norm
\[
  \mathcal N_{E/F}:\bigl(\tau_{\ge2}K(E)\bigr)_{hG}\longrightarrow \tau_{\ge2}K(F)
\]
is an equivalence.
Without complete splitting at the real places the same holds after inverting $2$.
\end{proposition}
\begin{proof}
Fix $\ell$ and $S$ as above.
The coefficient triangles comparing continuous Tate coefficients with reduction modulo $\ell^a$ are compatible with sections.
The equivalences of Lemma~\ref{lem:paper-derived-tate-norm} modulo $\ell^a$ therefore give the corresponding equivalence in the derived $\ell$-complete category.

We now pass to \'etale $K$-theory.
Write $\mathcal{KU}_\ell$ for the periodic \'etale $K$-theory sheaf, and consider the geometric Postnikov tower
\[
  \mathcal P_L^m=R\Gamma_\star(V_L,\tau_{\ge2m}\mathcal{KU}_\ell),
  \qquad
  \gr^m\mathcal P_L \simeq R\Gamma_\star(V_L,\Z_\ell(m))[2m].
\]
We interpret the continuous layers by inverse limits of their finite-coefficient versions.
The connectivity bound for $\mathcal P_L^m$ is $2m-2$.
Since $H^0(V_L,\Z_\ell(1))=0$ and relative $H^0$ vanishes as well, weights at most one make no contribution in degrees at least two.
Hence $\mathcal P_L^2\simeq\tau_{\ge2}R\Gamma_\star(V_L,\mathcal{KU}_\ell)$.
Homotopy orbits preserve the lower connectivity bound of the tails, so arbitrarily high layers cannot contribute in any fixed degree \cite{ClausenMathew}.
Applying Lemma~\ref{lem:paper-derived-tate-norm} to the finite-coefficient versions of the layers and inducting on finite truncations, we obtain $(\mathcal P_E^2)_{hG}\simeq\mathcal P_F^2$ in the derived $\ell$-complete category.

Next we return to algebraic $K$-theory.
For ordinary sections, we use the $\ell$-complete comparison $K(V_L)\to K^{\et}(V_L)$ in degrees at least two.
In the real-relative case, we restore the real terms.
After $S$ is enlarged, weak approximation supplies $S$-units $u$ realizing all real sign vectors.
The symbols $\{-1,u\}\in K_2(V_L)$ hit the generators of $\pi_2(K(L_v)^\wedge_2)\simeq\pi_2(KO^\wedge_2)=\Z/2$ at the specified real places, so the restriction $K(V_L)^\wedge_2\to\bigoplus_{v\ \mathrm{real}}KO^\wedge_2$ is surjective on $\pi_2$.
In the commutative square formed by this restriction, the restriction $R\Gamma(V_L,\mathcal{KU}_2)\to\bigoplus_{v\ \mathrm{real}}KO^\wedge_2$, and the comparison maps, the vertical maps are isomorphisms in degrees at least two, so the lower restriction is also surjective on $\pi_2$.
Hence the relative, absolute, and real terms still form fiber sequences after taking $2$-connective covers, and the algebraic relative term is identified with $\tau_{\ge2}R\Gamma_{\mathrm{rel}}(V_L,\mathcal{KU}_2)$ \cite[Theorems 7.12--7.13]{ClausenMathew}.
On $V_E$ the real terms are a finite sum of induced $G$-spectra and the derived norm is an equivalence on them.
The equivalence on the relative terms then gives it on the absolute terms \cite{RognesWeibel,ClausenMathew}.

Finally, we reduce these equivalences modulo $\ell$ and pass to the filtered colimit over $S$.
The $K$-theory, connective covers, orbits, and reduction modulo $\ell$ all commute with this colimit, giving an equivalence over the fields modulo every $\ell$.
Rationally the norm is an equivalence by the restriction-transfer identities and averaging.
The integral equivalence follows from the vanishing of the cofiber rationally and modulo each prime.
After two is inverted, only the odd primes and the rational comparison are needed.
\end{proof}

We next construct the transfer from homotopy orbits that appears in the splitting theorems.
Let $E/F$ be a finite Galois splitting field of $D$ with group $G$.
The group $G$ acts on the category of right $D_E$-modules by semilinear scalar change fixing $D$.
Restriction of scalars from $D_E$ to $D$ has natural isomorphisms with its composites with this action.
These isomorphisms preserve the underlying $D$-module and are coherent with group multiplication, so they give the transfer
\[
 \mathcal N^{D}_{E/F}:\bigl(\tau_{\ge2}K(D_E)\bigr)_{hG}\longrightarrow \tau_{\ge2}K(D)
\]
from homotopy orbits.
For $j\ge2$, its composite with the edge map $K_j(D_E)\to\pi_j\bigl((\tau_{\ge2}K(D_E))_{hG}\bigr)$ is restriction of scalars, that is, the splitting transfer $\theta_E$ under the identification $K(D_E)\simeq K(E)$ by Morita equivalence.

In both splitting theorems, we obtain a section of $\tau_{\ge2}g$ once its composite with the transfer from homotopy orbits is an equivalence.

\begin{lemma}\label{lem:paper-split-from-norm}
Let $F$ be a global field, let $g$ be the comparison of Theorem~\ref{thm:paper-global-comparison}, put $Z=\fib(g)$, and let $E/F$ be a finite Galois splitting field of $D$ with group $G$.
If $D$ is ramified at a real place, invert $2$ throughout.
Assume that $h=(\tau_{\ge2}g)\circ\mathcal N^{D}_{E/F}$ is an equivalence.
Then $\mathcal N^{D}_{E/F}\circ h^{-1}$ is a section of $\tau_{\ge2}g$, so that
\[
 \tau_{\ge2}K(D)\simeq\tau_{\ge2}K(F)\oplus\fib(\tau_{\ge2}g)
 \quad\text{with first component }\tau_{\ge2}g,
 \qquad
 \cofib(\mathcal N^{D}_{E/F})\simeq\fib(\tau_{\ge2}g),
\]
and $g_j$ is split surjective for every $j\ge2$.
If moreover $g_1$ is an isomorphism and $g_2$ is surjective, then $\fib(\tau_{\ge2}g)\simeq Z$, and $\pi_jZ=\ker g_j$ for $j\ge2$.
\end{lemma}
\begin{proof}
Since $(\tau_{\ge2}g)\circ\mathcal N^{D}_{E/F}\circ h^{-1}=hh^{-1}$ is the identity, the map $s=\mathcal N^{D}_{E/F}\circ h^{-1}$ is a section of $\tau_{\ge2}g$.
The section gives the direct-sum decomposition with first component $\tau_{\ge2}g$, and $g_j$ is split surjective for $j\ge2$.
Since $\mathcal N^{D}_{E/F}=s\circ h$ and $h$ is an equivalence, $\cofib(\mathcal N^{D}_{E/F})\simeq\cofib(s)$, and the cofiber of the section $s$ of $\tau_{\ge2}g$ is $\fib(\tau_{\ge2}g)$.
If $g_1$ is an isomorphism and $g_2$ is surjective, the long exact sequence of homotopy groups for $Z=\fib(g)$ gives the vanishing of $\pi_0Z$ and $\pi_1Z$, and the natural map $\fib(\tau_{\ge2}g)\to Z$ is an equivalence.
\end{proof}

We also need splitting fields in which prescribed real places split completely.

\begin{lemma}\label{lem:paper-galois-splitting-real}
Let $F$ be a global field, and let $S$ be a set of real places of $F$ at which $D$ splits.
Then there is a finite Galois splitting field $E/F$ of $D$ in which every place in $S$ splits completely.
\end{lemma}
\begin{proof}
If $S$ is empty, the Galois closure of a separable maximal subfield of $D$ suffices.
Otherwise $F$ is a number field, and we use weak approximation in $D$ to choose an element $\alpha$ that is close, at each place in $S$, to a matrix with distinct real eigenvalues.
The reduced characteristic polynomial of $\alpha$ is separable and has only real roots at each place in $S$.
Since $D$ is a division algebra, $F(\alpha)$ is a field, and it is a maximal subfield of degree $\deg D$ in which every place in $S$ splits completely.
The Galois closure of $F(\alpha)$ has the same property.
\end{proof}

\begin{theorem}\label{thm:paper-global-split}
Let $F$ be a number field and assume that $D$ splits at all real places, or else invert $2$ throughout.
Let $g$ be the comparison of Theorem~\ref{thm:paper-global-comparison}, and put $Z=\fib(g)$.
Let $E/F$ be a finite Galois splitting field of $D$ with group $G$, and if $D$ splits at all real places, assume that every real place of $F$ splits completely in $E$.
Then the composite $(\tau_{\ge2}g)\,\mathcal N^{D}_{E/F}$ is an equivalence.
Consequently,
\[
 \tau_{\ge2}K(D)\simeq\tau_{\ge2}K(F)\oplus Z
 \quad\text{with first component }\tau_{\ge2}g,
 \qquad Z\simeq\cofib(\mathcal N^{D}_{E/F}).
\]
For every $j\ge2$, $g_j$ is split surjective and $\pi_jZ=\ker g_j$.
In particular $Z\in\mathrm{Sp}_{\ge3}$.
\end{theorem}
\begin{proof}
Such a field $E$ exists by Lemma~\ref{lem:paper-galois-splitting-real}, applied to the set of all real places if $D$ splits at all of them and to the empty set otherwise.

We first compare the norms on the completed factors.
For $\ell\mid \deg D$, we restrict to $E$ the rank-$1$ class $x$ used in the comparison.
Since $x_E$ is restricted from $F$, Lemma~\ref{lem:paper-etale-rank-one} gives a coherently $G$-equivariant equivalence
\[
 m_{x_E}\colon L_{K(1),\ell}K(E)\xrightarrow{\ \sim\ }L_{K(1),\ell}K(D_E).
\]
Since $D_E$ is split, Morita equivalence and the Quillen--Lichtenbaum comparison for $E$ identify its $2$-connective cover with an equivariant equivalence
\[
 \alpha_E\colon(\tau_{\ge2}K(E))^\wedge_\ell\xrightarrow{\ \sim\ }(\tau_{\ge2}K(D_E))^\wedge_\ell.
\]
The projection formula of Lemma~\ref{lem:paper-etale-rank-one} gives
\[
   (\tau_{\ge2}g)^\wedge_\ell\,(\mathcal N^{D}_{E/F})^\wedge_\ell\,
           (\alpha_E)_{hG}=(\mathcal N_{E/F})^\wedge_\ell
\]
as maps of $K(F)^\wedge_\ell$-modules in the derived $\ell$-complete category.
Here the homotopy orbits are taken in the derived $\ell$-complete category, and reducing them modulo $\ell$ gives the same result as taking homotopy orbits of the reductions modulo $\ell$ before completion.
This applies at odd primes in the ramified real case as well.

After inverting $\deg D$, multiplication by the restriction to $E$ of $[D]/\deg D\in K_0(D)[1/\deg D]$ is a $G$-equivariant rank-$1$ equivalence $K(E)[1/\deg D]\to K(D_E)[1/\deg D]$.
By the projection formula, the homotopy orbits of this equivalence composed with $\mathcal N^{D}_{E/F}$ and then with $N_{D/F}/\deg D$ give $\mathcal N_{E/F}$.
Proposition~\ref{prop:paper-derived-norm} now proves that $(\tau_{\ge2}g)\mathcal N^{D}_{E/F}$ is an equivalence rationally and modulo all relevant primes, and hence is an equivalence with the stated coefficients.

Finally, $g_1$ and $g_2$ are isomorphisms by Lemma~\ref{lem:paper-global-low-degrees}, so Lemma~\ref{lem:paper-split-from-norm} yields both formulas and $\pi_jZ=\ker g_j$ for $j\ge2$.
Since $g_2$ is injective, $\pi_2Z=0$, and hence $Z\in\mathrm{Sp}_{\ge3}$.
On each completed factor the projection formula also gives the splitting as $K(F)^\wedge_\ell$-modules.
\end{proof}

\subsubsection{Splitting over global function fields by constant field transfer}
Over global function fields, we obtain the corresponding splitting (Theorem~\ref{thm:paper-global-function-split}) by transfer from a constant field extension.

Let $F$ be a global function field with constant field $k=\F_q$.
Let $f_v$ be the residue degree of a place $v$, and choose a finite constant extension $E=Fk'$ such that $s=[k':k]$ is a common multiple of the integers $\deg D\cdot f_v$ for the places $v$ where the local invariant of $D$ is nonzero.
At each such $v$, the local degree $s/\gcd(s,f_v)$ is divisible by $\deg D$, so $D_E$ has zero local invariants and splits by the local--global principle for the Brauer group.

\begin{lemma}\label{lem:paper-global-constant-norm}
For any finite constant extension $E/F$ with group $G=\Gal(E/F)$ and every $j\ge2$,
\[
 H_i(G,K_j(E))=0\quad(i>0),\qquad
 N_{E/F}:K_j(E)_G\xrightarrow{\sim}K_j(F).
\]
\end{lemma}
\begin{proof}
Let $\ell\ne p=\operatorname{char}k$.
We first show $H^2(F,\Q_\ell/\Z_\ell(m))=0$ for $m\ge2$.
Let $C/k$ be the smooth proper geometrically connected model of $F$.
By the Weil weights \cite{Deligne1974}, Frobenius minus one on $H^1(C_{\bar k},\Z_\ell(m))$ is invertible on the rationalized lattice.
Thus it is surjective on $H^1(C_{\bar k},\Q_\ell/\Z_\ell(m))$.
The finite field has cohomological dimension one, so the Hochschild--Serre edge map gives
\[
 H^2(C,\Q_\ell/\Z_\ell(m))\xrightarrow{\sim}H^0(k,H^2(C_{\bar k},\Q_\ell/\Z_\ell(m)))
 =H^0(k,\Q_\ell/\Z_\ell(m-1)).
\]
The Gysin map of a closed point $v$ becomes corestriction
\[
 H^0(k_v,\Q_\ell/\Z_\ell(m-1))\longrightarrow H^0(k,\Q_\ell/\Z_\ell(m-1)).
\]
Indeed, the trace of a geometric point class is one, and summing over conjugate points is corestriction.
For the arithmetic Frobenius $\Phi$ and $f=[k_v:k]$, corestriction is $1+\Phi+\cdots+\Phi^{f-1}$, which on $\Q_\ell/\Z_\ell(m-1)$ is multiplication by a nonzero integer and hence surjective.
The identity $(\Phi-1)(1+\Phi+\cdots+\Phi^{f-1})=\Phi^f-1$ implies surjectivity on invariants as well.
Also $H^1(k_v,\Q_\ell/\Z_\ell(m-1))=0$, so purity and localization give
\[
\begin{aligned}
 H^0(k_v,\Q_\ell/\Z_\ell(m-1))&\longrightarrow H^2(C,\Q_\ell/\Z_\ell(m))\\
 &\longrightarrow H^2(C-\{v\},\Q_\ell/\Z_\ell(m))\longrightarrow H^1(k_v,\Q_\ell/\Z_\ell(m-1)).
\end{aligned}
\]
Thus $H^2(C-\{v\},\Q_\ell/\Z_\ell(m))=0$, and the same argument on smaller open subsets of $C$ gives $H^2(F,\Q_\ell/\Z_\ell(m))=0$ by passage to the direct limit.

We next compute the torsion of $K_j(F)$ for $j\ge2$.
Harder's theorem and localization with closed support give $K_j(F)\otimes\Q=0$ for $j\ge2$ \cite[Corollary 2.3.2]{Harder}.
For an abelian group $M$ and a prime $\ell$, we write $M\{\ell\}$ for the $\ell$-primary torsion subgroup of $M$.
Since $[F:F^p]=p$, Geisser--Levine's calculation with $p$-coefficients shows that $K_j(F)\{p\}$ vanishes in this range \cite{GeisserLevine}.
Since $K_j(F)$ is torsion, the direct limit $K_j(F)\otimes(\Q_\ell/\Z_\ell)$ of the groups $K_j(F)/\ell^r$ vanishes, so the direct limit of the mod-$\ell^r$ Bockstein sequences gives $K_j(F;\Q_\ell/\Z_\ell)\simeq K_{j-1}(F)\{\ell\}$ for $j\ge3$.
We calculate the coefficient groups by the Quillen--Lichtenbaum comparison and \'etale descent \cite{ClausenMathew}.
Over $F\otimes_k\bar k$, which is a field since $k$ is algebraically closed in $F$, the \'etale $K$-theory coefficients are concentrated in even degrees, and the $H^2$ vanishing just proved leaves only $H^0(F,\Q_\ell/\Z_\ell(m))$ in total degree $2m$ and $H^1(F,\Q_\ell/\Z_\ell(m+1))$ in total degree $2m+1$.
Geometric connectedness gives $H^0(F,\Q_\ell/\Z_\ell(m))=H^0(k,\Q_\ell/\Z_\ell(m))$.
Compatibly with transfer, we have
\[
 \begin{aligned}
 K_{2m-1}(F)\{\ell\}&\simeq H^0(k,\Q_\ell/\Z_\ell(m))       &&(m\ge2),\\
 K_{2m}(F)\{\ell\}&\simeq H^1(F,\Q_\ell/\Z_\ell(m+1))      &&(m\ge1).
 \end{aligned}
\]
The same identifications hold over $E$.

Finally we verify the exactness of the periodic resolution for the cyclic group $G$.
In odd degree we set $M=\Q_\ell/\Z_\ell(m)$, and in even degree we set $M=H^1(F\otimes_k\bar k,\Q_\ell/\Z_\ell(m+1))$.
In both degrees, we have
\[
 K_j(F)\{\ell\}\simeq M^\Phi,\qquad K_j(E)\{\ell\}\simeq M^{\Phi^s},
\]
where $j=2m-1$ in odd degree, $j=2m$ in even degree, and $s=[E:F]$.
In even degree this follows from $H^1(k,\Q_\ell/\Z_\ell(m+1))=0$ and the Hochschild--Serre edge map.
Since $F\otimes_k\bar k$ has $\ell$-cohomological dimension at most one, vanishing of $H^2(F,\Q_\ell/\Z_\ell(m+1))$ and $H^2(E,\Q_\ell/\Z_\ell(m+1))$ implies that $\Phi-1$ and $\Phi^s-1$ are surjective on $M$.
In odd degree both are nonzero integer multiplications on $\Q_\ell/\Z_\ell(m)$ and thus surjective.
Put $N_s=1+\Phi+\cdots+\Phi^{s-1}$.
Since $\Phi^s-1$ is surjective and $N_s(\Phi-1)=(\Phi-1)N_s=\Phi^s-1$, the operator $N_s$ is surjective on $M$ and induces a surjection $N_s\colon M^{\Phi^s}\to M^\Phi$.
If $z\in M^{\Phi^s}$ and $N_sz=0$, we write $z=(\Phi-1)y$.
Then $(\Phi^s-1)y=N_sz=0$, so the kernel of $N_s\colon M^{\Phi^s}\to M^\Phi$ is $(\Phi-1)M^{\Phi^s}$.
These two checks give the exactness of the periodic resolution in positive degrees, so higher group homology vanishes and the norm on coinvariants is an isomorphism on $\ell$-primary torsion.
Summing over the primes gives the claim.
\end{proof}

\begin{theorem}\label{thm:paper-global-function-split}
Let $F$ be a global function field.
Let $g$ be the comparison of Theorem~\ref{thm:paper-global-comparison}, and put $Z=\fib(g)$.
Let $E=Fk'$ be a finite constant extension splitting $D$, with group $G=\Gal(E/F)$.
Then the composite $(\tau_{\ge2}g)\,\mathcal N^{D}_{E/F}$ is an equivalence.
Consequently,
\[
 \tau_{\ge2}K(D)\simeq\tau_{\ge2}K(F)\oplus Z
 \quad\text{with first component }\tau_{\ge2}g,
 \qquad Z\simeq\cofib(\mathcal N^{D}_{E/F}).
\]
For every $j\ge1$, $g_j$ is surjective, and for $j\ge2$,
\begin{equation}\label{eq:paper-global-function-quotient}
 \pi_jZ\simeq K_j(D)\big/\theta_E K_j(E)
  =K_j(D)\Big/\sum_{\substack{L/F\ \text{finite separable}\\ D_L\ \text{split}}}
                            \theta_L K_j(L).
\end{equation}
\end{theorem}
\begin{proof}
We first compute the homotopy orbits.
The orbit spectral sequence
\[
 H_i(G,\pi_j\tau_{\ge2}K(D_E))
 \ \Longrightarrow\ \pi_{i+j}\bigl((\tau_{\ge2}K(D_E))_{hG}\bigr)
\]
has a finite filtration in each total degree by connectivity.
Under the usual identification $\pi_jK(D_E)=K_j(E)$ by Morita equivalence, the action of each $\gamma\in G$ agrees with the field action.
Indeed, composing the semilinear action of $\gamma$ with the inverse of the field action gives an $E$-linear automorphism of $D_E\simeq M_{\deg D}(E)$, which is inner by Skolem--Noether and hence acts as the identity on $K$-theory.
By Lemma~\ref{lem:paper-global-constant-norm}, only the row $i=0$ remains, so the edge map identifies $\pi_j\bigl((\tau_{\ge2}K(D_E))_{hG}\bigr)$ with $K_j(E)_G$, and under this identification $\pi_j(\mathcal N^{D}_{E/F})$ is induced by the splitting transfer $\theta_E$.
The transfer identity $g_j\theta_E=N_{E/F}$ therefore shows that $(\tau_{\ge2}g)\mathcal N^{D}_{E/F}$ induces an isomorphism on every homotopy group.
By Lemma~\ref{lem:paper-global-low-degrees}, $g_1$ and $g_2$ are isomorphisms, so Lemma~\ref{lem:paper-split-from-norm} gives the direct-sum decomposition, $\cofib(\mathcal N^{D}_{E/F})\simeq Z$, and the surjectivity of $g_j$ for every $j\ge1$.

Finally we identify $\pi_jZ$.
For any finite separable splitting field $L/F$, the compositum $EL/L$ is a finite constant extension, so $N_{EL/L}\colon K_j(EL)\to K_j(L)$ is surjective for $j\ge2$ by Lemma~\ref{lem:paper-global-constant-norm}.
Transitivity of transfer gives $\theta_LK_j(L)\subseteq\theta_EK_j(E)$, and $E$ itself is a finite separable splitting field.
Together with the long exact sequence of the splitting, this gives~\eqref{eq:paper-global-function-quotient}.
\end{proof}

\subsubsection{The third \texorpdfstring{$K$}{K}-group and real signatures}
This subsubsection computes $K_3(D)$ with integral coefficients (Theorem~\ref{thm:global-k3}).
Together with the splittings above, it shows that the comparison fiber is $3$-connected (Corollary~\ref{cor:global-fibre-connectivity}).

Let $J_F$ be the image of $K_3^M(F)\to K_3(F)$.
For a number field $F$, the signatures at the real places give an isomorphism $J_F\simeq\bigoplus_{v\ \mathrm{real}}\Z/2$, and for a global function field, $J_F=0$ \cite{BassTate,MerkurjevSuslinK3}.
For a number field, write $J_{F}^{\mathrm{ram}}\subset J_F$ for the subgroup corresponding to the summands at the real places where $D$ is ramified, and $J_{F}^{\mathrm{sp}}$ for the subgroup corresponding to the summands at the real places where $D$ splits.
For a global function field, put $J_F^{\mathrm{ram}}=0$.
In any splitting field of $D$, every place above a real place at which $D$ is ramified is complex; hence restriction $K_3(F)\to K_3(L)$ kills $J_F^{\mathrm{ram}}$ for every splitting field $L$, and restriction is defined on $K_3(F)/J_F^{\mathrm{ram}}$.

\begin{theorem}\label{thm:global-k3}
If $F$ is a global function field with constant field $\F_q$, there is a unique isomorphism
\[
 K_3(D)\xrightarrow{\ \sim\ }K_3(F)
\]
compatible with restriction to $K_3(L)$ for every finite separable splitting field $L/F$; moreover, restriction from the constant field is an isomorphism $K_3(\F_q)\simeq K_3(F)$, and $K_3(\F_q)$ is cyclic of order $q^2-1$.
If $F$ is a number field, there is a unique isomorphism
\[
 K_3(D)\xrightarrow{\ \sim\ }K_3(F)/J_{F}^{\mathrm{ram}}
\]
compatible with restriction to $K_3(L)$ for every finite Galois splitting field $L/F$ in which every real place of $F$ at which $D$ splits is completely split.
These isomorphisms are compatible with restriction and transfer for finite extensions.
\end{theorem}
\begin{proof}
\emph{Step 1: reduction to weights two and three.}
By Theorem~\ref{thm:paper-global-weight-one}, $\CH^1(F,3;D)=0$, and the stages $\mathcal{F}^rK(D)$ with $r\ge4$ vanish in degrees at most three (over a global function field, the homotopy coniveau tower is available by Lemma~\ref{lem:app-twisted-tower}), so
\[
 \CH^3(F,3;D)\longrightarrow \pi_3 \mathcal{F}^2K(D)
   \longrightarrow \CH^2(F,3;D)\longrightarrow 0,\qquad
 \pi_3\mathcal{F}^2K(D)\twoheadrightarrow K_3(D).
\]
On the commutative side, $\CH^1(F,4)=0$ as well, so the commutative edge map $\pi_3\mathcal{F}^2K(F)\to K_3(F)$ is an isomorphism and the image of weight three is $J_F$ \cite[Corollary 6.1.4]{KahnLevine2010}.

\emph{Step 2: weight two.}
By Propositions~\ref{prop:global-two-chow-rows} and \ref{prop:global-real-weight-two}, the weight-two reduced norm is an isomorphism integrally.
For a finite separable extension $L/F$, transfer shows that the kernel of restriction on $K_3^{\mathrm{ind}}(F)$ is killed by the extension degree.
The coefficient sequence and $H^0_{\mot}(F,\Q(2))=0$ identify the torsion subgroup of $K_3^{\mathrm{ind}}(F)$ with $H^0(F,\Q/\Z(2))$, and restriction is inclusion of invariants \cite{Voevodsky}; over a global function field, the characteristic-primary part is treated in Step~5.
Hence restriction is injective on $K_3^{\mathrm{ind}}(F)$.
Since $\CH^2(F,3)\simeq K_3^{\mathrm{ind}}(F)$ and the weight-two reduced norm is an isomorphism, restriction of $\CH^2(F,3;D)$ to a splitting field is also injective, with the same image as on the commutative side.

\emph{Step 3: weight three.}
As in the proof of Theorem~\ref{thm:paper-global-weight-one} for $\CH^2(k,2;D)$, pushforward from a maximal subfield at each admissible closed point shows that the group $\CH^3(F,3;D)$ is generated by the transfers of $\CH^3(E,3)=K_3^M(E)$ over finite splitting fields $E/F$ \cite{NesterenkoSuslin}.
For global function fields, these groups $K_3^M$ vanish.
For number fields, we choose a finite Galois splitting field $L/F$ in which every real place of $F$ at which $D$ splits is completely split; it exists by Lemma~\ref{lem:paper-galois-splitting-real}.
Every real place of any finite splitting field $E/F$ lies over a split real place of $F$, so $EL/E$ splits completely at every real place.
Restriction duplicates the real signatures of $K_3^M$ on the real branches, and transfer adds them over the branches.
Thus $K_3^M(EL)\to K_3^M(E)$ is surjective, and transitivity of transfer gives
\[
 K_3^M(L)_{\Gal(L/F)}\twoheadrightarrow \CH^3(F,3;D)
      \xrightarrow{\rho} K_3^M(F).
\]
The composite is the usual norm, and it is an isomorphism onto $J_{F}^{\mathrm{sp}}$: each Galois orbit of real branches corresponds to one split real place, and the unit signature on one branch maps to the unit signature below; here $K_3^M(F)\simeq J_F$ by the real signatures \cite{BassTate}.
Since restriction duplicates signatures, $\CH^3(F,3;D)\to J_L$ is injective with image $J_L^{\Gal(L/F)}$.
The commutative side has kernel $J_{F}^{\mathrm{ram}}$ and the same image \cite[Chapter II, Theorem 2.1(3)]{BassTate}, \cite[Lemmas 6.6.3--6.6.4]{KahnLevine2010}.

\emph{Step 4: assembly.}
For a number field, let $L$ be as in Step~3; for a global function field, let $L$ be any finite separable splitting field.
If an element of $\pi_3\mathcal{F}^2K(D)$ restricts to zero in $K_3(L)$, its weight-two image is zero, so it comes from weight three, where restriction is also injective.
Since this map factors through $K_3(D)$, one obtains $\pi_3\mathcal{F}^2K(D)\simeq K_3(D)$.
Given $x\in K_3(D)$, choose $y\in K_3(F)$ whose image in $K_3^{\mathrm{ind}}(F)=\CH^2(F,3)$ is $\rho$ of the weight-two image of $x$.
For global function fields, $L/F$ is any finite separable splitting field and $J_L=0$, so matching the weight-two images already makes the restrictions of $x$ and $y$ to $K_3(L)$ agree.
For number fields, the difference of the restrictions of $x$ and $y$ to $K_3(L)$ lies in $J_L^{\Gal(L/F)}$, which is the image of weight three on both sides, so after correcting $y$ by an element of $J_F$ the two restrictions agree.
The same argument with the two sides exchanged shows that the images of $K_3(D)$ and $K_3(F)$ in $K_3(L)$ agree.
The kernel on the commutative side is $J_F^{\mathrm{ram}}$, which is zero for a global function field.
Thus $K_3(D)$ and $K_3(F)/J_F^{\mathrm{ram}}$ both inject into $K_3(L)$ with equal images, so there is exactly one isomorphism between them compatible with restriction to $K_3(L)$.

\emph{Step 5: global function fields.}
For global function fields, weight three vanishes.
Here $J_F=0$, so $K_3^{\mathrm{ind}}(F)=K_3(F)$, and its characteristic-primary torsion vanishes by Proposition~\ref{prop:app-frobenius}(2); for the other primary parts we use Beilinson--Lichtenbaum comparison.
Finite-coefficient comparison and Harder's vanishing give, for $\ell\ne\operatorname{char} F$,
\[
 \CH^2(F,3)_{(\ell)}
 = H^1_{\mot}(F,\Z_{(\ell)}(2))
 \simeq H^0(F,\Q_\ell/\Z_\ell(2)).
\]
The action on the right-hand side factors through the constant field, whose Frobenius acts by multiplication by $q^2$.
Hence $K_3(F)$ is cyclic of order $q^2-1$, and restriction from the constant field is an isomorphism \cite{Harder,QuillenFinite}.

\emph{Step 6: independence and functoriality.}
Independence of the splitting field follows by comparing restriction images in a compositum.
Restriction and transfer for arbitrary finite extensions descend to the signature quotients by duplicating and adding signatures at ramified real places.
For a finite extension $E/F$, compatibility with the isomorphisms constructed from restriction images follows from the base-change formula and Morita equivalence over each field component of $E\otimes_F L$.
In particular, $\theta_E$ from a splitting field $E/F$ corresponds to the composite $K_3(E)\xrightarrow{N_{E/F}}K_3(F)\to K_3(F)/J_{F}^{\mathrm{ram}}$ on the commutative side.
\end{proof}

\begin{corollary}\label{cor:global-k3-structure}
Let $F$ be a number field with $r_2$ complex places, let $r_{\mathrm{sp}}$ be the number of real places at which $D$ splits, and put $w_2=|H^0(F,\Q/\Z(2))|$.
Then $K_3(D)$ has free rank $r_2$, and its torsion subgroup is $\Z/w_2$ if $r_{\mathrm{sp}}=0$ and $\Z/(2w_2)\oplus(\Z/2)^{r_{\mathrm{sp}}-1}$ if $r_{\mathrm{sp}}>0$.
\end{corollary}
\begin{proof}
By Theorem~\ref{thm:global-k3}, $K_3(D)\simeq K_3(F)/J_F^{\mathrm{ram}}$, and the free rank of $K_3(F)$ is $r_2$ \cite{Borel}.
Since $J_F^{\mathrm{ram}}$ is finite, the torsion of $K_3(D)$ is $T/J_F^{\mathrm{ram}}$, where $T$ is the torsion of $K_3(F)$.
The group $T$ is an extension of $H^0(F,\Q/\Z(2))=\Z/w_2$ by $J_F$.
For a lift $z\in T$ of a generator (it is torsion since $J_F$ is finite), $w_2z=\{-1,-1,-1\}$ \cite[\S11]{MerkurjevSuslinK3}, the element of $J_F$ that is nonzero at every real place.
Hence $T/J_F^{\mathrm{ram}}$ is an extension of $\Z/w_2$ by $J_F^{\mathrm{sp}}$ in which $w_2\bar z$ is the vector with all components one; it is zero if there are no split real places and of order two otherwise, which gives the stated torsion.
\end{proof}

\begin{corollary}\label{cor:global-fibre-connectivity}
The fiber $Z$ of the comparison map $g$ belongs to $\mathrm{Sp}_{\ge4}$.
If $D$ is ramified at a real place, the same holds for the fiber of the comparison map with $2$ inverted.
\end{corollary}
\begin{proof}
By Theorems~\ref{thm:paper-global-split} and~\ref{thm:paper-global-function-split}, $g_3$ is surjective.
By Theorem~\ref{thm:global-k3}, its source and target are abstractly isomorphic over $\Z$, or over $\Z[1/2]$ under real ramification.
Both are finitely generated: for number fields, localization makes $K_3(F)$ an image of $K_3(\mathcal O_F)$, since $K_2$ of finite fields vanishes, and for global function fields $K_3(F)$ is finite by the same theorem.
A surjection between isomorphic finitely generated modules over these Noetherian rings is an isomorphism, so $g_3$ is an isomorphism.
By Lemma~\ref{lem:paper-global-low-degrees}, $g_1$ and $g_2$ are isomorphisms.
The splitting in degrees at least two therefore makes the homotopy groups of $Z$ vanish in degrees at most three.
\end{proof}

\begin{remark}
At each prime $\ell\ne\operatorname{char}F$ to which Proposition~\ref{prop:global-two-chow-rows} applies, the comparison in degree three also agrees with the weight-two reduced norm.
Indeed, the rank-$1$ trivialization acts as the identity on the geometric weight-two layer of \'etale $K$-theory, and the reduced norm on $\CH^2(F,3;D)_{(\ell)}$ is the comparison on that layer.
Let $L$ be a splitting field as in Theorem~\ref{thm:global-k3}.
If $\ell\mid\deg D$, restriction carries the completed comparison $m_x^{-1}$ to $m_{x_L}^{-1}$; the latter acts as the identity in positive degrees after split Morita identification by Lemma~\ref{lem:paper-global-product-vanishing}.
If $\ell\nmid\deg D$, flat base change gives the same compatibility because $N_{D_L/L}$ is $\deg D$ times the split Morita comparison.
Since $K_3(L)_{(\ell)}$ injects into its derived completion, the two restrictions agree in $K_3(L)_{(\ell)}$.
At these primes $J_F{}_{(\ell)}=0$, so the restriction description in Theorem~\ref{thm:global-k3} identifies $g_3$ with the weight-two reduced norm.
\end{remark}

When $D$ is ramified at a real place, the comparison of Theorem~\ref{thm:paper-global-comparison} is defined after inverting $2$.
Integrally, $K_1(D)$ is the subgroup of $F^\times$ consisting of the elements positive at every ramified real place.
The group $K_1(D)$ has no two-torsion corresponding to the element $-1$ of $K_1(F)$.
\begin{example}\label{ex:paper-global-k3-real-ramification}
Let $F=\Q$, $D=(-1,-1)$, and $E=\Q(i)$.
There is no integral map $K_3(D)\to K_3(\Q)$ carrying all splitting transfers to the corresponding field norms.
\end{example}
\begin{proof}
By Theorem~\ref{thm:global-k3}, $K_3(D)=\Z/24$ and $K_3(\Q)=\Z/48$, and $\theta_E$ corresponds to the composite $K_3(E)\xrightarrow{N_{E/\Q}}K_3(\Q)\to K_3(\Q)/J_{\Q}^{\mathrm{ram}}$.
Lift a generator of $K_3^{\mathrm{ind}}(\Q)=\Z/24$ to $z\in K_3(\Q)=\Z/48$ and put $y=\res_{E/\Q}z$.
The surjection $K_3(\Q)=\Z/48\to K_3^{\mathrm{ind}}(\Q)=\Z/24$ has kernel $J_{\Q}=\{0,24\}$, so $z$ generates $\Z/48$.
Then $N_{E/\Q}y=2z$ has order $24$, while its image $\theta_E(y)$ in the signature quotient $K_3(D)$ has order $12$.
Since $12\theta_E(y)=0$ but $12N_{E/\Q}y\ne0$, no homomorphism $K_3(D)\to K_3(\Q)$ carries $\theta_E(y)$ to $N_{E/\Q}y$ \cite[\S11]{MerkurjevSuslinK3}.
\end{proof}

\subsubsection{The basic comparison in prime index}
For a prime $\ell\ne\operatorname{char}F$ at which the $\ell$-primary index of $D$ is $\ell$, the comparison can be analyzed directly by comparing twisted and commutative cycle complexes.
This case, treated in Proposition~\ref{prop:global-squarefree-spectrum}, is an input for the degree bounds in the next subsubsection.

\begin{proposition}\label{prop:global-squarefree-spectrum}
For a prime $\ell\mid \deg D$ with $\ell\ne\operatorname{char}F$, and $\ell\ne2$ if $D$ is ramified at a real place, such that $v_\ell(\deg D)=1$, the map $g$ of Theorem~\ref{thm:paper-global-comparison} is an equivalence after $\ell$-localization.
Consequently, if the prime-to-characteristic part of the index is squarefree, $g$ is an equivalence in the coefficient range of that theorem.
\end{proposition}

\begin{proof}
Since the rational comparison holds, the integral comparison after $\ell$-localization reduces to the $\ell$-completed comparison.
We split the primary components of the Brauer class other than the $\ell$-primary component over a separable extension of degree prime to $\ell$.
After $\ell$-localization, restriction to this extension is split injective, since its composite with transfer is multiplication by the degree, and the same holds for the cofiber of the localization unit $K(D)\to L_{K(1),\ell}K(D)$.
Since $m_x$ is an equivalence, it suffices to show that this localization unit is an equivalence after completion in positive degrees, and by the split injectivity this may be checked after restriction.
We are therefore reduced to the case where the division algebra in the Brauer class has degree $\ell$.
In positive characteristic, both composites of restriction and transfer for a finite purely inseparable extension are multiplication by its degree, which is prime to $\ell$.
Passing to the direct limit, we therefore move the calculation modulo $\ell^r$ to the perfect closure, as required by \cite[Theorem 8.2.2]{KahnLevine2010}, which is used below.
The \'etale side is unchanged by purely inseparable extension.

We first compare the twisted and commutative cycle complexes in positive weight.
By \cite[Theorem 8.2.2]{KahnLevine2010}, the reduced cycle map in positive weight $q$ is an isomorphism in motivic degrees $s<q$.
In degree $q$ it is injective, with cokernel the image of
\[
 K_q^M(F)\xrightarrow{[D]\cup-}
 H^{q+2}_{\et}(F,\Z/\ell(q+1)).
\]
Since $q+2\ge3$, the target $H^{q+2}_{\et}(F,\Z/\ell(q+1))$ vanishes when the $\ell$-cohomological dimension is at most two.
For $\ell=2$ over a number field with real places, restriction of the target to all real closures is injective, and $[D]=0$ at each real closure, so the image is again zero.
In degrees $s>q$, both twisted and commutative terms vanish because there are no cycles in simplicial degree $2q-s<q$.
Thus the positive-weight map is a quasi-isomorphism of complexes over a field.
It remains a quasi-isomorphism modulo $\ell^r$ and commutes with restriction to real closures \cite[Theorem 6.5.3]{KahnLevine2010}, \cite{Jannsen,Voevodsky}.

We compare the homotopy coniveau tower modulo $\ell^r$ with the tower of global sections of its \'etale hypersheafification.
The layer comparison is the map from the twisted motivic complex over a field to \'etale cohomology with coefficients $\Z/\ell^r(q)$, placed in homotopy degree $2q-s$.
This follows from the Kahn--Levine identification of the layers and from the fact that \'etale sheafification of the reduced cycle map identifies the coefficient sheaves.
By the quasi-isomorphism just proved, the comparison cofiber for $q>0$ identifies with the commutative one.
By Beilinson--Lichtenbaum comparison and vanishing of motivic cohomology in degrees $>q$, its cohomology occurs only for $s>q$, where it equals $H^s_{\et}(F,\Z/\ell^r(q))$.
For $q\ge2$, these cohomology groups vanish in cohomological dimension two.
In the real-relative case, restriction to the real closures is also an isomorphism in degrees $s\ge3$, so the real-relative comparison cofiber vanishes.
For $q=1$, both the absolute and real-relative cofibers lie in homotopy degrees $2-s\le0$ \cite[Theorem 6.5.5, \S6.9]{KahnLevine2010}.

Convergence of the towers is controlled by the following degree estimate.
The term with supports in codimension at least $q$ has no supports in simplicial degree below $q$, so the $q$th stage of the homotopy coniveau tower is concentrated in degrees at least $q$, and so are its geometric stalks.
The stalk of a layer modulo $\ell^r$ is a finite-coefficient Tate sheaf in degree $2q$.
In any fixed degree, the geometric homotopy group is computed using sufficiently many but finitely many layers, and is therefore a finite $\ell$-primary coefficient group.
Consequently, after \'etale hypersheafification, taking global sections or real-relative sections lowers connectivity by at most the cohomological-dimension bound two.
The tails at sufficiently high stages do not contribute to a fixed total degree.
Layerwise comparison therefore shows that $\pi_j$ of the comparison cofiber for the entire tower vanishes for $j\ge2$ and that $\pi_1$ of this cofiber injects into $\pi_1$ of the weight-zero comparison cofiber.

It remains to check weight zero.
The algebraic layer places the rank lattice $\deg D\cdot\Z$ with coefficients modulo $\ell^r$ in degree zero, and the comparison on $\pi_0$ is
\[
 \deg D\cdot\Z/\ell^r(\deg D\cdot\Z)\longrightarrow\Z/\ell^r,\qquad [\deg D]\longmapsto [\deg D].
\]
The \'etale layer has no positive homotopy groups, so $\pi_1$ of the weight-zero comparison cofiber is the kernel of this map.
In the real-relative case, $D$ splits over each real closure, and all local weight-zero comparison cofibers lie in degrees at most zero.
The same kernel occurs in $\pi_1$ of the relative fiber.
Under the coefficient transition, this kernel at level $r+v_\ell(\deg D)$ maps to zero at level $r$.
Thus the inverse limit of $\pi_1$ of the entire comparison cofiber is zero.
Since $\pi_2$ of the comparison cofiber vanishes at every finite stage, the Milnor exact sequence also makes $\pi_1$ of the derived-completed comparison cofiber zero.
This proves that the completed comparison is an isomorphism in every positive degree.
In the real-relative case, this proves the comparison for the relative terms.
Over the real closures $D$ is split, so Morita equivalence and the Quillen--Lichtenbaum comparison for the real closures give the comparison for the real terms.
If the relative and real terms are isomorphic in positive degrees, the long exact sequence gives isomorphisms of the absolute terms in degrees at least two; in degree one we use that $g_1=\Nrd_1$ is an isomorphism (Lemma~\ref{lem:paper-global-low-degrees}).
Applying the arithmetic square shows that $g$ is an equivalence after $\ell$-localization.
For the second assertion, the map $g[1/\deg D]=N_{D/F}/\deg D$ is an equivalence because $i_DN_{D/F}=N_{D/F}i_D=(\deg D)^2$.
If $p=\operatorname{char}F$ divides $\deg D$, Proposition~\ref{prop:app-frobenius} gives $K_i(D;\Z/p)=K_i(F;\Z/p)=0$ for $i\ge2$, so the $p$-completed comparison is the completion of $g_1=\Nrd_1$, and $g_1$ is an isomorphism by Lemma~\ref{lem:paper-global-low-degrees}.
Together with the first assertion at the other primes dividing $\deg D$, the arithmetic square for the primes dividing $\deg D$ shows that $g$ is an equivalence.
\end{proof}

\subsubsection{Prime decomposition and degree bounds for the comparison fiber}
We describe the primary factors of the comparison fiber using the localization fibers $B_\ell(D)$ (Theorem~\ref{thm:paper-global-fiber-bounds}).

For $\ell\mid \deg D$ with $\ell\ne\operatorname{char}F$, write
\[
 e_\ell=\ell^{v_\ell(\deg D)},\qquad
 B_\ell(D)=\fib\!\left(K(D)_{(\ell)}\longrightarrow L_{E(1),\ell}K(D)\right).
\]

\begin{theorem}\label{thm:paper-global-fiber-bounds}
Let $F$ be a global field; if $D$ is ramified at a real place, invert $2$ throughout and exclude $\ell=2$.
Let $g\colon\tau_{\ge1}K(D)\to\tau_{\ge1}K(F)$ be the comparison of Theorem~\ref{thm:paper-global-comparison}, and put $Z=\fib(g)$.
Then the fiber decomposes as
\begin{equation}\label{eq:paper-global-fiber-primes}
 Z\simeq
 \prod_{\substack{\ell\mid \deg D,\ \ell\ne\operatorname{char}F\\
                  \ell\ne2\ \text{if $D$ is ramified at a real place}}}
       \tau_{\ge2}B_\ell(D).
\end{equation}
For each factor,
\begin{equation}\label{eq:paper-global-fiber-amplitude}
 e_\ell\,\pi_j B_\ell(D)=0\quad (j\ge1),
 \qquad
 B_\ell(D)\in\mathrm{Sp}_{\le e_\ell-3}.
\end{equation}
The factor $\tau_{\ge2}B_\ell(D)$ is zero if $e_\ell=\ell$.
\end{theorem}

\begin{proof}
The functor $L_{E(1),\ell}$ is the smashing localization combining rationalization and height one.
Since $\ell$ is invertible over $F$, the $\mathrm{THH}(F)$-module structure makes multiplication by $\ell$ an equivalence on $\mathrm{THH}$ of algebras and schemes over $F$.
It is an equivalence of cyclotomic spectra and hence also on $\mathrm{TC}$, so $\mathrm{TC}_{(\ell)}$ is rational.
In the definition of Selmer $K$-theory, the map from $\mathrm{TC}_{(\ell)}$ to its $L_{E(1),\ell}$ localization is therefore an equivalence, which gives an equivalence $K^{\mathrm{Sel}}_{(\ell)}\simeq L_{E(1),\ell}K$ compatible with the comparison maps \cite[Definition 6.1]{ClausenMathew}.

We first bound the exponent of the factors.
We choose a finite separable splitting field $E/F$ of degree $\deg D$.
The virtual $\ell$-cohomological dimension of $E$ is at most two, so \cite[Theorem 6.18]{ClausenMathew} gives $B_\ell(E)\in\mathrm{Sp}_{\le-2}$.
By split Morita equivalence, we have $B_\ell(D_E)\simeq B_\ell(E)$.
We have $N_{D_E/D}\res_{E/F}=\deg D$ on both $K(D)$ and $L_{E(1),\ell}K(D)$.
Thus $\deg D$ kills $\pi_jB_\ell(D)$ for $j\ge1$.
These groups are $\ell$-local, so their exponent divides $e_\ell$.

For the upper degree bound, we split the prime-to-$\ell$ Brauer component over a finite separable extension $F'/F$ of degree prime to $\ell$.
This does not change the $\ell$-primary index, since $\ind(D_{F'})\mid\deg D\mid[F':F]\,\ind(D_{F'})$.
Restriction and transfer divided by the extension degree exhibit $\pi_jB_\ell(D)$ as a direct summand of $\pi_jB_\ell(D_{F'})$.
Thus by Morita equivalence it suffices to treat a division algebra of degree $e_\ell$.
For $X=\SB(D)$, we have the Quillen decomposition \cite[Section 8]{QuillenHigher}
\[
 K(X)\simeq\bigoplus_{i=0}^{e_\ell-1}K(D^{\otimes i}).
\]
Applying this decomposition to the localization unit and its fiber, we realize $B_\ell(D)$ as a summand of $B_\ell(X)$.
The residue fields of points of $X$ have virtual $\ell$-cohomological dimension at most $2+\dim X=e_\ell+1$.
The same theorem of Clausen and Mathew gives $B_\ell(X)\in\mathrm{Sp}_{\le e_\ell-3}$ and hence \eqref{eq:paper-global-fiber-amplitude}.

Next we take the derived $\ell$-completion of $g$.
The construction in Theorem~\ref{thm:paper-global-comparison} uses $L_{K(1),\ell}K(F)\simeq L_{K(1),\ell}K(D)$, so the completed comparison fiber in degrees at least two is
\[
 \tau_{\ge2}\fib\!\left(
   (\tau_{\ge1}K(D))^\wedge_\ell
      \longrightarrow \tau_{\ge1}L_{K(1),\ell}K(D)
 \right).
\]
Replacing $\tau_{\ge1}K(D)$ by $K(D)$ leaves this unchanged, since the completion of $\tau_{\le0}K(D)$ is concentrated in degrees at most one by $(-)^\wedge_\ell=\varprojlim_r(-)/\ell^r$ and the Moore coefficient sequence.
Removing $\tau_{\ge1}$ in the target leaves the same $2$-connective fiber.
The height-one arithmetic square gives $(L_{E(1),\ell}K(D))^\wedge_\ell\simeq L_{K(1),\ell}K(D)$, so the completed comparison fiber in this range is $\tau_{\ge2}B_\ell(D)^\wedge_\ell$.

The bounded exponent above implies that $\tau_{\ge1}B_\ell(D)$ is derived $\ell$-complete.
The Bockstein torsion terms have transition maps given by multiplication by $\ell$ and form a pro-zero system.
The quotient terms stabilize in each degree, so the Milnor sequence recovers the homotopy groups after completion.
Since the completion of $\tau_{\le0}B_\ell(D)$ is concentrated in degrees at most one, $\tau_{\ge2}B_\ell(D)^\wedge_\ell\simeq\tau_{\ge2}B_\ell(D)$.
On the other hand, $g$ is an isomorphism in degrees one and two by Lemma~\ref{lem:paper-global-low-degrees}, so $Z$ is $2$-connective and its completion is too.
We obtain $Z^\wedge_\ell\simeq\tau_{\ge2}B_\ell(D)$.

The construction of $g[1/\deg D]$ gives $Z[1/\deg D]=0$.
In characteristic $p$, the completed positive-degree comparison is an equivalence at $p$ by the rank-$1$ Frobenius comparison of Proposition~\ref{prop:app-frobenius} and the degree-one reduced norm.
The arithmetic square of $Z$ therefore leaves only the completed factors in \eqref{eq:paper-global-fiber-primes}.
When $e_\ell=\ell$, restrict to $F'$ as above and pass by Morita equivalence to the division algebra $D'$ in the Brauer class of $D_{F'}$, which has degree $\ell$.
Proposition~\ref{prop:global-squarefree-spectrum}, applied to $D'$ over $F'$, shows that the comparison for $D'$ is an equivalence after $\ell$-localization.
Its completed comparison fiber in degrees at least two is $\tau_{\ge2}B_\ell(D')$, as identified above, so $\tau_{\ge2}B_\ell(D_{F'})\simeq\tau_{\ge2}B_\ell(D')=0$.
Restriction and transfer, whose composite is an $\ell$-local unit, show that $\tau_{\ge2}B_\ell(D)$ is a summand of this zero spectrum.
\end{proof}

\begin{corollary}\label{cor:paper-global-fiber-finite}
Let $F$ be a number field, and let $\ell$ be a prime dividing $\deg D$, with $\ell\ne2$ if $D$ is ramified at a real place.
Then every homotopy group of $\tau_{\ge2}B_\ell(D)$ is a finite $\ell$-group.
\end{corollary}
\begin{proof}
We choose a maximal $\mathcal O_F$-order $\Lambda\subset D$ and a nonempty open subscheme $U\subset\Spec\mathcal O_F$ such that $\mathcal A=\Lambda|_U$ is Azumaya and $U$ contains no place above $\ell$.
The groups $K_j(\Lambda)$ are finitely generated \cite{Kuku1982}.
Maximal orders are hereditary and the support terms at the deleted places are $K$-theories of finite residue division rings, so localization gives finite generation of $K_j(\mathcal A)$ as well \cite[Section 5]{QuillenHigher}.
The support term for $K(\mathcal A)\to K(D)$ is the direct sum of $K(k_v)$ over the closed points $v$ of $U$.
Each $B_\ell(k_v)$ is concentrated in degrees at most $-2$, and height-one localization commutes with direct sums and localization sequences \cite{QuillenFinite,ClausenMathew}.
Thus $B_\ell(\mathcal A)\to B_\ell(D)$ is an isomorphism in degrees at least two.
The equivalence $Z^\wedge_\ell\simeq\tau_{\ge2}B_\ell(D)$ obtained in the proof of Theorem~\ref{thm:paper-global-fiber-bounds} (and $Z_{(\ell)}\to Z^\wedge_\ell$ is an equivalence, since the other factors of the finite product have bounded exponent prime to $\ell$) commutes, after $\ell$-completion, with the canonical maps to $(\tau_{\ge2}K(D))^\wedge_\ell$, the common source of the completed comparison and of the completed localization unit.
For $j\ge2$, the map $\pi_jB_\ell(D)\to\pi_jB_\ell(D)^\wedge_\ell$ is an isomorphism by the bounded exponent in \eqref{eq:paper-global-fiber-amplitude}.
Since the completed comparison splits by Theorem~\ref{thm:paper-global-split}, the composite $\pi_jB_\ell(D)\to K_j(D)_{(\ell)}\to\pi_jK(D)^\wedge_\ell$ is injective for $j\ge2$, and hence so is $\pi_jB_\ell(D)\to K_j(D)_{(\ell)}$.
The latter map factors through $\pi_jB_\ell(\mathcal A)\to K_j(\mathcal A)_{(\ell)}$, exhibiting $\pi_jB_\ell(D)$ as a subgroup of a finitely generated $\Z_{(\ell)}$-module.
By \eqref{eq:paper-global-fiber-amplitude}, the group $\pi_jB_\ell(D)$ is finite.
\end{proof}

\begin{corollary}\label{cor:paper-global-fiber-prime-four}
If $e_\ell=\ell$ or $e_\ell=4$ for every $\ell\mid \deg D$ with $\ell\ne\operatorname{char}F$, then $g$ is an equivalence (with the coefficient convention of Theorem~\ref{thm:paper-global-fiber-bounds}).
Over a global function field, we have $\pi_jZ=0$ for $j<4$, and the $\ell$-primary part of $\pi_jZ$ can survive only in the range
\[
 4\le j\le e_\ell-3.
\]
\end{corollary}
\begin{proof}
If $e_\ell=4$ then $e_\ell-3=1$ and $\tau_{\ge2}B_\ell(D)=0$.
The case $e_\ell=\ell$ is the last assertion of Theorem~\ref{thm:paper-global-fiber-bounds}.
Over a global function field, $\pi_jZ=0$ for $j<4$ by Corollary~\ref{cor:global-fibre-connectivity}.
We combine this with \eqref{eq:paper-global-fiber-amplitude}.
\end{proof}

\begin{remark}
Equations~\eqref{eq:paper-global-function-quotient} and \eqref{eq:paper-global-fiber-primes} describe the same fiber by transfers and by height-one localization, respectively.
The latter gives the prime decomposition and an exponent bound for $K_j(D)/\theta_EK_j(E)$ over a global function field.
\end{remark}

\subsubsection{Invariance of the fiber under powers of the Brauer class}
We use Adams operations and the prime decomposition~\eqref{eq:paper-global-fiber-primes} to compare the fibers for different powers of the Brauer class.

\begin{proposition}
\label{prop:paper-brauer-power-fiber}
Let $F$ be a global field and let $D$ be a central division algebra over $F$.
Let $a\ge1$ be an integer with $\gcd(a,\per D)=1$.
Let $g_D$ and $g_{D^{\otimes a}}$ be the comparisons of Theorem~\ref{thm:paper-global-comparison} for $D$ and for the division algebra in the Brauer class of $D^{\otimes a}$, and put $Z_D=\fib(g_D)$ and $Z_{D^{\otimes a}}=\fib(g_{D^{\otimes a}})$.
For number fields with real ramification, the comparison fibers are taken with $2$ inverted, as in that theorem.
There is a noncanonical equivalence of spectra
\[
 Z_D\simeq Z_{D^{\otimes a}}.
\]
Consequently, two central simple algebras whose classes generate the same cyclic subgroup of $\Br(F)$ have equivalent comparison fibers under the same coefficient convention.
Moreover, if nonnegative integers $i,j$ satisfy $\gcd(i,\per D)=\gcd(j,\per D)$, then $Z_{D^{\otimes i}}\simeq Z_{D^{\otimes j}}$, where $D^{\otimes0}=F$.
\end{proposition}
The proof of Proposition~\ref{prop:paper-brauer-power-fiber} is given after Lemmas~\ref{lem:paper-brauer-power-correspondence}--\ref{lem:paper-brauer-power-descent} below.

For the auxiliary comparison, let $F$ be any field and let $D$ be a central division algebra over $F$ with $\deg D\ge2$, and take $a\ge2$.
Write $D^i=D^{\otimes i}$ and $X=\SB(D^{\mathrm{op}})$.
We regard the universal right $D^{\mathrm{op}}$-ideal as a left $D$-bundle, and let $\phi_i:K(D^i)\to K(X)$ be the maps defined by its tensor powers.
The Quillen decomposition is
\[
 \phi=\bigoplus_{i=0}^{\deg D-1}\phi_i:
 P=\bigoplus_{i=0}^{\deg D-1}K(D^i)\xrightarrow{\sim}K(X).
\]
The negative universal bundles give a semiorthogonal decomposition ordered by $i=\deg D-1,\ldots,0$; reorder the $K$-theory direct sum in increasing order.
We take the inclusions and projections to be exact functors on the differential graded enhancements of the categories of perfect complexes \cite[\S8]{QuillenHigher}, \cite[Theorem 4.1]{Bernardara2009}.
The same universal bundle defines $\phi_i$ for all $i\ge0$, giving
\[
 \phi^{(a)}=\bigoplus_{i=0}^{\deg D-1}\phi_{ai}:
 P_a=\bigoplus_{i=0}^{\deg D-1}K(D^{ai})\longrightarrow K(X).
\]

\begin{lemma}\label{lem:paper-brauer-power-correspondence}
If $\gcd(a,\per D)=1$, then $\phi^{(a)}[1/a]\colon P_a[1/a]\to K(X)[1/a]$ is an equivalence.
\end{lemma}
\begin{proof}
Since $\per D\mid\deg D$ and multiplication by $a$ permutes $\Z/(\per D)$, the classes $[D^i]$ and $[D^{ai}]$ for $0\le i<\deg D$ agree as multisets, so permuting components and applying Morita equivalences gives $P_a\simeq P$.
Composing the projections in the Quillen decomposition with $\phi_{ai}$, we obtain derived tensor products with perfect bimodule complexes between the components.
We use the $K_0$-classes of these bimodule complexes as correspondences, with composition given by derived tensor product.
Taking Morita ranks after splitting gives a faithful representation of the resulting endomorphism ring of correspondences as
\[
 \Pi=\left\{(n_{ij})\in M_{\deg D}(\Z):
 n_{ij}\in\ind\bigl((D^i)^{\mathrm{op}}\otimes_F D^j\bigr)\Z\right\}
\]
with ordinary matrix multiplication as composition.
By additivity of $K$-theory, equalities of correspondences yield equalities of homotopy classes of maps of spectra \cite[Proposition 2.25(i), \S10, (10.1)]{TabuadaVanDenBerghSeparable}.

After splitting, put $t=[\mathcal O(-1)]$.
The maps $\phi$ and $\phi^{(a)}$ compare the two lists $1,t,\ldots,t^{\deg D-1}$ and $1,t^a,\ldots,t^{a(\deg D-1)}$ in $\Z[t]/(1-t)^{\deg D}$.
Writing $h=1-t$, we have $1-t^a=ah+O(h^2)$, so the degree-$r$ factor in the $h$-filtration is $a^r$.
The matrix comparing the two lists therefore has determinant
\[
                    \pm a^{\deg D\,(\deg D-1)/2}.
\]
Because $\gcd(a,\per D)=1$, the algebras $D^i$ and $D^{ai}$ have the same splitting fields and hence the same index.
Through the chosen permutation and the Morita equivalences $P_a\simeq P$, the map $\phi^{-1}\phi^{(a)}$ becomes an endomorphism of $P$ whose matrix in $\Pi$ is the split $K_0$-matrix above up to a permutation of components, so its determinant is $\pm a^{\deg D(\deg D-1)/2}$.
The Cayley--Hamilton identity expresses its inverse as a polynomial in it with coefficients in $\Z[1/a]$, so the inverse belongs to $\Pi[1/a]$.
Realizing this inverse in $K$-theory shows that $\phi^{(a)}[1/a]$ is an equivalence.
\end{proof}

\begin{lemma}\label{lem:paper-brauer-power-adams}
The stable Adams operation $\psi_X^a$ is an equivalence on $K(X)[1/a]$.
\end{lemma}
\begin{proof}
Let $KGL$ be the motivic $\PP^1$-spectrum representing algebraic $K$-theory, let $f_q$ be the $q$-effective cover, and let $s_q=\cofib(f_{q+1}\to f_q)$ be the $q$th slice.
Write $\operatorname{Fil}^qK(X)=(f_qKGL)(X)$; in particular, $\operatorname{Fil}^0K(X)\simeq K(X)$.
The unit $\eta:\mathbf1[1/a]\to KGL[1/a]$ induces an isomorphism on the zeroth slice, and $\psi^a\eta=\eta$ implies $s_0\psi^a=1$.
Since $\psi^a(b)=ab$ for the Bott element $b$, we have $s_r\psi^a=a^r$ for every $r$.
Put $m=\dim X$.
The $q$th stage of the motivic filtration $\operatorname{Fil}^\bullet K(X)[1/a]$ belongs to $\mathrm{Sp}_{\ge q-m}$.
A fixed $\pi_n$ can be computed in the finite quotient $\operatorname{Fil}^0/\operatorname{Fil}^q$ for $q>n+m$, and $a^r$ is a unit on each of its slices.
Thus $\psi_X^a$ induces an isomorphism on every $\pi_n$.
For the stable operation and its action on the Bott element, see \cite[Example 9.6]{AnnalaHoyoisIwasa}; the unit equivalence on slices and the convergence estimate are \cite[Corollaries 9.5 and 9.7]{BachmannElmantoMorrow}.
Integer localization commutes with effective covers and slices because the effective subcategories have compact generators.
\end{proof}

\begin{lemma}\label{lem:paper-brauer-power-descent}
For every $i\ge0$ there are Adams operations $\psi_i^a\colon K_n(D^i)\to K_n(D^{ai})$, $n\ge0$, such that
\[
 (\psi_X^a)_*(\phi_i)_*=(\phi_{ai})_*\psi_i^a\qquad\text{on }K_n(D^i)[1/a].
\]
\end{lemma}
\begin{proof}
We choose a finite Galois splitting field $L/F$, put $G=\Gal(L/F)$, and choose $\gamma$-semilinear isomorphisms $\xi_\gamma$ for $\gamma\in G$ realizing the descent of $D_L=\operatorname{End}_L(W)$, with $\xi_1=1$.
Write $\xi_\gamma\xi_\delta=c_{\gamma,\delta}\xi_{\gamma\delta}$.
The crossed product of the resulting normalized $2$-cocycle $c$ acts on $W$ and has commutant $D$, so the Brauer class of $c$ is $-[D]$.
Let $\mathcal M_i$ be the exact category of finite-dimensional $L$-vector spaces $V$ equipped with $\gamma$-semilinear isomorphisms $\zeta_\gamma:V\to V$, for $\gamma\in G$, satisfying
\[
                     \zeta_\gamma\zeta_\delta=c_{\gamma,\delta}^{i}\zeta_{\gamma\delta}.
\]
Morphisms are $L$-linear maps commuting with the $\zeta_\gamma$.
The category $\mathcal M_i$ is the category of left modules over the crossed product of $c^i$.
The opposite algebra of this crossed product has Brauer class $i[D]$, so $\mathcal M_i$ is Morita equivalent to the category of right $D^i$-modules.
Tensor products add the exponents of the multipliers, and $a$th exterior powers of flags multiply them by $a$.
After splitting, a right $D^i$-module has the form $V\otimes_L(W^*)^{\otimes i}$.
Since $(W^*)^{\otimes i}$ has multiplier $c^{-i}$, the multiplier of $V$ is $c^i$.
The variety $X_L$ parametrizes lines in $W^*$, so $\mathcal O(-1)$ has multiplier $c^{-1}$, which cancels in $V\otimes_L\mathcal O(-i)$.
This ordinary descent defines $\phi_i$.
The $\zeta_\gamma$ are compatible with the natural isomorphisms for tensor products and exterior powers of flags, so the natural isomorphisms descend to natural isomorphisms on the original exact categories.
This system of graded exact categories gives exterior-power operations $\lambda_i^a$ and Adams operations $\psi_i^a:K_n(D^i)\to K_n(D^{ai})$ on the groups \cite[Definition 1.1, Definition 2.1, Proposition 2.2, (2.8), Theorem 4.1]{KockZanchettaPower}.
The Newton polynomial has total tensor degree $a$, and in positive degrees $\psi_i^a=(-1)^{a-1}a\lambda_i^a$; the target of $\psi_i^a$ is therefore $K_n(D^{ai})$.
Applying the functorial construction of exterior-power operations to the natural isomorphisms above gives the displayed identity for the Adams operation on $K_n(X)$ defined by exterior powers.

It remains to identify this Adams operation with the stable operation $\psi_X^a$ on $K_n(X)[1/a]$.
Precomposing the stable operation with $KGL\to KGL[1/a]$, and composing the operation defined by the Newton polynomial in exterior powers with the map to $[1/a]$-coefficients, we obtain two pointed maps $\Z\times\mathrm{Gr}\to\Omega^\infty KGL[1/a]$.
Both send the class of a line bundle to its $a$th power, so by the splitting principle they induce the same natural transformation on $K_0$ of smooth $F$-schemes.
The proof of \cite[Lemma 1.2.10]{Riou2010} that restriction on $K_1$ of finite Grassmannians is surjective remains valid after inverting $a$, so the target satisfies the condition of \cite[Theorem 1.1.6]{Riou2010}, and the two maps agree.
Hence they agree on every $K_n\to K_n[1/a]$, and by the universal property of localization of coefficients they agree on $K_n(X)[1/a]$ \cite[Example 9.6]{AnnalaHoyoisIwasa}.
\end{proof}

\begin{proof}[Proof of Proposition~\ref{prop:paper-brauer-power-fiber}]
The assertion is immediate if $a=1$ or $\deg D=1$, so we assume that $a\ge2$ and $\deg D\ge2$, and we use the notation above.
By Lemmas~\ref{lem:paper-brauer-power-correspondence}, \ref{lem:paper-brauer-power-adams} and~\ref{lem:paper-brauer-power-descent}, the equivalence of spectra
\[
 (\phi^{(a)})^{-1}\psi_X^a\phi:P[1/a]\xrightarrow{\sim}P_a[1/a]
\]
induces, on all homotopy groups, a diagonal matrix preserving the index $i$.
Each diagonal component is an isomorphism of groups, so the map of spectra defined using inclusion and projection,
\[
 \Psi_D^{(a)}=\operatorname{pr}_1(\phi^{(a)})^{-1}\psi_X^a\phi\operatorname{in}_1:
 K(D)[1/a]\longrightarrow K(D^{\otimes a})[1/a]
\]
is an equivalence.

Now assume that $F$ is a global field.
If $\ell\mid \deg D$ and $\ell\ne\operatorname{char}F$, then $\ell\nmid a$, since $\per D$ and $\deg D$ have the same prime divisors.
We localize the equivalence above at $\ell$ and use naturality of the localization unit for $L_{E(1),\ell}$ to obtain $B_\ell(D)\simeq B_\ell(D^{\otimes a})$.
Unit powers preserve the index, so combining the finitely many factors of~\eqref{eq:paper-global-fiber-primes} gives $Z_D\simeq Z_{D^{\otimes a}}$ without inverting $a$.
If there is real ramification, invert $2$ as in Theorem~\ref{thm:paper-global-comparison}.

Since generators of the same cyclic subgroup of $\Br(F)$ are related by a unit power, the first consequence follows.
If $\gcd(i,\per D)=\gcd(j,\per D)$, then $i[D]$ and $j[D]$ generate the same cyclic subgroup, so the second consequence also follows.
\end{proof}

\subsubsection{Examples with vanishing comparison fiber at arbitrary index}
Example~\ref{ex:paper-global-fiber-rational} gives algebras of arbitrary index whose comparison fiber vanishes.

\begin{example}\label{ex:paper-global-fiber-rational}
Let $k$ be a finite field and put $F=k(t)$.
Let $k_n/k$ be the extension of degree $n$, and let $\sigma$ be its Frobenius.
Then
\[
 D=(k_n(t)/k(t),\sigma,t)
\]
has index $n$, and the comparison $g$ in positive degrees is an equivalence.
\end{example}
\begin{proof}
The local invariant at $t=0$ has order $n$, so the cyclic algebra $D$ has index $n$.
Put $\mathcal A=k_n[u,u^{-1};\sigma]$ and identify its center with $k[t,t^{-1}]$ by $t=u^n$.
After base change to $k_n$, the coefficient algebra decomposes into diagonal idempotents, and $u$ becomes a cyclic permutation matrix whose full cycle multiplies by $t$.
These matrices generate the matrix algebra, so $\mathcal A$ is Azumaya over $k[t,t^{-1}]$.

By the twisted Bass--Heller--Swan decomposition \cite[Corollary 5.3]{KirsteinKremer} and the vanishing of the Nil terms \cite[Corollary 4.2]{KirsteinKremer}, which follows from regularity of the coefficient field $k_n$, we have
\[
 K(\mathcal A)\simeq\cofib\!\left(1-\sigma_*:K(k_n)\longrightarrow K(k_n)\right).
\]
Height-one localization is exact, and taking fibers is an exact functor in a stable category, so the same decomposition holds for $B_\ell$.
The bound $B_\ell(k_n)\in\mathrm{Sp}_{\le-2}$ for finite fields \cite[Theorem 6.18]{ClausenMathew} gives $B_\ell(\mathcal A)\in\mathrm{Sp}_{\le-1}$.
In the localization sequence for $K(\mathcal A)\to K(D)$, the support term is the direct sum of $K(k_v)$ over the closed points $v$ of $\G_{m,k}$.
Indeed, the Azumaya algebra splits over each finite residue field $k_v$, so Morita equivalence and d\'evissage apply \cite[\S5]{QuillenHigher}.
Since $L_{E(1),\ell}$ commutes with direct sums, applying $B_\ell$ to this localization sequence gives a fiber sequence whose support term is the direct sum of the spectra $B_\ell(k_v)$.
Each $B_\ell(k_v)$ is in degrees at most $-2$, so $B_\ell(D)$ is in degrees at most $-1$.
All factors in~\eqref{eq:paper-global-fiber-primes} therefore vanish.
\end{proof}

\begin{remark}
Quillen's finite-field calculation gives $K_{2m-1}(\F_q)=\Z/(q^m-1)$ \cite{QuillenFinite}, so every connective cover of the $K$-theory of a finite field is unbounded above.
Moreover, every nonzero upper truncation $\tau_{\le N}\tau_{\ge1}K(\F_q)$ has a nonzero homotopy group in degree at most three.
By Theorem~\ref{thm:paper-global-fiber-bounds}, the comparison fiber is bounded above, and by Corollary~\ref{cor:global-fibre-connectivity}, it lies in degrees at least four for a global function field.
It remains open to describe the homotopy groups and Postnikov invariants of the comparison fiber in this finite range of degrees for global function fields.
\end{remark}

\subsubsection{Eilenberg--Mac Lane decompositions over number fields}\label{subsec:paper-global-fiber-gem}
We prove the odd-primary part of Theorem~\ref{thm:paper-global-fiber-gem} here and use the Bott element for the $2$-primary part in the next two subsubsections.

A finite generalized Eilenberg--Mac Lane spectrum (finite GEM) is a finite direct sum of suspensions of Eilenberg--Mac Lane spectra.
Let $F$ be a number field, let $D$ be a central division algebra over $F$, and let $Z$ be the fiber of the comparison in Theorem~\ref{thm:paper-global-comparison}.
If there is real ramification, take the comparison and $Z$ with $2$ inverted.

\begin{theorem}\label{thm:paper-global-fiber-gem}
For each prime $\ell\mid\deg D$, put $e_\ell=\ell^{v_\ell(\deg D)}$.
There is a noncanonical equivalence of spectra
\[
 Z \simeq
 \bigoplus_{\ell\mid\deg D}\ \bigoplus_{j=4}^{e_\ell-3}
           \Sigma^j H(\pi_jZ_{(\ell)}).
\]
If there is real ramification, the sum omits $\ell=2$; an empty degree range contributes zero.
Each $\pi_jZ_{(\ell)}$ is a finite $\ell$-group and, by the splitting in Theorem~\ref{thm:paper-global-split}, is the kernel of $K_j(D)_{(\ell)}\to K_j(F)_{(\ell)}$.
\end{theorem}

We prove Theorem~\ref{thm:paper-global-fiber-gem} after the following lemmas.
Fix a prime $\ell$ as in the theorem and put $R=K(F)_{(\ell)}$ and $\mathsf{X}=\tau_{\ge2}B_\ell(D)$.

\begin{lemma}\label{lem:paper-gem-positive-action}
There is a finite layer $E$ of the cyclotomic $\Z_\ell$-extension of $F$, with $G=\Gal(E/F)$, such that the following hold.
\begin{enumerate}
\item The action map $\tau_{\ge2}R\otimes_R\mathsf{X}\to\mathsf{X}$ is nullhomotopic as a map of $R$-modules.
\item The map
\[
 (\tau_{\ge1}K(E)_{(\ell)})_{hG}\otimes_R\mathsf{X}\longrightarrow\mathsf{X},
 \qquad y\otimes z\longmapsto\mathcal N_{E/F}(y)\cdot z,
\]
is nullhomotopic as a map of $R$-modules.
\end{enumerate}
\end{lemma}
\begin{proof}
By Theorem~\ref{thm:paper-global-fiber-bounds}, $\mathsf X\simeq Z_{(\ell)}$ as spectra.
We use the $R$-module structure on $\mathsf X$ induced by the localization fiber defining $B_\ell(D)$.
The localization map is a module map, and restriction, transfer, and their projection formula also induce maps and identities on fibers.

We take $F_\infty$ to be the compositum of $F$ with the cyclotomic $\Z_\ell$-extension of $\Q$, so every real place splits completely.
We take a sufficiently high finite layer $E$ in the cyclotomic $\Z_\ell$-extension $F_\infty/F$, and put $G=\Gal(E/F)$.
The decomposition group of every finite place is open in $\Gal(F_\infty/F)$.
Indeed, away from $\ell$, the Frobenius of the place acts through its absolute norm, whose image in the $\Z_\ell$-quotient of $\Z_\ell^\times$ is nonzero; above $\ell$, the ramification indices in the local cyclotomic extension are unbounded.
We may therefore choose $E$ so that, at each of the finitely many places where the $\ell$-primary Brauer class of $D$ is nontrivial, the local degree is divisible by the order of the local invariant of that class.
The invariants of this primary component vanish at real places, so Brauer--Hasse--Noether implies that $\ind(D_E)$ is prime to $\ell$.
Applying the argument that proves the first relation of \eqref{eq:paper-global-fiber-amplitude} to $D_E$ over $E$, with a separable splitting field of degree $\ind(D_E)$, we obtain $\ind(D_E)\,\pi_jB_\ell(D_E)=0$ for $j\ge1$.
Since $\ind(D_E)$ is prime to $\ell$ and $B_\ell(D_E)$ is $\ell$-local, we have $\tau_{\ge2}B_\ell(D_E)=0$.

Since $R$ is connective and $\mathsf X=\tau_{\ge2}B_\ell(D)$, by adjunction for connective covers of $R$-modules, a map of $R$-modules from a $2$-connective $R$-module to $\mathsf X$ is nullhomotopic once its composite with $\mathsf X\to B_\ell(D)$ is.
Restriction $\mathsf{X}\to B_\ell(D_E)$ is nullhomotopic as a map of $G$-equivariant $R$-modules, since $\mathsf{X}$ lies in degrees at least $4$ and $B_\ell(D_E)$ in degrees at most $1$, so the $t$-structure on equivariant modules makes the mapping space contractible.
By the projection formula, the map
\[
 (\tau_{\ge2}K(E)_{(\ell)})_{hG}\otimes_R \mathsf{X}\ \longrightarrow\ B_\ell(D),
 \qquad y\otimes z\ \longmapsto\
       N_{D_E/D}(y\cdot\res_{E/F}z)
\]
is nullhomotopic and equals the map $y\otimes z\mapsto\mathcal N_{E/F}(y)\cdot z$ given by the action through the derived norm.
Using the equivalence $(\tau_{\ge2}K(E)_{(\ell)})_{hG}\simeq\tau_{\ge2}R$ of Proposition~\ref{prop:paper-derived-norm} and factorization of the connective tensor product through the truncation, we see that the composite of the action map in~(1) with $\mathsf X\to B_\ell(D)$ is nullhomotopic.
The source of the action map in~(1) is $6$-connective, since $\mathsf X$ lies in degrees at least $4$, so the principle above gives~(1).
The same projection formula holds with $\tau_{\ge2}K(E)_{(\ell)}$ replaced by $\tau_{\ge1}K(E)_{(\ell)}$, so the composite of the map in~(2) with $\mathsf{X}\to B_\ell(D)$ is nullhomotopic as a map of $R$-modules.
The source of the map in~(2) is $5$-connective, since $\mathsf{X}$ lies in degrees at least $4$.
The same principle gives~(2).
\end{proof}

\begin{lemma}\label{lem:paper-gem-retract}
A retract of a finite GEM is a finite GEM.
Moreover, for every prime $\ell$, a bounded-above and bounded-below $H\Z_{(\ell)}$-module splits as a module into the direct sum of the suspended Eilenberg--Mac Lane spectra of its homotopy groups.
\end{lemma}
\begin{proof}
Order the terms by increasing degree and combine terms of the same degree.
Maps from higher to lower degrees vanish by the $t$-structure, so the ideal of off-diagonal terms in the endomorphism ring is nilpotent.
Let $\varepsilon$ be the idempotent defining the retract and $\varepsilon_0$ its diagonal part.
Then $u=\varepsilon\varepsilon_0+(1-\varepsilon)(1-\varepsilon_0)$ is congruent to $1$ modulo the off-diagonal ideal, hence invertible, and satisfies $\varepsilon u=u\varepsilon_0$.
Thus $\varepsilon=u\varepsilon_0u^{-1}$ is conjugate to $\varepsilon_0$.
The diagonal idempotents specify direct summands of the individual coefficient groups.
For the second assertion, $\Z_{(\ell)}$ has global dimension $1$, so the Ext groups of degree at least $2$ containing the Postnikov extensions vanish.
\end{proof}

\begin{proof}[The odd-primary part of Theorem~\ref{thm:paper-global-fiber-gem}]
Tensor $\tau_{\ge2}R\to R\to\tau_{\le1}R$ with $\mathsf{X}$.
By Lemma~\ref{lem:paper-gem-positive-action}(1), the map $\mathsf{X}\to \tau_{\le1}R\otimes_R \mathsf{X}$ is a split monomorphism of $R$-modules.
For odd $\ell$, one has $\tau_{\le1}\mathbb S_{(\ell)}\simeq H\Z_{(\ell)}$.
Truncating the unit map gives $H\Z_{(\ell)}\to\tau_{\le1}R$, so $\tau_{\le1}R\otimes_R \mathsf{X}$ is an $H\Z_{(\ell)}$-module.
We truncate the entire retract diagram to the finite range of degrees supporting the homotopy groups of $\mathsf{X}$.
Lemma~\ref{lem:paper-gem-retract} shows that $\mathsf{X}$ is also a finite GEM.
\end{proof}

\subsubsection{The Bott element and the \texorpdfstring{$2$}{2}-primary decomposition}
The $2$-primary part rests on two lemmas: Lemma~\ref{lem:paper-gem-bott} on the Bott element and Lemma~\ref{lem:paper-gem-mod-two} on the decomposition modulo two.

\begin{lemma}\label{lem:paper-gem-bott}
Let $\ell=2$.
Let $\widetilde\eta\in\pi_2(\mathbb S/2)$ be an element whose Bockstein is the Hopf element $\eta$, and let $b\in K_2(F;\Z/2)$ be its image under the unit map.
The maps $\Sigma^2 \mathsf{X}\to \mathsf{X}/2$ induced by $b$ and $\Sigma \mathsf{X}\to \mathsf{X}$ induced by $\eta$ are nullhomotopic as maps of $R$-modules.
\end{lemma}
\begin{proof}
For the finite cyclotomic layer $E/F$ of Lemma~\ref{lem:paper-gem-positive-action}, $-1\in N_{E/F}E^\times$.
Indeed, at each finite place the local reciprocity map $F_v^\times\to\Gal(F_\infty/F)$ kills finite-order elements, since its target is a torsion-free subgroup of $\Z_2$.
The real places split completely.
Thus $-1$ is a local norm at every place and is a global norm by the Hasse norm theorem for cyclic extensions \cite[\S1]{HasseNorm1931}.

With integral coefficients, put $Y=(\tau_{\ge1}K(E))_{hG}$.
The derived norm in degrees at least $2$ makes $\pi_2Y\to K_2(F)$ surjective.
In degree $1$, $\pi_1Y=(E^\times)_G$, and cyclic Hilbert 90 makes the field norm $\pi_1Y\to F^\times$ injective.
The preimage of $-1$ therefore lies in $\pi_1Y[2]$.
Under the sphere unit, modeled by permutations of finite sets, the transposition representing $\eta$ becomes a matrix of determinant $-1$, so $\eta$ maps to $-1\in K_1(F)$ and the Bockstein of $b$ is $-1$.
We lift the preimage of $-1$ through the Moore coefficient sequence
\[
 0\longrightarrow \pi_2Y/2\longrightarrow\pi_2(Y/2)
   \longrightarrow\pi_1Y[2]\longrightarrow0.
\]
The norm of this lift differs from $b$ by a class in $K_2(F)/2$, and surjectivity of $\pi_2Y\to K_2(F)$ lets us correct the lift, so $b$ lifts to $\pi_2(Y/2)$.
We apply Lemma~\ref{lem:paper-gem-positive-action}(2) to $Y_{(2)}$ and reduce modulo $2$.
It follows that the action of $b$ is zero.
Composing with the Bockstein shows that the action of $\eta$ is zero as well.
\end{proof}

\begin{lemma}\label{lem:paper-gem-mod-two}
The spectrum $\mathsf{X}/2$ is a finite direct sum of suspensions of $H\F_2$.
\end{lemma}
\begin{proof}
Put $\mathsf{W}=\tau_{\le1}R\otimes_R(\mathsf{X}/2)$.
Reducing the nullhomotopy of Lemma~\ref{lem:paper-gem-positive-action}(1) modulo $2$ makes $\mathsf{X}/2$ a retract of $\mathsf{W}$.
If $\iota$ denotes the inclusion of the bottom cell of the Moore spectrum and $\partial$ the projection to its top cell, then $2\operatorname{id}_{\mathbb S/2}=\iota\eta\partial$ \cite[Proposition 4]{SchwedeMoore}.
Indeed, a self-map whose composite with $\iota$ is zero comes from the order-two group $\pi_1(\mathbb S/2)$, and its nonzero element is $\iota\eta\partial$.
The inequality $2\operatorname{id}_{\mathbb S/2}\ne0$ follows from the action of $Sq^2$ on the bottom degree of $\mathbb S/2\wedge\mathbb S/2$ and the Cartan formula.
Lemma~\ref{lem:paper-gem-bott} therefore gives $2\operatorname{id}_{\mathsf{X}/2}=0$.

We regard $\mathsf{W}$ as a $ko_{(2)}$-module through the ring maps
\[
 ko_{(2)}\longrightarrow \tau_{\le1}ko_{(2)}
    \simeq \tau_{\le1}\mathbb S_{(2)}\longrightarrow\tau_{\le1}R.
\]
The middle equivalence comes from the isomorphisms induced by the sphere unit on $\pi_0,\pi_1$.
The action of $\eta$ on $\mathsf{W}$ is obtained from its action on $\mathsf{X}$ in Lemma~\ref{lem:paper-gem-bott} by base change to $\tau_{\le1}R$ and reduction modulo $2$, so it is zero.
Wood's sequence $\Sigma ko\xrightarrow{\eta}ko\to ku$ makes $\mathsf{W}$ a retract of $\mathsf{V}=ku_{(2)}\otimes_{ko_{(2)}}\mathsf{W}$.
Here the connective version of Wood's sequence follows from the two-cell description of $\Sigma^{-2}\mathbb{CP}^2$ in the periodic version and comparison with complexification on connective covers.
In degrees $0,1$, use that $\pi_0ko\to\pi_0ku$ is an isomorphism and $\eta:\pi_0ko\to\pi_1ko$ is surjective \cite[Theorem 9.7]{MathewNaumannNoel2017}.

The image $b_0\in\pi_2(ko/2)$ of $\widetilde\eta$ maps under complexification to the nonzero element of $\pi_2(ku/2)=\F_2$, namely the reduction of the complex Bott element $\beta\in\pi_2ku$.
Indeed, in Wood's sequence reduced modulo $2$, the map $\eta\colon\pi_0(ko/2)\to\pi_1(ko/2)$ is an isomorphism, so $\pi_2(ko/2)\to\pi_2(ku/2)$ is surjective with kernel $\eta\pi_1(ko/2)$, generated by the reduction of the integral class $\eta^2$.
The Bockstein vanishes on this kernel but takes $b_0$ to $\eta\ne0$, so $b_0$ does not lie in the kernel.

Since $b_0$ comes from the sphere, its action $\Sigma^2\mathsf{W}\to\mathsf{W}/2$ is obtained by base change and reduction modulo $2$ from the action in Lemma~\ref{lem:paper-gem-bott}.
After base change to $ku_{(2)}$, the composite $\Sigma^2\mathsf{V}\xrightarrow{\beta}\mathsf{V}\to \mathsf{V}/2$ is therefore zero.
Since $2\operatorname{id}_{\mathsf{V}}=0$ (the nullhomotopy of $2\operatorname{id}_{\mathsf X/2}$ is $R$-linear, so base change gives a $ku_{(2)}$-linear one), the map $\mathsf{V}\to \mathsf{V}/2$ is a split monomorphism of $ku_{(2)}$-modules, so $\beta_{\mathsf{V}}=0$.
As $ku_{(2)}/\beta\simeq H\Z_{(2)}$, the spectrum $\mathsf{V}$ is a retract of $H\Z_{(2)}\otimes_{ku_{(2)}}\mathsf{V}$.
Consequently, $\mathsf{X}/2$ is a retract of an $H\Z_{(2)}$-module.
We truncate the entire retract diagram to the bounded degree range of $\mathsf{X}/2$.
Lemma~\ref{lem:paper-gem-retract} shows that $\mathsf{X}/2$ is a finite GEM.
Finally, $2\operatorname{id}_{\mathsf{X}/2}=0$ implies that every coefficient group is an $\F_2$-vector space.
\end{proof}

\subsubsection{From finite-coefficient to integral decompositions}
It remains to pass from decompositions with finite coefficients to integral ones, which completes the $2$-primary part of Theorem~\ref{thm:paper-global-fiber-gem}.

\begin{lemma}\label{lem:paper-gem-operations}
Let $M,M'$ be finite $\ell$-groups.
If $n\ge2$, then $\ell[HM,\Sigma^nHM']=0$.
\end{lemma}
\begin{proof}
Richter's account of Kochman's calculation shows that $H_q(H\Z;\Z_{(\ell)})$ is annihilated by $\ell$ for $q>0$ \cite[\S6.1, p.~297]{Richter2006}.
The cofiber of the unit $\mathbb S\to H\Z$ lies in degrees at least $2$, and smashing with $H\Z_{(\ell)}$ preserves this connectivity.
Thus $H_0(H\Z;\Z_{(\ell)})=\Z_{(\ell)}$ and $H_1(H\Z;\Z_{(\ell)})=0$.

We may assume $M=\Z/\ell^a$.
By Lemma~\ref{lem:paper-gem-retract}, each finite-degree truncation of $H\Z_{(\ell)}\wedge H\Z$ splits into suspended Eilenberg--Mac Lane spectra.
In the cofiber of multiplication by $\ell^a$, the free summand in degree zero gives $H(\Z/\ell^a)$, and each summand $\Sigma^qHA$ in positive degree gives $\Sigma^qHA\oplus\Sigma^{q+1}HA$, since $\ell A=0$.
Thus $H_0(HM;\Z_{(\ell)})=\Z/\ell^a$, $H_1(HM;\Z_{(\ell)})=0$, and $H_q(HM;\Z_{(\ell)})$ is annihilated by $\ell$ for $q\ge2$; in particular, no extension in positive degrees produces $\Z/\ell^2$.
We compute $[HM,\Sigma^n HM']$ by the extension--restriction adjunction for $H\Z_{(\ell)}$-modules.
Its summands are
\[
 \operatorname{Hom}_{\Z_{(\ell)}}(H_n(HM;\Z_{(\ell)}),M'),
 \qquad
 \operatorname{Ext}^1_{\Z_{(\ell)}}(H_{n-1}(HM;\Z_{(\ell)}),M').
\]
Here we again use the splitting afforded by global dimension $1$.
For $n\ge2$, both summands are annihilated by $\ell$.
\end{proof}

\begin{lemma}\label{lem:paper-gem-detection}
Let $T$ be a bounded-above and bounded-below spectrum such that every $\pi_iT$ is a finite $\ell$-group.
If $T/\ell$ is a direct sum of suspensions of $H\F_\ell$, then $T$ is a finite GEM.
\end{lemma}
\begin{proof}
If a map from a direct sum of suspensions of $H\F_\ell$ is surjective on $\pi_*$, then its target has the same form and the map admits a section.
Indeed, lift a homogeneous basis of the target to the source and realize it by maps from factors of the same degree; this gives an equivalence onto the target.

Put $T_m=\tau_{\le m}T$.
The map $T/\ell\to T_m/\ell$ is surjective on every homotopy group: it is an isomorphism in degrees at most $m$, and in degree $m+1$ the Bockstein gives a surjection onto $\pi_mT[\ell]$.
Thus $T_m/\ell$ is also such a direct sum.

Inductively assume that $T_{m-1}$ is a finite GEM.
Put $M'=\pi_mT$ and write the Postnikov invariant as $k:T_{m-1}\to\Sigma^{m+1}HM'$.
The composites of the maps from $T/\ell$ show that $T_m/\ell\to T_{m-1}/\ell$ is likewise surjective on all homotopy groups.
By the observation above it admits a section, so $k/\ell=0$.
By naturality of reduction modulo $\ell$, the composite of $k$ with the reduction $\Sigma^{m+1}HM'\to\Sigma^{m+1}HM'/\ell$ factors through $k/\ell=0$, so the coefficient cofiber sequence makes $k$ a multiple of $\ell$.
The degree difference from each summand of $T_{m-1}$ to $\Sigma^{m+1}HM'$ is at least $2$, so Lemma~\ref{lem:paper-gem-operations} gives $k=0$.
\end{proof}

\begin{proof}[Completion of the proof of Theorem~\ref{thm:paper-global-fiber-gem}]
For each prime in the decomposition of Theorem~\ref{thm:paper-global-fiber-bounds}, the groups $\pi_jZ_{(\ell)}$ are finite $\ell$-groups supported in degrees $4\le j\le e_\ell-3$ by \eqref{eq:paper-global-fiber-amplitude} and Corollaries~\ref{cor:global-fibre-connectivity} and~\ref{cor:paper-global-fiber-finite}.
When the $2$-primary factor occurs, Lemmas~\ref{lem:paper-gem-mod-two} and~\ref{lem:paper-gem-detection} applied to $T=\mathsf X\simeq Z_{(2)}$ give the $2$-primary part.
\end{proof}

\subsubsection{The first remaining group}
By Corollary~\ref{cor:global-fibre-connectivity}, the first homotopy group of the comparison fiber that can be nonzero is $\pi_4Z$.
Over a number field, we describe $\pi_4Z_{(\ell)}$ in Proposition~\ref{prop:global-fibre-pi4} when $\ell$ is odd or $F$ has no real places.

\begin{proposition}\label{prop:global-fibre-pi4}
Let $F$ be a number field, let $\ell$ be a prime, and assume either that $\ell$ is odd or that $F$ has no real places.
There is a natural short exact sequence
\[
\begin{split}
 0\longrightarrow&
 \coker\bigl[\CH^1(F,5;D)_{(\ell)}
      \xrightarrow{d_2}\CH^2(F,4;D)_{(\ell)}\bigr]\\
 \longrightarrow&\,\pi_4Z_{(\ell)}
 \longrightarrow\CH^1(F,4;D)_{(\ell)}
 \longrightarrow 0.
\end{split}
\]
Here $d_2$ is the differential in the twisted motivic spectral sequence.
\end{proposition}
\begin{proof}
All groups in this proof are localized at $\ell$.
We index the spectral sequence of the homotopy coniveau tower $\mathcal{F}^rK(D)$ by weight $r$ and total degree $m$; then $d_s$ goes from $(r,m)$ to $(r+s-1,m-1)$.
Weight four in total degree four vanishes by Proposition~\ref{prop:global-two-chow-rows} and the Bass--Tate calculation.
Weight zero and weights at least five also vanish.

First compare weight three with the actual comparison map $g_4$.
Since $\pi_4\mathcal{F}^4K(D)=\pi_3\mathcal{F}^4K(D)=0$, one has
\[
 \pi_4 \mathcal{F}^3K(D)\simeq\CH^3(F,4;D).
\]
On the commutative side, $H^0_{\mot}(F,\Z_{(\ell)}(2))$ and $H^{-1}_{\mot}(F,\Z_{(\ell)}(2))$ vanish by vanishing of the rational groups and Beilinson--Lichtenbaum comparison with divisible coefficients.
Thus weight two vanishes in total degrees four and five.
Weight one also vanishes in these total degrees, so the commutative edge map $\CH^3(F,4)\to K_4(F)$ is an isomorphism.

With finite coefficients, Morita equivalence over $F^{\mathrm{sep}}$ and the motivic--\'etale comparison identify the geometric weight-$r$ layer of the \'etale sheafified homotopy coniveau tower $\mathcal{F}^rK(D)$ with $\Z/\ell^a(r)[2r]$.
The layer comparison is given by reduced rank.
Since these layers are concentrated in degree $2r$ and the connectivity of the tails of the support tower tends to infinity, the sheafified $r$th stage is the geometric $2r$-connective cover.
Multiplication by the rank-$1$ class therefore preserves the stages and acts as the identity on the weight-three layer.
Since $\cd_\ell F\le2$, the stages of weight at least four are $8$-connective at geometric points, $6$-connective on sections and $5$-connective after the derived inverse limit, so their $\pi_4$ and $\pi_3$ vanish as needed in the long exact sequence.
Consequently, the target of the weight-three comparison on $\pi_4$ is $H^2_{\mathrm{cont}}(F,\Z_\ell(3))$.
Thus, after derived completion, the composite $\pi_4\mathcal{F}^3K(D)\to K_4(D)\xrightarrow{g_4}K_4(F)$ equals the composite
\[
 \pi_4\mathcal{F}^3K(D)\simeq\CH^3(F,4;D)\longrightarrow\CH^3(F,4)\longrightarrow K_4(F),
\]
where $\CH^3(F,4;D)\to\CH^3(F,4)$ is the weight-three reduced norm and $\CH^3(F,4)\to K_4(F)$ is the commutative edge map \cite[\S6.9]{KahnLevine2010}, \cite{ClausenMathew}.

For an $\ell$-local spectrum, the map to its derived completion has rational fiber, so its kernel on homotopy groups is divisible.
By Lemma~\ref{lem:paper-global-no-divisible}, $K_4(F)_{(\ell)}$ therefore injects into $\pi_4(K(F)^\wedge_\ell)$, and the comparison above holds with $\Z_{(\ell)}$-coefficients.
Write $\mathcal{F}^3K_4(D)$ for the image of $\pi_4\mathcal{F}^3K(D)$ in $K_4(D)$.
Proposition~\ref{prop:global-two-chow-rows} with $r=3$ implies that $g_4$ maps $\mathcal{F}^3K_4(D)$ isomorphically onto $K_4(F)$.
In particular, all of weight three survives in $K_4(D)$, and $\pi_4Z_{(\ell)}\simeq K_4(D)_{(\ell)}/\mathcal{F}^3K_4(D)_{(\ell)}$.

Next we consider the differential from weight one,
\[
 d_2:\CH^1(F,4;D)\longrightarrow\CH^2(F,3;D).
\]
The comparison map from its source goes to the \'etale weight-one term $H^{-2}_{\et}(F,\Z/\ell^a(1))$, which vanishes because it is in negative cohomological degree.
The comparison on the target is injective: the weight-two reduced norm is an isomorphism by Proposition~\ref{prop:global-two-chow-rows}.
The commutative group $\CH^2(F,3)_{(\ell)}$ is finitely generated over $\Z_{(\ell)}$, so it is detected by its quotients modulo $\ell^a$, and by the coefficient sequence and Beilinson--Lichtenbaum comparison each quotient $\CH^2(F,3)_{(\ell)}/\ell^a$ injects into $H^1_{\et}(F,\Z/\ell^a(2))$.
Since the comparison with the \'etale spectral sequence commutes with $d_2$, the differential is zero.
The target of the next differential, $\CH^3(F,3;D)$, vanishes by Proposition~\ref{prop:global-two-chow-rows} with $r=3$ and \cite[Chapter II, Theorem 2.1(3)]{BassTate}, and higher targets vanish by the dimension condition on supports.
Differentials entering weight one also vanish, leaving $\CH^1(F,4;D)$ in weight one.
Differentials out of weight two vanish for the same reason, and the only incoming differential is $d_2$ from total degree five in weight one.
The resulting two-step filtration in weights one and two gives the displayed short exact sequence \cite[Corollary 6.1.4, \S6.9]{KahnLevine2010}.
\end{proof}

\begin{problem}\label{prob:global-fibre-vanishing}
Do all the finite groups $\pi_jZ_{(\ell)}$ in Theorem~\ref{thm:paper-global-fiber-gem} vanish over number fields?
A nonzero group requires $e_\ell\ge7$, since otherwise the range $4\le j\le e_\ell-3$ is empty, and $e_\ell\ne\ell$, since Theorem~\ref{thm:paper-global-fiber-bounds} gives $\tau_{\ge2}B_\ell(D)=0$ for $e_\ell=\ell$; so the smallest possible values of $e_\ell$ are $8=2^3$ and $9=3^2$.
In degree four, at odd primes or over a field without real places, Proposition~\ref{prop:global-fibre-pi4} reduces this question to conditions on the twisted higher Chow groups in that proposition.
\end{problem}

\section{Reduced norms over two-dimensional local fields}\label{sec:two-local}
This section proves that the reduced norm on $K_2$ is an isomorphism over two-dimensional local fields.
Surjectivity follows from Kato's existence theorem, and localization computes $\SK_2$ in tame situations (Section~\ref{sec:second-two-tame}).
The vanishing of $\SK_2$ in the remaining wild case is proved by an induction on the index (Section~\ref{sec:second-two-wild}).
\subsection{Norms from splitting fields and tame localization}\label{sec:second-two-tame}
In this subsection, we first show that, over a two-dimensional local field, the field norms from splitting fields generate $K_2(F)$ (Theorem~\ref{thm:second-two-surj}); in particular, the reduced norm on $K_2$ is surjective.
Under the hypotheses stated below, we then give localization formulas for $\SK_2(D)$ when $D/F$ is tamely ramified: in first residue characteristic zero (Proposition~\ref{prop:second-residue-zero}) and, for $\ell$-primary $D$ with $\ell$ different from the first residue characteristic, after $\ell$-localization (Proposition~\ref{prop:second-residue-zero-ell}).
Along the way we prove that $\SK_1(D)=0$ over two-dimensional local fields (Corollary~\ref{cor:review-two-local-sk1-zero}).
These results reduce the vanishing of $\SK_2$ over two-dimensional local fields to the wild case of Section~\ref{sec:second-two-wild}.

Let $F=F_h$ be an $h$-dimensional local field, with the fields $F_{h-1},\ldots,F_0$ as in the Conventions.
We write the characteristic sequence in the order $\operatorname{char} F_h,\ldots,\operatorname{char} F_1$.
For a central division algebra $D$ over $F$ and $\alpha=[D]\in\Br(F)$, let $N_2(\alpha)\subseteq K_2(F)$ be the norm subgroup of Definition~\ref{def:paper-transfer}.
We identify $K_2^M(F)$ with $K_2(F)$ by Matsumoto's theorem \cite{Matsumoto}.
The Morita-normalized splitting transfer $\theta_E$ satisfies $\Nrd_2\theta_E=N_{E/F}$, and hence $N_2(\alpha)\subset\Nrd_2K_2(D)$.

\begin{theorem}\label{thm:second-two-surj}
Let $F$ be a two-dimensional local field, let $D$ be a central division algebra over $F$, and put $\alpha=[D]$.
Then $N_2(\alpha)=K_2(F)$.
Consequently, $C_2(D/F)=0$.
\end{theorem}
\begin{proof}
We apply \cite[Proposition 5]{KatoExistence} to $\alpha\in H_m^2(F)=\Br(F)[m]$ with $m=\per(\alpha)$, using logarithmic coefficients for the characteristic-primary part.
Taking $n=q=r=2$ gives $q+r>n+1$, so the Milnor $K_2$-norms from finite extensions killing $\alpha$ generate $K_2(F)$.
\end{proof}

For the tame localization formulas we need the vanishing of $\SK_1$ over two-dimensional local fields, which we deduce from the following factorization of norm-one elements.

\begin{lemma}\label{lem:second-sl1-two-generators}
Let $F$ be a two-dimensional local field and $D$ a central division algebra over $F$.
Every element of $\operatorname{SL}_1(D)(F)$ is a product of two elements of $\operatorname{SL}_1(D)(F)$, each of which generates a separable maximal subfield of $D$.
\end{lemma}
\begin{proof}
The group $\operatorname{SL}_1(D)$ is smooth and geometrically connected, and open balls in \'etale coordinates for the nontrivial Henselian valuation on $F$ show that its $F$-points are Zariski dense.
The regular semisimple elements form a nonempty open subset $\Omega$, and an element of $\Omega(F)$ generates a separable maximal subfield.
For $g\in\operatorname{SL}_1(D)(F)$, the open set $\Omega\cap g\Omega^{-1}$ is nonempty, so by density it has an $F$-point $h$.
Then $g=h\cdot(h^{-1}g)$ with $h,h^{-1}g\in\Omega(F)$.
\end{proof}

\begin{corollary}
\label{cor:review-two-local-sk1-zero}
If $F$ is a two-dimensional local field and $D$ is a central division algebra over $F$, then $\SK_1(D)=0$.
\end{corollary}
\begin{proof}
Period equals index over a two-dimensional local field and over each of its finite extensions \cite[Corollary 5.1, Proposition 5.2]{Chipchakov2022}.
Moreover, Kato's duality and norm-image formula say that, for a Brauer class $\alpha$ of prime-power order, the subgroup $\alpha\cup F^\times\subset H^3(F,\Q/\Z(2))\simeq\Q/\Z$ has order $\per(\alpha)$, and the kernel of the cup product with $\alpha$ on $F^\times$ is the reduced-norm image.
Under the invariant map, restriction to a finite extension is multiplication by its degree.
The same formula with logarithmic coefficients applies to the characteristic-primary component \cite[\S3.1, Theorem 3(2)]{KatoLocalII}, \cite[Theorem 2]{KatoExistence}, \cite[Example 1.1(4)]{HuWu2020}.

The group $\SK_1(D)$ is killed by $\deg D$.
For the $\ell$-primary part, transfer makes restriction injective along an extension of degree prime to $\ell$ splitting the other Brauer components.
Morita equivalence then reduces to $\deg D=\ell^r$, and we induct on $r$ from the split case $r=0$.

By Lemma~\ref{lem:second-sl1-two-generators}, it is enough to treat a norm-one element $a$ that generates a separable maximal subfield.

Let $\widetilde L$ be the normal closure of $L=F(a)$ over $F$, let $G=\Gal(\widetilde L/F)$, and let $H$ be the subgroup fixing $L$.
We choose an $\ell$-Sylow subgroup $P$ of $G$ containing an $\ell$-Sylow subgroup of $H$, and pass to the extension $\widetilde L^P/F$ of degree prime to $\ell$.
The index and the degree of $L$ are unchanged, and restriction on $\SK_1$ remains injective.
A maximal subgroup of the $\ell$-group $P$ containing $H\cap P$ shows that, after this base change, $L$ contains a cyclic intermediate field $T/F$ of degree $\ell$.
We again denote the base field by $F$.

Put $B=Z_D(T)$ and $\alpha=[D]$.
Then $\deg B=\ell^{r-1}$, and period equals index, so $\per B=\ell^{r-1}$.
The element $b=\Nrd_{B/T}(a)$ satisfies $N_{T/F}b=1$; for a generator $\sigma$ of $\Gal(T/F)$, Hilbert 90 gives $b=c/\sigma(c)$ with $c\in T^\times$.
The order of $\alpha_T\cup c$ divides $\ell^{r-1}$, and the restriction of $\alpha\cup F^\times$ is also the subgroup of $\Q/\Z$ of order $\ell^{r-1}$.
We may therefore choose $u\in F^\times$ with $\alpha_T\cup(cu)=0$, and the norm-image formula then gives $y\in B^\times$ with $\Nrd_{B/T}(y)=cu$.
We choose $v\in D^\times$ such that conjugation by $v$ induces $\sigma$ on $T$.
Then
\[
 \Nrd_{B/T}([y,v])=c/\sigma(c)=b.
\]
Thus $a[y,v]^{-1}$ is a norm-one element of $B$, and its class in $K_1(B)$ vanishes by induction.
Since $[y,v]$ also vanishes in $K_1(D)$, so does the class of $a$.
We descend from the prime-to-$\ell$ extension and use the preceding factorization into two elements to conclude.
\end{proof}

We also need two lemmas, which are used again in Section~\ref{sec:second-three-coker}: surjectivity of reduction on $K_2$ for Azumaya algebras with division residue algebra over complete discrete valuation rings, and the vanishing of the cup product of the Brauer class with the image of the reduced norm.

\begin{lemma}\label{lem:second-three-azumaya-reduction}
Let $R$ be a complete discrete valuation ring with uniformizer $\pi$, and let $\widetilde D$ be an Azumaya algebra over $R$ whose residue algebra $\widetilde D/\pi\widetilde D$ is a division algebra.
Then $K_2(\widetilde D)\to K_2(\widetilde D/\pi\widetilde D)$ is surjective.
\end{lemma}
\begin{proof}
Put $I=\pi\widetilde D$, and carry out the argument of Lemma~\ref{lem:app-deep-unit-norm} over $R$.
Since $\widetilde D$ is Azumaya, its reduced trace pairing is perfect.
For units $s_i$ lifting a basis of the residue division algebra, the differential of $(x_i)\mapsto\prod_i[s_i,x_i]$ surjects onto $\ker\operatorname{Trd}$ on the residue, so the map to the smooth group $\operatorname{SL}_1(\widetilde D)$ is smooth at the identity.
Hensel's lemma expresses every norm-one principal unit as a product $\prod_i[s_i,x_i]$ with $x_i\in1+I$.
The relative group $K_1(\widetilde D,I)$ is generated by units in $1+I$, and, by the relative elementary relation $[1+bx]=[1+xb]$ used in the proof of Lemma~\ref{lem:app-orbit-stable-image}, the class of each $[s_i,x_i]$ is zero in $K_1(\widetilde D,I)$.
If a unit represents an element of the kernel of $K_1(\widetilde D,I)\to K_1(\widetilde D)$, then it has reduced norm one, so this map is injective.
The relative exact sequence $K_2(\widetilde D)\to K_2(\widetilde D/I)\to K_1(\widetilde D,I)\to K_1(\widetilde D)$ then makes the map $K_2(\widetilde D)\to K_2(\widetilde D/I)$ surjective.
\end{proof}

\begin{lemma}\label{lem:second-three-containment}
Let $K$ be a field, let $A$ be a central simple algebra over $K$, and put $\alpha=[A]$.
\begin{enumerate}
\item If $\operatorname{char}K=0$, then $\alpha\cup\Nrd_2(x)=0$ in $H^4(K,\Q/\Z(3))$ for every $x\in K_2(A)$.
If $\operatorname{char}K=p>0$, the same holds for the prime-to-$p$ components.
\item If $\operatorname{char}K=p>0$ and $A$ is a division algebra, then the $p$-primary component of $\alpha\cup\Nrd_2(x)$ vanishes for every $x\in K_2(A)$.
\end{enumerate}
\end{lemma}
\begin{proof}
(1) In characteristic zero, we use the \'etale motivic descent of Kahn--Levine \cite[\S6.9]{KahnLevine2010}.
The weight-two edge map of this descent is the reduced norm: this follows by passing to a splitting Severi--Brauer function field and using injectivity of restriction on $K_2(K)$.
By the module structure and the formula $d_2^A(1)=\alpha$ in \cite[Proposition 6.9.1]{KahnLevine2010}, the image of the edge map lies in the kernel of $d_2^A(x)=\alpha\cup x$.
In positive characteristic with invertible coefficients, the same calculation applies because passage to the perfect closure leaves the relevant Galois cohomology unchanged; here we also use the compatibility of $\Nrd_2$ with base change.

(2) Choose a complete Cohen discrete valuation ring $R$ with residue field $K$ and uniformizer $\pi$, and an Azumaya algebra $\widetilde D/R$ whose residue is $A$ \cite[Expos\'e XXIV, Proposition 8.1]{SGA3}.
By Lemma~\ref{lem:second-three-azumaya-reduction}, we can lift $x\in K_2(A)$ to $\widetilde x\in K_2(\widetilde D)$; let $y\in K_2(\operatorname{Frac}R)$ be its generic-point reduced norm.
Then $y$ comes from $K_2(R)$ and reduces to $\Nrd_2(x)$.
Indeed, after restriction to the discrete valuation ring of the generic point of the special fiber of the relative Severi--Brauer variety, where the algebra splits, the boundary of $y$ vanishes; since restriction on residue multiplicative groups is injective, the boundary of $y$ itself vanishes, and $y$ lifts to $K_2(R)$ by localization.
There the reduction on $K_2$ agrees with $s_\pi=\partial(\{\pi\}\cdot-)$ at the generic point, so two integral lifts with the same generic class have the same reduction; since $K_2(K)\to K_2(K(\SB(A)))$ is injective \cite[Corollaries 5.6--5.7]{Suslin1987}, the reduction of $y$ is $\Nrd_2(x)$.
The class $y$ has zero boundary, so it is a sum of unit symbols in $K_2(\operatorname{Frac}R)$.
Indeed, by Matsumoto's theorem it suffices to treat the symbols $\{\pi,1-\pi c\}$ with $c\in R$.
If $c\in R^\times$, the Steinberg relation gives $\{\pi,1-\pi c\}=-\{c,1-\pi c\}$.
If $c\in\pi R$, we write $1-\pi c=(1-\pi(1+c-\pi c))/(1-\pi)$ and apply the unit case to both factors.
Kato's unramified map $H^r(K)\hookrightarrow H^r(\operatorname{Frac}R)$ sends Artin--Schreier--Witt characters to unramified characters and sends products with unit symbols to products with arbitrary unit lifts \cite[\S5.1, pp.~55--56]{Kurihara2000}.
Thus $\alpha\cup\Nrd_2(x)$ maps to $[\widetilde D_{\operatorname{Frac}R}]\cup y$.
Part~(1), applied over the field $\operatorname{Frac}R$ of characteristic zero, makes this product zero.
\end{proof}

\begin{proposition}\label{prop:second-residue-zero}
Let $k$ be a local field of dimension $r\ge1$ with $\operatorname{char}k=0$, and let $F=k((t))$.
Let $B$ be the residue division algebra of a central division algebra $D/F$ (tamely ramified, since $\operatorname{char}k=0$), and let $L$ be the center of $B$.
The group $G=\Gal(L/k)$ is cyclic.
If $\SK_1(B)=0$, then
\[
 \SK_2(D)\simeq
 \ker\!\left(K_2(B)_G
   \xrightarrow{N_{L/k}\Nrd_{B/L,2}}K_2(k)\right).
\]
If $\Nrd_{B/L,2}$ is surjective, then $\SK_2(B)_G\to\SK_2(D)$ is surjective.
\end{proposition}

We keep the notation of Proposition~\ref{prop:second-residue-zero}, whose proof follows the next two lemmas.
In the Witt decomposition $[D]=a+(\chi,t)$ with $a\in\Br(k)$ and $\chi\in H^1(k,\Q/\Z)$, we have $L=k(\chi)$ and $[B]=a_L$.
Let $R$ be the valuation ring of $D$.
By \cite[Theorem 1]{Zheglov2004} and its proof, the residue map $R\to B$ has a ring section $\iota\colon B\to R$, the coefficient embedding, which extends the fixed coefficient field $k\subset F=k((t))$; this justifies the identification $L\otimes_kF=L((t))$ used below.
Since $B$ is a separable $k$-algebra, the coefficient embedding over the fixed $k$ is obtained by lifting through the square-zero kernels of $R/J^{s+1}\to R/J^s$, where $J=J(R)$, and passing to the inverse limit.
Normalize $w$ on $D$ by $w(D^\times)=\Z$, so $w(t)=e$, where $e=[L:k]$.
Under the tame isomorphism $\Gamma_D/\Gamma_F\simeq\Gal(L/k)$, let $\sigma_L$ be the generator corresponding to the value $1$.
Since $[B]=a_L$ is restricted from $k$, the automorphism $\sigma_L$ lifts to a semilinear automorphism $\sigma$ of $B$.
Extending $\sigma$ to $B\otimes_kF$ and applying Skolem--Noether to the two resulting $F$-embeddings of $B\otimes_kF$ into $D$, we obtain $u\in D^\times$ with $u\,\iota(b)=\iota(\sigma(b))\,u$ for $b\in B$.
The action on the residue center gives $w(u)\equiv1\bmod e$, so after multiplying by a power of the central element $t$ we may assume that $u$ is a uniformizer of $D$.
With this uniformizer, $R=B[[u;\sigma]]$.
We write $\sigma_*$ for the automorphism of $K_*(B)$ induced by $\sigma$; since $\sigma^{[L:k]}$ is inner, $\sigma_*$ defines an action of $G$.

\begin{lemma}\label{lem:second-residue-zero-localization}
In this setup, let $i_*\colon K_*(B)\to K_*(R)$ be the map from closed support, and put $V_2=\ker(K_2(R)\to K_2(B))$.
Then $i_*=\iota_*\circ(1-(\sigma^{-1})_*)$, the group $V_2$ is uniquely divisible, and the localization sequence takes the form
\[
 0\longrightarrow K_2(B)_G\oplus V_2
 \xrightarrow{\ \iota_*+\mathrm{incl}\ } K_2(D)\xrightarrow{\partial_D}K_1(B)^G
 \longrightarrow0.
\]
\end{lemma}
\begin{proof}
The right $R$-module $B$ has the projective resolution $0\to uR\to R\to B\to0$ of length one.
The right linear isomorphism $R\to uR$, $r\mapsto ur$, twists the left $B$-action by $\sigma^{-1}$, since $\iota(b)ur=u\,\iota(\sigma^{-1}(b))r$.
Additivity therefore gives $i_*=\iota_*\circ(1-(\sigma^{-1})_*)$.
Since $\iota$ splits reduction, $K_2(R)=\iota_*K_2(B)\oplus V_2$, and $\iota_*$ is injective on $K_*(B)$.
Hence $\coker(i_*\colon K_2(B)\to K_2(R))=K_2(B)_G\oplus V_2$ and $\ker(i_*\colon K_1(B)\to K_1(R))=K_1(B)^G$.
The localization sequence $K_2(B)\to K_2(R)\to K_2(D)\to K_1(B)\to K_1(R)$ now gives the displayed sequence.
Finite-coefficient rigidity in residue characteristic zero makes $V_2$ uniquely divisible \cite[\S3, Corollary 2 to Theorem 2]{SuslinYufryakov}.
\end{proof}

\begin{lemma}\label{lem:second-residue-zero-norms}
In this setup, let $U\subset D$ be an unramified maximal subfield.
Then
\[
 \partial_F\circ\Nrd_{D,2}=\Nrd_{B/L,1}\circ\partial_D\quad\text{inside }\bar U^\times,
\]
and, for $x\in K_2(B)$,
\[
 \Nrd_{D,2}(\iota_*x)=\res_{F/k}N_{L/k}\Nrd_{B/L,2}(x).
\]
Moreover, $\res_{F/k}\colon K_2(k)\to K_2(F)$ is split injective.
\end{lemma}
\begin{proof}
Such a field $U$ exists: lift a separable maximal subfield $K$ of $B$ to $\iota(K)((t))$.
Let $\mathcal O_U$ be the valuation ring of $U$ and $\bar U$ its residue field, a maximal subfield of $B$ containing $L$.
Restriction of scalars from $R$ to $\mathcal O_U$ is compatible with the localization sequences, and after the split Morita equivalence it is $\res_{U/F}\Nrd_{D,2}$ at the generic point and $\res_{\bar U/L}\Nrd_{B/L,1}$ on closed support.
Since $e(U/F)=1$, the field-side boundary $K_2(U)\to K_1(\bar U)$ composed with $\res_{U/F}$ is $\partial_F$ followed by restriction to the residue field $\bar U$, and the resulting localization diagram gives the first identity.

For the second identity, put $M=L((t))$ and $C=Z_D(M)$.
By Skolem--Noether, $\sigma^{[L:k]}=\operatorname{Ad}(\beta)$ on $B$ for some $\beta\in B^\times$, so $\beta^{-1}u^{[L:k]}$ commutes with $B$ and the centralizer of $L$ in $R=B[[u;\sigma]]$ is $B[[\beta^{-1}u^{[L:k]}]]$.
Hence $C=B\otimes_LM$: it contains this algebra, and the double centralizer theorem shows that the degrees agree.
We use the Morita $(C,D_M)$-bimodule $D$.
After this Morita equivalence, restriction of scalars from $D_M$ to $D$ is extension of scalars along $C\to D$, so the composite of extension of scalars along $C\to D$ with $\Nrd_{D,2}$ is $N_{M/F}\Nrd_{C/M,2}$ \cite[Proposition 6.6.1(1)]{KahnLevine2010}.
Therefore
\[
 \begin{aligned}
 \Nrd_{D,2}(\iota_*x)
 &=N_{M/F}\Nrd_{C/M,2}(x)\\
 &=N_{M/F}\res_{M/L}\Nrd_{B/L,2}(x)\\
 &=\res_{F/k}N_{L/k}\Nrd_{B/L,2}(x)
 \end{aligned}
\]
by this formula, compatibility of reduced norms with base change, and the base-change formula for $L\otimes_kF=M$.
Finally, specialization at $t$ splits $\res_{F/k}$.
\end{proof}

\begin{proof}[Proof of Proposition~\ref{prop:second-residue-zero}]
Let $\xi\in\SK_2(D)$.
By Lemma~\ref{lem:second-residue-zero-norms}, $\Nrd_{B/L,1}(\partial_D\xi)=\partial_F(1)=1$, and $\SK_1(B)=0$ gives $\partial_D\xi=0$.
By Lemma~\ref{lem:second-residue-zero-localization}, $\xi$ therefore lies in $K_2(B)_G\oplus V_2$.
Its $V_2$-component vanishes because $\deg D$ kills $\SK_2(D)$ and $V_2$ is uniquely divisible.
Thus $\SK_2(D)$ is the kernel of the reduced norm on $K_2(B)_G$, which Lemma~\ref{lem:second-residue-zero-norms} and the injectivity of $\res_{F/k}$ identify with the displayed kernel.

For the final assertion, suppose that $\Nrd_{B/L,2}$ is surjective; we first show that $\SK_1(B)=0$ holds automatically, so that the displayed identification is available.
Since $L/k$ is finite, $L$ is a local field of dimension $r$, and $\operatorname{char}L=0$.
If $r\le2$, then $\SK_1(B)=0$ by Theorem~\ref{thm:paper-local-norm} for $r=1$ and by Corollary~\ref{cor:review-two-local-sk1-zero} for $r=2$.
If $r\ge3$, then Lemma~\ref{lem:second-three-containment}(1) and the surjectivity of $\Nrd_{B/L,2}$ give $[B]\cup K_2(L)=0$.
Writing a Milnor symbol of length $r-1$ as the product of its first two entries with the rest, we obtain $[B]\cup K_{r-1}^M(L)=0$, and Kato's top-degree duality gives $[B]=0$ \cite[Theorem 2]{KatoExistence}; hence $B=L$ and again $\SK_1(B)=0$.
We now use cyclic $K_2$ Hilbert 90 \cite[Theorem 14.1]{MerkurjevSuslin1982}: the norm $N_{L/k}\colon K_2(L)_G\to K_2(k)$ is injective.
Hence the kernel equals $\ker(K_2(B)_G\to K_2(L)_G)$, and right exactness of coinvariants shows that this kernel is the image of $\SK_2(B)_G$.
\end{proof}

When the first residue characteristic is $p>0$, the same kernel formula holds after $\ell$-localization for $\ell\ne p$.
In the next proposition, $F$ is a complete discretely valued field with residue field $k$ of characteristic $p$, $D$ is an $\ell$-primary central division algebra over $F$ with residue division algebra $B$, and $L$ is the center of $B$.
Since $D/F$ is tamely ramified ($\ell\ne p$), the conjugation $z\mapsto u^{-1}zu$ by a uniformizer $u$ of $D$ restricts to a generator of $\Gal(L/k)$.
Another normalized choice $u'$ gives a unit $u'u^{-1}$ of the valuation ring $\mathcal O_D$ of $D$, whose residue acts by an inner automorphism, so the induced action $\bar\sigma_*$ of $\Gal(L/k)$ on $K_2(B)$ does not depend on the choice.
This is the action in the statement.

\begin{proposition}\label{prop:second-residue-zero-ell}
Let $F$ be a complete discretely valued field whose residue field $k$ is a local field of dimension at most two with $\operatorname{char}k=p>0$.
Let $\ell\ne p$ be a prime, let $D$ be an $\ell$-primary central division algebra over $F$ (so $D/F$ is tamely ramified), let $B$ be its residue division algebra, and let $L$ be the center of $B$.
\begin{enumerate}
\item If $\SK_1(B)=0$, then
\[
 \SK_2(D)_{(\ell)}\simeq
 \ker\!\left(K_2(B)_{\Gal(L/k)}
   \xrightarrow{N_{L/k}\Nrd_{B/L,2}}K_2(k)\right)_{(\ell)}.
\]
\item If $k$ is a nonarchimedean local field, then $\SK_2(D)_{(\ell)}=0$.
\end{enumerate}
\end{proposition}
The proof of Proposition~\ref{prop:second-residue-zero-ell} is given after Lemma~\ref{lem:second-residue-zero-ell-integral} below.

Let $\mathcal O_D$ be the valuation ring of $D$, so that $B=\mathcal O_D/J(\mathcal O_D)$.
By the structure of tame division algebras, the ramification index $e$ of $D/F$ equals $[L:k]$ and divides $\deg D$, so it is a power of $\ell$.
Let $u$ be a uniformizer of $D$, let $\pi$ be a uniformizer of $F$, and write $\sigma_*$ and $\bar\sigma_*$ for the actions of conjugation $z\mapsto u^{-1}zu$ on $K_2(\mathcal O_D)$ and on $K_2(B)$.
Let $r\colon K_2(\mathcal O_D)\to K_2(B)$ be the reduction and $i_*\colon K_2(B)\to K_2(\mathcal O_D)$ the map from closed support.
There is no coefficient field in general, so we lift through an unramified centralizer.
Let $M/F$ be the unramified lift of $L/k$, and let $\mathcal O_0$ be a maximal order in the centralizer of $M$ in $D$.
Then $\mathcal O_0$ is an Azumaya algebra over $\mathcal O_M$ with residue algebra $B$.

\begin{lemma}\label{lem:second-residue-zero-ell-integral}
In the situation of Proposition~\ref{prop:second-residue-zero-ell}, the reduction $r$ is surjective, $\ker r$ is uniquely $\ell$-divisible, and
\[
 \im i_*\cap\ker r=0.
\]
\end{lemma}
\begin{proof}
Lemma~\ref{lem:second-three-azumaya-reduction}, applied to the Azumaya algebra $\mathcal O_0$ over the complete discrete valuation ring $\mathcal O_M$, whose residue algebra is the division algebra $B$, shows that $K_2(\mathcal O_0)\to K_2(B)$ is surjective; hence so is $r$.

By the rigidity of Suslin and Yufryakov \cite[\S3, Theorem 2 and Corollary 2]{SuslinYufryakov}, the fiber
\[
 V=\fib\bigl(K(\mathcal O_D)\to K(B)\bigr)
\]
is uniquely $\ell$-divisible.
The image $C=\im\bigl(K_3(B)\to\pi_2V\bigr)$ is $\ell$-torsion-free.
Moreover, $C/\ell^a$ is a quotient of $K_3(B)/\ell^a$, which injects into $K_3(B;\Z/\ell^a)$.
At each residue step, rigidity and localization give $|K_i(B;\Z/\ell^a)|\le|K_i(\bar B;\Z/\ell^a)|\cdot|K_{i-1}(\bar B;\Z/\ell^a)|$ for the residue algebra $\bar B$.
Since $k$ has dimension at most two, starting from $K_3(B)$ takes at most two residue steps and leaves only the groups $K_i(\F_q;\Z/\ell^a)$ with $i\ge1$.
Their orders are bounded independently of $a$, and hence so is the order of $C/\ell^a$.
Since $C$ is $\ell$-torsion-free, $|C/\ell^a|=|C/\ell|^a$, so the uniform bound forces $C/\ell=0$, and $\ker r\simeq\pi_2V/C$ is uniquely $\ell$-divisible.

The localization sequence gives
\[
 i_*\circ r=1-\sigma_*,\qquad r\circ i_*=1-\bar\sigma_*,
\]
where the first identity follows from additivity for the sequence of bimodules $J(\mathcal O_D)\to\mathcal O_D\to B$, and the second comes from the derived tensor product with the residue algebra, whose $\operatorname{Tor}_0$ and $\operatorname{Tor}_1$ give the two terms.
If $y=i_*b$ lies in $\ker r$, then $\bar\sigma_*b=b$, and hence $\sigma_*y=y$.
Since $\pi u^{-e}\in\mathcal O_D^\times$ and $\pi$ is central, conjugation by $u^{-e}$ is conjugation by a unit of $\mathcal O_D$, so $\sigma^e$ is inner on $\mathcal O_D$ and $\sigma_*^e=1$.
Since $r$ is surjective, the first identity shows that $1+\sigma_*+\cdots+\sigma_*^{e-1}$ kills $\im i_*$; applied to the $\sigma_*$-invariant element $y$, this gives $e\,y=0$.
Since $e$ is a power of $\ell$ and $\ker r$ is uniquely $\ell$-divisible, $y=0$.
\end{proof}

\begin{proof}[Proof of Proposition~\ref{prop:second-residue-zero-ell}]
We compare boundaries by restriction of scalars to the unramified maximal subfield $U\subset D$ obtained by lifting a separable maximal subfield of $B$ through $\mathcal O_0$; the resulting diagram commutes on all of $K_2(D)$, which we also use below for arbitrary integral lifts.
After Morita equivalence, restriction of scalars to the valuation ring $\mathcal O_U$ of $U$ is $\res_{U/F}\Nrd_{D,2}$ at the generic point and $\res_{\bar U/L}\Nrd_{B/L,1}$ on closed support, where $\bar U$ is the residue field of $U$.
Since $e(U/F)=1$, the field-side boundary is restriction to the residue field, so the boundary of $\SK_2(D)$ lies in $\SK_1(B)$.
Let $s_\pi=\partial([\pi]\cdot-)\colon K_2(F)\to K_2(k)$ be the specialization.
On lifts from $\mathcal O_0$, the composite $K_2(\mathcal O_0)\to K_2(D)\xrightarrow{\Nrd_2}K_2(F)\xrightarrow{s_\pi}K_2(k)$ equals the reduction $K_2(\mathcal O_0)\to K_2(B)$ followed by $N_{L/k}\Nrd_{B/L,2}$.
This is computed through the centralizer of $M$, at the generic point of the special fiber of the relative Severi--Brauer scheme where the algebra splits, as in the proof of Lemma~\ref{lem:second-three-containment}(2), together with the compatibility of specialization and transfer along the unramified extension $M/F$ for the common uniformizer $\pi$.
Let $y\in\ker r$.
After restriction to $U$ and Morita equivalence, the reduction of $y$ is computed from the filtration of $\mathcal O_D/\pi\mathcal O_D$ by the $J^i$ $(0\le i<e)$, where $J=J(\mathcal O_D)$, whose graded pieces are twists of $B$; hence it is a sum of twists of $r(y)=0$, and $s_\pi(\res_U\Nrd_2y)=0$.
Taking the norm from $U$ gives $\deg D\cdot s_\pi(\Nrd_2y)=0$, and since $\ker r$ is $\ell$-divisible and $\deg D$ is a power of $\ell$, $s_\pi(\Nrd_2y)=0$.
As $K_2(\mathcal O_0)\to K_2(B)$ is surjective, every class of $K_2(\mathcal O_D)$ is a lift from $\mathcal O_0$ plus an element of $\ker r$, so the formula above holds on all of $K_2(\mathcal O_D)$ after $\ell$-localization.

(1) Suppose that $\SK_1(B)=0$.
By localization, an element of $\SK_2(D)$, whose boundary vanishes, lifts to $K_2(\mathcal O_D)$.
By Lemma~\ref{lem:second-residue-zero-ell-integral}, reduction gives an exact sequence
\[
 0\longrightarrow\ker r\longrightarrow\coker i_*\longrightarrow K_2(B)_{\Gal(L/k)}\longrightarrow0,
\]
and $\coker i_*$ embeds into $K_2(D)$ by localization.
Since $\ker r$ is uniquely $\ell$-divisible and $\deg D$ kills $\SK_2(D)$, the reduction is injective on $\SK_2(D)_{(\ell)}$.
Hence $\SK_2(D)_{(\ell)}$ lies in the $\ell$-localized kernel of $N_{L/k}\Nrd_{B/L,2}$ on $K_2(B)_{\Gal(L/k)}$.
Conversely, working in the $\ell$-localized groups, so that the lift and the correction below are taken $\ell$-locally, the reduced norm of a lift of an element of the residue norm kernel lies in
\[
 J'_F=\ker\partial_F\cap\ker s_\pi\subseteq K_2(F),
\]
and, since $s_\pi$ agrees with the reduction $r_F\colon K_2(\mathcal O_F)\to K_2(k)$ on $\ker\partial_F$ by the projection formula, $J'_F$ is the image of $\ker r_F$.
On the commutative side, $i_*r_F=0$ and $r_F$ is surjective, so $i_*=0$ and $K_2(\mathcal O_F)\to K_2(F)$ is injective; hence $J'_F\simeq\ker r_F$, which is uniquely $\ell$-divisible by the same rigidity and boundedness.
Scalar extension followed by $\Nrd_{D,2}$ is multiplication by $\deg D$ on $J'_F$, so we may divide this norm by $\deg D$ in $J'_F$, extend scalars, and subtract to obtain a lift with zero reduced norm.
This correction comes from a class in $K_2(\mathcal O_F)$ with zero reduction, so it does not change the residue class in $K_2(B)_{\Gal(L/k)}$.
This gives the reverse comparison after $\ell$-localization.

(2) If $k$ is a nonarchimedean local field, then so is $L$, and Theorem~\ref{thm:paper-local-norm} gives $\SK_1(B)=0$ and the isomorphism $\Nrd_{B/L,2}$.
Cyclic $K_2$ Hilbert 90 \cite[Theorem 14.1]{MerkurjevSuslin1982} makes $N_{L/k}$ injective on $K_2(L)_{\Gal(L/k)}$, so the kernel in (1) vanishes.
\end{proof}

\subsection{Integral induction in wild ramification}\label{sec:second-two-wild}
The object of this subsection is the following theorem, which completes the computation of $K_2$ of division algebras over two-dimensional local fields.

\begin{theorem}\label{thm:second-two-vanishing}
Let $F$ be a two-dimensional local field.
For every central division algebra $D$ over $F$, one has, integrally,
\[
 \SK_2(D)=0,\qquad \Nrd_2\colon K_2(D)\xrightarrow{\ \sim\ }K_2(F).
\]
\end{theorem}

The proof of Theorem~\ref{thm:second-two-vanishing} is given at the end of this subsection, after that of Theorem~\ref{thm:second-two-wild}.

\begin{theorem}\label{thm:second-two-wild}
Let $k$ be a nonarchimedean local field of characteristic $p>0$, and let $F$ be a complete discretely valued field with first residue field $k$.
The field $F$ may have characteristic $p$ or zero.
For every $p$-primary central division algebra $D/F$, one has $\SK_2(D)=0$.
\end{theorem}

We outline the proof of Theorem~\ref{thm:second-two-wild}.
It is an induction on the minimal unramified index.
For a finite Galois unramified extension $F'/F$ realizing the minimal index, the kernel of restriction on $\SK_2$ is contained in the corestriction of $\SK_2(D')$, where $D'$ is the division algebra in the Brauer class of $D_{F'}$; hence it suffices to prove $\SK_2(D')=0$, and we may assume that $\deg D=\ind(D_{F^{\mathrm{nr}}})=p^q$.
If $q=0$, then $D$ is split.

If $q>0$, after adjoining $\mu_p$ if necessary and passing to a further finite Galois unramified extension, the cyclic cover lemma gives a cyclic extension $T/F$ of degree $p$ inside $D$ whose centralizer $B$ has minimal unramified index $p^{q-1}$.
This is where the induction hypothesis enters: applied to $B$, it gives $\SK_2(B)=0$, and together with the isomorphism of the reduced norm on splitting transfers it gives $K_2(B)=I_T(B)$.
This handles the generic images of the corner components of the orbits of the relative cover, and the lemma on orbit images handles the other components, so the generic image of the whole $\pi_2$ of the orbits lies in $I_F(D)$.

Finally, for $\xi\in\SK_2(D)$, we lift $\xi$ to an integral class with zero residue, lift its trace through the surjection on $\pi_2$ of degree-zero $\mathrm{TC}$, and remove the residue component.
At a fixed finite truncation we represent the remaining relative orbit class by an integral class and, in the mixed-characteristic homogeneous-unit type, add back an integral representative of the class subtracted with the residue.
The generic image of the resulting class lies in $I_F(D)$; its difference from the original lift vanishes at that truncation, so the generic image of the difference lies in $I_F(D)$ by congruence separation (Lemma~\ref{lem:app-congruence-separation}).
Hence $\xi\in I_F(D)$, and since $\Nrd_2$ is injective on $I_F(D)$, we conclude that $\xi=0$.

\subsubsection{Unramified descent and minimal index}
We first reduce Theorem~\ref{thm:second-two-wild} to division algebras whose index does not drop over the maximal unramified extension.
Recall from Definition~\ref{def:paper-transfer} the subgroup $I_F(D)=\sum_E\theta_E K_2(E)$, with $E/F$ ranging over finite splitting fields.

\begin{proposition}
\label{prop:review-splitting-transfer-norm-iso}
For a two-dimensional local field $F$ and any central simple algebra $D/F$, the map $\Nrd_2\colon I_F(D)\xrightarrow{\sim}K_2(F)$ is an isomorphism.
\end{proposition}
\begin{proof}
If $\ell\ne\operatorname{char}F$, then $\cd_\ell F\le3$, so Corollary~\ref{cor:review-small-cd-splitting-transfer} applies.
If $\operatorname{char}F=p$, then $[F:F^p]=p^2$, and Lemma~\ref{lem:review-p-rank-two-frobenius-norm} handles the remaining $p$-component.
Since $\Nrd_2\colon I_F(D)\to K_2(F)$ is an isomorphism after localization at every prime, the integral map is also an isomorphism.
\end{proof}

\begin{lemma}
\label{lem:review-minimal-order-k1-support}
Let $F$ be a two-dimensional local field whose first residue field $k$ has characteristic $p$, and let $D/F$ be a central division algebra of degree $n=p^q$.
Assume $\ind D_{F^{\mathrm{nr}}}=n$, and set $P=\mathcal O_D$ and $\Delta=P/J(P)$.
Then the ramification index is $n$, the residue field is $\Delta=k^{1/n}$, and $\per D_{F^{\mathrm{nr}}}=n$.
Moreover, the sequence
\[
 0\longrightarrow\Delta^\times/k^\times
 \longrightarrow K_1(P)\longrightarrow
 \Nrd_1(P^\times)\longrightarrow0
\]
is exact.
For a finite unramified extension $F'/F$, the map $K_1(P)\to K_1(\mathcal O_{D_{F'}})$ is injective, as is the map to the direct limit over all finite unramified extensions.
\end{lemma}
\begin{proof}
At every finite unramified stage, which is again a two-dimensional local field, the index stays $n$, and so does the period by the equality of period and index over every finite extension \cite[Corollary 5.1, Proposition 5.2]{Chipchakov2022}.
A smaller positive multiple of the Brauer class vanishing over $F^{\mathrm{nr}}$ would vanish at a finite stage by descent of its splitting data, so $\per D_{F^{\mathrm{nr}}}=n$.
If the residue center had a nontrivial separable part, or the residue algebra had a noncommutative division-algebra part, a finite separable residue extension splitting that part would produce a nontrivial idempotent.
Henselianity lifts that idempotent to the order, contradicting the fact that the generic fiber is a division algebra.
Thus $\Delta/k$ is purely inseparable.
Since $[k:k^p]=p$, the extension is simple; lifting a primitive element gives $f=[\Delta:k]\le n$.
Lifting a generator of the value group similarly gives $e\le n$ for the ramification index.
The defectless dimension formula $n^2=ef$ now yields $e=f=n$, and hence $\Delta=k^{1/n}$ \cite{WadsworthValuation}.

Let $i_*:K_1(\Delta)\to K_1(P)$ be the support map.
For $x\in K_2(D)$, put $c=\partial_Dx$.
For a separable maximal subfield $E\subset D$ and the split Morita equivalence $\mu_E$, put $\rho_E=\mu_E^{-1}\res_{E/F}\colon K_2(D)\to K_2(E)$.
Then $\theta_E\circ\rho_E=n$, $N_{E/F}\circ\rho_E=N_{D/F}$, and $\Nrd_2\theta_E=N_{E/F}$ \cite[Proposition 6.6.1(1)]{KahnLevine2010}, so $N_{D/F}=\Nrd_2\theta_E\rho_E=n\Nrd_2$.
The norm $N_{D/F}=n\Nrd_2$ satisfies $\partial_FN_{D/F}=N_{\Delta/k}\partial_D$, and $N_{\Delta/k}(c)=c^n$ because $\Delta/k$ is purely inseparable of degree $n$, so
\[
 c^n=(\partial_F\Nrd_2x)^n.
\]
Since $\Delta^\times$ has no $p$-torsion, it follows that $c=\partial_F\Nrd_2x\in k^\times$.
The inclusion $k^\times\subseteq\partial_DK_2(D)$ follows from Theorem~\ref{thm:second-two-surj} and surjectivity of the boundary for the commutative field, so $\ker i_*=k^\times$.
Since $\SK_1(D)=0$ by Corollary~\ref{cor:review-two-local-sk1-zero}, we have $\ker[K_1(P)\to\Nrd_1(P^\times)]=\im i_*$, which proves the displayed exact sequence.

For a finite unramified extension with residue extension $k'/k$, the new residue field is $\Delta'=\Delta k'$.
The boundary comparison alone gives $\ker i_*'\subseteq(k')^\times$.
If $z\in K_1(P)$ becomes zero after base change, then $\Nrd_1(z)=1$, so $z=i_*c$ for some $c\in\Delta^\times/k^\times$.
Vanishing after base change forces $c\in(k')^\times$.
The linear disjointness of purely inseparable and separable extensions gives $\Delta\cap k'=k$, hence $z=0$.
Injectivity of the map to the direct limit follows from continuity of $K$-theory.
\end{proof}

\begin{lemma}
\label{lem:review-geometric-bockstein-detection}
In the situation of the preceding lemma, let $P^{\mathrm{nr}}$ be the union of all finite unramified base changes of $P$, and let $\pi$ be a uniformizer of $\mathcal O_F$.
There is $j_0$ such that, for all $a\ge1$ and $j\ge j_0$,
\[
 K_1(P)[p^a]\longrightarrow K_1(P^{\mathrm{nr}}/\pi^jP^{\mathrm{nr}})
\]
is injective.
Consequently, if a class in $K_2(P;\Z/p^a)$ vanishes on this geometric truncation, its Bockstein is zero.
\end{lemma}
\begin{proof}
We use the injectivity in Lemma~\ref{lem:app-deep-unit-norm}, uniform over unramified extensions.
The hypothesis $\ind D_{F^{\mathrm{nr}}}=n$ says that every $D_L$ is a division algebra of the same degree and that $P_L$ is obtained by base change from the same integral lattice.
Thus there is $j_1$ such that, for every finite unramified $L/F$ and $j\ge j_1$, the map $\Nrd_1$ is injective on the image of $1+\pi^jP_L$ in $K_1(P_L)$.

In mixed characteristic put $e=v_F(p)$.
For $j>e/(p-1)$, the group $1+\pi^j\mathcal O_L$ contains no nontrivial root of unity of $p$-power order.
Indeed, if $\zeta\ne1$ has order $p^r$, then $v_L(\zeta-1)=e/(p^{r-1}(p-1))$, and $e$ is constant under unramified extension.
In equal characteristic the multiplicative group has no $p$-torsion.
Choosing $j_0$ sufficiently large therefore makes the $\pi^j$-congruence kernel in every $K_1(P_L)$ free of $p$-primary torsion.
Here lifting general linear groups shows that this congruence kernel equals $\im(1+\pi^jP_L\to K_1(P_L))$.

If $z\in K_1(P)[p^a]$ vanishes on the geometric truncation, it already vanishes at a finite stage $L/F$.
Then $z_L$ is a $p$-primary torsion element of the preceding congruence kernel, hence is zero.
The unramified injectivity in Lemma~\ref{lem:review-minimal-order-k1-support} gives $z=0$.
The last assertion follows by naturality of the Bockstein in the Moore coefficient sequence.
\end{proof}

\begin{lemma}
\label{lem:review-geometric-moore-integral-lift}
Let $F$ be a two-dimensional local field whose first residue field $k$ has characteristic $p$, and let $B/F$ be a central division algebra with $n=p^q=\deg B=\ind B_{F^{\mathrm{nr}}}$.
Let $R=\mathcal O_F$, $C=\mathcal O_B$, $\Delta=C/J(C)$, and let $\pi$ be a uniformizer of $R$.
Let $Y_C=K(C,\Delta)$ and $Y_s=K(C/\pi^sC,\Delta)$.
For every $a,j\ge1$, there is $b\ge j$ such that, for every geometrically zero class in $\pi_2(Y_b/p^a)$, the image of that class in $Y_j/p^a$ is represented by an integral element of $\pi_2Y_C$.
\end{lemma}
\begin{proof}
We lift a generator of $\Delta/k$ to $w\in C^\times$.
Since $[\Delta:k]=\deg B=n$ by Lemma~\ref{lem:review-minimal-order-k1-support}, the residue degree of $E=F(w)$ is at least $n$.
The bound on the degree of a maximal subfield gives $[E:F]=n$ and ramification index one.
The finite free $R$-algebra $R[w]$ has residue field $\Delta$ and maximal ideal $\pi R[w]$, so it is a discrete valuation ring, and hence $R[w]=\mathcal O_E\subset C$.
Since both entries of a Matsumoto symbol in $K_2(\Delta)$ lift to units of this commutative order, the map $K_2(C)\to K_2(\Delta)$ is surjective, and the same lifts work in every truncation.
Thus
\[
\begin{gathered}
 \pi_1Y_C=M,\qquad M=\ker[K_1(C)\to\Delta^\times],\\
 \pi_1Y_s=M/\mathcal U_s(C),\qquad \mathcal U_s(C)=\im(1+\pi^sC\to K_1(C)).
\end{gathered}
\]
By Corollary~\ref{cor:app-integral-cut-stable}, for sufficiently deep $r\ge j$ the image of $\pi_2Y_r$ in $\pi_2Y_j$ equals the image of $\pi_2Y_C$.
Increase $r$ so that $r\ge j_0$ for the integer $j_0$ of Lemma~\ref{lem:review-geometric-bockstein-detection} and so that, by Lemma~\ref{lem:app-deep-unit-norm}, the reduced norm is injective on $\mathcal U_r(C)$.

For a sufficiently deep $b$, write the Bockstein of the given class as $\beta_b\in(M/\mathcal U_b(C))[p^a]$ and lift it to $\beta\in M$.
Then $p^a\beta\in \mathcal U_b(C)$.
Put $t=\Nrd_1(\beta)$, so $t^{p^a}\in1+\pi^bR$.
In mixed characteristic, if $b$ is sufficiently large, Hensel's lemma provides a $p^a$-th root $\eta$ of $t^{p^a}$ in $1+\pi^{r+c}R$, where $c$ satisfies $1+\pi^{r+c}R\subset\Nrd_1(\mathcal U_r(C))$.
In equal characteristic, $v_F(t^{p^a}-1)=p^a v_F(t-1)$, so $b\ge p^a(r+c)$ allows us to take $\eta=t\in1+\pi^{r+c}R$.
We choose $\gamma\in \mathcal U_r(C)$ with $\Nrd_1(\gamma)=\eta$.
The element $p^a(\beta-\gamma)$ belongs to $\mathcal U_r(C)$ and has reduced norm one, so it is zero.
Thus $\beta_0=\beta-\gamma\in M[p^a]$ represents the image of $\beta_b$ at depth $r$.

Since the given Moore class is geometrically zero, $\beta_0$ vanishes on the geometric truncation at depth $r$.
Lemma~\ref{lem:review-geometric-bockstein-detection} gives $\beta_0=0$.
The image of the given class in $Y_r/p^a$ therefore has zero Bockstein and comes from some $u_r\in\pi_2Y_r$.
By the choice of $r$ via Corollary~\ref{cor:app-integral-cut-stable}, the image of $u_r$ in $Y_j$ is represented by an element of $\pi_2Y_C$.
\end{proof}

Let $F'/F$ be a finite Galois unramified extension of a two-dimensional local field $F$, let $D/F$ be a central division algebra, and let $D'$ be the division algebra in the Brauer class of $D_{F'}$.
Put $R=\mathcal O_F$, $R'=\mathcal O_{F'}$, $P=\mathcal O_D$ and $C=\mathcal O_{D'}$.
The order $P\otimes_RR'$ is hereditary.
Indeed, finite \'etale base change is separable, so every base-changed module is a direct summand of the module obtained by restricting and inducing it again; the bound of one on projective dimension is preserved.
After rectangular diagonal blocks are removed by Morita equivalence, the lattice-chain classification of hereditary orders therefore identifies $P\otimes_RR'$ with an Iwahori order with corner $C$ \cite[Chapter 9]{ReinerMaximalOrders}.

\begin{lemma}
\label{lem:review-corner-h-one-transfer}
Let $F$ be a two-dimensional local field whose first residue field $k$ has characteristic $p$, and let $D/F$ be a central division algebra.
Let $R=\mathcal O_F$, $P=\mathcal O_D$, and let $\pi$ be a uniformizer of $R$.
Let $a\ge1$.
Let $F'/F$ be a finite Galois unramified extension, put $R'=\mathcal O_{F'}$, and let $k'/k$ be the residue extension.
Let $D'$ be the division algebra in the Brauer class of $D_{F'}$, and put $C=\mathcal O_{D'}$.
After removing rectangular diagonal blocks by Morita equivalence, we write $P\otimes_RR'$ as an Iwahori order $H$ with corner $C$.
Let
\[
 \mathcal X_s=K(P/\pi^sP,P/J(P))/p^a,\qquad
 \mathcal Y_s=K(C/\pi^sC,C/J(C))/p^a.
\]
These may also be regarded as sheaves on the finite \'etale site, as in Lemma~\ref{lem:app-nilpotent-etale-descent}.
Suppose that classes $y_s\in\pi_2\mathcal X_s(k)$ at this fixed coefficient $p^a$ are compatible for all $s$ and have zero image at the geometric stalk.
For every $j$ and every prescribed $b_0\ge j$, there are $b\ge b_0$ and a geometrically zero $z_b\in\pi_2\mathcal Y_b(k')$ such that $y_j$ is the image, at depth $j$, of the transfer of $z_b$ from the corner.
\end{lemma}
\begin{proof}
Put $V_s=\pi_3(\mathcal X_s)_{\bar k}$ and $W_s=\pi_3(\mathcal Y_s)_{\bar k}$.
They are discrete modules for $G_k$ and $G_{k'}$, respectively.
Corner induction and finite \'etale norms give a coefficient map
\[
 c_{b,j}\colon\operatorname{Ind}_{G_{k'}}^{G_k}W_b\longrightarrow V_j.
\]
At the geometric stalk, it is the sum of the corner maps indexed by the embeddings $k'\hookrightarrow\bar k$.
Denote its image by $L_{b,j}\subset V_j$.

By Lemma~\ref{lem:app-iwahori-derived-cut}, the transition maps on the relative corner cofibers are nullhomotopic through a factorization via the Koszul algebra defined over $k'$, so the nullhomotopies respect $G_{k'}$.
For fixed $b\ge b_0$, we may therefore take $r$ sufficiently deep that
\[
 V_r|_{G_{k'}}\longrightarrow W_b\longrightarrow V_j|_{G_{k'}}
\]
equals $V_r\to V_j$.
Indeed, if the composite from $\mathcal X_r$ to the corner cofiber of $\mathcal X_b$ is nullhomotopic, lift it to $\mathcal Y_b$ and then pass to depth $j$.
Thus $\im(V_r\to V_j)\subseteq L_{b,j}$.

By Lemma~\ref{lem:app-nilpotent-etale-descent}, $y_r$ comes from some $h_r\in H^1(k,V_r)$.
The image of $h_r$ in $H^1(k,V_j)$ comes from $H^1(k,L_{b,j})$.
The map $\operatorname{Ind}W_b\to L_{b,j}$ is surjective by definition; its kernel is a discrete $p$-primary module, so $H^2(k,\ker(\operatorname{Ind}W_b\to L_{b,j}))=0$.
Consequently,
\[
 H^1(k,\operatorname{Ind}_{G_{k'}}^{G_k}W_b)
 \twoheadrightarrow H^1(k,L_{b,j}).
\]
By Shapiro's lemma, we have $H^1(k,\operatorname{Ind}_{G_{k'}}^{G_k}W_b)\simeq H^1(k',W_b)$, which gives a geometrically zero $z_b\in\pi_2\mathcal Y_b(k')$ through the edge map of Lemma~\ref{lem:app-nilpotent-etale-descent}.
At the geometric stalk, transfer is the sum over the preceding embeddings.
Compatibility with the descent edge map therefore identifies the transfer of $z_b$ at depth $j$ with $y_j$.
\end{proof}

\begin{lemma}\label{lem:second-two-inertial-reduction}
Let $F$ and $D$ be as in Lemma~\ref{lem:review-corner-h-one-transfer}.
Let $F'/F$ be a finite Galois unramified extension with $\ind D_{F'}=n$, where $n=\ind D_{F^{\mathrm{nr}}}$, and let $D'$ be the division algebra in the Brauer class of $D_{F'}$.
Then
\[
 \ker[\SK_2(D)\longrightarrow\SK_2(D')]
 \ \subseteq\ \cor_{F'/F}\SK_2(D').
\]
We do not assume that $[F':F]$ is prime to $p$.
\end{lemma}
\begin{proof}
Put $R=\mathcal O_F$, $P=\mathcal O_D$, $C=\mathcal O_{D'}$, and $R'=\mathcal O_{F'}$.
As explained before Lemma~\ref{lem:review-corner-h-one-transfer}, after rectangular diagonal blocks are removed by Morita equivalence, $P\otimes_RR'$ is an Iwahori order with corner $C$.
An element $\xi$ of $\ker[\SK_2(D)\to\SK_2(D')]$ is $p$-primary torsion, since the prime-to-$p$ parts of $\SK_2$ vanish: for each prime $\ell\ne p$, split the other primary components of $[D]$ over a finite separable extension of degree prime to $\ell$, which detects the $\ell$-primary part by restriction and transfer, and apply Proposition~\ref{prop:second-residue-zero-ell}(2) to the division algebra representing the result.
The residue algebra $P/J(P)$ is a division algebra whose center $l$ is a nonarchimedean local field, so Theorem~\ref{thm:paper-local-norm} gives reduced-norm isomorphisms $K_1(P/J(P))\simeq l^\times$ and $K_2(P/J(P))\simeq K_2(l)$.
The localization boundary of $\xi$ vanishes since $l^\times$ has no $p$-torsion, so $\xi$ lifts to $K_2(P)$.
The support image is $p$-divisible, since $K_2(l)$ is uniquely $p$-divisible \cite{GeisserLevine}.
Adjusting the lift by a support class therefore makes a $p$-power multiple of the adjusted lift zero, so the integral lift itself may be chosen $p$-primary torsion.
Since residue $K_2$ has no $p$-torsion, the adjusted lift has zero reduction, so we choose a relative integral lift $y$ of $\xi$.
Let $\pi$ be a uniformizer of $R$, and put $\Delta'=C/J(C)$.
Let $y'$ be the image of $y$ in $K_2(C)$ after base change and the corner equivalence (Lemma~\ref{lem:app-iwahori-derived-cut}).
Since $\xi|_{D'}=0$, the localization sequence gives $y'\in i_*K_2(\Delta')$, and the same holds after every further finite unramified base change.
Every finite separable extension $\Delta''$ of $\Delta'$ is a nonarchimedean local field of characteristic $p$ with $[\Delta'':\Delta''^p]=p$, so Geisser--Levine gives $K_2(\Delta'';\Z/p^a)=K_3(\Delta'';\Z/p^a)=0$.
The first vanishing kills the image of the closed-support classes in $K_2(C;\Z/p^a)$, and the second makes $\pi_2$ of the relative group inject into the absolute one, so the relative Moore class also vanishes.
Hence, for each fixed $a$, the compatible family has zero image at the geometric stalk, as required in Lemma~\ref{lem:review-corner-h-one-transfer}.

Fix a finite separable splitting field $L_0/F$, and choose a depth $j$ at which Lemma~\ref{lem:app-congruence-separation} holds for $L_0$.
We choose $p^a$ annihilating $K_2(P/\pi^jP,P/J(P))$, which exists because this relative group is $p$-primary of bounded exponent (Lemma~\ref{lem:app-nilpotent-trace}).
Let $b_0$ be a depth at which Lemma~\ref{lem:review-geometric-moore-integral-lift} can be applied to the corner $C$, and apply Lemma~\ref{lem:review-corner-h-one-transfer}.
The lemma applies because $\deg D'=n$ is a power of $p$: the residue field of $F^{\mathrm{nr}}$ is separably closed and $\cd_\ell(F^{\mathrm{nr}})\le1$ for $\ell\ne p$, so $\Br(F^{\mathrm{nr}})$ is $p$-primary; moreover, $D'$ has minimal unramified index $n$, since $(F')^{\mathrm{nr}}=F^{\mathrm{nr}}$.
At depth $j$, the image of the geometrically zero class $z_b$ is represented by an integral $z\in K_2(C,C/J(C))$.
Thus $y$ and the transfer of $z$ from the corner agree in $\pi_2\mathcal X_j$, with $\mathcal X_j$ as in Lemma~\ref{lem:review-corner-h-one-transfer}.
Since $p^a$ kills $K_2(P/\pi^jP,P/J(P))$, the Moore sequence injects this integral group into $\pi_2\mathcal X_j$, so the equality holds integrally at depth $j$.
By Lemma~\ref{lem:app-congruence-separation}, we have $\xi-\cor_{F'/F}(z|_{D'})\in I_F(D)$.
By Proposition~\ref{prop:review-splitting-transfer-norm-iso}, $K_2(D')=I_{F'}(D')\oplus\SK_2(D')$.
Write $z|_{D'}=z_I+z_0$ with $z_I\in I_{F'}(D')$ and $z_0\in\SK_2(D')$.
By transitivity of transfer, $\cor_{F'/F}I_{F'}(D')\subseteq I_F(D)$, so $\xi-\cor_{F'/F}(z_0)\in I_F(D)$.
Since $\Nrd_2\circ\cor_{F'/F}=N_{F'/F}\circ\Nrd_2$, the same element lies in $\SK_2(D)$.
As $I_F(D)\cap\SK_2(D)=0$ by the injectivity of $\Nrd_2$ on $I_F(D)$ (Proposition~\ref{prop:review-splitting-transfer-norm-iso}), we obtain $\xi=\cor_{F'/F}(z_0)$, and hence
\[
 \xi\in\cor_{F'/F}\SK_2(D').
\]
\end{proof}

By Lemma~\ref{lem:second-two-inertial-reduction}, it suffices to induct on $q$ for a division algebra achieving the minimal unramified index $n=p^q$.
We replace the division algebra accordingly and fix $R=\mathcal O_F$, $P=\mathcal O_D$, a uniformizer $\pi$ of $F$, and $J_D=J(P)$.
Lemma~\ref{lem:review-minimal-order-k1-support} gives ramification index $n$, residue field $\Delta=k^{1/n}$, and period $n$ over the maximal unramified field.

\subsubsection{Two cyclic covers that lower the index}
To run the induction on the minimal unramified index, we construct, after a suitable base extension, a cyclic extension $T/F$ of degree $p$ inside $D$ whose centralizer has smaller minimal unramified index.
There are two types of such extensions, according to the ramification of $T/F$.

\begin{lemma}\label{lem:second-two-cyclic-cover}
Let $k$ be a nonarchimedean local field of characteristic $p$, let $F$ be a complete discretely valued field with first residue field $k$, and let $D/F$ be a central division algebra with $\deg D=\ind D_{F^{\mathrm{nr}}}=n=p^q>1$.
Put $m=n/p$, $P=\mathcal O_D$ and $\Delta=P/J(P)$.
After adjoining $\mu_p$ if necessary and passing to a finite Galois unramified extension, and renaming the resulting fields and orders again $F$, $k$, $D$, $P$, $\Delta$ and $\pi$, there are a cyclic extension $T/F$ of degree $p$, a generator $\sigma$ of $\Gal(T/F)$ and $x\in D^\times$ with the following properties.
\begin{enumerate}
\item $T\subset D$, and conjugation by $x$ induces $\sigma$ on $T$.
\item $B=Z_D(T)$ is a division algebra with $\deg B=\ind B_{T^{\mathrm{nr}}}=m$; it has ramification index $m$, its residue field $\Delta_0$ is the field of $m$-th roots of the residue field of $T$, and $\SK_1(B)=0$.
\item $m[D]=(T/F,\sigma,\Nrd_{B/T,1}(x^p))$, where $\sigma$ is the generator in the chosen cyclic normal form.
\item Either $T/F$ is totally ramified and $x$ is a unit whose residue generates $\Delta$ over $k$ \textup{(}the first type, or homogeneous-unit type\textup{)}, or $T/F$ has ramification index one and purely inseparable residue degree $p$, and $x$ is a uniformizer of $D$ \textup{(}the second type, or homogeneous-uniformizer type\textup{)}.
\item With $C=\mathcal O_B$, one has $P=\bigoplus_{j=0}^{p-1} Cx^j$.
\end{enumerate}
\end{lemma}
\begin{proof}
Adjoining $\mu_p$ is an extension of degree prime to $p$, so it preserves the $p$-primary index, and restriction to it detects the $p$-primary part of $\SK_2$.
It also preserves the minimal unramified index, since the maximal unramified extension after adjoining is $F(\mu_p)F^{\mathrm{nr}}$, of degree prime to $p$ over $F^{\mathrm{nr}}$.
The finite Galois unramified stage is not assumed to have prime degree; we return from it by the inclusion through transfer in Lemma~\ref{lem:second-two-inertial-reduction}.
Modulo classes killed by the maximal unramified extension, the top component $m[D]$ has one of the following presentations using a cyclic extension $T/F$ of degree $p$:
\begin{enumerate}
\item $T/F$ is totally ramified, and the coefficient is a unit whose residue is not a $p$-th power.
\item $T/F$ has ramification index one and purely inseparable residue degree $p$, and the coefficient is a uniformizer.
\end{enumerate}
In equal characteristic, the leading pole order in the Artin--Schreier--Cartier presentation gives these two types; see the proof of \cite[\S4, Lemma 5, pp.~338--339]{KatoLocalI}.
In mixed characteristic, depth zero and positive depths prime to $p$ give the first type, and positive depths divisible by $p$ give the second \cite[\S2, Proposition 2(i)--(iii), pp.~321--322, and \S4, p.~339]{KatoLocalI}; at depth zero, Kato's symbol $\{f,\pi\}$ is used in the form $\{\pi,f^{-1}\}$, which gives the totally ramified Kummer field $T=F(\pi^{1/p})$.
Indeed, at an intermediate depth $s$, write $z^p=1+c$ with $v_F(c)=s<pv_F(p)/(p-1)$.
If $p\nmid s$, then $v(z-1)=s/p$.
If $s=ph$, then $Y=(z-1)/\pi^h$ satisfies $\overline Y^{\,p}=\overline{c/\pi^s}$: the intermediate binomial coefficients are divisible by $p$, and the corresponding terms vanish modulo $\pi$ because $v_F(p)-(p-1)h>0$.
In the second type, $\overline{c/\pi^s}\notin k^p$, so $\overline Y$ is a purely inseparable residue root.
At the terminal depth $pv_F(p)/(p-1)$, we obtain a separable Artin--Schreier equation by subtracting $1$ from a $p$-th root and dividing by $\zeta_p-1$, and hence that layer is killed unramifiedly.

In both types $T^{\mathrm{nr}}=TF^{\mathrm{nr}}$, since the residue extension of $T/F$ is trivial or purely inseparable, so the separable residue extensions come from $k$.
Both types of $T/F$ have no nontrivial unramified subextension, so $[T^{\mathrm{nr}}:F^{\mathrm{nr}}]=p$.
After removing the difference that vanishes over the maximal unramified extension, $m[D]$ is a cyclic class for $T$, hence vanishes over $T^{\mathrm{nr}}$ and over some finite stage $TF'$.
The upper bound is
\[
 \ind D_{T^{\mathrm{nr}}}\le\ind D_{TF'}=\per D_{TF'}\le m,
\]
where $F'/F$ is a finite unramified extension such that $m[D]$ splits over $TF'$, and the equality holds over the two-dimensional local field $TF'$ \cite[Corollary 5.1, Proposition 5.2]{Chipchakov2022}.
The lower bound $\ind D_{T^{\mathrm{nr}}}\ge\per D_{T^{\mathrm{nr}}}\ge m$ follows from $\per D_{F^{\mathrm{nr}}}=n$ and $\cor_{T^{\mathrm{nr}}/F^{\mathrm{nr}}}\res_{T^{\mathrm{nr}}/F^{\mathrm{nr}}}=p$ on $\Br(F^{\mathrm{nr}})$.
Hence the minimal unramified index after restriction to $T$ is $m$.
Since $T^{\mathrm{nr}}=TF^{\mathrm{nr}}$, Lemma~\ref{lem:second-two-inertial-reduction} lets us pass to a further finite unramified stage of $F$ at which $\ind D_T=m$ and the difference from the normal form vanishes.
The index of $D$ remains $n$, and since $\ind D/\ind D_T=n/m=[T:F]$, the field $T$ embeds into $D$.
Then $B=Z_D(T)$ is a division algebra of minimal unramified index $m$, and Lemma~\ref{lem:review-minimal-order-k1-support} applied to $B/T$ shows that $B$ has ramification index $m$ and residue field equal to the full field of $m$-th roots of the residue field of $T$.
We have $\per B=\ind B$, and Corollary~\ref{cor:review-two-local-sk1-zero} gives $\SK_1(B)=0$.

For a generator $\sigma$ of $\Gal(T/F)$ and $x\in D^\times$ implementing $\sigma$ by conjugation,
\[
 m[D]=(T/F,\sigma,\Nrd_{B/T,1}(x^p)).
\]
This identity is obtained in the corner of $D^{\otimes_F m}$ corresponding to $T^{\otimes_F m}\to T$.
The algebra $B^{\otimes_T m}$ is split.
For the standard representation $V$ of $\operatorname{GL}_m$, scalars act trivially on $V^{\otimes m}\otimes\det(V)^{-1}$, so $g\mapsto g^{\otimes m}/\det g$ is a representation of $\operatorname{PGL}_m$ and applies functorially to the cocycle and to the $\sigma$-semilinear action; hence the carry in the semilinear action on that corner is $(x^p)^{\otimes m}/\Nrd_{B/T,1}(x^p)$.
If the semilinear action is $P\sigma$, replacing the corner generator $x'$ by $P^{-1}x'$ changes it to componentwise $\sigma$, and the remaining coefficient is $(P^{-1}x')^p=\Nrd_{B/T,1}(x^p)$.

In the first type, correction by $B^\times$ makes $x$ a unit.
Residues of norms of units in a totally ramified extension are $p$-th powers, and $\overline{\Nrd_{B/T,1}(x^p)}=\bar x^{pm}=\bar x^n$ because the residue field $k^{1/m}$ of $B$ is purely inseparable of degree $m$ over $k$.
Hence $\bar x^n\notin k^p$, and $\bar x$ generates $\Delta$.
In the second type, $\Nrd_{B/T,1}(x^p)/\pi$ is a norm from $T/F$, whose valuation lies in $p\Z$, and $v_F(\Nrd_{B/T,1}(x^p))=v_D(x)$, with the valuation of a uniformizer of $D$ normalized to one; hence $v_D(x)\equiv1\bmod p$.
Since $v_D(B^\times)=p\Z$, we may take $x$ to be a uniformizer.
Put $C=\mathcal O_B$ and $\Delta_0=C/J(C)$.
In both types, $P=\bigoplus_{j=0}^{p-1} Cx^j$.
In the first type, after cancelling the minimum valuation of the coefficients by an element of $B$, the independence of $1,\bar x,\ldots,\bar x^{p-1}$ over $\Delta_0$ gives the decomposition; in the second, the valuations of the terms are distinct modulo $p$.
\end{proof}

We pass to the extension of Lemma~\ref{lem:second-two-cyclic-cover}, rename the resulting fields and orders again $F$, $k$, $D$, $P$, $\Delta$ and $\pi$, and fix $T$, $\sigma$, $x$, $B$, $C$ and $\Delta_0$ as in that lemma.
Induction on the minimal unramified index gives $\SK_2(B)=0$.
Together with the isomorphism $\Nrd_2|_{I_T(B)}\colon I_T(B)\to K_2(T)$ of Proposition~\ref{prop:review-splitting-transfer-norm-iso}, this yields $K_2(B)=I_T(B)$.

\subsubsection{The degree-zero relative trace}
The cyclic extension of Lemma~\ref{lem:second-two-cyclic-cover} makes $P$ a $\Z/p$-graded order.
In this subsubsection, we show that the map from the $p$-completed $\mathrm{TC}$ of the total-degree-zero part of $\mathrm{THH}(P)$ to $\mathrm{TC}(P)^\wedge_p$ is surjective on $\pi_2$, and we then remove the residue contribution.
Let $G=\Z/p$.
For a $G$-graded ring $S$, its cover $\mathcal Q(S)$ has vertices $i\in G$, morphisms $S_{j-i}$, and the strict $G$-action that simultaneously translates the vertices.
We write $K(\mathcal Q(S))$ for the $K$-theory of the category of perfect right $\mathcal Q(S)$-modules.
We write $\mathbb T$ for the circle group.

\begin{lemma}
\label{lem:review-cyclic-zero-tc}
Let $A$ be a finite abelian group with $pA=0$, and let $S=\bigoplus_{g\in A}S_g$ be an $A$-graded discrete ring.
Let
\[
 X=\mathrm{THH}(S)=\bigoplus_{g\in A}X_g
\]
be the total-degree decomposition of the cyclic bar construction.
This decomposition preserves the $\mathbb T$-action, and $X_0$ is a $p$-cyclotomic subobject.
We call $\mathrm{TC}(X_0;p)^\wedge_p$ the degree-zero $\mathrm{TC}$ of $S$; we do not claim that it is a direct summand of $\mathrm{TC}(S)^\wedge_p$.
With this notation there is a fiber sequence
\[
 \mathrm{TC}(X_0;p)^\wedge_p\longrightarrow
 \mathrm{TC}(S)^\wedge_p\longrightarrow
 \left(\bigoplus_{g\ne0}\Sigma(X_g)_{h\mathbb T}\right)^\wedge_p.
\]
If $\pi_1(((X_g)_{h\mathbb T})^\wedge_p)=0$ for all $g\ne0$, the first map is therefore surjective in degree two.

For the graded cover $\mathcal Q_A(S)$ with vertices $A$ and morphisms $S_{j-i}$ from $i$ to $j$, there is a natural equivalence
\[
 \left(\mathrm{TC}(\mathcal Q_A(S);p)_{hA}\right)^\wedge_p
 \simeq\mathrm{TC}(X_0;p)^\wedge_p.
\]
On the left, $A$ acts on $\mathcal Q_A(S)$ by vertex translation.
The same assertions hold for relative fibers of maps of graded rings or covers.
\end{lemma}
\begin{proof}
Total degree in the cyclic bar construction is preserved by face maps given by multiplication, degeneracies given by the unit, and cyclic permutation.
The cyclotomic Frobenius from the $p$-fold Tate diagonal repeats a simplex of total degree $g$ exactly $p$ times, sending it to degree $pg=0$.
Thus $X_0$ is a sub-$p$-cyclotomic spectrum, and the cofiber $Q=X/X_0$, with underlying spectrum $\bigoplus_{g\ne0}X_g$, has zero Frobenius.

Since $Q$ has zero Frobenius, the Nikolaus--Scholze formula for bounded-below $p$-cyclotomic spectra identifies $\mathrm{TC}(Q;p)^\wedge_p$ with the fiber of the canonical map $(Q^{h\mathbb T})^\wedge_p\to(Q^{t\mathbb T})^\wedge_p$.
By the circle norm fiber sequence, we have $\mathrm{TC}(Q;p)^\wedge_p\simeq(\Sigma Q_{h\mathbb T})^\wedge_p$.
The functor $\mathrm{TC}(-;p)$ is exact between stable categories; applying it to $X_0\to X\to Q$ gives the stated sequence.
Since circle homotopy orbits commute with the finite direct sum $\bigoplus_{g\ne0}X_g$, we have $Q_{h\mathbb T}\simeq\bigoplus_{g\ne0}(X_g)_{h\mathbb T}$.
The degree-two assertion follows from the long exact sequence \cite[Proposition II.1.9, Lemma II.4.2]{NikolausScholzeTC}.

In each simplicial degree, the lifts of a total-degree-zero simplex in the cyclic bar construction form a free $A$-set indexed by the starting vertex.
Taking homotopy orbits under translation of the starting vertex gives $\mathrm{THH}(\mathcal Q_A(S))_{hA}\simeq X_0$ via the map that forgets the vertices.
For every prime $q$, the $q$-fold Tate diagonal repeats a cycle $q$ times from the same starting vertex and commutes with translation of vertices, so the $A$-action and this map are cyclotomic \cite[Corollary II.1.7, \S III.2]{NikolausScholzeTC}.
Since the forgetful functor from cyclotomic spectra preserves colimits, the computation of free $A$-orbits in each cyclic degree gives an equivalence of cyclotomic spectra.
The absolute $\mathrm{THH}$ is connective, and $\mathrm{TC}/p$ preserves all colimits of connective cyclotomic spectra \cite[Theorem 2.7]{ClausenMathewMorrow}, so it may be applied to the $A$-bar construction.
Induction using the coefficient exact sequence gives the same assertion with Moore coefficients $p^a$, and passing to the inverse limit over $a$ yields the displayed equivalence in the $p$-complete category.
The relative version follows from the two absolute versions, because finite fibers commute with homotopy orbits and $\mathrm{TC}$ is exact.
\end{proof}

\begin{lemma}\label{lem:review-finite-pbasis-differentials}
Let $R$ be a complete discrete valuation ring with uniformizer $\pi$ and residue field $k$ of characteristic $p$.
Let $\{\bar u_i\}$ be a finite $p$-basis of $k$, and choose lifts $u_i\in R$.
In mixed characteristic, $\Omega^1_{R/\Z}\otimes_R R/p^aR$ is generated by $d\pi$ and the $du_i$ for every $a\ge1$.
Moreover, for $f\in\pi^tR$,
\[
 df\in\pi^{t-1}R\,d\pi+
       \sum_i\pi^tR\,du_i
 \quad\pmod {p^a\Omega^1_{R/\Z}}.
\]
In equal characteristic, choose a coefficient field and take the $u_i$ in it.
Then the $u_i$ together with $\pi$ form a $p$-basis of the ring $R$, and $\Omega^1_{R/\F_p}$ is free on $d\pi$ and the $du_i$.
\end{lemma}
\begin{proof}
In mixed characteristic, put $e=v_R(p)$ and choose $s\ge a$ and $M\ge ea+1$.
By lifting the residue field's $p^s$-basis expansion one digit at a time, we express every $r\in R$ as
\[
 r=\sum_{t=0}^{M-1}\pi^t
      \sum_{0\le I_i<p^s}b_{t,I}^{p^s}u^I+\pi^Mc
 \qquad(b_{t,I},c\in R).
\]
For any derivation $d$ of $R$ into a module killed by $p^a$, we have $d(b^{p^s})=0$ and $d(\pi^Mc)=0$.
Thus $dr$ lies in the submodule generated by the stated finite collection of differentials.
For $f=\pi^tr$, the Leibniz rule gives $df=t\pi^{t-1}r\,d\pi+\pi^t\,dr$, which proves the depth assertion.
This argument applies to universal differentials, so no finite-type hypothesis on $R$ is needed.

In equal characteristic, write $R=k[[\pi]]$.
Using the $p$-basis expansion of coefficients and the decomposition of power-series exponents according to their remainders modulo $p$, we obtain
\[
 R=\bigoplus_{0\le I_i<p,\ 0\le j<p}
                  R^p u^I\pi^j.
\]
This is precisely the definition of a ring $p$-basis, proving the assertion about differentials.
\end{proof}

\begin{lemma}
\label{lem:review-finite-free-pro-kernel}
Let $R$ be a discrete valuation ring with uniformizer $\varpi$, and let $L,L'$ be finite free $R$-modules.
If $f\colon L\to L'$ is injective, there is $c\ge0$ such that, for all $j\ge0$ and all $s\ge j+c$,
\[
 \ker(f\bmod\varpi^s)\longrightarrow L/\varpi^jL
\]
is zero.
\end{lemma}
\begin{proof}
We clear denominators in a left inverse over the fraction field and obtain $g\colon L'\to L$ and $c\ge0$ with $gf=\varpi^c\mathrm{id}_L$.
If $a\in L$ satisfies $f(a)\in\varpi^sL'$, then $\varpi^ca\in\varpi^sL$ and hence $a\in\varpi^{s-c}L\subseteq\varpi^jL$.
\end{proof}

\begin{lemma}
\label{lem:review-nonzero-hc-one}
For each nonzero degree $g$ of the $\Z/p$-graded maximal order $P$ constructed above, one has, in mixed characteristic,
\[
 \{HC_1(P/p^aP\,/\,\Z/p^a)_g\}_{a\ge1}=0
 \quad\text{in the pro category}
\]
and, in equal characteristic, $HC_1(P/\F_p)_g=0$.
\end{lemma}
\begin{proof}
We normalize the valuation on the center by $v_D(\pi)=n$.
The ring $\gr D$ is commutative.
Indeed, since $\Delta/k$ is purely inseparable, conjugation by a homogeneous element induces the identity $k$-automorphism of $\Delta$.
The quotient of value groups is cyclic, and every homogeneous component is the product of $\Delta$ and a power of one uniformizer; these commute.
We calculate cyclic one-forms in $\Omega^1(P)_\natural=\Omega^1(P)/[P,\Omega^1(P)]$ and then subtract actual exact differentials.
The boundary is $b(y\,dz)=[y,z]$.

We begin with the homogeneous-unit type.
Let $\eta$ be a degree-zero uniformizer and $x$ a degree-one unit.
Residue and value-group bases give the $R$-basis $\eta^i x^j$ ($0\le i,j<n$) of $P$.
We choose the lift $u$ of the central $p$-basis element as the constant term of the minimal polynomial of $x$.
Then
\[
 \varphi_\eta(\eta)=\eta^n+\pi\sum_{i<n}b_i\eta^i=0,
 \quad b_0\in R^\times,\qquad
 \varphi_x(x)=x^n+\sum_{0<j<n,\,p\mid j}c_jx^j-u=0,
 \quad c_j\in\pi R.
\]
Only exponents divisible by $p$ occur in $\varphi_x$, here and in the homogeneous-uniformizer type below: if $f$ is the minimal polynomial of $x^p$ over $F$, then $f(X^p)$ kills $x$ and has degree $[F(x):F]$, so it is the minimal polynomial of $x$.
Lemma~\ref{lem:review-finite-pbasis-differentials} expresses every coefficient differential on a finite truncation in terms of $d\pi,du$.

For $y$ of nonzero degree, we apply $y\,d\varphi_\eta(\eta)=0$.
Let $\psi(y)$ be the coefficient of $d\eta$ after cyclic permutation.
Then $[\psi(y),\eta]=[y,\varphi_\eta(\eta)]=0$.
The field $F(\eta)$ is a self-centralizing maximal subfield in degree zero, so $\psi(y)=0$.
The relation therefore involves only $d\pi$ and $du$.
The coefficient operator of $d\pi$ is multiplication by the unit $\sum_i b_i\eta^i$ plus $\pi$ times an integral operator, hence is invertible.
The coefficient of $du$ is divisible by $\pi$.
A $d\pi$ coefficient of depth $r$ can therefore be replaced by a $du$ coefficient of depth at least $r+n$.

Next we use $y\,d\varphi_x(x)=0$ and eliminate $d\pi$ as above.
The coefficient operator of $du$ is $-1$ plus $\pi$ times an integral operator, hence is invertible.
For the coefficient of $dx\,x^{-1}$, each cyclic sum on the associated graded ring is multiplication by the corresponding exponent $j$.
Every such $j$ is divisible by $p$, so this coefficient is at least one step deeper than the input.
Thus a $du$ coefficient of depth $r$ can be replaced by a $dx\,x^{-1}$ coefficient of depth at least $r+1$; for $d\pi$, the gain is at least $n+1$.

We take the exact differential of $h=f\eta^i x^j$, where $f\in\pi^tR$ and $j\equiv g\pmod p$.
After scalar differentials are eliminated by the two relations above, the coefficient of $dx\,x^{-1}$ is $j h$ modulo depth $v_D(h)+1$.
Indeed, the $d\pi$-part of $df$ loses at most $n$ in depth, but its replacement gains $n+1$.
The $du$-part loses no depth and gains one step.
The cyclic sum in $d(x^j)$ is multiplication by $j$ in the commutative associated graded ring.
Since $j\equiv g\ne0$ is a unit, subtracting such actual exact differentials depth by depth eliminates the $dx\,x^{-1}$ coefficient.
Modulo $p^a$ this stops after finitely many steps, leaving a form $z\,d\eta$.

We apply the same procedure to the homogeneous-uniformizer type.
Let $w$ be a degree-zero unit and $x$ a degree-one uniformizer.
The basis is now $w^i x^j$ ($0\le i,j<n$), with relations
\[
 \varphi_w(w)=w^n+\sum_{i<n}c_iw^i-u=0,
 \quad c_i\in\pi R,\qquad
 \varphi_x(x)=x^n+\pi\sum_{i<m}b_i x^{pi}=0,
 \quad b_0\in R^\times.
\]
For an input $y$ of nonzero degree, the $dw$ coefficient of $y\,d\varphi_w(w)$ commutes with $w$.
Since $F(w)$ is a self-centralizing maximal subfield in degree zero, this coefficient vanishes integrally as well.
On the other hand, the $d\pi$ coefficient of $y\,d\varphi_x(x)$ is multiplication by a unit plus $\pi$ times an integral operator.
The $dx$ coefficient of $y\,d\varphi_x(x)$ increases the depth of $y$ by at least $n$.
In the highest term, $x^{n-1}$ contributes $n-1$, and the cyclic sum of length $n$, divisible by $p$, contributes one further step.
The lower terms satisfy the same bound because of their explicit factor $\pi$ and the cyclic sums for exponents $pi$.

These two relations therefore eliminate $d\pi,du$ simultaneously, moving a scalar-differential coefficient of depth $r$ to a $dx$ coefficient of depth at least $r+n$.
The leading $dx$ coefficient of the exact differential of $h=fw^i x^j$ is $j h x^{-1}$.
The $d\pi$-part of $df$ loses at most $n$, so after elimination it is at least one step deeper than the leading depth $v_D(h)-1$.
Again $j\equiv g\ne0$ allows successive elimination of $dx$ by exact differentials, leaving $z\,dw$.

Thus, for both types, every $HC_1$ class on a finite truncation is represented respectively by
\[
 z\,d\eta,\quad [z,\eta]=0\pmod {p^aP};
 \qquad z\,dw,\quad [z,w]=0\pmod {p^aP}.
\]
The coefficient $z$ belongs to the finite free module $P_g/p^aP_g$.
The maps $[\,\cdot\,,\eta]$ and $[\,\cdot\,,w]$ are injective on $P_g$, because the self-centralizing subfields lie in degree zero.
By Lemma~\ref{lem:review-finite-free-pro-kernel}, for any fixed truncation, these kernels at a sufficiently deep stage have zero image in that truncation.
The inverse system of $HC_1$ groups is therefore pro-zero.

In equal characteristic, $\Omega^1_{R/\F_p}$ is finite free and the cyclic one-forms form a finitely generated, hence complete, $R$-module.
The primitives chosen successively above have depths tending to infinity, so they converge.
Continuity of differentiation lets us subtract their sum as one exact differential.
Integral injectivity of the remaining boundary makes $HC_1$ itself zero.
\end{proof}

This lemma supplies the vanishing needed for comparison with absolute degree-zero $\mathrm{TC}$.

\begin{lemma}\label{lem:second-two-degree-zero-surjective}
For the $\Z/p$-graded maximal order $P$ of Lemma~\ref{lem:second-two-cyclic-cover}, with $X=\mathrm{THH}(P)$, the map
\[
 \mathrm{TC}(X_0;p)^\wedge_p\longrightarrow\mathrm{TC}(P)^\wedge_p
\]
of Lemma~\ref{lem:review-cyclic-zero-tc} is surjective on $\pi_2$.
\end{lemma}
\begin{proof}
For each $g\ne0$, we compute the completion of the circle orbits $Z_g=(\mathrm{THH}(P)_g)_{h\mathbb T}$.
In mixed characteristic, $P$ is flat over $\Z$, and low-degree linearization and base change give
\[
 \pi_1(Z_g/p^a)\simeq
 HC_1(P/p^aP\,/\,\Z/p^a)_g.
\]
Here we use the degree-one comparison of circle orbits for $\mathrm{THH}(P)\to HH(P/\Z)$.
In the term of simplicial degree $s$ of the cyclic bar construction, the map from the topological tensor product to the tensor product over $\Z$ is an isomorphism on $\pi_0$ and $\pi_1$: since $P$ is $\Z$-flat, the classical terms are concentrated in degree zero, and the topological ones have $\pi_0=P^{\otimes(s+1)}$ and $\pi_1=0$.
This low-degree comparison is preserved by Moore reduction, realization and circle orbits, and the classical reduction modulo $p^a$ is the cyclic complex of $P/p^aP$.
Lemma~\ref{lem:review-nonzero-hc-one} therefore makes the inverse system on the left pro-zero.
In equal characteristic, the fiber of $\mathrm{THH}(\F_p)\to H\F_p$ lies in degrees at least two, so $\pi_1Z_g\simeq HC_1(P/\F_p)_g=0$ by Lemma~\ref{lem:review-nonzero-hc-one}.
The Moore sequence identifies the transitions on $\pi_1(Z_g/p^a)=\pi_0Z_g[p^a]$ with multiplication by $p$.
Since $p\pi_0Z_g=0$, these transitions are zero, so the inverse system is pro-zero in this case as well.

In either case, the inverse system $\{\pi_1Z_g/p^a\}_a$ has surjective transitions and embeds into the preceding pro-zero system, so every term is zero.
Thus $\pi_1Z_g$ is $p$-divisible, and the maps $\pi_1Z_g[p^{a+1}]\xrightarrow{p}\pi_1Z_g[p^a]$ are surjective.
Both outer terms of the Moore exact sequence
\[
 0\to\pi_2Z_g/p^a\to\pi_2(Z_g/p^a)\to \pi_1Z_g[p^a]\to0
\]
have surjective transitions, so the middle term satisfies Mittag--Leffler as well.
The Milnor exact sequence and $\varprojlim_a\pi_1(Z_g/p^a)=0$ give $\pi_1((Z_g)^\wedge_p)=0$.
By Lemma~\ref{lem:review-cyclic-zero-tc}, the $\mathrm{TC}$ of the cofiber of the degree-zero subobject is $(\bigoplus_{g\ne0}\Sigma Z_g)^\wedge_p$, whose degree two vanishes.
This gives surjectivity in degree two of the map from absolute degree-zero $\mathrm{TC}$ to $\mathrm{TC}(P)^\wedge_p$.
\end{proof}

The lift provided by Lemma~\ref{lem:second-two-degree-zero-surjective} need not vanish in the residue orbits $(\mathrm{TC}(\bar H)^\wedge_p)_{hG}$, where $\bar H$ is the semisimple residue of the cover $H=\mathcal Q(P)$; removing this component makes the remaining class relative, so that Lemma~\ref{lem:app-orbit-stable-image} applies.
In the homogeneous-uniformizer type, $\mathcal Q(P)$ is the Iwahori order with $C$ on and above the diagonal and $J(C)$ below it.
The group $G$ regularly permutes the $p$ copies of $\Delta$ in the semisimple residue of $\mathcal Q(P)$, so the residue orbit is $\mathrm{TC}(\Delta)^\wedge_p$.
The groups $\pi_2\mathrm{TC}(\Delta)^\wedge_p$ and $\pi_3\mathrm{TC}(\Delta)^\wedge_p$ vanish because $[\Delta:\Delta^p]=p$ \cite[Theorem 6.5]{ClausenMathewMorrow}, \cite[Theorem 8.1]{GeisserLevine}.
In the homogeneous-unit type, put $E=F(x)$.
Then $E/F$ is a maximal subfield of ramification index one and residue field $\Delta$.
The ring $\mathcal O_E=\bigoplus_{i<p}\mathcal O_{F(x^p)} x^i$ defines a graded subcover.
In equal characteristic, choose a coefficient field $\Delta\subset\mathcal O_E$ containing $x$, where $\bar x$ is a $p$-basis of $\Delta$; since $E^p\subset F(x^p)$, the field $\Delta_0=\Delta^p$ lies in the degree-zero part $\mathcal O_{F(x^p)}$, so the residue cover has a graded section.
Once the cover is forgotten, the subtracted residue term is zero because $\pi_2\mathrm{TC}(\Delta)^\wedge_p=0$.

We now treat mixed characteristic.
Put $S=\mathcal O_{F(x^p)}$, $l=\Delta_0$, and $v=x^p\in S^\times$.
The reduction of the commutative subcover is the same graded cover of $\Delta=\bigoplus_{i<p}\Delta_0\bar x^i$ as the reduction of $\mathcal Q(P)$.
Multiplication by $x$ identifies successive vertices with the corner.
On $\operatorname{Perf}(S)$ the transported action of $G=C_p$ is the identity for each group element.
For $0\le i,j<p$, writing $[i+j]$ for the representative in $\{0,\ldots,p-1\}$, its composition isomorphism is multiplication by
\[
 x^{i+j-[i+j]}=v^{\lfloor(i+j)/p\rfloor}.
\]
The action after reduction is given by the same formula with $\bar v\notin l^p$.
We show that reduction on the $p$-completed homotopy orbits of $\mathrm{TC}(S)$ is surjective on $\pi_2$.
For this we compare the trace with the determinant.
For a finite separable $l'/l$, let $S_{l'}$ be the corresponding finite \'etale local $S$-algebra.
On the small \'etale site of $l$, let $\mathcal D_S$ and $\mathcal D_l$ be the Picard spectra of graded lines over the algebras $S_{l'}$ and over the fields $l'$, respectively, with the symmetry on ranks $r,r'$ given by $(-1)^{rr'}$.
For $\Lambda\in\{S,l\}$, put $\mathcal D(\Lambda)=R\Gamma(l,\mathcal D_\Lambda)^\wedge_p$.
The group $G$ acts on these spectra through the action above, with the composition isomorphisms $v_\Lambda^{\lfloor(i+j)/p\rfloor}$, which act on determinants of rank $n$ by $v_\Lambda^{n\lfloor(i+j)/p\rfloor}$, where $v_S=v$ and $v_l=\bar v$.

\begin{lemma}\label{lem:second-two-determinant-comparison}
In this setup, the determinant and the cyclotomic trace give a map of reduction arrows of $G$-spectra
\[
 \bigl(\mathrm{TC}(S)^\wedge_p\to\mathrm{TC}(l)^\wedge_p\bigr)
 \longrightarrow\bigl(\mathcal D(S)\to\mathcal D(l)\bigr).
\]
The map $\mathrm{TC}(l)^\wedge_p\to\mathcal D(l)$ is an equivalence, and the map $\bigl((\mathrm{TC}(S)^\wedge_p)_{hG}\bigr)^\wedge_p\to\bigl(\mathcal D(S)_{hG}\bigr)^\wedge_p$ is surjective on $\pi_2$.
\end{lemma}
\begin{proof}
We first take the signed determinant of perfect complexes and then take Moore coefficients.
This determinant identifies the zero and first homotopy sheaves, and a scalar automorphism $v$ on a complex of rank $n$ induces multiplication by $v^n$ on its determinant \cite[Proposition 5.3 and \S12]{BhattScholzeWitt}.
The finite \'etale base-change functors preserve the scalar composition isomorphisms.
The cyclotomic trace is natural for these functors and for their natural isomorphisms \cite[Theorem 7.4]{BlumbergGepnerTabuada}, and so is the determinant \cite[Proposition 5.3]{BhattScholzeWitt}.
Over the closed \'etale site, the cyclotomic trace identifies the sheafified Postnikov towers of $K/p^a$ and $\mathrm{TC}/p^a$ by strict Henselian comparison.
For each $\Lambda\in\{S,l\}$, the global sections of the latter tower compute $\mathrm{TC}(\Lambda)/p^a$ by proper base change for the identity of $\operatorname{Spec}\Lambda$, since $\Lambda$ is Henselian along $(p)$ \cite[Theorems 6.1 and 6.3]{ClausenMathewMorrow}.
Inverting the trace equivalence, composing with the determinant, and taking the inverse limit gives the map of reduction arrows of $G$-spectra in the statement.
For $l$, the determinant is an equivalence with finite Moore coefficients on the geometric stalk.
Indeed, a separable extension of $l$ has $p$-rank one, so logarithmic differentials in degree at least two vanish.
The norm-residue isomorphism and the comparison with Quillen $K$-theory then give $K_i(-;\Z/p^a)=0$ there for $i\ge2$, and the units have no $p$-torsion \cite[Corollary (2.8)]{BlochKato1986}, \cite[Theorem 8.1]{GeisserLevine}.
It follows by descent that $\mathrm{TC}(l)^\wedge_p\to\mathcal D(l)$ is an equivalence.

Next we bound the comparison fiber.
Write $U$ for the sheaf of principal units $1+\mathfrak m_{S_{l'}}$.
Its higher cohomology on $l$ vanishes.
At a finite Galois stage the successive quotients of the principal-unit filtration are additive residue modules, whose higher cohomology vanishes by the normal basis theorem.
The invariants of the finite congruence quotients form a surjective tower, so the vanishing passes to the inverse limit there, and continuous cochains on the geometric stalk are the filtered union of the finite-stage cochains.
Also $U(S)$ has bounded $p$-power torsion, since a primitive $p^a$-th root of unity satisfies
$v_S(\zeta-1)=v_S(p)/(p^{a-1}(p-1))$.
The rank map is unchanged by reduction, and the scalar coherences act trivially on objects of rank zero.
Thus the fiber of $\mathcal D(S)\to\mathcal D(l)$ is $\Sigma H\widehat{U(S)}_p$ with trivial $G$-action.
The relative group $K_1(S,l)$ is $U(S)$ because $K_2(S)\to K_2(l)$ is surjective by lifting Matsumoto symbols, and $K(S,l)$ is $1$-connective.
The determinant induces this identification on $\pi_1$ of the relative spectra.
The Henselian comparison gives $K(S,l)^\wedge_p\simeq \mathrm{TC}(S,l)^\wedge_p$ \cite[Theorem 4.36]{ClausenMathewMorrow}.
The relative comparison fiber is therefore $2$-connective, also after completion since the bottom Moore homotopy groups form a surjective tower.
The absolute comparison fiber over $S$ has the same bound, since the comparison over $l$ is an equivalence.
Homotopy orbits and subsequent $p$-completion preserve this bound, so the comparison over $S$ is surjective on $\pi_2$ of the completed orbits.
For a bounded-below $p$-complete spectrum, the homotopy orbits are already $p$-complete, since the $G$-bar construction has $p$-complete finite skeleta and the cofibers of the skeleta in its realization have increasing connectivity.
\end{proof}

We compute the reduction map on these determinant orbits.

\begin{lemma}\label{lem:second-two-determinant-orbits}
In this setup, assume that $l$ is a nonarchimedean local field.
Then the reduction map
\[
 \pi_2\bigl((\mathcal D(S)_{hG})^\wedge_p\bigr)\longrightarrow\pi_2\bigl((\mathcal D(l)_{hG})^\wedge_p\bigr)
\]
is surjective.
\end{lemma}
\begin{proof}
Put
\[
 \mathcal A=R\Gamma(l,H\underline{\Z})^\wedge_p,\qquad
 \mathcal M_\Lambda=R\Gamma(l,H\mathbf G_{m,\Lambda})^\wedge_p\quad(\Lambda\in\{S,l\}).
\]
The rank and unit terms in $\Sigma\mathcal M_\Lambda\to\mathcal D(\Lambda)\to\mathcal A$ have trivial action.
Write $\kappa_\Lambda$ for the connecting map with constant action and $\iota_{v_\Lambda}\colon\mathcal A\to\mathcal M_\Lambda$ for the map induced by $n\mapsto v_\Lambda^n$, where $v_S=v$ and $v_l=\bar v$.
The connecting map for the weak action is
\[
 \kappa_\Lambda+(\Sigma^2\iota_{v_\Lambda})\circ e_G(\mathcal A),
 \qquad e_G(\mathcal A)\colon\mathcal A\longrightarrow\Sigma^2\mathcal A,
\]
where $e_G(\mathcal A)$ is the tensor product of the standard double extension $H\Z_p\to\Sigma^2H\Z_p$, given by the class $e_G\in H^2(G,\Z)$ represented by $(g^i,g^j)\mapsto\lfloor(i+j)/p\rfloor$, with the rank module $\mathcal A$ with trivial action.
Indeed, before taking sections, the space of extensions with identified endpoints is the mapping space from $H\underline{\Z}$ to $\Sigma^2H\mathbf G_m$.
It is $2$-truncated, and the two actions agree on vertices and edges.
Their difference on each $2$-simplex is the scalar determinant $v_\Lambda^{n\lfloor(i+j)/p\rfloor}$.
The same calculation in the category of reduction arrows gives the constant-action comparison of the sign terms and the compatible carry terms.
For $p=2$, the residue sign is trivial because $-1=1$ in $l$, so $\kappa_l=0$.
For odd $p$, tensor square of the line gives a map to the unsigned graded-line spectrum that is an equivalence after $p$-completion, and again $\kappa_l=0$.
We retain $\kappa_S$ as a map of spectra.

We now determine the degree-two orbits.
The two nonzero homotopy groups of $\mathcal A$ are $\Z_p$ in degree zero and $X_l=H^1_{\mathrm{cont}}(l,\Z_p)$ in degree $-1$.
Those of $\mathcal M_l$ are $\widehat{l^\times}_p$ and $T_p\operatorname{Br}(l)=\Z_p$ in the same degrees.
These descriptions follow from $\operatorname{cd}_p(l)=1$, Hilbert~90 and the injectivity of the $p^a$-power map on the \'etale unit sheaf of $l$.
Reduction on $\mathcal M_\Lambda$ is surjective on $\pi_0$ and an isomorphism on $\pi_{-1}$ because its fiber is $H\widehat{U(S)}_p$.
The map $\iota_{\bar v}$ in degree zero sends $1$ to a primitive element of $\widehat{l^\times}_p$ because $\bar v\notin l^p$.
In degree $-1$ it is cup product with $\bar v$.
The local invariant pairing makes this map surjective modulo $p$, and $\operatorname{cd}_p(l)=1$ lifts degree-one characters to $\Z_p$-characters, so the map to $T_p\operatorname{Br}(l)$ is surjective \cite[\S5]{Kurihara2000}.

Put $X_l^0=\ker(X_l\xrightarrow{\cup\bar v}T_p\operatorname{Br}(l))$ and $U^0=\widehat{l^\times}_p/\Z_p\bar v$.
The modules $\widehat{l^\times}_p$ and $X_l$ are $p$-complete and $p$-torsion-free.
For the units this follows from the injectivity of multiplication by $p$ on $l^\times$, and for $X_l$ from the surjective reduction maps on $H^1(l,\Z/p^a)$ and the coefficient exact sequence.
A $p$-complete, $p$-torsion-free $\Z_p$-module is the completion of a free module obtained by lifting a basis modulo $p$.
Maps from the resulting Eilenberg--Mac Lane spectrum in derived $p$-complete $H\Z_p$-modules can thus be computed on those generators before completion.
We can therefore split the two homotopy degrees of $\mathcal A$ and $\mathcal M_l$ and split off the primitive summands in the two degreewise maps of $\iota_{\bar v}$.
Since $e_G(\mathcal A)$ is natural for the chosen decomposition of complete $H\Z_p$-modules, on the residue side, where $\kappa_l=0$, the carry survives only on the corresponding two summands.
Each of the two matched rank and unit summands of the carry is the extension represented by
\[
 0\longrightarrow\Z_p\xrightarrow{\ 1+g+\cdots+g^{p-1}\ }\Z_p[G]
 \xrightarrow{\ g-1\ }\Z_p[G]\longrightarrow\Z_p\longrightarrow0.
\]
Consequently its middle complex is free over $\Z_p[G]$ in degrees $1,0$, and its derived coinvariants are $\Z_p\xrightarrow{0}\Z_p$ in these degrees.
The complex for the degree $-1$ pairing is its shift by $-1$, with derived coinvariants in degrees $0,-1$.
Neither contributes in degree two.
The remaining summands have trivial action, so the degree-two determinant orbits are generated by $H_1(G,U^0)$ and $H_3(G,X_l^0)$.
The unit contribution lifts from $\Sigma\mathcal M_S$ by the surjectivity of unit reduction.

For the character contribution, we restrict the reduction arrow of extensions to the rank summand $\Sigma^{-1}HX_l^0$ of $\mathcal A$.
The restriction of $\iota_v$ is null in derived $p$-complete $H\Z_p$-modules because it is zero on $\pi_{-1}$ and the source is a shifted completed free module.
We choose a nullhomotopy over $S$ and reduce it to obtain the compatible nullhomotopy over $l$.
The carry then disappears on this reduction arrow, leaving the constant-action comparison, including the sign $\kappa_S$ when $p=2$.
Its relative boundary from the character summand lies in the orbits of the constant map
\[
 \Sigma^{-1}HX_l^0\longrightarrow\Sigma^2H\widehat{U(S)}_p .
\]
The contribution from $H_3(G,X_l^0)$ belongs to the reduced summand under $\Sigma^\infty BG_+\simeq\mathbb S\vee\Sigma^\infty BG$.
The constant map preserves this splitting, and $\pi_2$ of the reduced target is $\widetilde H_0(BG,\widehat{U(S)}_p)=0$.
Thus the character contribution lifts as well.
Changing its chosen splitting over $l$ changes its image by unit classes, which already lift.
\end{proof}

By Lemmas~\ref{lem:second-two-determinant-comparison} and~\ref{lem:second-two-determinant-orbits}, the comparison with $\mathrm{TC}$ now lifts every degree-two residue-orbit class from the commutative subcover.
We subtract such a lift to remove the residue image of the absolute degree-zero class.
Once the cover is forgotten, the subtracted term lies in $\pi_2\mathrm{TC}(\mathcal O_E)^\wedge_p$; since $\pi_2\mathrm{TC}(\Delta)^\wedge_p=0$, it lifts to a relative class.
By Corollary~\ref{cor:app-commutative-cut-representative}, the image of this relative class in the prescribed truncation is represented by an integral element of $K_2(\mathcal O_E,\Delta)$ whose generic-point image lies in the image of $\theta_E$.

\subsubsection{Integral transfer from orders and vanishing}
We now complete the proof of Theorem~\ref{thm:second-two-wild}.
Under the inductive hypothesis $K_2(B)=I_T(B)$, we will show in Lemma~\ref{lem:second-two-orbit-image} below that the image of the degree-two homotopy-orbit classes of the relative cover lies in $I_F(D)$.
After representing the image of the relative $\mathrm{TC}$ orbit class at a finite truncation by the trace of an integral class and applying the congruence separation, the injectivity of $\Nrd_2$ on $I_F(D)$ closes the induction.

We consider the relative spectrum $Y=K(H,\bar H)$ of the cover $H=\mathcal Q(P)$, where $\bar H$ is the residue cover used above.
The relative Morita comparison with the corner makes $Y$ $1$-connective, with $\pi_1Y=K_1(C,\Delta_0)$ and action $\tau=\operatorname{Ad}(x)$ of $G$.
In the homogeneous-uniformizer type, the relative $K_1$ classes of $a\in1+J(C)$ placed at different corners agree: for $i<j$ put $b=e_{ij}$ and $d=(a-1)e_{ji}\in J_H=\ker(H\to\bar H)$ and use $[1+bd]=[1+db]$; hence after vertex translation the remaining action is conjugation by the change of coordinates.
Since $\deg B=\ind B_{T^{\mathrm{nr}}}=m$, lifting Matsumoto symbols to a commutative order as in the proof of Lemma~\ref{lem:review-geometric-moore-integral-lift} makes $K_2(C)\to K_2(\Delta_0)$ surjective, so $\pi_1Y=M=\ker[K_1(C)\to\Delta_0^\times]$.
For the corner comparison in the homogeneous-uniformizer type, use the two-term derived Iwahori quotient of Lemma~\ref{lem:app-iwahori-derived-cut}.
Transitions on the corner cofibers of finite truncations vanish in one step.
Hence there is an exact sequence
\[
 (\pi_2Y)_G\longrightarrow \pi_2(Y_{hG})
    \longrightarrow H_1(G,M)\longrightarrow 0.
\]
The generic-point transfer of the left-hand image lies in $I_F(D)$ because $K_2(B)=I_T(B)$.
The right-hand edge can also be realized by integral subcovers, as follows.

\begin{lemma}
\label{lem:review-generic-source-edge-generation}
Let $Y\to Z$ be an equivariant map of spectra for a finite group $G$, with $Y$ $1$-connective and $Z$ connective.
Let $S_\lambda\to Z$ be maps of connective $G$-spectra, let $q\colon\pi_2(Z_{hG})\to \Lambda$ be a homomorphism, and let $I\subset \Lambda$ be a subgroup.
Assume that $q$ sends the image of $(\pi_2Z)_G$ and the entire image of every $\mathcal F_1\pi_2((S_\lambda)_{hG})$ into $I$.
If the image in $H_1(G,\pi_1Z)$ of the edge of $u\in\pi_2(Y_{hG})$ belongs to $\sum_\lambda\im[H_1(G,\pi_1S_\lambda)\to H_1(G,\pi_1Z)]$, then $q$ sends the image of $u$ in $\pi_2(Z_{hG})$ into $I$.
Here $\mathcal F_a$ denotes the filtration of the homotopy-orbit homology spectral sequence.
\end{lemma}
\begin{proof}
For a connective spectrum $X$, the first quotient in degree two is
\[
 \mathcal F_1/\mathcal F_0=E^\infty_{1,1}
 =H_1(G,\pi_1X)/\im
   [d_2\colon H_3(G,\pi_0X)\to H_1(G,\pi_1X)],
\]
and $\mathcal F_0\pi_2(X_{hG})=\im[(\pi_2X)_G\to\pi_2(X_{hG})]$.
Since $Y$ is $1$-connective, $u\in\mathcal F_1\pi_2(Y_{hG})$.
We express its edge image as a finite sum of classes from the specified source $H_1$ groups.
We lift each class in $E^\infty_{1,1}$ to the source $\mathcal F_1\pi_2$.
After their images in $Z$ are subtracted, naturality shows that the remainder lies in $\mathcal F_0\pi_2(Z_{hG})$.
The hypotheses on $q$ now apply.
\end{proof}

As above, we write $R=\mathcal O_F$, $P=\mathcal O_D$, $J_D=J(P)$, $C=\mathcal O_B$, and $J_C=J(C)$.
For $b\in C^\times$, write
\[
 w_b=bx,\qquad E_b=F(w_b).
\]
In the homogeneous-unit type, $\bar b\bar x\notin\Delta_0$, so $\bar b\bar x$ generates $\Delta/k$ and $E_b/F$ is a maximal subfield of residue degree $n$ and ramification index one.
In the homogeneous-uniformizer type, $v_D(w_b)=1$, so $E_b/F$ is a totally ramified maximal subfield of degree $n$.
Since $w_b$ acts nontrivially by conjugation on $T$, we have $w_b\notin B$, but $w_b^p\in B$.
By the Kummer presentation $X^p-w_b^p$ (the purely inseparable presentation in equal characteristic), we have $[E_b:F(w_b^p)]=p$ and $[F(w_b^p):F]=m$.
Moreover, $T\cap E_b=F$, so $T(w_b^p)/T$ is a maximal subfield of $B$.
A residue basis or an Eisenstein basis, respectively, gives
\[
 \mathcal O_{E_b}=\bigoplus_{j=0}^{p-1}\mathcal O_{F(w_b^p)}w_b^j.
\]
This is the commutative graded subcover used below.
For a Brauer class $\beta$ over a two-dimensional local field $L$, put $\chi_\beta(z)=\inv_L(\beta\cup z)$.

\begin{lemma}
\label{lem:review-homogeneous-unit-norm-generation}
With the notation above, let $U_D=1+J_D$, $U_C=1+J_C$, $U_F=1+\mathfrak m_F$, and $V=\Nrd_{B/T}(U_C)$.
For any valued field $L$, write $U(L)=1+\mathfrak m_L$ and $U_T=U(T)$.
In the homogeneous-unit type, the ratios
\[
 \frac{N_{E_b/F}f(bx)}{N_{F(x)/F}f(x)}
\]
for $b\in1+J_C$ and $f(X)=1+hX$ ($h\in R$), together with $N_{T/F}V$, generate $\ker\chi_{[D]}\cap U_F$.
In the homogeneous-uniformizer type, the groups $N_{E_b/F}U(E_b)$ ($b\in C^\times$), together with $N_{T/F}V$, generate the same group.
One can choose $u\in U_F$ such that these generators, together with $u^m=N_{F(x^p)/F}u$, generate $V^G=\ker\chi_{p[D]}\cap U_F$.
\end{lemma}
\begin{proof}
First approximate principal units of $D$ one valuation layer at a time.
In the homogeneous-unit type, let $s$ be a uniformizer of $C$, take $c\in C$, and put $b=1+cs^r$.
The leading coefficient of $f(bx)/f(x)$ at depth $r$ is
\[
 \bar c\,\frac{\bar x f'(\bar x)}{f(\bar x)}.
\]
This follows by expanding the terms in which $b-1$ occurs exactly once.
We choose distinct nonzero elements $h_1,\ldots,h_p$ of the infinite field $k$, choose lifts $\widetilde h_i\in R$, and set $f_i=1+\widetilde h_iX$.
If $\sum_i c_i h_i\bar x/(1+h_i\bar x)=0$ with $c_i\in\Delta_0$, divide by $\bar x$ and clear denominators.
Then the polynomial $\sum_i c_i h_i\prod_{j\ne i}(1+h_jX)\in\Delta_0[X]$ has degree at most $p-1$ and vanishes at $\bar x$.
Since $[\Delta_0(\bar x):\Delta_0]=p$, this polynomial is zero.
Substitution of $X=-1/h_i$ gives $c_i=0$ for every $i$.
These logarithmic derivatives therefore form a basis of $\Delta_0(\bar x)$ over $\Delta_0$.
Finite products of the ratios thus realize every leading coefficient of a principal unit.

In the homogeneous-uniformizer type, we use $x^r$ as reference.
If $p\mid r$, principal units of the corner $C$ give all coefficients of that layer.
If $p\nmid r$, the coefficient of $1+a(bx)^r$ is $\bar a\bar b^r$.
Here $\Delta_0=\Delta$; choose $e_0$ with $\Delta^{p^{e_0}}\subset k$ and integers $s_0,t_0$ with $rs_0+p^{e_0}t_0=1$.
Every $z\in\Delta^\times$ can then be written $z=(z^{p^{e_0}})^{t_0}(z^{s_0})^r$.
We lift $\bar a=(z^{p^{e_0}})^{t_0}\in k^\times$ to the center $R$ and $\bar b=z^{s_0}\in\Delta^\times$ to $C^\times$.
The single principal unit $1+a(bx)^r$ realizes the leading coefficient $z$.
Hence every layer is realized in this type as well.

After finitely many layers, the norm of the deep remainder lies in $N_{T/F}V$.
Indeed, by separability of $T/F$ and nondegeneracy of the reduced trace, $\Tr_{T/F}\operatorname{Trd}_{B/T}$ is a nonzero $F$-linear functional, so after scaling by a power of the fixed uniformizer $\pi$ of $F$ we can choose $z\in J_C$ with $\Tr_{T/F}\operatorname{Trd}_{B/T}(z)\ne0$.
The polynomial
\[
 g(X)=N_{T/F}\Nrd_{B/T}(1+Xz)
\]
satisfies $g(0)=1$ and $g'(0)\ne0$, so by Hensel's lemma its image contains all sufficiently deep principal units of $F$.
Thus finitely many layer corrections and one deep norm correction suffice, without using an infinite product.
This proves the reduced-norm generation in both types.

Kato's norm-image formula and the purely inseparable residue extension give
\[
 \Nrd_{D/F}(U_D)=\ker\chi_{[D]}\cap U_F,
 \qquad V=\ker\chi_{[B]}\cap U_T.
\]
For example, if the reduced norm of a unit has residue one, that residue is the purely inseparable degree power of the residue of the unit, so the unit itself has residue one.
The projection formula and restriction formula for invariants give $V^G=\ker\chi_{p[D]}\cap U_F$.

Finally, Kato duality shows that the image of $\chi_{[D]}$ on $U_F$ has order $n$: the annihilator of $U_F$ is $\ker(\Br(F)\to\Br(F^{\mathrm{nr}}))$.
Indeed, the pairings with the two Witt components $\Br(k)$ and $H^1(k,\Q/\Z)$ give, respectively, the finite-order continuous characters of the valuation and residue-unit components of $F^\times/U_F$.
Ordinary local reciprocity and Kato duality show that these are all such characters.
Minimality of the index gives $m[D]\notin\ker(\Br(F)\to\Br(F^{\mathrm{nr}}))$, so the image on $U_F$ has order $n$ \cite[Theorem 2]{KatoExistence}.
We choose $u\in U_F$ such that $\chi_{[D]}(u)$ has order $n$.
Then $u^m=N_{F(x^p)/F}u$ belongs to $\ker\chi_{p[D]}\cap U_F$, and $\chi_{[D]}(u^m)$ has order $p$.
Since $\chi_{p[D]}=p\chi_{[D]}$ and $\chi_{[D]}(U_F)\subset\Q/\Z$ is cyclic of order $n$, the map $\chi_{[D]}$ identifies
\[
 (\ker\chi_{p[D]}\cap U_F)/(\ker\chi_{[D]}\cap U_F)
\]
with the subgroup of order $p$ of $\chi_{[D]}(U_F)$.
This quotient therefore has order $p$ and is generated by the class of $u^m$, so the group $\ker\chi_{[D]}\cap U_F$ and the element $u^m$ generate $V^G$.
\end{proof}

\begin{lemma}\label{lem:second-two-orbit-image}
In the situation of Lemma~\ref{lem:second-two-cyclic-cover}, assume that $K_2(B)=I_T(B)$, and let $Y=K(H,\bar H)$ be the relative spectrum of the cover $H=\mathcal Q(P)$ as above.
Then the map $\pi_2(Y_{hG})\to K_2(D)$, which forgets the cover through $K(B)_{hG}$ and passes to the generic point, satisfies
\[
 \im\bigl[\pi_2(Y_{hG})\longrightarrow K_2(D)\bigr]\subseteq I_F(D).
\]
\end{lemma}
\begin{proof}
We describe integral sources for the generators in Lemma~\ref{lem:review-homogeneous-unit-norm-generation}.
In the homogeneous-uniformizer type, use the relative Iwahori cover of $E_b$.
The relative corner of this cover is $\mathcal O_{F(w_b^p)}$, with degree-one group $U(F(w_b^p))$.
The element $w_b$ fixes $F(w_b^p)$, and $\operatorname{Ad}(b)$ acts trivially on $K_1(C)$; hence the class in $M$ of $a_b\in U(F(w_b^p))$ is $G$-invariant.
We take $a_b=N_{E_b/F(w_b^p)}z$ with $z\in U(E_b)$, so that by transitivity its reduced norm is $N_{F(w_b^p)/F}a_b=N_{E_b/F}z$.
The action and its carry are taken from the actual graded subring $\mathcal O_{E_b}=\bigoplus_{j=0}^{p-1}\mathcal O_{F(w_b^p)}w_b^j$, using $w_b^p$.

In the homogeneous-unit type, put $a_b=N_{E_b/F(w_b^p)}f(bx)$ and $a_1=N_{F(x)/F(x^p)}f(x)$.
Here $b\in1+J_C$, so $bx$ and $x$ have the same reduction, and the $p$-power norm of $\Delta/\Delta_0$ gives $\bar a_b=\bar a_1$.
Set $R_b=\mathcal O_{E_b}\times_\Delta\mathcal O_{F(x)}$ and $S_b=\tau_{\ge0}K(\mathcal Q(R_b))$.
Both projections are graded and have the same composite to the residue, so their difference, with its canonical residue nullhomotopy, gives a $G$-equivariant map $S_b\to Y$.
The corner unit pair $(a_b,a_1)$ is a $G$-invariant class of $\pi_1S_b$ whose image in $M$ is $[a_b]-[a_1]$.
Indeed, the degree-one unit $r=(bx,x)$ commutes with the degree-zero corner, and vertex translation fixes this unit pair.
The carry of the weak action is $r^p$.
For the additional class $u$, use the commutative relative cover of $F(x)$ and its corner principal unit $u\in U(F(x^p))$.

We apply Lemma~\ref{lem:review-generic-source-edge-generation} with these sources, $Y\to K(B)$, and $I_F(D)\subset K_2(D)$.
The map that forgets source orbits factors, as a spectrum map, through $K(E_b)$; at the generic point it is extension of scalars from $E_b$ to $D$, that is, the normalized transfer $\theta_{E_b}$ given by the $(E_b,D)$-bimodule $D$.
For the source $S_b$, it factors through the difference of the two maps at the generic point $R_b\otimes_RF=E_b\times F(x)$.
Thus the \emph{entire} image of the first filtration step lies in $I_F(D)$, through $\theta_{E_b}$ or $\theta_{E_b}-\theta_{F(x)}$, respectively.
The image of the corner $K_2(B)$ also lies in $I_F(D)$ by the inductive hypothesis $K_2(B)=I_T(B)$ and transitivity of transfer.
On the other hand, since $G$ is cyclic, we have natural isomorphisms $H_1(G,M)\simeq M^G/N_GM$ and $H_1(G,K_1(B))\simeq K_1(B)^G/N_GK_1(B)$.
Since $\SK_1(B)=0$, the reduced norm $\Nrd_{B/T}\colon K_1(B)\to\Nrd_{B/T}(B^\times)$ is a $G$-equivariant isomorphism, and $\Nrd_{B/T}(M)=V$ because the image of $U_C$ in $K_1(C)$ lies in $M$ and $\Nrd_{B/T}(M)\subseteq U_T\cap\ker\chi_{[B]}=V$.
Hence $\Nrd_{B/T}$ maps the image of $H_1(G,M)$ in $H_1(G,K_1(B))$ onto the image of $\Nrd_{B/T}(M^G)\subseteq V^G$ in $\Nrd_{B/T}(B^\times)^G/N_{T/F}\Nrd_{B/T}(B^\times)$.
By Lemma~\ref{lem:review-homogeneous-unit-norm-generation} and $N_{T/F}V=\Nrd_{B/T}(N_GM)$, the generic-point image of $H_1(G,M)$ is therefore generated by these source edges.
We therefore obtain the asserted inclusion in both types.
\end{proof}

\begin{proof}[Completion of the proof of Theorem~\ref{thm:second-two-wild}]
The case $n=1$ is split.
Assume $n>1$, and let $\xi\in\SK_2(D)$.
For $c=\partial_D\xi$, the boundary comparison in the proof of Lemma~\ref{lem:review-minimal-order-k1-support} gives $c=\partial_F\Nrd_2\xi=1$, so $\xi$ lifts to some $y\in K_2(P)$.
We first show that the support map $i_*\colon K_2(\Delta)\to K_2(P)$ vanishes.
Lifting Matsumoto symbols through a maximal subfield of ramification index one makes the reduction $r\colon K_2(P)\to K_2(\Delta)$ surjective.
For a uniformizer $x_0$ of $D$, put $\alpha=\operatorname{Ad}(x_0^{-1})_*$ on $K_2(P)$.
Additivity for the bimodule sequence $J_D\to P\to\Delta$ gives $i_*r=1-\alpha$, since $J_D=x_0P$ and the left action on $J_D$ is twisted by $\operatorname{Ad}(x_0^{-1})$, as in the proof of Lemma~\ref{lem:second-residue-zero-ell-integral}.
Conjugation by $x_0$ fixes the purely inseparable residue, so $r\alpha=r$ and $(1-\alpha)^2=i_*r(1-\alpha)=i_*r-i_*r\alpha=0$.
Since $x_0^n\in\pi P^\times$ and $\pi$ is central, $\operatorname{Ad}(x_0^n)$ is conjugation by a unit of $P$, so $1=\alpha^n=(1-(1-\alpha))^n=1-n(1-\alpha)$.
Thus $ni_*r=n(1-\alpha)=0$, and surjectivity of $r$ gives $ni_*=0$, so $i_*$ from the uniquely $p$-divisible group $K_2(\Delta)$ is zero.
By the localization sequence $K_2(\Delta)\xrightarrow{i_*}K_2(P)\to K_2(D)$, the map $K_2(P)\to K_2(D)$ is therefore injective.
Hence $n\xi=0$ gives $ny=0$, and $r(y)=0$ because $K_2(\Delta)$ is uniquely $p$-divisible.

Fix a finite separable splitting field $L_0/F$, and choose a depth $j$ at which Lemma~\ref{lem:app-congruence-separation} holds for $L_0$.
We lift the trace of $y$ by the surjection in degree two from absolute degree-zero $\mathrm{TC}$ to $\mathrm{TC}(P)^\wedge_p$ of Lemma~\ref{lem:second-two-degree-zero-surjective} and remove the residue component of the lift.
The order of the steps matters: in the homogeneous-unit type the residue orbit is removed by the $\mathrm{TC}$ class of the commutative subcover (in the homogeneous-uniformizer type, $\pi_2$ of the residue orbit already vanishes), the image with the cover forgotten is lifted to relative $\mathrm{TC}$ and compared there using $\pi_3\mathrm{TC}(\Delta)^\wedge_p=0$, and only then is it represented, at the prescribed depth $j$, by an integral class.
In the mixed-characteristic homogeneous-unit type, the image of the subtracted $\mathrm{TC}$ class with the cover forgotten is represented at depth $j$ by an integral class from $\mathcal O_E$, already compared on the residue and on the prescribed truncation.
By the equivalence in Lemma~\ref{lem:review-cyclic-zero-tc}, the remaining relative orbit class is $t_{\mathrm{rel}}\in\pi_2((\mathrm{TC}(H,\bar H)^\wedge_p)_{hG})$.
Here the outer $p$-completion in that lemma can be removed, since $\mathrm{TC}(H,\bar H)^\wedge_p$ is bounded below (the absolute terms lie in $\mathrm{Sp}_{\ge-1}$ and the relative term in $\mathrm{Sp}_{\ge-2}$) and $G$ is finite.
We then send $t_{\mathrm{rel}}$ to depth $j$: by the last assertion of Lemma~\ref{lem:app-orbit-stable-image}, its image at depth $j$ is represented by the trace of an integral element of $\pi_2(Y_{hG})$.
The generic-point image of that integral element lies in $I_F(D)$ by Lemma~\ref{lem:second-two-orbit-image}.
To compare with $y$, note that $r(y)=0$, so $y$ lifts to a relative class $y_{\mathrm{rel}}$.
In relative $\mathrm{TC}$ we have
\[
 \operatorname{tr}(y_{\mathrm{rel}})=f_*(t_{\mathrm{rel}})+b_{\mathrm{rel}},
\]
where $f_*$ forgets the cover and $b_{\mathrm{rel}}$ is the relative correction added back in the mixed-characteristic homogeneous-unit type (zero otherwise).
The absolute equality gives this relative one because $\pi_3\mathrm{TC}(\Delta)^\wedge_p=0$, so that the map from $\pi_2$ of relative $\mathrm{TC}$ to $\pi_2$ of absolute $\mathrm{TC}$ is injective.
We send this agreement to the prescribed finite truncation, where the relative trace comparison for nilpotent ideals (Lemma~\ref{lem:app-nilpotent-trace}) gives the agreement with integral coefficients.
Adding back the integral representatives, namely the representative from the cover and, in the mixed-characteristic homogeneous-unit type, the subtracted representative from $\mathcal O_E$, produces $z\in K_2(P)$ with $z|_D\in I_F(D)$ and $y-z\in\ker[K_2(P)\to K_2(P/\pi^jP)]$.
By Lemma~\ref{lem:app-congruence-separation}, we have $(y-z)|_D\in\theta_{L_0}K_2(L_0)$.
Thus $\xi\in I_F(D)$, and injectivity in Proposition~\ref{prop:review-splitting-transfer-norm-iso} gives $\xi=0$.
This closes the induction, and unramified descent in Lemma~\ref{lem:second-two-inertial-reduction} gives the result for the original $F,D$.
\end{proof}

We reduce Theorem~\ref{thm:second-two-vanishing} to Theorem~\ref{thm:second-two-wild}.

\begin{proof}[Proof of Theorem~\ref{thm:second-two-vanishing}]
Surjectivity is proved in Theorem~\ref{thm:second-two-surj}.
When the first residue field $k$ has characteristic zero, write $F=k((t))$ by Cohen's structure theorem, and let $B$ be the residue division algebra of $D$, with center $L$.
Since $L$ is a nonarchimedean local field, Theorem~\ref{thm:paper-local-norm} gives $\SK_1(B)=0$ and the isomorphism $\Nrd_{B/L,2}\colon K_2(B)\xrightarrow{\sim}K_2(L)$.
This isomorphism respects the action of $\Gal(L/k)$, so Proposition~\ref{prop:second-residue-zero} and cyclic Hilbert 90 for $K_2$ \cite[Theorem 14.1]{MerkurjevSuslin1982} give $\SK_2(D)=0$.
Suppose now that $k$ has characteristic $p>0$; then $k$ is a nonarchimedean local field.
For a prime $\ell\ne p$, an extension of degree prime to $\ell$ splitting the other primary Brauer components reduces the $\ell$-primary part of $\SK_2(D)$, by restriction and transfer, to that of an $\ell$-primary division algebra, and Proposition~\ref{prop:second-residue-zero-ell}(2) shows that it vanishes.
After an extension of degree prime to $p$ splitting the prime-to-$p$ Brauer components, restriction and transfer show that the $p$-primary part of $\SK_2$ is detected over that extension.
By Morita equivalence, the remaining cases are therefore exactly those of Theorem~\ref{thm:second-two-wild}.
\end{proof}

\section{Reduced norms over three-dimensional local fields}\label{sec:three-local}
This section treats three-dimensional local fields.
We compute the image and the cokernel of the reduced norm on $K_2$ in every characteristic (Section~\ref{sec:second-three-coker}).
We then show that $\SK_2$ vanishes when the first residue field has characteristic zero (Section~\ref{sec:second-three-charzero}), and we discuss the vanishing of $\SK_2$ in general (Section~\ref{sec:paper-further-results}).
\subsection{The cokernel of the reduced norm}\label{sec:second-three-coker}
In this subsection, we compute the image and the cokernel of the reduced norm on $K_2$ over three-dimensional local fields (Theorem~\ref{thm:second-three-coker}).
For a three-dimensional local field $F$, we write $\inv_F\colon H^4(F,\Q/\Z(3))\xrightarrow{\sim}\Q/\Z$ for Kato's isomorphism \cite[Theorem 2]{KatoExistence}, as in the Conventions.
For a three-dimensional local field $F$ and a central division algebra $D$ over $F$ with $\alpha=[D]$, duality reduces the computation to showing that a Brauer class that vanishes over every finite splitting field of $\alpha$ lies in the subgroup generated by $\alpha$.

For a field $F$ and $z\in\Br(F)$, let $B_z(F)$ denote the common kernel in $\Br(F)$ of restriction to all finite splitting fields of $z$.
For the primary decomposition $z=\sum_\ell z_\ell$, one has $B_z(F)=\bigoplus_\ell B_{z_\ell}(F)$.
Indeed, $B_z(F)$ is killed by $\ind(z)$.
Over a splitting field of $z_\ell$, we split the other components in degree prime to $\ell$ by successively adjoining separable maximal subfields for them; restriction and corestriction then show that the $\ell$-component lies in $B_{z_\ell}(F)$.
The reverse inclusion follows from the definition.
The \emph{Witt part} is the subgroup of Brauer classes killed by the maximal unramified extension with respect to the top discrete valuation of $F$, whose residue field is the first residue field $k$.
A class in the Witt part has a Witt decomposition $i(a)+(\chi,\pi)$ with $a\in\Br(k)$ and $\chi\in H^1(k,\Q/\Z)$; we call it purely unramified if $\chi=0$.
In first residue characteristic $p$, we call a $p$-primary class wild if it does not belong to the Witt part.
For a two-dimensional local field $F$, the equality $B_z(F)=\langle z\rangle$ follows from Kato's degree-one norm theorem and duality, since the sum of $K_1$ norms from all splitting fields is the reduced-norm image of the division algebra in the Brauer class $z$ \cite[\S3.1, Theorem 3(2)]{KatoLocalII}, \cite[Theorem 2]{KatoExistence}.
Here, for the division algebra $D$ in the class $z$, the group $K_1(D)$ is generated by transfers from separable splitting fields: elements generating separable maximal subfields form a nonempty open subset of $D^\times$, and every element is a ratio of two of them.

\begin{lemma}\label{lem:second-three-norm-quotient}
Let $F$ be a local field of dimension $h\ge2$, let $\ell$ be a prime, and let $L/F$ be a cyclic extension of degree $\ell$.
Then $F^\times/N_{L/F}L^\times$ is an $\F_\ell$-vector space of dimension at least two.
\end{lemma}
\begin{proof}
We use $v_F(N_{L/F}x)=f(L/F)\,v_L(x)$.
The quotient is killed by $\ell$.
Let $k$ be the first residue field.
If the extension is unramified, there are the valuation quotient $\Z/\ell$ and the residue norm quotient of units.
For the cyclic residue extension $l/k$, if the residue norm quotient of units vanished, the projection formula would make the Milnor norm $K_{h-1}^M(l)\to K_{h-1}^M(k)$ surjective, whereas its cokernel is $\Gal(l/k)\simeq\Z/\ell$ by higher local class field theory \cite[\S0, p.\ 168]{KatoExistence}.
If the extension is totally ramified, residues of unit norms lie in $k^{\times\ell}$.
A preimage of a unit norm is a unit, so the images of units in the norm quotient surject onto $k^\times/k^{\times\ell}$.
When $\ell\ne\operatorname{char}k$, tame cyclic total ramification implies $\mu_\ell\subset k$.
Reduction preserves the orders of roots of unity at invertible primes, and at a mixed-characteristic step these orders are bounded by $v(\zeta_{\ell^r}-1)=v(\ell)/(\ell^{r-1}(\ell-1))$.
Thus $k$ has a root of unity of maximal $\ell$-power order.
Such a root and a uniformizer give independent classes in $k^\times/k^{\times\ell}$.
If $\ell=\operatorname{char}k$, the classes of a uniformizer $w$ and $1+w$ are independent.
Indeed, valuation first removes the exponent of $w$.
For a nonzero power of $1+w$ with exponent less than $\ell$, subtracting one gives valuation one, but subtracting one from an $\ell$-th power gives valuation divisible by $\ell$.
The remaining case has ramification index one and purely inseparable residue degree $\ell$ (here $\ell=\operatorname{char}k$).
If a norm is a first principal unit, the residue of its preimage is one, so the preimage is also a first principal unit.
The norm of the preimage lies in the second principal units: for a common uniformizer $\pi$ one has $N_{L/F}(1+\pi a)\equiv1+\pi\Tr_{l/k}(\bar a)\pmod{\pi^2}$, and $\Tr_{l/k}=0$ for the purely inseparable residue extension $l/k$.
Thus $N_{L/F}L^\times\cap U_F^1\subseteq U_F^2$, so the image of the principal units in the norm quotient surjects onto $U_F^1/U_F^2\simeq(k,+)$.
Since $k$ is infinite, the dimension is at least two.
These cases exhaust the possibilities because a complete discretely valued field satisfies $[L:F]=e(L/F)f(L/F)$.
\end{proof}

\begin{lemma}\label{lem:second-three-cyclic-family}
Let $\ell$ be a prime, let $F$ be a local field of dimension at least two, and let $\eta\in\Br(F)\{\ell\}$ and $\theta\in H^1(F,\Q_\ell/\Z_\ell)$.
If $\theta$ has order $n$, then
\[
 \bigcap_{u\in F^\times}\langle\eta+(\theta,u)\rangle
   =\langle n\eta\rangle .
\]
\end{lemma}
\begin{proof}
The case $n=1$ is clear, so assume $n>1$.
By Lemma~\ref{lem:second-three-norm-quotient} and the description of cyclic relative Brauer groups applied to the character $(n/\ell)\theta$ of order $\ell$, the group $(n/\ell)(\theta,F^\times)$ also has dimension at least two.
We may therefore choose $\delta=(\theta,u)$ with $(n/\ell)\delta\notin\langle\eta\rangle$.
Then $\delta$ has order $n$ and $\langle\delta\rangle\cap\langle\eta\rangle=0$, since every nonzero subgroup of the cyclic $\ell$-group $\langle\delta\rangle$ contains $(n/\ell)\delta$; hence
\[
 \langle\eta\rangle\cap\langle\eta+\delta\rangle
     =\langle n\eta\rangle .
\]
The terms for $u=1$ and this $u$ give one inclusion, and the other follows because every $(\theta,u)$ is killed by $n$.
\end{proof}

\begin{lemma}\label{lem:second-three-witt-common-kernel}
Let $F$ be a complete discretely valued field whose first residue field $k$ is an $h$-dimensional local field with $h\ge2$.
For a prime $\ell$, assume $B_z(k)=\langle z\rangle$ for every $z\in\Br(k)\{\ell\}$.
If $\alpha\in\Br(F)\{\ell\}$ is split by a finite unramified extension, then $B_\alpha(F)=\langle\alpha\rangle$.
\end{lemma}
\begin{proof}
Restriction and corestriction through a separable maximal subfield make $B_\alpha(F)$ $\ell$-primary.
Every $\beta$ in this group is killed by the same finite unramified splitting field as $\alpha$, so write the Witt decompositions
\[
 \alpha=i(a)+(\chi,\pi),\qquad\beta=i(b)+(\psi,\pi).
\]
Here $\pi$ is a uniformizer, $i$ is unramified lifting, $a,b\in\Br(k)\{\ell\}$, and $\chi,\psi\in H^1(k,\Q_\ell/\Z_\ell)$.
This decomposition of unramifiedly split classes holds at every prime and respects residue-field restriction for finite extensions of ramification index one \cite[\S3, (3.3)--(3.4), Theorem 3.2]{AljadeffSonnWadsworth}.
Put $m=\operatorname{ord}\chi$ and $L=k(\chi)$, with $m=1,L=k$ when $\chi=0$.
Every finite $M/k$ lifts to a finite extension $F_M/F$ of ramification index one and residue field $M$: lift monic minimal polynomials successively along a tower of simple extensions.
The lifted ring is a complete discrete valuation ring, since its residue is a field and its maximal ideal is generated by $\pi$.
For every finite splitting field $M/L$ of $a_L$, both $a_M$ and $\chi_M$ vanish, since $M$ splits $a_L$ and contains $L$, so $\alpha_{F_M}=0$ by the Witt decomposition.
Since $\beta\in B_\alpha(F)$, also $\beta_{F_M}=0$, and uniqueness of the Witt decomposition over $F_M$ gives $\psi_M=0$.
By \cite[Proposition 5]{KatoExistence} with $n=q=h$ and $r=2$, we have $\sum_M N_{M/L}K_h^M(M)=K_h^M(L)$, where $M$ runs over these splitting fields.
The projection formula and duality on the $H^1(L,\Q_\ell/\Z_\ell)$ side therefore give $\psi_L=0$, hence $\psi=c\chi$ for some integer $c$ \cite[Proposition 5, Theorem 2]{KatoExistence}.
Replacing $\beta$ by $\beta-c\alpha$, we may assume $\beta=i(b)$.

We take any $u\in k^\times$ and any finite splitting field $E/k$ of $z=a+(\chi,u)$.
We adjoin $t^m=\pi/\tilde u$ to the lift $F_E$, where $\tilde u$ is a unit lift of $u$.
For $m>1$ this is an Eisenstein extension, so the residue remains $E$.
Over this field, writing $\pi=\tilde u\,t^m$, we have $(\chi,\pi)=(\chi,\tilde u)+m(\chi,t)=i((\chi,u))$ because $m\chi=0$; hence $\alpha=i(z_E)=0$ there, and so $i(b_E)=0$.
Unramified lifting is injective, so $b_E=0$.
The hypothesis and Lemma~\ref{lem:second-three-cyclic-family}, applied over $k$ with $\eta=a$ and $\theta=\chi$ of order $n=m$, give
\[
 b\in\bigcap_{u\in k^\times}\langle a+(\chi,u)\rangle
   =\langle ma\rangle,\qquad i(ma)=m\alpha .
\]
Thus the original $\beta$ is an integral multiple of $\alpha$.
\end{proof}

\begin{theorem}\label{thm:second-three-coker}
Let $F$ be a three-dimensional local field, let $D$ be a central division algebra over $F$, and let $\alpha$ be the Brauer class of $D$.
Then
\[
\begin{aligned}
 N_2(\alpha)=\Nrd_2K_2(D)
       &=\ker(\alpha\cup-:K_2(F)\to H^4(F,\Q/\Z(3))),\\
 C_2(D/F)&\xrightarrow[\ \inv_F(\alpha\cup-)\ ]{\ \sim\ }
       \tfrac{1}{\per(\alpha)}\Z/\Z .
\end{aligned}
\]
\end{theorem}

We prove Theorem~\ref{thm:second-three-coker} at the end of this subsection, after a sequence of lemmas.
For a three-dimensional local field $F$ and $\alpha\in\Br(F)$, we write $B_\alpha$ for $B_\alpha(F)$.
Lemma~\ref{lem:second-three-duality} reduces the theorem to the equality $B_\alpha=\langle\alpha\rangle$, and Lemma~\ref{lem:second-three-containment} in Section~\ref{sec:second-two-tame} places the image of the reduced norm in the kernel of the cup product.
For a wild class of period $p$, the remaining lemmas pass to an extension over which the relevant classes stay nonzero and $\alpha$ has index $p$ (Lemma~\ref{lem:second-three-wild-descent}), and they give the norm image there (Lemmas~\ref{lem:second-three-kahn-levine} and~\ref{lem:second-three-char-p-index-p}).
In the proof of the theorem, classes split by an unramified extension are treated first by Lemma~\ref{lem:second-three-witt-common-kernel}, and the equality $B_\alpha=\langle\alpha\rangle$ over $F$ for wild $p$-primary classes of higher period is proved by induction.

\begin{lemma}\label{lem:second-three-duality}
Let $F$ be a three-dimensional local field and $\alpha\in\Br(F)$.
Then $N_2(\alpha)$ is an open subgroup of finite index in $K_2(F)$, and Kato duality gives
\[
 (K_2(F)/N_2(\alpha))^\vee=B_\alpha .
\]
Moreover, $B_\alpha=\langle\alpha\rangle$ if and only if $N_2(\alpha)=\ker(\alpha\cup-)$, and in this case $\inv_F(\alpha\cup-)$ induces an isomorphism $K_2(F)/N_2(\alpha)\simeq\tfrac1{\per(\alpha)}\Z/\Z$.
\end{lemma}
\begin{proof}
The subgroup $N_2(\alpha)$ of $K_2(F)$ is open and of finite index \cite[Theorems 2--3]{KatoExistence}, so all characters of $K_2(F)/N_2(\alpha)$ are continuous.
Under Kato duality, a continuous character of $K_2(F)$ has the form $x\mapsto\inv_F(\beta\cup x)$ for a unique Brauer class $\beta$.
For each finite splitting field $E/F$, separable or not, the projection formula and the compatibility of the top invariant with corestriction show that this character vanishes on $N_{E/F}K_2(E)$ if and only if $\beta_E=0$ \cite[\S3.2, Lemma 1(2) and Proposition 1(2)(ii)]{KatoLocalII}.
Thus the dual of the finite quotient $K_2(F)/N_2(\alpha)$ is $B_\alpha$.

The class $\alpha$ generates $B_\alpha$ if and only if its character is faithful.
Since $\inv_F$ is an isomorphism, this holds if and only if $N_2(\alpha)=\ker(\alpha\cup-)$, and then the image of the character is the subgroup of order $\per(\alpha)$.
\end{proof}

\begin{lemma}\label{lem:second-three-kahn-levine}
Let $B$ be a central simple algebra over a field $k$ whose index $e$ is squarefree and prime to $\operatorname{char}k$.
Then
\[
 K_2(B)=\sum_{\substack{E/k\ \mathrm{finite\ separable}\\ B_E\ \mathrm{split}}}\theta_EK_2(E),
\]
and $\Nrd_2K_2(B)=N_2([B])$ is the kernel of
\[
 [B]\cup-\colon K_2(k)\longrightarrow H^4\bigl(k,(\Z/e)(3)\bigr).
\]
\end{lemma}
\begin{proof}
By \cite[Proposition 6.6.1(2)]{KahnLevine2010}, $K_2(B)$ is generated by the images of $\theta_E$ for finite splitting fields $E/k$.
If a finite splitting field $E/k$ is not separable, let $E_s$ be its maximal separable subextension; the index of $B_{E_s}$ divides both $e$, which is prime to $\operatorname{char}k$, and $[E:E_s]$, a power of $\operatorname{char}k$, so $E_s$ splits $B$ and $\theta_E=\theta_{E_s}N_{E/E_s}$.
This gives the displayed generation, and consequently $\Nrd_2K_2(B)=N_2([B])$.

The inclusion of $N_2([B])$ in the kernel of $[B]\cup-$ follows from the projection formula.
For the reverse inclusion, suppose first that $k$ is perfect.
Then \cite[Theorem 8.2.2(4)]{KahnLevine2010}, which uses the Beilinson--Lichtenbaum comparison in weights up to three \cite{Voevodsky}, identifies $N_2([B])$ with the kernel of $[B]\cup-$ after the exponential characteristic is inverted.
Both quotients of $K_2(k)$, by $N_2([B])$ and by the kernel of $[B]\cup-$, are killed by $e$: the first because $N_{E/k}\res_{E/k}=e$ for a separable maximal subfield $E$ of degree $e$ of the division algebra in the class of $B$, and the second because the coefficients are $\Z/e$.
Since $e$ is prime to $\operatorname{char}k$, inverting the exponential characteristic does not change the comparison, so the identification holds integrally.
If $\operatorname{char}k=q>0$ and $k$ is not perfect, let $z$ lie in the kernel of $[B]\cup-$.
Over the perfect closure, $z$ lies in the norm subgroup; we descend the resulting norm presentation to a finite purely inseparable stage and transfer it back.
This puts a $q$-power multiple of $z$ into $N_2([B])$, and since $q$ is prime to $e$, the element $z$ itself lies in $N_2([B])$.
\end{proof}

\begin{lemma}\label{lem:second-three-coefficients}
Let $F'$ be a three-dimensional local field and $p$ a prime.
The natural map
\[
 H^4(F',\Z/p(3))\longrightarrow H^4(F',\Q_p/\Z_p(3))
\]
is injective, and the top invariant identifies $H^4(F',\Q_p/\Z_p(3))$ with $\Q_p/\Z_p$.
If $p=\operatorname{char}F'$, the coefficients are logarithmic.
\end{lemma}
\begin{proof}
By the coefficient exact sequence, the kernel is the cokernel of multiplication by $p$ on $H^3(F',\Q_p/\Z_p(3))$.
The norm residue isomorphism in degree three identifies this group with the $p$-divisible group $K_3^M(F')\otimes\Q_p/\Z_p$: for $p\ne\operatorname{char}F'$ by \cite[Theorems 6.16 and 6.17]{Voevodsky}, and for $p=\operatorname{char}F'$ by \cite[Corollary (2.8)]{BlochKato1986}.
For $p=\operatorname{char}F'$, successive Kato residue maps followed by the trace on the finite last residue field identify $H^4(F',\Q_p/\Z_p(3))$ with $\Q_p/\Z_p$, and the top invariant does the same in general \cite[\S5.2]{Kurihara2000}.
\end{proof}

\begin{lemma}\label{lem:second-three-char-p-index-p}
Let $F'$ be a three-dimensional local field of characteristic $p$, let $B$ be a central simple algebra over $F'$ of index $p$, and put $\beta=[B]$.
Then
\[
 N_2(\beta)
 =\ker\bigl(\beta\cup-\colon K_2(F')
       \longrightarrow H^4(F',\Q_p/\Z_p(3))\bigr).
\]
\end{lemma}
\begin{proof}
The inclusion from left to right follows from the projection formula in Galois cohomology \cite[\S2 (C3)]{Otabe2023}, applied with finite coefficients and passed to the direct limit: for any finite splitting field $E/F'$ and any $w\in K_2(E)$,
\[
 \beta\cup N_{E/F'}w
 =\cor_{E/F'}(\beta_E\cup w)=0.
\]
For the reverse inclusion, take $z\in K_2(F')$ with $\beta\cup z=0$.
By Lemma~\ref{lem:second-three-coefficients}, $H^4(F',\Z/p(3))\simeq\Z/p$ is identified with the subgroup of $H^4(F',\Q_p/\Z_p(3))$ killed by $p$.
Let $\overline c=\beta\cup z\in H^4(F',\Z/p(3))$ be the product of $\beta\in\Br(F')[p]$ with the norm-residue symbol of $z$.
Its image in the direct-limit coefficient group is zero, so $\overline c=0$ by this identification.

We choose a complete Cohen discrete valuation ring $R$ of mixed characteristic $(0,p)$ with residue field $F'$, uniformizer $\pi$, and fraction field $\mathcal K$.
The Henselian Brauer-group isomorphism \cite[IV, Corollary 2.13]{MilneEtale}
\[
 \Br(R)\xrightarrow{\ \sim\ }\Br(F')
\]
and the fixed-degree lifting theorem for $\operatorname{PGL}_p$-torsors over Henselian local rings \cite[Expos\'e XXIV, Proposition 8.1]{SGA3} provide an Azumaya algebra $\widetilde B$ over $R$ of degree $p$ whose special fiber is isomorphic to the division algebra Brauer equivalent to $B$.
Its generic fiber $\widetilde B_{\mathcal K}$ also has index $p$: the index divides its degree, and the Brauer class cannot vanish over $\mathcal K$ because $\Br(R)\to\Br(\mathcal K)$ is injective \cite[IV, Corollary 2.6]{MilneEtale}.
By Matsumoto's theorem, we write $z$ as a finite sum of Milnor symbols; lifting their entries to units of $R$ gives
\[
 \widetilde z\in K_2^M(R),\qquad
 \operatorname{sp}(\widetilde z)=z.
\]

We use Kato's unramified lifting map with finite $p$-coefficients
\[
 \iota_{F'/\mathcal K}^q\colon
 H^q(F',\Z/p(q-1))
 \longrightarrow H^q(\mathcal K,\mu_p^{\otimes(q-1)})
\]
for $q=2,4$.
It is defined by lifting Artin--Schreier characters and the unit entries of symbols \cite[Theorem~3]{KatoGaloisCohomology1982}, \cite[\S5.1]{Kurihara2000}.
In degree two, we write $\beta$ as a sum of cyclic symbols and lift their Artin--Schreier characters to finite \'etale $R$-algebras and their units to $R^\times$.
The sum of the lifted cyclic classes extends over $R$ with special fiber $\beta$, so $\Br(R)\simeq\Br(F')$ gives $\iota_{F'/\mathcal K}^2(\beta)=[\widetilde B_{\mathcal K}]$.
In degree four, Kato's explicit product formula gives
\[
 \iota_{F'/\mathcal K}^4(\beta\cup z)
  =[\widetilde B_{\mathcal K}]\cup\widetilde z_{\mathcal K}
  \quad\text{in }H^4(\mathcal K,\mu_p^{\otimes3}).
\]
The left-hand side is zero because $\overline c=0$, so the finite-coefficient cup product over $\mathcal K$ vanishes.

Since $\mathcal K$ has characteristic zero and $\widetilde B_{\mathcal K}$ has index $p$, Lemma~\ref{lem:second-three-kahn-levine} gives finite separable splitting fields $L_j/\mathcal K$ of $\widetilde B_{\mathcal K}$ and elements $w_j\in K_2(L_j)$ such that
\[
 \widetilde z_{\mathcal K}=\sum_jN_{L_j/\mathcal K}(w_j)
    \quad\text{in }K_2(\mathcal K).
\]
The integral closure $R_j$ of $R$ in $L_j$ is a complete discrete valuation ring; let $F'_j$ be its residue field.
Since $\Br(R_j)\to\Br(L_j)$ is injective and $\widetilde B_{L_j}$ splits, each finite residue extension $F'_j/F'$ splits $B$.
The residue extensions need not be separable, which is harmless because the definition of $N_2(\beta)$ includes all finite splitting fields.

Write $\partial_{\mathcal K}\colon K_3(\mathcal K)\to K_2(F')$ and $\partial_{L_j}\colon K_3(L_j)\to K_2(F'_j)$ for the boundary maps in Quillen's localization sequences, and identify $K_2$ of the residue fields with Milnor $K_2$ by Matsumoto's theorem; under this identification the Quillen and Milnor norms agree.
Since $\widetilde z$ comes from $K_2(R)$, the $K(R)$-linearity of the boundary and $\partial_{\mathcal K}\{\pi\}=1$ give the first equality below.
The projection formula in Quillen $K$-theory, which expresses the compatibility of tensor products with restriction of scalars, and the compatibility of norms with boundaries then give
\[
 \begin{aligned}
 z
 &=\partial_{\mathcal K}\bigl(\{\pi\}\cdot\widetilde z_{\mathcal K}\bigr)\\
 &=\sum_j\partial_{\mathcal K}N_{L_j/\mathcal K}
             \bigl(\res_{L_j/\mathcal K}\{\pi\}\cdot w_j\bigr)\\
 &=\sum_jN_{F'_j/F'}\,\partial_{L_j}
             \bigl(\res_{L_j/\mathcal K}\{\pi\}\cdot w_j\bigr).
 \end{aligned}
\]
The compatibility of norms with boundaries in the last equality follows from Quillen's localization sequence and d\'evissage \cite[\S5, Theorems 4 and 5]{QuillenHigher}: restriction of scalars along the finite free extension $R_j/R$ preserves modules of finite length, and on closed support it is restriction of scalars from $F'_j$ to $F'$.
Hence the boundary commutes with norms even when $F'_j/F'$ is inseparable, and no ramification factor appears.
Since each $F'_j$ splits $B$, every term of the final sum lies in $N_2(\beta)$, and hence $z\in N_2(\beta)$.
\end{proof}

The remaining lemmas treat wild classes of period $p$ by index descent.
In both, $F$ is a three-dimensional local field with first residue characteristic $p$ and first residue field $k$, $\alpha\in\Br(F)[p]$ is wild, and $S\subset\Br(F)[p]\setminus\{0\}$ is a finite set containing $\alpha$ with at most one purely unramified class.
In mixed characteristic we first adjoin $\mu_p$ and, if necessary, pass to a totally ramified extension of degree prime to $p$; these restrictions are injective on $p$-torsion and preserve wildness, and they create no new purely unramified classes, since the character component of a Witt class becomes $e\cdot\res\chi$ with $e$ prime to $p$.
We choose the degree of the totally ramified extension to achieve the depth inequalities used below.
After these preparatory extensions we again write $F$, $\alpha$ and $S$ for the field and the classes, and only then choose $k$, the uniformizer $\pi$ and the residue $p$-basis.

\begin{lemma}\label{lem:second-three-lambda}
In this situation, there are a residue $p$-basis $(a,b)$ of $k$ and a unit $C$ of $F$ with residue $a$ such that, for $F_i=F(C^{1/p^i})$ and the completion $L^\wedge$ of $L=\bigcup_iF_i$, every member of $S$ is nonzero over $L^\wedge$ and $\alpha_{L^\wedge}$ is wild.
The residue field $k_\infty$ of $L^\wedge$ satisfies $[k_\infty:k_\infty^p]=p$, with $p$-basis $b$.
\end{lemma}
\begin{proof}
For $\gamma\in\Br(F)[p]$ we use the leading term of $\gamma$ for Kato's filtration in equal characteristic (indexed by the pole order $i$) and for the Bloch--Kato filtration in mixed characteristic (indexed by the depth $j$ of $U^j$).
In equal characteristic for pole order $i>0$, and in mixed characteristic for depth $0<j<c$, its graded pieces are $\Omega^1_k$ if $p$ does not divide the index and $\Omega^1_k/Z^1_k\oplus k/k^p$ if $p$ divides it; in mixed characteristic, the piece in degree zero is $\nu^2(k)\oplus\nu^1(k)$.
The terminal layer, pole order zero in equal characteristic and depth $c$ (defined below) in mixed characteristic, is the Witt part.
The components below are those of this leading term, and a class is wild exactly when its leading term lies outside the terminal layer.
For a $p$-basis $(a,b)$ of $k$, $\partial_a,\partial_b$ are the derivations dual to $da,db$.

We perfect one direction of a residue $p$-basis $(a,b)$, keeping every member of $S$ nonzero and keeping $\alpha$ wild.
The forbidden directions for $da$ are the nonzero one-form components of $S$, the differentials $dw$ of nonzero components $w\in k/k^p$, and the logarithmic-differential directions of nonzero degree-zero $\nu^1(k)$ components; there are finitely many of them.
A nonzero character component in the Witt decomposition remains nonzero under purely inseparable extension.
If $S$ contains no purely unramified class, we choose $a$ directly so that $da$ avoids the forbidden directions.
Otherwise, let $\gamma\in\Br(k)[p]$ be the purely unramified class in $S$.
Two-dimensional local duality gives $a_0$ with $\gamma\cup a_0\ne0$.
Then $a_0\notin k^p$, so $da_0$ extends to a basis of $\Omega^1_k$.
We choose $db_0$ independent of $da_0$, and vary $a=a_0+t b_0$ with sufficiently deep $t\in k^p$.
The cup-product character kills some principal-unit group, so $\gamma\cup a\ne0$ is preserved.
Meanwhile, $da=da_0+t\,db_0$ takes infinitely many directions and can avoid the finitely many forbidden ones.
For $k_i=k(a^{1/p^i})$, the projection formula and $N_{k_i/k}(a^{1/p^i})=a$ give $\gamma_{k_i}\ne0$ at every stage, and hence over the direct limit $k_\infty$.
We choose $b$ to complete $a$ to a $p$-basis.

We lift $a$ to a unit $C$, and put $F_i=F(C^{1/p^i})$.
The polynomial degree gives $[F_i:F]\le p^i=[k_i:k]$, and the residue field contains $k_i$, so $F_i/F$ has ramification index one and residue field $k_i$.
The residue $k_\infty$ of the completion $F_\infty^\wedge=L^\wedge$ of $\bigcup_i F_i$ has $p$-basis $b$.
In particular, $\ker(\Omega^1_k\to\Omega^1_{k_\infty})=k\,da$.
We combine this with the freedom to vary $C$ modulo $\pi^2$.

Suppose first that $F$ has equal characteristic, and write $F=k((\pi))$.
Kato's leading term at pole order $i>0$ is measured by $\Omega^1_k$ if $p\nmid i$, and by $\Omega^1_k/Z^1_k\oplus k/k^p$ if $p\mid i$, where $Z^1_k=\ker d$.
In the first case, a nonzero form remains nonzero if the direction of $da$ avoids that form; in the second case, a nonzero second component $w\bmod k^p$ remains nonzero if the direction of $da$ avoids $dw$.
If the second component vanishes, take a representative $\omega=f\,da+g\,db$ of a nonzero class in $\Omega^1_k/Z^1_k$; then $d\omega\ne0$.
We set $C=a+\lambda\pi$ for $\lambda\in k^p$ and choose a new coefficient field of $F_\infty^\wedge$ for its residue field $k_\infty$ containing the prescribed lift of $b$; let $a^\sharp$ be the lift of $a$ in it.
Then $C$ and $a^\sharp$ have the same residue and both have all $p$-power roots in $F_\infty^\wedge$, so if they differed, the positive valuation of $C/a^\sharp-1$ would be divisible by every $p^n$; hence $C=a^\sharp$.
To first order, the difference between the old and new lifts of any $u\in k$ is $-\lambda\pi\partial_a u$, by the binomial expansion of $u=\sum u_{ij}^p a^i b^j$.
Since $dC=0$ in $F_\infty^\wedge$, we have
\[
 \omega\equiv g\,db-\lambda\pi
       \bigl((\partial_a g)\,db+f\,d\log\pi\bigr)
       \pmod{\pi^2}.
\]
In this case, where $p\mid i$, the second component vanishes and the representative satisfies $d\omega\ne0$, the closed form lowered by Cartier to pole order $i/p$ is the restricted form $g\,db$; it is closed over the coefficient field, which has a one-element $p$-basis.
At pole order $i-1$, eliminating the $d\log\pi$ term with the exact differential $d(\pi^{-(i-1)}\lambda f)$ makes the coefficient of $\lambda$ equal to $(\partial_b f-\partial_a g)\,db$, which is nonzero since $d\omega\ne0$.
The coefficient of the image of the restricted closed form $g\,db$ under the Cartier operator is given by its residue in the new coefficient field, which does not involve $\lambda$.
This also applies when $i/p=i-1$, which happens only for $i=p=2$.

Suppose next that $F$ has mixed characteristic, so that the preliminary extensions above have been made.
Put $e=v_F(p)$ and $c=pe/(p-1)$.
The terminal term $U^c$ of the Bloch--Kato unit filtration modulo $p$ corresponds to the Witt part: with $\tau=\zeta_p-1$, the residue of $\bigl((1+\tau y)^p-1\bigr)/\tau^p$ is $y^p-y$.
Thus the Kummer character of $1+(\zeta_p-1)^pu$ is the unramified lift of the Artin--Schreier character of $\bar u$, and taking a unit and a uniformizer in the second slot generates the two Witt components \cite[(5.1)--(5.3)]{BlochKato1986}.
For $0<j<c$, the graded pieces are $\Omega^1_k$ if $p\nmid j$ and $\Omega^1_k/Z^1_k\oplus k/k^p$ if $p\mid j$.
The degree-zero piece is $\nu^2(k)\oplus\nu^1(k)$, where $\nu^r(k)=K_r^M(k)/p$ is identified by logarithmic differentials \cite[(5.1)--(5.3)]{BlochKato1986}.
Scaling $e,c,j$ by the degree of the preparatory totally ramified extension, we have arranged $j+1<c$ for each relevant positive depth, and $e+j/p-j\ge2$ when $p\mid j$; we have also arranged $e=v_F(p)\ge2$, which is used below even when $S$ has no class of positive depth.
For the change $\pi'^d=\pi$ with $p\nmid d$, we have $d\log\pi=d\cdot d\log\pi'$, so the nonzero graded pieces are preserved while these depth inequalities are achieved.
The choice of $da$ avoiding the forbidden directions preserves nonzero one-forms (for $p\nmid j$), components in $k/k^p$, and nonzero degree-zero $\nu^1(k)$ components.

For the remaining cases, take integral lifts $\widetilde a,\widetilde b,\widetilde\lambda$ of $a,b,\lambda$ and put $C=\widetilde a+\widetilde\lambda\pi$.
Fix a lift $\widetilde u$ of $u=\sum x_{rs}^p a^r b^s$, and set $\widetilde u_0=\sum \widetilde x_{rs}^pC^r\widetilde b^s$.
Then
\[
 \overline{(\widetilde u-\widetilde u_0)/\pi}=z_{\widetilde u}-\lambda\,\partial_a u,
\]
where $z_{\widetilde u}$ is the value of the left-hand side at $C=\widetilde a$.
For positive depth $j$ divisible by $p$ with zero second component, lift $\omega=f\,da+g\,db$ by $\{1+\pi^j\widetilde f \widetilde a,\widetilde a\}+ \{1+\pi^j\widetilde g \widetilde b,\widetilde b\}$.
In the first term, $C$ in $\widetilde a=C(1-\widetilde\lambda\pi/C)$ becomes a $p$-th power.
We apply the first-order Steinberg relation $\{1+x,1+y\}\equiv\{1+xy,-x^{-1}\}\bmod U^{j+2}$ \cite[(4.1.1)]{BlochKato1986} to the first symbol with $x=\pi^j\widetilde f\widetilde a$ and $y=-\widetilde\lambda\pi/C$.
Then $\overline{xy/\pi^{j+1}}=-\lambda f$, and since $p\mid j$ and $da=0$ after restriction, the term of degree $j+1$ in $-\lambda f\,d\log(-x^{-1})$ is $\lambda\,df$.
The $\widetilde u_0$-part of the coefficient of the second term can be written modulo $\pi^2$ as $\sum_{r<p}h_r^p \widetilde b^r$.
Since $e+j/p-j\ge2$, we have $1+\pi^jh_0^p\equiv(1+\pi^{j/p}h_0)^p\pmod{U^{j+2}}$, so the $r=0$ part is a $p$-th power modulo $U^{j+2}$.
For $r>0$, multiplying the symbol by $r$ and using Steinberg moves it to a symbol with second slot $-(\pi^j h_r^p)^{-1}$, and the same vanishing follows.
The new second slot of the second symbol is a $p$-th power up to sign, and the sign for $p=2$ vanishes in this degree since $-1\in U_F^e$ and $e\ge2$.
Thus $-\lambda(\partial_a g)\,db$ is the $\lambda$-dependent part of the second term; a term independent of $\lambda$ may also occur, which is why at most one value of $\lambda$ per class is excluded below.
Hence the coefficient of $\lambda$ is the nonzero form $(\partial_b f-\partial_a g)\,db$.

At degree zero, if the $\nu^1(k)$ component vanishes, lift the nonzero $\omega\in\nu^2(k)$ by a sum of unit symbols.
The ring $\mathcal O_{F_\infty^\wedge}/(\pi^2)$ has characteristic $p$.
We consider the section sending an element of $k_\infty^p$ to the $p$-th power of a lift of its root; modulo $\pi^2$ in characteristic $p$, this $p$-th power does not depend on the lift and is compatible with $p$-th powers of $b$, so the section is a ring homomorphism.
Using the power basis in $b$, we extend this section to a section of $k_\infty$.
Symbols along this section vanish modulo $(p,U^2)$, since $K_2^M(k_\infty)/p=0$.
Taking the difference in the preceding first-order expansion, we find that the coefficient of $\lambda$ for an original symbol $\{\widetilde u,\widetilde v\}$ is
\[
 -(\partial_a u/u)\,d\log v+(\partial_a v/v)\,d\log u.
\]
For the sum it is therefore $-\iota_{\partial_a}\omega$; a nonzero $da\wedge db$ component maps to a nonzero $db$ component.

In every case, contributions already lying in the next layer come from residue restriction and are independent of $\lambda$.
Avoiding at most one $\lambda$ for each class keeps all members of $S$ nonzero and $\alpha$ wild.
\end{proof}

The following density statement is used in the proof of the next lemma.

\begin{lemma}\label{lem:paper-henselian-density}
Let $H$ be a Henselian field for a rank-$1$ valuation, and let $\widehat H$ be its completion.
On any smooth finite-type $H$-scheme, the $H$-points are dense in the $\widehat H$-points.
\end{lemma}
\begin{proof}
Locally, a smooth scheme is standard \'etale over $\mathbb A^n_H$: it is given by $g(a,y)=0$ with $\partial_yg(a,y)\ne0$, where $g\in H[a_1,\dots,a_n][y]$ is monic in $y$.
Let $(a,b)$ be an $\widehat H$-point.
After translating and rescaling, we may assume that the coefficients of $g(a,-)$ and the point lie in the valuation ring.
Since $H$ is dense in $\widehat H$, we may choose $a'\in H^n$ and $b'\in H$ so close to $a$ and $b$ that $|g(a',b')|<|\partial_yg(a',b')|^2$.
Newton's form of Hensel's lemma in the Henselian rank-one field $H$ then gives a root $b''\in H$ of $g(a',-)$ with $|b''-b'|\le|g(a',b')|/|\partial_yg(a',b')|$, which can be made arbitrarily small, and $(a',b'')$ is an $H$-point close to $(a,b)$.
Shrinking the neighborhood keeps the open condition $\partial_yg\ne0$.
\end{proof}

\begin{lemma}\label{lem:second-three-wild-descent}
In this situation, there is a finite extension $F'/F$, namely a stage $F_i$ of the tower in Lemma~\ref{lem:second-three-lambda}, such that every member of $S$ is nonzero over $F'$, $\alpha_{F'}$ is wild, and $\ind(\alpha_{F'})=p$.
\end{lemma}
\begin{proof}
Let $(a,b)$, $C$, $F_i$ and $L$ be as in Lemma~\ref{lem:second-three-lambda}.
The residue field $k_\infty$ of $L^\wedge$ satisfies $[k_\infty:k_\infty^p]=p$, so the wild class $\alpha_{L^\wedge}$ of period $p$ has index $p$ by \cite[\S4, Lemma 5, p.~337]{KatoLocalI}, whose hypotheses (a $p$-torsion class not split by the maximal unramified extension, and residue $p$-rank one) hold here.
Let $A$ be a central simple algebra over $F$ with class $\alpha$, and let $X=\SB_p(A_L)$ be the generalized Severi--Brauer variety parametrizing right ideals of reduced dimension $p$.
Since $X$ is smooth and $L$ is Henselian and dense in $L^\wedge$, Lemma~\ref{lem:paper-henselian-density} and $X(L^\wedge)\ne\varnothing$ give $X(L)\ne\varnothing$.
Finite presentation descends this point to a finite stage $F_i$, where $\ind(\alpha_{F_i})\mid p$.
All members of $S$ are nonzero over $L^\wedge$, hence over $F_i$; in particular, $\ind(\alpha_{F_i})=p$.
If $\alpha_{F_i}$ belonged to the Witt part, so would its restriction $\alpha_{L^\wedge}$, because a finite unramified extension splitting a Witt class lifts along the base change of ramification index one; this contradicts Lemma~\ref{lem:second-three-lambda}.
Hence $\alpha_{F_i}$ is wild.
\end{proof}

\begin{proof}[Proof of Theorem~\ref{thm:second-three-coker}]
By Lemma~\ref{lem:second-three-duality}, $(K_2(F)/N_2(\alpha))^\vee=B_\alpha$.
The splitting transfers give $N_2(\alpha)\subset\Nrd_2K_2(D)$, and Lemma~\ref{lem:second-three-containment} gives $\Nrd_2K_2(D)\subset\ker(\alpha\cup-)$.
By definition, $B_\alpha$ contains $\langle\alpha\rangle$, and it remains to prove equality.
Since $B_\alpha=\bigoplus_\ell B_{\alpha_\ell}$ and $\langle\alpha\rangle=\bigoplus_\ell\langle\alpha_\ell\rangle$, it suffices to do so for $\ell$-primary $\alpha$.
A class in the Witt part is split by a finite unramified extension, so Lemma~\ref{lem:second-three-witt-common-kernel} applies with $h=2$, its hypothesis being the two-dimensional case recalled at the beginning of this subsection; this covers all $\ell$-primary classes with $\ell$ different from the first residue characteristic, since they lie in the Witt part.

Suppose that $\alpha$ is wild of period $p$ and that $\beta\in B_\alpha[p]\setminus\langle\alpha\rangle$.
Put $S=\{\alpha\}\cup\{\beta-j\alpha:j\in\F_p\}$.
Every member of $S$ is nonzero, and $S$ contains at most one purely unramified class: the difference of two would be a nonzero multiple of the wild class $\alpha$.
Lemma~\ref{lem:second-three-wild-descent} gives a finite extension $F'=F_i$ of the original field $F$, through the preparatory extensions, over which every member of $S$ is nonzero, $\alpha_{F_i}$ is wild, and $\ind(\alpha_{F_i})=p$.
We claim that $N_2(\alpha_{F_i})=\ker(\alpha_{F_i}\cup-)$.
In characteristic zero, apply Lemma~\ref{lem:second-three-kahn-levine} with $k=F_i$ and $e=p$, and use Lemma~\ref{lem:second-three-coefficients} to pass from $\Z/p(3)$ to $\Q/\Z(3)$ coefficients.
In characteristic $p$, apply Lemma~\ref{lem:second-three-char-p-index-p} with $F'=F_i$ and $\beta=\alpha_{F_i}$.
By the equivalence in Lemma~\ref{lem:second-three-duality} over $F_i$, the claim gives $B_{\alpha_{F_i}}=\langle\alpha_{F_i}\rangle$.
Since $\beta_{F_i}\in B_{\alpha_{F_i}}$, there is $j\in\F_p$ such that the element $\beta-j\alpha$ of $S$ vanishes over $F_i$, contradicting the nonvanishing of the members of $S$.

Next put $\per(\alpha)=p^r$ and $\theta=p^{r-1}\alpha$.
For any finite extension $E/F$ splitting $\theta$, let $E_s$ be its maximal separable subextension.
In the Galois closure $E^{\mathrm{gal}}_s$ of $E_s/F$, choose a Sylow $p$-subgroup $P$ of $\Gal(E^{\mathrm{gal}}_s/F)$ containing a Sylow $p$-subgroup of $\Gal(E^{\mathrm{gal}}_s/E_s)$, and put $F''=(E^{\mathrm{gal}}_s)^P$.
Both $[F'':F]$ and $[EF'':E]$ are prime to $p$; for the latter, $[E_sF'':E_s]=[H:H\cap P]$ with $H=\Gal(E^{\mathrm{gal}}_s/E_s)$, and $E/E_s$ is purely inseparable.
A subgroup chain of the $p$-group gives a tower of degree-$p$ extensions for $E_sF''/F''$; the purely inseparable extension $EF''/E_sF''$ also has a degree-$p$ tower.
Let $M'/M$ be the first step where $\theta$ vanishes.
We have $\per(\alpha_M)=p^r$ and $\per(\alpha_{M'})=p^{r-1}$.
Restriction in degree $p$ lowers order by at most a factor of $p$: for a separable step this follows from corestriction, and for a purely inseparable step $M'/M$ of degree $p$, the composite $M\hookrightarrow M'\xrightarrow{x\mapsto x^p}M$ is the Frobenius, which acts on $\Br(M)\{p\}$ as multiplication by $p$.
For $r\ge2$, we assume by induction that $B_{\alpha'}=\langle\alpha'\rangle$ for every wild class $\alpha'$ with $p^{r-1}\alpha'=0$ over a three-dimensional local field.
By Lemma~\ref{lem:second-three-witt-common-kernel}, the same equality holds for every $p$-primary class $\alpha'$ in the Witt part over such a field.
Still for $r\ge2$, we take $\beta\in B_\alpha[p]$.
Every finite extension of $M'$ splitting $\alpha_{M'}$ also splits $\alpha$ over $F$; hence restriction sends $B_\alpha$ into $B_{\alpha_{M'}}$, and $\beta_{M'}\in B_{\alpha_{M'}}[p]$.
Induction gives $\beta_{M'}=j p^{r-2}\alpha_{M'}$ for some integer $j$.
If $j\not\equiv0\pmod p$, then $\beta_M-jp^{r-2}\alpha_M$ has order $p^2$ and cannot vanish in a degree-$p$ extension.
Thus $\beta_{M'}=0$, hence $\beta_{EF''}=0$.
As $[EF'':E]$ is prime to $p$, also $\beta_E=0$, so $\beta\in B_\theta[p]$.
Moreover, $B_\theta[p]=\langle\theta\rangle$, by Lemma~\ref{lem:second-three-witt-common-kernel} if $\theta$ lies in the Witt part, and by the period-$p$ case if $\theta$ is wild.
For $r=1$, we have already proved $B_\alpha[p]=\langle\alpha\rangle$.
In both cases, $B_\alpha$ is finite and $B_\alpha[p]$ has $\F_p$-dimension one, so $B_\alpha$ is cyclic.
If a generator $\delta$ had order $p^{r+q}$, $q>0$, then $\alpha$ would be a unit multiple of $p^q\delta$.
At the first vanishing step $M'/M$ for $\theta$, $\alpha_M$ has order $p^r$, so $\delta_M$ has order $p^{r+q}$.
Consequently, $\delta_{M'}$ has order at least $p^{r+q-1}$.
By the same restriction property as above, every finite extension of $M'$ splitting $\alpha_{M'}$ splits $\alpha$, so $\delta_{M'}\in B_{\alpha_{M'}}$.
But for $r\ge2$, by induction, $B_{\alpha_{M'}}=\langle\alpha_{M'}\rangle$ is killed by $p^{r-1}$.
For $r=1$, one has $\alpha_{M'}=0$ and hence $B_{\alpha_{M'}}=0$.
Both cases are contradictions.
Therefore $B_\alpha=\langle\alpha\rangle$, so Lemma~\ref{lem:second-three-duality} gives $N_2(\alpha)=\ker(\alpha\cup-)$.
With the inclusions $N_2(\alpha)\subseteq\Nrd_2K_2(D)\subseteq\ker(\alpha\cup-)$ shown above, this gives $N_2(\alpha)=\Nrd_2K_2(D)=\ker(\alpha\cup-)$, and the same lemma identifies $C_2(D/F)=K_2(F)/N_2(\alpha)$ with $\tfrac1{\per(\alpha)}\Z/\Z$ via $\inv_F(\alpha\cup-)$, which proves the theorem.
\end{proof}

\subsection{Vanishing of the kernel in first residue characteristic zero}\label{sec:second-three-charzero}
When the first residue field has characteristic zero, the results of Sections~\ref{sec:second-two-tame} and~\ref{sec:second-two-wild} also give the vanishing of the kernel of the reduced norm.

\begin{corollary}\label{cor:second-three-charzero}
For every three-dimensional local field $F$ whose first residue field has characteristic zero, and every central division algebra $D$ over $F$, one has $\SK_2(D)=0$.
This includes both characteristic sequences $(0,0,0)$ and $(0,0,p)$.
\end{corollary}
\begin{proof}
By Cohen's structure theorem, write $F=k((t))$, where $k$ is the first residue field of $F$.
Let $B$ be the first residue division algebra of $D$, and put $L=Z(B)$.
The center $L$ is a two-dimensional local field, so Corollary~\ref{cor:review-two-local-sk1-zero} gives $\SK_1(B)=0$, which is the hypothesis of Proposition~\ref{prop:second-residue-zero}.
By that proposition, $L/k$ is cyclic and
\[
 \SK_2(D)\simeq\ker\bigl[K_2(B)_{\Gal(L/k)}
  \xrightarrow{\ N_{L/k}\Nrd_{B/L,2}\ }K_2(k)\bigr].
\]
The equivariant map $\Nrd_{B/L,2}\colon K_2(B)\to K_2(L)$ is an isomorphism by Theorem~\ref{thm:second-two-vanishing}, and cyclic Hilbert 90 for $K_2$ \cite[Theorem 14.1]{MerkurjevSuslin1982} makes $N_{L/k}$ injective on $K_2(L)_{\Gal(L/k)}$.
Thus the kernel is zero.
\end{proof}

The following observation on the unramified case $D=B((t))$, used in Remark~\ref{rem:paper-sk2-dimension-four}, holds in every characteristic.

\begin{proposition}\label{prop:paper-laurent-sk2}
Let $B$ be a central division algebra over a field $k$, and put $F=k((t))$ and $D=B((t))$.
Then the map $\SK_1(B)\to\SK_2(D)$, $x\mapsto x\cdot\{t\}$, is well defined and injective.
\end{proposition}
\begin{proof}
Multiplication by $\{t\}$ commutes with reduced norms in every characteristic: $K_1(B)$ is generated by splitting transfers, and for $x=\theta_E(e)$ the projection formula and $\Nrd_2\theta_{E((t))}=N_{E((t))/F}$ give $\Nrd_{D,2}(x\cdot\{t\})=N_{E/k}(e)\cdot\{t\}=\Nrd_{B/k,1}(x)\cdot\{t\}$.
Hence multiplication by $\{t\}$ maps $\SK_1(B)$ into $\SK_2(D)$, and it is injective since the localization boundary $K_2(D)\to K_1(B)$ sends $x\cdot\{t\}$ to $\pm x$.
The same computation applies with $E(t)$ in place of $E((t))$.
\end{proof}

\subsection{Special cases, further results, and open problems}\label{sec:paper-further-results}
This subsection discusses the vanishing of $\SK_2$ beyond Corollary~\ref{cor:second-three-charzero}.
We state conjectures for three-dimensional local fields and for two- and three-dimensional global fields, record a special case, and collect further observations on their scope.
\begin{conjecture}\label{conj:paper-three-local-sk2}
For every three-dimensional local field $F$ and every central division algebra $D$ over $F$, one has $\SK_2(D)=0$.
\end{conjecture}

The conjecture is suggested by the vanishing of $\SK_2(D)$ over local fields of dimension at most two: in dimension zero $D$ is commutative by Wedderburn's theorem, and in dimensions one and two this is Theorems~\ref{thm:paper-local-norm} and~\ref{thm:second-two-vanishing}.
In dimension at least four, the analogous statement fails (Remark~\ref{rem:paper-sk2-dimension-four} below).
In dimension three, Corollary~\ref{cor:second-three-charzero} proves the conjecture when the first residue field has characteristic zero.
The next remark records a mixed-characteristic case with positive first residue characteristic.
The general case is open.

We record several further results, some of which are stated without proof.
We include these statements to make their precise scope available for future work.

\begin{remark}
Let $F$ be a three-dimensional local field with first residue field $k$ of characteristic $p>0$ and with characteristic sequence $(0,p,p)$.
For a uniformizer $\pi$ and a unit $u$ with residue $\bar u\notin k^p$, one has $\SK_2(D)=0$ for every central division algebra $D$ over $F$ of period $p^2$ split by $F(\pi^{1/p},u^{1/p})$.
We omit its proof.
This result is not used elsewhere in the paper.
\end{remark}

\begin{remark}\label{rem:paper-sk2-dimension-four}
The analogous vanishing statement fails in every dimension at least four.
Let $\ell$ be a prime and let $k_0$ be a one-dimensional local field with two linearly disjoint cyclic extensions of degree $\ell$; such fields exist in characteristic zero and in positive characteristic.
Platonov constructed a central division algebra $B$ of index $\ell^2$ over $k_0((x))((y))$ with $\SK_1(B)\simeq\Z/\ell$ \cite[Theorems 4.7 and 5.9]{Platonov1976}.
Over the four-dimensional local field $k_0((x))((y))((t))$, Proposition~\ref{prop:paper-laurent-sk2} shows that $D=B((t))$ has $\SK_2(D)\supseteq\SK_1(B)\cdot\{t\}\ne0$.
For dimension $d>4$, adjoin $d-4$ further central Laurent variables $s$ one at a time: restriction maps $\SK_2$ into $\SK_2$, and the specialization $\partial(\{s\}\cdot-)$ sends the restriction of a class back to that class up to sign, so a nonzero class stays nonzero.
These examples delimit the scope of Conjecture~\ref{conj:paper-three-local-sk2}.
\end{remark}

\begin{remark}\label{rem:paper-wild-cokernel}
The reduced-norm cokernel $C_2(D/F)=\coker(\Nrd_2)$ is cyclic of order $\per[D]$ in dimension three by Theorem~\ref{thm:second-three-coker}.
It is killed by $\deg D$, so it is $p$-primary when $\deg D$ is a power of the characteristic $p$ of the last residue field.
Over iterated Laurent series fields of dimension at least four, this $p$-primary cokernel can be finite or uncountable; we give one example of each, without proof.
\begin{enumerate}
\item Let $F=\Q_2((t_2))((t_3))((t_4))$, a four-dimensional local field with last residue field $\F_2$.
Let $\chi$ be the character of the base-changed extension $F(i)/F$ of the ramified quadratic extension $\Q_2(i)/\Q_2$, and let $D$ be the quaternion division algebra over $F$ with class $(\chi,t_2)$.
The Brauer residue of $[D]$ at $t_2$, computed over $\Q_2((t_2))$ before the outer Laurent variables $t_3,t_4$ are added, is a nonzero character, so $D$ has index two over $\Q_2((t_2))$, and adding Laurent variables preserves this.
Then $C_2(D/F)\simeq(\Z/2)^3$.
\item Let $k=\F_q((u))((v))$, with $q$ a power of a prime $p$, and let $L=k(y)$ with $y^p-v^{p-1}y=u$.
The residue polynomial $X^p-u$ is irreducible, so $L/k$ has ramification index one and purely inseparable residue extension of degree $p$, and the roots $y+cv$ with $c\in\F_p$ make $L/k$ cyclic of degree $p$.
Let $D$ be the base change of the cyclic algebra $(L((t))/k((t)),\sigma,t)$ to the four-dimensional local field $F=k((t))((s))$.
Over $k((t))$, the Brauer residue at $t$ is a nonzero cyclic character, so the cyclic algebra has index $p$.
Adjoining the Laurent variable $s$ preserves the index, so $D$ has index $p$.
Then $C_2(D/F)$ contains $k^\times/N_{L/k}L^\times$ as a direct summand, and this group has the cardinality of $\F_q((u))$; in particular, $C_2(D/F)$ is uncountable.
\end{enumerate}
\end{remark}

\begin{conjecture}\label{conj:paper-global-sk2}
Let $F$ be a two- or three-dimensional global field, that is, a finitely generated field of transcendence degree $h\in\{2,3\}$ over a finite field or of transcendence degree $h-1$ over a number field.
Then $\SK_2(D)=0$ for every central division algebra $D$ over $F$.
\end{conjecture}

\begin{remark}
An $h$-dimensional global field is a finitely generated field of transcendence degree $h$ over a finite field or of transcendence degree $h-1$ over a number field.
We have also proved the following implication for a two-dimensional global field $F$ and a central division algebra $D$ over $F$; we omit the proof:
\[
v_2(\deg D)\le3,\qquad v_3(\deg D)\le2,\qquad
v_\ell(\deg D)\le1\;(\ell\ge5)
\quad\Longrightarrow\quad\SK_2(D)=0.
\]
The cases of indices $8$ and $9$, including characteristics $2$ and $3$, involve Adams operations and Bott lifting.
There are central division algebras with nonzero $\SK_2$ over four-dimensional global fields in both zero and positive characteristic.
Take Platonov's algebra $B$ with $\SK_1(B)\ne0$ over a purely transcendental extension $k_0(x,y)$ of a global field $k_0$ \cite[Theorem 6.2]{Platonov1976}, and choose $0\ne\xi\in\SK_1(B)$.
The projection formula, as in the proof of Proposition~\ref{prop:paper-laurent-sk2}, gives $\xi\cdot\{t\}\in\SK_2(B(t))$, and its image in $\SK_2(B((t)))$ is nonzero by that proposition.
These results are consistent with Conjecture~\ref{conj:paper-global-sk2}, which predicts the vanishing of $\SK_2(D)$ at arbitrary index over two- and three-dimensional global fields.
\end{remark}

\section{Locally compact modules and the Artin map}\label{sec:lca}
In this section, we compare the $K$-theory of locally compact $D$-modules with that of locally compact $F$-modules and relate this comparison to the Artin map.
The motivation comes from Clausen's description of the Artin map of class field theory through the $K$-theory of locally compact $F$-modules \cite{ClausenArtin}, and from its noncommutative analogue by Braunling--Henrard--van Roosmalen \cite{BHRAdeles}, which describes the id\`ele class group of $D$ modulo the id\`eles of local reduced norm one through $K_1$ of locally compact $D$-modules.
Braunling, Henrard and van Roosmalen conjecture that $K_2(\operatorname{LCA}_A)\cong K_2(\operatorname{LCA}_F)$ for every finite-dimensional simple $\Q$-algebra $A$ with center $F$ \cite[Conjecture 10.2]{BHRAdeles}.
For $A=D$, Theorem~\ref{thm:paper-lca-k2} below gives this isomorphism with integral coefficients, also when $D$ is ramified at a real place.

Our comparison shows that, in positive degrees, the fiber of the comparison $c^{\mathrm{lc}}_D$ is exactly the suspension $\Sigma Z$ of the global comparison fiber of Section~\ref{sec:global-spectrum}, a finite sum of shifted Eilenberg--Mac Lane spectra over number fields, with $2$ inverted under real ramification, and that this fiber is invisible to the Artin map.
On $\pi_1$ the compatibility gives the formula for the pro-$p$ Artin map on id\`ele classes: local reduced norms, followed by the global Artin map and the projection to the maximal pro-$p$ quotient, with $p$ odd under real ramification.

The $K$-theory of locally compact modules fits into a localization sequence between $K(D)$ and the restricted product of the local $K$-theories, and a comparison square built from a common integral model of $D$ induces a comparison $c^{\mathrm{lc}}_D$ between the cofibers.
If the global comparison of Section~\ref{sec:global-spectrum} is kept fixed, this square may have to be corrected.

Throughout this section, $F$ is a number field and, by Morita equivalence, $D$ is a central division algebra over $F$.
We work with integral coefficients when $D$ splits at every real place.
If $D$ is ramified at a real place, all spectra and groups in this section, except in Theorem~\ref{thm:paper-lca-k2}, are taken with $\Z[1/2]$ coefficients; the restricted products below are then formed integrally and $2$ is inverted afterwards, not place by place.
We construct $c^{\mathrm{lc}}_D$ (Theorem~\ref{thm:paper-lca-comparison}), identify its fiber with the suspension of the global comparison fiber (Theorem~\ref{thm:paper-lca-fiber-splitting}), prove the integral isomorphism on $K_2$ (Theorem~\ref{thm:paper-lca-k2}), and show that it is compatible with Clausen's Artin map (Theorem~\ref{thm:paper-lca-artin}).

\subsection{Locally compact modules and localization in positive degrees}
\label{subsec:paper-lca-localization}
We first set up the localization fiber sequences for locally compact modules over $D$ and over $F$ in positive degrees.
The input is the theorem of Braunling--Henrard--van Roosmalen, which we also need for second-countable modules.

For $R=D,F$, let $\operatorname{LCA}_R$ denote the category of locally compact right modules over the discrete ring $R$, with continuous $R$-linear maps.
We give it the standard exact structure in which closed embeddings are admissible monomorphisms and open surjections are admissible epimorphisms.
Write $\operatorname{AB}_R=\operatorname{LCA}_{R,\mathrm{ab}}$ for its full exact subcategory of adelic blocks.
Concretely, $\operatorname{AB}_R$ is the filtered $2$-colimit of the following categories.
For a maximal order $\mathfrak R$ in $R$ and a finite set $S$ of finite places, put
\[
 \mathcal J_R^S=
 \prod_{v\mid\infty}\operatorname{proj}(R_v)
 \times\prod_{v\in S}\operatorname{proj}(R_v)
 \times\prod_{v\notin S,\ v\nmid\infty}
                       \operatorname{proj}(\mathfrak R_v).
\]
For $S\subset S'$, the transition functor is extension of scalars at the newly added places.
This description is independent of the maximal order.
We use the theorem of Braunling--Henrard--van Roosmalen in the following form \cite[Definition 4.3, Proposition 4.5, Theorems 6.1, 6.2 and 6.13]{BHRAdeles}.

\begin{theorem}\label{thm:paper-lca-source}
There is a natural fiber sequence
\begin{equation}\label{eq:paper-lca-bhr}
 K(R)\xrightarrow{\lambda_R}K(\operatorname{AB}_R)
       \longrightarrow K(\operatorname{LCA}_R).
\end{equation}
The first map is induced by $P\mapsto P\otimes_R\mathbb A_R$.
For the middle term, there is an equivalence of spectra
\begin{equation}\label{eq:paper-lca-restricted-spectrum}
 K(\operatorname{AB}_R)\simeq
 \varinjlim_S\left(
 \prod_{v\mid\infty}K(R_v)\times\prod_{v\in S}K(R_v)
 \times\prod_{v\notin S,\ v\nmid\infty}K(\mathfrak R_v)\right).
\end{equation}
In particular, $K_j(\operatorname{AB}_R)=\prod'_vK_j(R_v)$, where the restriction is that the $v$-component lies in the image of $K_j(\mathfrak R_v)$ at all but finitely many finite places $v$.
\end{theorem}

The restricted product in~\eqref{eq:paper-lca-restricted-spectrum} is formed in spectra, using the product theorem for nonconnective $K$-theory of exact categories and the compatibility of $K$-theory with filtered colimits.

\begin{lemma}\label{lem:paper-lca-countable-localization}
Let $R\in\{D,F\}$, let $\mathcal C=\operatorname{LCA}_{R,\aleph}$ be the full exact subcategory of second-countable objects of $\operatorname{LCA}_R$, and let $\mathcal C_{\mathrm c}$ be its compact subcategory.
\begin{enumerate}
\item Every adelic block lies in $\mathcal C$.
\item The localization sequences of \cite[Proposition 6.3 and Corollary 6.6]{BHRAdeles}, for the quotient by $\mathcal C_{\mathrm c}$ and for the image $\mathcal V\simeq\operatorname{AB}_R$ of the adelic blocks in the exact hull $(\mathcal C/\mathcal C_{\mathrm c})^{\mathrm{ex}}$, hold with $\operatorname{LCA}_R$ replaced by $\mathcal C$.
\item The discrete quotient functor $M\mapsto M/Q(M)$ of \cite[Proposition 2.11]{BHRAdeles} induces an equivalence
\[
 (\mathcal C/\mathcal C_{\mathrm c})^{\mathrm{ex}}/\mathcal V
 \simeq\operatorname{Mod}_\aleph(R)/\operatorname{mod}(R),
\]
whose inverse is induced by the inclusion of countable discrete modules; here $\operatorname{Mod}_\aleph(R)$ is the category of countable discrete modules and $\operatorname{mod}(R)$ is the category of finitely generated modules.
\end{enumerate}
\end{lemma}
\begin{proof}
Second countability of locally compact groups is preserved under passage to closed subgroups, Hausdorff quotients, and finite iterated extensions.
A discrete quotient of a second-countable group is countable.
Since $R$ itself is countable, every finitely generated $R$-module is countable.
Every adelic block is second countable, since it is a restricted product of finite-dimensional local components over countably many places, so restricting to $\mathcal C$ does not change the category of adelic blocks.

We now explain how the localization proof of Braunling--Henrard--van Roosmalen restricts to this subcategory.
Let $M_0$ be an object of $\operatorname{LCA}_R$, and choose an order $\Lambda\subset R$.
As in the proof of \cite[Proposition~2.11]{BHRAdeles}, choose a conflation in $\operatorname{LCA}_\Lambda$
\[
 V\oplus C\hookrightarrow M_0\twoheadrightarrow M_0/(V\oplus C),
\]
where $V$ is a vector $R$-module, $C$ is a compact $\Lambda$-module that is open in a closed complement of $V$, and $M_0/(V\oplus C)$ is discrete, so that $V\oplus C$ is open in $M_0$; the summand $C$ need not be an $R$-module, but its saturation below is.
The split torsion pair of \cite[Proposition 2.11]{BHRAdeles} is constructed from the open $R$-submodule $Q(M_0)=\bigcup_{n\ge1}n^{-1}(V\oplus C)$ and the discrete quotient $M_0/Q(M_0)$.
If $M_0\in\mathcal C$, then $Q(M_0)$ is second countable and $M_0/Q(M_0)$ is countable.
The decomposition of \cite[Lemma 2.12]{BHRAdeles} splits the quasi-adelic part into a compact part and an adelic part by using only closed submodules and quotients, so both parts stay in $\mathcal C$.
Compact $R$-modules and adelic blocks are injective in the whole category and remain injective after restriction.

Consequently, we can apply to this subcategory the localization by the compact subcategory in \cite[Proposition 6.3]{BHRAdeles} and the full faithfulness of adelic blocks in \cite[Lemma 6.4]{BHRAdeles}.
Indeed, the quotient $\mathcal C/\mathcal C_{\mathrm c}$ and its exact hull are built intrinsically.
The category $\mathcal C$ has closed kernels, Hausdorff cokernels and finite pushouts and pullbacks, $\mathcal C_{\mathrm c}$ is a Serre subcategory, and images of maps from compact objects are closed and compact, so the first localization of \cite[\S5]{BHRAdeles} applies within $\mathcal C$; the arguments of \cite[Lemma 6.4, Proposition 6.5]{BHRAdeles} are then carried out with roofs in this quotient.
Let $\mathcal V$ denote the image of the adelic blocks in the exact hull $(\mathcal C/\mathcal C_{\mathrm c})^{\mathrm{ex}}$; then $\mathcal V\simeq\operatorname{AB}_R$.
The verification of percolation in \cite[Proposition 6.5]{BHRAdeles} is unchanged.
The main conditions of percolation use the fact that an admissible image of an adelic block is quasi-adelic and becomes an adelic block after removing its compact part, together with the injectivity of adelic blocks.
The remaining conditions involve only finitely many kernels, cokernels, and pushouts.
We therefore obtain the localization sequence of \cite[Corollary 6.6]{BHRAdeles} in this subcategory.

For~(3), we follow \cite[Proposition 6.8 and Corollary 6.9]{BHRAdeles}.
By \cite[Lemma 2.18]{BHRAdeles}, the discrete quotient functor $M\mapsto M/Q(M)$ is exact up to a finitely generated defect, so it induces an exact functor to $\operatorname{Mod}_\aleph(R)/\operatorname{mod}(R)$.
In the other direction, the inclusion of countable discrete modules kills $\operatorname{mod}(R)$ in $(\mathcal C/\mathcal C_{\mathrm c})^{\mathrm{ex}}/\mathcal V$, because for finitely generated projective $P$ the cokernel of $P\to P\otimes_R\mathbb A_R$ is compact.
The two functors are mutually inverse: a countable discrete module is its own discrete quotient, and $M\to M/Q(M)$ becomes an isomorphism once the compact and adelic parts of $Q(M)$ are removed.
None of these operations requires uncountable objects.
\end{proof}

\begin{lemma}\label{lem:paper-lca-countability}
Let $R\in\{D,F\}$.
The inclusion of the subcategory $\operatorname{LCA}_{R,\aleph}$ of Lemma~\ref{lem:paper-lca-countable-localization} induces an equivalence
\[
 K(\operatorname{LCA}_{R,\aleph})\xrightarrow{\ \sim\ }
 K(\operatorname{LCA}_R).
\]
Theorem~\ref{thm:paper-lca-source} also holds with $\operatorname{LCA}_R$ replaced by $\operatorname{LCA}_{R,\aleph}$.
\end{lemma}
\begin{proof}
Put $\mathcal C=\operatorname{LCA}_{R,\aleph}$, and write $\mathcal C_{\mathrm c}$ for its compact subcategory.
By Lemma~\ref{lem:paper-lca-countable-localization}, the localization of Braunling--Henrard--van Roosmalen and its discrete quotient restrict to $\mathcal C$.

Since $\mathcal C_{\mathrm c}$ is closed under countable products and $\operatorname{Mod}_\aleph(R)$ under countable sums, the Eilenberg swindle gives
\[
 K(\mathcal C_{\mathrm c})=0,
 \qquad K(\operatorname{Mod}_\aleph(R))=0.
\]
From the two quotient identifications and the localization sequences of Lemma~\ref{lem:paper-lca-countable-localization}, we obtain the same bicartesian square as in the proof of \cite[Theorem 6.2]{BHRAdeles}:
\[
\begin{array}{ccc}
 K(\operatorname{mod}(R))&\longrightarrow&K(\operatorname{AB}_R)\\
 \downarrow&&\downarrow\\
 K(\operatorname{Mod}_\aleph(R))&\longrightarrow&K(\mathcal C)
\end{array}
\]
The top arrow is induced by $P\mapsto P\otimes_R\mathbb A_R$.
The homotopy making the square commute comes from the exact sequence with a diagonal discrete lattice and a compact quotient
\[
 P\longrightarrow P\otimes_R\mathbb A_R
 \longrightarrow (P\otimes_R\mathbb A_R)/P.
\]
Since $K(\operatorname{Mod}_\aleph(R))=0$, the sequence $K(R)\to K(\operatorname{AB}_R)\to K(\mathcal C)$ is a fiber sequence.
The inclusion $\mathcal C\subset\operatorname{LCA}_R$ is compatible with the exact sequence above and induces a map of fiber sequences
\[
\begin{array}{ccccc}
 K(R)&\longrightarrow&K(\operatorname{AB}_R)&\longrightarrow&K(\mathcal C)\\
 \|&&\|&&\downarrow\\
 K(R)&\longrightarrow&K(\operatorname{AB}_R)&\longrightarrow&K(\operatorname{LCA}_R)
\end{array}
\]
whose bottom row is~\eqref{eq:paper-lca-bhr}.
Therefore the right vertical map is also an equivalence.
\end{proof}

Set
\[
 X_R=\tau_{\ge1}K(R),\qquad
 B_R=\tau_{\ge1}K(\operatorname{AB}_R),\qquad
 Y_R=\tau_{\ge1}K(\operatorname{LCA}_R).
\]
The map $K_0(R)\to K_0(\operatorname{AB}_R)$ is injective.
Indeed, if $R_v=M_{r_v}(\Delta_v)$ with $\Delta_v$ a division algebra, the positive generator of $K_0(R)$ maps to $r_v$ times the positive generator of $K_0(R_v)$.
Thus
\begin{equation}\label{eq:paper-lca-positive-row}
 X_R\xrightarrow{\lambda_R}B_R\longrightarrow Y_R
\end{equation}
is also a fiber sequence.
The map $K_1(R)\to K_1(\operatorname{AB}_R)$ is injective as well: for $R=D$, this follows from $\SK_1(D)=0$~\cite{Wang1950} and $F^\times\hookrightarrow F_v^\times$; for $R=F$, it follows from the injectivity of the diagonal map.
Hence applying $\tau_{\ge2}$ to all three terms of~\eqref{eq:paper-lca-positive-row} again gives a fiber sequence \cite[Construction 8.2 and Remark 8.3]{BHRAdeles}.

\subsection{Local--global comparison from a common integral model}
\label{subsec:paper-lca-comparison}
The two localization sequences are compared by means of a common integral model of $D$.
The difficulty is that the global comparison and the local comparisons are constructed separately, so the square relating them must be produced at the level of spectra and not only on homotopy groups.

\begin{lemma}
\label{lem:paper-lca-common-model}
Let $S_0$ be a finite set of finite places containing the finite places at which $D$ is ramified and the places above primes dividing $\deg D$, large enough that $D$ extends to an Azumaya algebra $\mathcal A$ over $U_0=\Spec\os_{F,S_0}$.
The geometric rank of a class over $U_0$ means its geometric rank at the generic point.
Let $\ell$ be a prime dividing $\deg D$, with $\ell\ne2$ if there is real ramification.
At a ramified real place $v$, where $2$ is inverted, the class $[\mathbb H]/2$ serves as a rank-$1$ class; here $\mathbb H$ is the division algebra in the Brauer class of $D_v$.
\begin{enumerate}
\item One can choose
\[
 x_{\ell,U_0}\in\pi_0L_{K(1),\ell}K(\mathcal A)
\]
of geometric rank $1$, with the prescribed rank-$1$ Morita class at every split real place.
\item For $x=x_{\ell,U_0}$, the multiplication map
\[
 m_x\colon L_{K(1),\ell}K(U_0)\longrightarrow L_{K(1),\ell}K(\mathcal A)
\]
is an equivalence.
\end{enumerate}
\end{lemma}
\begin{proof}
Apply the lifting of Proposition~\ref{prop:global-rank} over $U_0$, where $\ell$ is invertible.
Restriction of cohomology with finite coefficients to the real places is an isomorphism in degrees at least three \cite[I, Theorem 4.10(c)]{MilneADT}.
For degree two, choose a place $v\mid\ell$.
By local duality $H^2(F_v,M)\to H^0(F,M^D)^\vee$ is surjective, so by Poitou--Tate the obstruction for any tuple at the real places can be cancelled at $v$, which gives surjectivity onto the real components in degree two \cite[I, Corollary 4.16]{MilneADT}.
The long exact sequence of the real-relative fiber therefore shows that the real-relative sections functor has cohomological dimension at most two.
The only geometric homotopy sheaves with $\Z/\ell^a$ coefficients are $\Z/\ell^a(b)$ in even degrees $2b$.
If no values at real places need to be prescribed, the obstruction from the rank edge map lies in $H^3(U_0,\Z/\ell^a(1))=0$.
In the real-relative case, compute the obstruction to the prescribed tuple of real rank-$1$ classes using relative sections.
The only term remaining in total degree $-1$ is the weight-zero $H^1_{\mathrm{rel}}$, and we write $\delta$ for the second arrow in
\[
 H^0(U_0,\Z/\ell^a)\longrightarrow
 \prod_{v\ \mathrm{real}}H^0(F_v,\Z/\ell^a)
 \longrightarrow H^1_{\mathrm{rel}}(U_0,\Z/\ell^a).
\]
The obstruction $\delta((1)_v)$ vanishes because $(1)_v$ is the restriction of the diagonal element $1$.
Thus for each $a$ there is an element with the required rank and prescribed real values.

Let $E_a$ be the periodic \'etale hypersheaf $V\mapsto(L_{K(1),\ell}K(\mathcal A|_V))/\ell^a$; since $\ell$ is invertible on $U_0$, its global sections are $(L_{K(1),\ell}K(\mathcal A))/\ell^a$.
The rank is a section of a constant sheaf on the connected scheme $U_0$, so rank one at the generic point means rank one at every geometric point.
The fiber of the map from $\Gamma(U_0,\tau_{\ge0}E_a)$ to the pullback prescribing the global rank and real components of the same rank is $R\Gamma_{\mathrm{rel}}(U_0,\tau_{\ge2}E_a)$; by relative cohomological dimension at most two, its $\pi_{-1}$ vanishes and its $\pi_0$ is $H^2_{\mathrm{rel}}(U_0,\Z/\ell^a(1))$.
Thus any discrepancy with the previous choice under coefficient reduction $\ell^{a+1}\to\ell^a$ can be corrected by the weight-one $H^2$ (the relative $H^2$ in the real-relative case).
The next $H^3$ in the coefficient sequence vanishes, so reduction on this finite $H^2$ is surjective, and we can choose compatible lifts inductively.
Finiteness of the finite-coefficient groups over a finite $S_0$ gives the Mittag--Leffler condition, so the Milnor $\varprojlim^1$ vanishes, and $x_{\ell,U_0}$ is the inverse limit of this compatible system.

For $v\notin S_0$, the algebra $\mathcal A_v$ splits and $v\nmid\ell$.
Rigidity and the calculation over finite residue fields give
\[
 \pi_0L_{K(1),\ell}K(\os_v)=\Z_\ell.
\]
Hence $x_{\ell,U_0}|_v$ is precisely the rank-$1$ Morita class.
The same holds at infinite places: by prescription at split real places and by the rank calculation at complex places.
At a ramified real place, for odd $\ell$, geometric rank identifies $\pi_0$ of the height-one localization with $\Z_\ell$, and the rank-$1$ class in this group is $[\mathbb H]/2$.
We have $m_{[\mathbb H]/2}^{-1}=N_{\mathbb H/\mathbb R}/2$.
Thus, with $x=x_{\ell,U_0}$, the class $\beta_\ell=N_{\mathcal A/U_0}(x)-\deg D$ used below also restricts to zero at every infinite place.
Since $\mathcal A$ is \'etale locally a matrix algebra over $U_0$ and $\ell$ is invertible on $U_0$, the proof of Lemma~\ref{lem:paper-etale-rank-one} applies over $U_0$ and proves~(2).
The multiplication map $m_x$ and its inverse are constructed by restriction from the common model, so we obtain at the same time the homotopies expressing compatibility with restriction and the coherence of their composites.
\end{proof}

The single element $x_{\ell,U_0}$ and its base changes allow the completed comparisons over the global field, its completions at finite places, and the integer rings at $v\notin S_0$ to be chosen compatibly with restriction maps.
At $v\notin S_0$ and at split infinite places, these are the usual Morita comparisons.
At ramified real places, Morita invariance identifies the comparison with the normalized norm $N_{\mathbb H/\mathbb R}/2$; in the coordinates of the central simple algebra $D_v$ of degree $\deg D$ this is $N_{D_v/F_v}/\deg D$, since $N_{D_v/F_v}=(\deg D/2)N_{\mathbb H/\mathbb R}$.

We also fix Morita duals of these classes.
Let $\ell$ be a prime, assumed odd if $D$ is ramified at a real place, and put $L_\ell=L_{K(1),\ell}$.
If $\ell\mid\deg D$, let $x$ be the class $x_{\ell,U_0}$ of Lemma~\ref{lem:paper-lca-common-model}; if $\ell\nmid\deg D$, put $x=[\mathcal A]/\deg D$, which is defined $\ell$-locally and has geometric rank $1$.
For the Morita pairing
\[
 \mu\colon L_\ell K(\mathcal A)\otimes_{L_\ell K(U_0)}L_\ell K(\mathcal A^{\mathrm{op}})\longrightarrow L_\ell K(U_0),
\]
the map $\mu(x,-)$ is an equivalence, and we put $y=\mu(x,-)^{-1}(1)$, the \emph{Morita dual} of $x$.
If $\ell\mid\deg D$, then $\ell$ is invertible on $U_0$; \'etale locally on $U_0$ the tensor product of a simple right module and a simple left module is free of rank one, and $\mu(x,-)$ becomes multiplication by a unit, so the assertion follows by \'etale hyperdescent.
If $\ell\nmid\deg D$, we use the free-module argument instead.
The underlying vector bundle of $\mathcal A$ has trivial determinant, since \'etale locally $\mathcal A\simeq\operatorname{End}(E)$ and $\det\operatorname{End}(E)\simeq\os$ canonically, a trivialization that descends because $\mathrm{PGL}_n$ has no nontrivial characters.
With $d=\deg D$, let $i$ and $N$ denote extension of scalars to $\mathcal A$ and the norm.
Then $N\circ i$ is multiplication by $[\mathcal A]$ and $i\circ N$ by $[\mathcal A^\vee]$, the latter by the Morita equivalence $\mathcal A^e\simeq\operatorname{End}(\mathcal A)$; by the triviality of the determinant, both are $d^2$ in $K_0(U_0)$.
Hence $i/d$ and $N/d$ are mutually inverse, and $\mu([\mathcal A]/d,[\mathcal A^{\mathrm{op}}]/d)=[\mathcal A]/d^2=1$; on the opposite side $m_{[\mathcal A^{\mathrm{op}}]/d}$ is also an equivalence, so $\mu(x,-)$ is an equivalence and this element is the Morita dual, that is, $y=[\mathcal A^{\mathrm{op}}]/d$.
In both cases, $\mu(x,y)=1$, the map $m_y\colon L_\ell K(U_0)\to L_\ell K(\mathcal A^{\mathrm{op}})$ is an equivalence, and $\mu(-,y)\simeq m_x^{-1}$.
Over $F$, its completions, and the integer rings, we use the base changes of $x$ and $y$.

Let $g_D\colon X_D\to X_F$ be the comparison of Theorem~\ref{thm:paper-global-comparison} whose completed components are $m_x^{-1}$ for the base changes to $F$ of the classes $x_{\ell,U_0}$ ($\ell\mid\deg D$) and whose component with $\deg D$ inverted is $N_{D/F}/\deg D$.
The only remaining choice is the gluing homotopy in the arithmetic square, which we fix; it need not be compatible with the local comparisons, which is why the square below may have to be corrected.
Write $\nu_{j,D_v}$ for the isomorphisms of Theorem~\ref{thm:paper-local-norm} for $D_v$ at finite places $v$.
At infinite places, we use the same notation for the Morita comparison above, or for $N_{\mathbb H/\mathbb R}/2$ at a ramified real place.

\begin{theorem}
\label{thm:paper-lca-comparison}
One can choose a homotopy-commutative diagram
\begin{equation}\label{eq:paper-lca-comparison-square}
\begin{array}{ccccc}
 X_D&\xrightarrow{\lambda_D}&B_D&\longrightarrow&Y_D\\
 {\scriptstyle g_D}\downarrow&&
 {\scriptstyle h_D}\downarrow\,\simeq&&{\scriptstyle c^{\mathrm{lc}}_D}\downarrow\\
 X_F&\xrightarrow{\lambda_F}&B_F&\longrightarrow&Y_F
\end{array}
\end{equation}
whose rows are fiber sequences and whose middle vertical map is an equivalence.
Moreover, in every degree $j\ge1$, one can arrange that
\begin{equation}\label{eq:paper-lca-all-positive-norms}
 (h_D)_j=\prod'_v\nu_{j,D_v}.
\end{equation}
The dependence of $h_D$ on the choice of $g_D$ is described in Remark~\ref{rem:paper-lca-comparison-choices}.
\end{theorem}
The proof of Theorem~\ref{thm:paper-lca-comparison} is given after Lemma~\ref{lem:paper-lca-regulator-retract} below.

\begin{remark}\label{rem:paper-lca-comparison-choices}
If the global comparison is chosen simultaneously from the common model, $h_D$ itself can be taken to be the restricted product of the local spectrum comparisons.
If this gluing homotopy of $g_D$ is kept fixed, a correction preserving~\eqref{eq:paper-lca-all-positive-norms} may be needed; it is constructed in the proof below.
\end{remark}

We first construct, from the common integral model, a global comparison $g_c$ compatible with the local comparisons (Lemmas~\ref{lem:paper-local-d-arithmetic}--\ref{lem:paper-lca-compatible-gluing}), and then use a rational splitting (Lemma~\ref{lem:paper-lca-regulator-retract}) to return from $g_c$ to the given $g_D$.
Throughout, commutativity of squares of spectra is established by explicit homotopies and is not deduced from equality on homotopy groups.

For a positive integer $n$ and a spectrum $T$, put
\[
 T^\wedge_n=\prod_{\ell\mid n}T^\wedge_\ell.
\]
The arithmetic square for this finite set of primes is
\begin{equation}\label{eq:paper-finite-d-arithmetic}
 T\simeq T[1/n]\mathop{\times}_{T^\wedge_n[1/n]}T^\wedge_n,
\end{equation}
as one verifies by rationalization and comparison modulo each prime; the spectrum $T^\wedge_n[1/n]$ is rational.

\begin{lemma}
\label{lem:paper-local-d-arithmetic}
Let $k$ be a $p$-adic local field, let $D'$ be a central simple algebra over $k$ of degree $e$, not necessarily a division algebra, and let $n$ be a positive multiple of $e$.
One can choose the equivalence of Theorem~\ref{thm:paper-local-spectrum} so that it induces all the isomorphisms of Theorem~\ref{thm:paper-local-norm}, using a prescribed geometric rank-$1$ class at each $\ell\mid n$ and keeping the component with $n$ inverted equal to $N_{D'/k}/e$.
\end{lemma}
\begin{proof}
Put $X_{D'}=\tau_{\ge1}K(D')$ and $X_k=\tau_{\ge1}K(k)$.
If $D'\simeq M_r(\Delta)$, then under Morita equivalence $N_{D'/k}=r\,N_{\Delta/k}$, so Theorem~\ref{thm:paper-local-norm} gives $\pi_jN_{D'/k}=r\deg\Delta\,\nu_{j,\Delta}=e\,\nu_{j,D'}$ in all positive degrees; thus the factor is the degree $e$ of $D'$, not its index.

The prescribed completed comparisons and $N_{D'/k}$ agree on the overlap after dividing by $e$, by the projection formula and the vanishing of products with classes of geometric rank $0$ (Lemma~\ref{lem:paper-rank-zero-product}).
We choose a comparison homotopy in the category of rational spectra and use~\eqref{eq:paper-finite-d-arithmetic} to obtain an initial equivalence $g^0\colon X_{D'}\to X_k$.
The odd-degree calculation and the splitting-transfer and Bockstein argument in the proof of Theorem~\ref{thm:paper-local-spectrum} show that $g^0$ agrees with the prescribed $\nu_{j,D'}$ in odd degrees and on torsion in integral homotopy.
The components with $n$ inverted agree, both being $N_{D'/k}/e$, and in positive even degrees the difference on the completed components factors through a rational group into a finite group, hence vanishes.
Any remaining discrepancy is confined to positive even degrees $j$ and factors through the torsion-free quotient $K_j(D')/K_j(D')_{\mathrm{tors}}$, which is a rational vector space.

By exactness of the arithmetic fiber sequence, the target of this discrepancy is identified with
\[
 \begin{gathered}
 V_j/(\Gamma_{n,j}+\Lambda_j),\qquad V_j=\pi_{j+1}(X_k)^\wedge_n[1/n],\\
 \Gamma_{n,j}=\im\bigl(\pi_{j+1}X_k[1/n]\to V_j\bigr),\\
 \Lambda_j=\im\bigl(\pi_{j+1}(X_k)^\wedge_n\to V_j\bigr).
 \end{gathered}
\]
Since $j+1\ge3$ is odd, the completed groups at $\ell\ne p$ are finite, but $\pi_{j+1}(X_k)^\wedge_p$ is a finitely generated $\Z_p$-module.
Thus $V_j=0$ if $p\nmid n$, but if $p\mid n$, then
\[
 V_j=\pi_{j+1}(X_k)^\wedge_p\otimes_{\Z_p}\Q_p,\qquad
 \Lambda_j=\im(\pi_{j+1}(X_k)^\wedge_p\to V_j)=\pi_{j+1}(X_k)^\wedge_p/(\pi_{j+1}(X_k)^\wedge_p)_{\mathrm{tors}}
\]
and $\Lambda_j$ is a full $\Z_p$-lattice in $V_j$.
Since $\Gamma_{n,j}$ is $p$-divisible, Lemma~\ref{lem:paper-rational-correction}(3) gives
\[
 \operatorname{Ext}^1_{\Z}(\Q,\Gamma_{n,j}+\Lambda_j)=0.
\]
Lemma~\ref{lem:paper-rational-correction}(1), applied to~\eqref{eq:paper-finite-d-arithmetic}, now corrects $g^0$ to an equivalence $g$ that induces the prescribed isomorphisms, and the correction preserves the prescribed components.
\end{proof}

\begin{lemma}
\label{lem:paper-arithmetic-arrow-gluing}
Suppose given maps $a_R\colon X_R\to B'_R$ ($R=D,F$), an equivalence $h'\colon B'_D\to B'_F$, maps $u_{\deg D}\colon X_D[1/\deg D]\to X_F[1/\deg D]$ and $u_\ell\colon(X_D)^\wedge_\ell\to (X_F)^\wedge_\ell$ for each prime $\ell\mid\deg D$, and homotopies of squares compatible with $h'$.
We take the difference between the two comparisons on the overlap, including these homotopies of squares, and obtain a class
\[
 \Xi\in
 \bigl[X_D[1/\deg D],\fib
       \{(X_F)^\wedge_{\deg D}[1/\deg D]\to (B'_F)^\wedge_{\deg D}[1/\deg D]\}\bigr].
\]
If $\Xi=0$, one can choose $g_c\colon X_D\to X_F$ with the prescribed arithmetic components and homotopies of squares.
\end{lemma}
\begin{proof}
Apply~\eqref{eq:paper-finite-d-arithmetic} objectwise in the arrow category $\mathrm{Ar}=\operatorname{Fun}(\Delta^1,\mathrm{Sp})$.
Limits in $\mathrm{Ar}$ are computed objectwise, and the mapping spaces of $\mathrm{Ar}$ are given by the homotopy pullback
\[
 \operatorname{Map}_{\mathrm{Ar}}(a_D,a_F)
 \simeq\operatorname{Map}(X_D,X_F)
 \mathop{\times}_{\operatorname{Map}(X_D,B'_F)}
                         \operatorname{Map}(B'_D,B'_F).
\]
We fix the component $B'_D\to B'_F$ to be $h'$; the resulting homotopy fiber, if nonempty, is a torsor for $\operatorname{Map}(X_D,\fib(a_F))$.
By taking the difference between the two comparisons on the overlap, we pass from this torsor without a chosen base point to the zero fiber and obtain the class $\Xi$ above.
A nullhomotopy of $\Xi$ therefore supplies the homotopy needed to glue the two arithmetic components; we choose one and take the pullback, which preserves the prescribed component homotopies.
\end{proof}

We fix the following choices, which are used in Lemmas~\ref{lem:paper-lca-relative-cohomology}--\ref{lem:paper-lca-compatible-gluing}.
For each $v\in S_0$, apply Lemma~\ref{lem:paper-local-d-arithmetic} to $D'=D_v$ with $n=\deg D'=\deg D$; thus the prescribed primes are all $\ell\mid\deg D$, even when the local index is smaller.
Fix the completed components to be those determined by $x_{\ell,U_0}|_v$, and take $N_{D_v/F_v}/\deg D$ on the component with $\deg D$ inverted.
If there is real ramification, the completed components are taken at odd primes $\ell$ only.
We use the integral Morita comparisons at the remaining places and $N_{D_v/F_v}/\deg D$ at ramified real places, which is $N_{\mathbb H/\mathbb R}/2$ under Morita invariance.
The comparison homotopies from integer rings to fields are obtained by base change from the choices over the integer rings.
The products at each finite stage and their transition maps therefore give an equivalence $h_0\colon B_D\to B_F$.

We write $B'_R$ for the factor of $B_R$ obtained by omitting the finitely many field components in $S_0$,
\[
 B'_R=\tau_{\ge1}\varinjlim_{S\supseteq S_0}\left(
 \prod_{v\mid\infty}K(R_v)\times\prod_{v\in S\setminus S_0}K(R_v)
 \times\prod_{v\notin S,\ v\nmid\infty}K(\mathfrak R_v)\right).
\]
We write $a_R\colon X_R\to B'_R$ for the diagonal map and $h'\colon B'_D\to B'_F$ for the comparison between these factors.
With $h'$ fixed, we glue the global comparison using $u_{\deg D}=N_{D/F}/\deg D$ and the $u_\ell$ obtained from $m_x^{-1}$ on the common model.
On the integer-ring factors, the field factors at finite places $v\notin S_0$, and the split infinite-place factors of $B'_D$, the norm for the degree-$\deg D$ matrix algebra is $\deg D$ times the Morita comparison.
At a ramified real place it is $(\deg D/2)N_{\mathbb H/\mathbb R}$.
Thus the equality
\[
 h'[1/\deg D]\,a_D[1/\deg D]\simeq a_F[1/\deg D]\,(N_{D/F}/\deg D)
\]
of maps of spectra follows from natural isomorphisms of these exact functors.
In what follows, exclude $\ell=2$ when there is real ramification.
We take the homotopies of the squares relating $u_{\deg D}$ and $u_\ell$ to $h'$ from the projection formula for the tensor product and from base change of the common model.

For each $\ell\mid \deg D$, the element
\[
 \beta_\ell=N_{\mathcal A/U_0}(x_{\ell,U_0})-\deg D\in\pi_0L_{K(1),\ell}K(U_0)
\]
vanishes on every local component of $B'_F$.
We take these nullhomotopies from the two specified square homotopies: on each integer-ring component, their difference, multiplied by $\deg D$, precomposed with $m_x$ and evaluated at the unit, is a nullhomotopy of $\beta_\ell$, and on the field components at finite places we use its base change; at the infinite places we restrict the same two squares directly along $U_0\to F_v$.
The comparison multiplied by $\deg D$ lies on the $\ell$-completed side because of the natural norm identifications $(\deg D)u_{\deg D}=N$ and $(\deg D)h'_v\simeq N_v$; we subtract the rank-$1$ side from the norm side.
This gives nullhomotopies compatible with changes of the finite set $S\supset S_0$.
Thus $\beta_\ell$ lifts to an element $\beta_{\ell,\mathrm{rel},S}$ of
\[
 \mathcal F_\ell^S=\fib\!\left(
 L_{K(1),\ell}K(U_S)\longrightarrow \Pi_\ell^S\right),
\]
where
\[
 \Pi_\ell^S=
 \prod_{v\mid\infty}L_{K(1),\ell}K(F_v)
 \times\prod_{v\in S\setminus S_0}L_{K(1),\ell}K(F_v)
 \times\prod_{v\notin S,\ v\nmid\infty}
                         L_{K(1),\ell}K(\os_v).
\]
The relative weight-zero $H^0$ vanishes by injectivity of the diagonal restriction to the infinite places.
The odd geometric layers also vanish, so this lift has descent filtration at least $2$.

\begin{lemma}\label{lem:paper-lca-relative-cohomology}
Let $S\supseteq S_0$ be a finite set of finite places, where all places above $\ell$ lie in $S_0$, and assume $\ell\ne2$ if there is real ramification.
Define the relative continuous cohomology $H^i_{\mathrm{rel},S}(\Z_\ell(t))$ by the fiber of the map from the global complex over $U_S$ to the product of the local complexes in $\Pi_\ell^S$, in which each infinite place occurs exactly once.
Then
\[
 H^i_{\mathrm{rel},S}(\Z_\ell(t))=0\quad(i>3),\qquad
 H^3_{\mathrm{rel},S}(\Z_\ell(t))\text{ is finite}\quad(t\ge2).
\]
\end{lemma}
\begin{proof}
On the global side for finite $S$, inverse limits can be taken using finiteness of finite-coefficient groups; on the local side, derived inverse limits commute with products.
We use Milne's global and local duality.
The $H^i$ of the integer-ring factors vanish for $i\ge2$, and for $v\in S\setminus S_0$ one has $v\nmid\ell$, so
\[
 H^2(F_v,\Z_\ell(t))=
 \Z_\ell/(q_v^{t-1}-1)\Z_\ell\quad(t\ge2)
\]
is finite.
Restriction from global cohomology to the real places is an isomorphism in degrees at least three \cite[I, Theorem 4.10(c)]{MilneADT}.
The relative long exact sequence therefore expresses relative $H^3$ as the cokernel of a map into the direct sum of the finitely many finite local groups $H^2$ above and the finite groups at the real places, and higher relative degrees vanish \cite[I, Theorem 4.10 and Corollary 4.16]{MilneADT}.
\end{proof}

Here and below, the \emph{descent filtration} is the filtration induced by the relative geometric Postnikov tower; its $s$th graded piece in degree $j$ is a subquotient (the $E_\infty$-term) of $H^s_{\mathrm{rel},S}(\Z_\ell(t))$ with $2t-s=j$.
By Lemma~\ref{lem:paper-lca-relative-cohomology}, this tower converges in each fixed degree, and filtrations four and higher contribute nothing.

\begin{lemma}\label{lem:paper-lca-obstruction-vanishing}
For the common-model classes $x_{\ell,U_0}$ for all primes $\ell\mid\deg D$, with $\ell\ne2$ if there is real ramification, the comparison $h'$, the relative lifts $\beta_{\ell,\mathrm{rel},S}$ and the homotopies fixed above, the obstruction class $\Xi$ of Lemma~\ref{lem:paper-arithmetic-arrow-gluing} vanishes.
\end{lemma}
\begin{proof}
To compute $\Xi$ on the common model, fix $\ell$ and put $x=x_{\ell,U_0}$.
On the overlap, with $d=\deg D$, we have $d(u_d-u_\ell)\simeq(m_{N_{D/F}(x)}-d)\,m_x^{-1}=m_{\beta_\ell}m_x^{-1}$ by the projection formula, and the nullhomotopies of $\beta_\ell$ fixed above make this a map into the relative fiber.
Over finite $U_S$, the projection formula $N_{\mathcal A|_{U_S}/U_S}m_x\simeq m_{N(x)}$ leads us to consider the relative map on the common model
\begin{equation}\label{eq:paper-lca-relative-defect}
 m_{\beta_{\ell,\mathrm{rel},S}}m_x^{-1}.
\end{equation}
We use this expression as a map from $\tau_{\ge1}K(\mathcal A|_{U_S})$, so it factors through $\tau_{\ge1}\mathcal F_\ell^S$; by the choice of the nullhomotopies of $\beta_\ell$, it represents $d\,\Xi_\ell$ as a map into the relative fiber, not merely as a formula, where $\Xi_\ell$ is the $\ell$-component of $\Xi$, and multiplying by $d^{-1}$ on the overlap gives $\Xi_\ell$.
Since finite Moore coefficients commute with products and filtered colimits, we compare the common-model diagrams with finite coefficients in positive degrees for each finite $S$, take the colimit in the completed category, and then pass to the overlap.
Each global class used to test the obstruction is computed by returning to the relative spectrum $\mathcal F_\ell^S$ at a finite $S$ before this passage.
Since the target is rational, the source of $\Xi$ can be replaced by $X_{D,\Q}$.

A global class of positive odd degree $j=2m-1$ comes from some $U_S$.
Since $j$ is odd, the geometric layer in degree $j$ vanishes, so the completed image of this class has descent filtration at least $1$.
The map $m_x^{-1}$ is a map of \'etale modules and preserves the descent filtration.
The relative fiber $\mathcal F_\ell^S$ is a module over the global side.
For finite $S$, the class $\beta_{\ell,\mathrm{rel},S}$ lies in relative geometric Postnikov stage at least $2$, and after $m_x^{-1}$ the class of odd degree $j$ lifts to absolute geometric stage at least $j+1$.
Their product in~\eqref{eq:paper-lca-relative-defect} therefore lies in relative stage at least $j+3$, that is, in descent filtration at least $3$, and the only possible term remaining in degree $j=2m-1$ is $H^3_{\mathrm{rel},S}(\Z_\ell(m+1))$.
This group is finite because $m+1\ge2$, so it vanishes after rationalization.
In positive even degrees the source is $K_{2m}(D)\otimes\Q=0$, since $N_{D/F}i_D=i_DN_{D/F}=(\deg D)^2$ gives $K_{2m}(D)\otimes\Q\cong K_{2m}(F)\otimes\Q$, which vanishes by \cite{Borel}.
Integral classes generate the rationalized homotopy groups, so $\Xi$ is zero on all homotopy groups, and hence $\Xi=0$ because maps between rational spectra are determined by their effect on homotopy groups.
\end{proof}

\begin{lemma}
\label{lem:paper-lca-compatible-gluing}
For the choices in Lemma~\ref{lem:paper-lca-common-model}, one can choose a comparison $g_c\colon X_D\to X_F$ and an equivalence $h_0\colon B_D\to B_F$ given by the restricted product of the local comparisons such that
\[
 h_0\lambda_D\simeq\lambda_Fg_c.
\]
The maps $g_c$ and $g_D$ have the same completed components and the same component with $\deg D$ inverted.
\end{lemma}
\begin{proof}[Proof of Lemma~\ref{lem:paper-lca-compatible-gluing}]
The equivalence $h_0$ was constructed above.
By Lemma~\ref{lem:paper-lca-obstruction-vanishing}, Lemma~\ref{lem:paper-arithmetic-arrow-gluing} gives $g_c$, and the gluing preserves $h'$ and the prescribed homotopies of completed squares.

Finally, we restore the finitely many places $v\in S_0$.
For $X_{F_v}=\tau_{\ge1}K(F_v)$, the two localization maps agree on the completed components and the component with $\deg D$ inverted.
The obstruction to lifting their difference, together with the homotopies on those prescribed components, lies in
\[
 [X_D,\Omega (X_{F_v})^\wedge_{\deg D}[1/\deg D]]=[X_{D,\Q},\Omega (X_{F_v})^\wedge_{\deg D}[1/\deg D]].
\]
The rationalized source has homotopy only in positive odd degrees, and $(X_{F_v})^\wedge_{\deg D}[1/\deg D]$ has zero homotopy in positive even degrees.
Therefore the group above vanishes; equivalently,
\[
 \pi_1\operatorname{Map}\bigl(X_D,(X_{F_v})^\wedge_{\deg D}[1/\deg D]\bigr)=0.
\]
Hence the difference on the overlap between the homotopies prescribed on the two arithmetic components can be removed, and the places in $S_0$ can be restored while retaining the prescribed projection-formula homotopies of the common model.
\end{proof}

\begin{lemma}\label{lem:paper-lca-regulator-retract}
The map $\lambda_{F,\Q}\colon X_{F,\Q}\to B_{F,\Q}$ is a split monomorphism.
\end{lemma}
\begin{proof}
In degree one, $F^\times\to\mathbb A_F^\times$ is injective and remains so after rationalization.
In positive even degrees, $K_{2m}(F)_\Q=0$.
Let $r_1$ and $r_2$ be the numbers of real and complex places of $F$, respectively.
For $m\ge2$, put $r_m=r_1+r_2$ if $m$ is odd, and $r_m=r_2$ if $m$ is even.
In odd degrees $2m-1\ge3$, the Borel regulator through the infinite-place terms,
\[
 K_{2m-1}(F)\longrightarrow
 \prod_{v\mid\infty}K_{2m-1}(F_v)\longrightarrow\mathbb R^{r_m},
\]
is injective on the rational global $K$-group \cite{Borel}.
Hence $\lambda_{F,\Q}$ is injective on every homotopy group.
Since the category of rational spectra is equivalent to the category of complexes over $\Q$ and short exact sequences of $\Q$-vector spaces split, we can choose
\[
 \sigma_F\colon B_{F,\Q}\longrightarrow X_{F,\Q},\qquad
 \sigma_F\lambda_{F,\Q}\simeq\mathrm{id}.
\]
\end{proof}

\begin{proof}[Proof of Theorem~\ref{thm:paper-lca-comparison}]
We take $g_c$ and $h_0$ from Lemma~\ref{lem:paper-lca-compatible-gluing}.
Since $g_c-g_D$ vanishes on each component of the arithmetic square, we choose $b_0\colon X_D\to\Omega (X_F)^\wedge_{\deg D}[1/\deg D]$ with $g_c-g_D=\partial_F b_0$, where $\partial_F$ is the boundary in the arithmetic fiber sequence.
The target is rational, so $b_0$ factors through $X_{D,\Q}$.
Let $f$ be the composite of this factorization with $\partial_F$.
Then
\[
 g_c-g_D=fq_D,\qquad
 f\colon X_{D,\Q}\longrightarrow X_F,\qquad q_Ff\simeq0,
\]
where $q_D$ and $q_F$ are the rationalization maps.
Here $q_Ff\simeq0$ because the composite $\Omega (X_F)^\wedge_{\deg D}[1/\deg D]\xrightarrow{\partial_F}X_F\to X_F[1/\deg D]$ vanishes and $q_F$ factors through $X_F\to X_F[1/\deg D]$.
Since $g_D[1/\deg D]=N_{D/F}/\deg D$ and $N_{D/F}i_D=i_DN_{D/F}=(\deg D)^2$, the map $g_{D,\Q}$ is an equivalence, and putting
\[
 t=f(g_{D,\Q})^{-1}\colon X_{F,\Q}\longrightarrow X_F,
 \qquad \varepsilon=tq_F
\]
gives $q_Ft=0$, $\varepsilon^2=0$, and $g_c=(1+\varepsilon)g_D$.

Using $\sigma_F$ from Lemma~\ref{lem:paper-lca-regulator-retract}, we define
\[
 \zeta=\lambda_Ft\,\sigma_Fq_{B_F},\qquad b=1+\zeta.
\]
The identities $\sigma_F\lambda_{F,\Q}=1$ and $q_Ft=0$ give
\[
 \zeta^2=0,\qquad
 b^{-1}=1-\zeta,\qquad
 b\lambda_F=\lambda_F(1+\varepsilon).
\]
Here we choose $b_0$ together with the prescribed nullhomotopies of $g_c-g_D$ on the arithmetic components, that is, in the fiber of the corresponding mapping spaces; these data give $b\lambda_F\simeq\lambda_F(1+\varepsilon)$ and hence the square for $h_D$.
The correction pair $(\varepsilon,\zeta)$ is the map of arrow categories induced by $(t,\lambda_Ft\sigma_F)$ and $\sigma_F\lambda_{F,\Q}\simeq1$, so it factors through $\lambda_{F,\Q}$; in particular, for every prime $p$, the height-one localization $L_{K(1),p}$ kills the correction together with the homotopy of the square.
Consequently, $h_D=b^{-1}h_0$ is an equivalence, and
\[
 h_D\lambda_D\simeq b^{-1}\lambda_Fg_c
 \simeq b^{-1}\lambda_F(1+\varepsilon)g_D
 \simeq\lambda_Fg_D.
\]
We use this homotopy to define the map $c^{\mathrm{lc}}_D$ between cofibers.

Finally, we check that this correction does not change the norms on homotopy groups.
For $j\ge1$, the source of $t$ is rational, so the image of $t_*$ in $K_j(F)$, or in $K_j(F)[1/2]$ under real ramification, is divisible and hence zero by Lemma~\ref{lem:paper-global-no-divisible}.
Thus $\zeta_*=0$ and $b_*=1$, so $h_D$ satisfies the same formula~\eqref{eq:paper-lca-all-positive-norms} as $h_0$.
Moreover, $\zeta$ factors through a rational spectrum and hence vanishes with all finite coefficients, so the completed components are unchanged as well.
\end{proof}

\begin{remark}\label{rem:paper-lca-choice}
Theorem~\ref{thm:paper-lca-comparison} is an existence result, and its construction involves choices.
We do not deduce that the correction maps $t$ or $\zeta$ vanish as maps of spectra merely because they vanish on homotopy groups.
On the other hand, since $g_c=(1+\varepsilon)g_D$, returning from $g_c$ to $g_D$ changes the comparison by a self-equivalence of the target and therefore preserves the comparison fiber and the splitting in degrees at least two.
When expressing $m_x^{-1}$ on a periodic completed component by the Morita pairing $\mu$, we use the Morita dual $y$ of $x$ fixed after Lemma~\ref{lem:paper-lca-common-model} and its base changes.
This normalizes the comparison map more strongly than merely prescribing geometric rank.
The homotopies on the completed components of~\eqref{eq:paper-lca-comparison-square} are also chosen to be the projection-formula homotopies associated with base change of the common model and evaluation against $y$.
Corrections factoring through rational spectra vanish under height-one localization, so the height-one localization of~\eqref{eq:paper-lca-comparison-square} is preserved after correction.
The map $c^{\mathrm{lc}}_D$ is the cofiber map of the specified square.
Independence of the choice of commutativity homotopy is not asserted.
\end{remark}

\begin{lemma}
\label{lem:paper-lca-periodic-evaluation}
Let $p$ be a prime, assumed odd if there is real ramification.
Let $x$ and $y$ be the class at $p$ on the common model and its Morita dual, fixed after Lemma~\ref{lem:paper-lca-common-model}.
Put $L=L_{K(1),p}$, and let $h_0$ be the equivalence of Lemma~\ref{lem:paper-lca-compatible-gluing}.
Under the equivalences $LB_R\simeq L\Theta_{R,p}$ of the proof, $Lh_0$ is the $L$-localization of the evaluation map $\Theta_{D,p}\to\Theta_{F,p}$ given componentwise by $\mu(-,y|_v)\simeq m_{x|_v}^{-1}$.
Consequently $L$ applied to~\eqref{eq:paper-lca-comparison-square} is the map of localized rows given by evaluation against $y$, and $Lc^{\mathrm{lc}}_D$ is the induced map on cofibers.
\end{lemma}
\begin{proof}
First observe that $LT=0$ for every bounded-above spectrum $T$.
We choose a finite type $1$ spectrum with a positive-degree $v_1$ self-map \cite{HopkinsSmith1998}.
Smashing its telescope with $T$ gives a contractible spectrum, since in each fixed degree the terms of the resulting sequential colimit eventually vanish.
The smash product of $K(1)$ with this finite spectrum is a finite direct sum of suspensions of nonzero $K(1)$, and the self-map acts on it by an equivalence.
Hence $K(1)\otimes T=0$ and therefore $LT=0$.
This applies in particular to the part removed by taking a connective cover.

Let $\mathfrak R$ be a maximal order in the common integral model, and define, in spectra,
\[
 \Theta_{R,p}=\varinjlim_S\left(
 \prod_{v\mid\infty}LK(R_v)\times\prod_{v\in S}LK(R_v)
 \times\prod_{v\notin S,\ v\nmid\infty}LK(\mathfrak R_v)\right).
\]
We do not assume that this object is already local.
The natural comparison $B_R\to \Theta_{R,p}$ induces isomorphisms on mod-$p$ homotopy groups in degrees $\ge2$.
Indeed, the finite-coefficient comparisons over nonarchimedean local fields, integer rings, $\mathbb R$, and $\mathbb C$ all hold in this uniform range.
For the nonarchimedean field components, use Lemma~\ref{lem:paper-local-compare}.
Restricting the colimit to $S\supset S_0$ makes every algebra in an integer-ring component split.
At $v\nmid p$, use Henselian rigidity and the finite-residue-field calculation.
At $v\mid p$, the localization sequence
\[
 K(k_v)\longrightarrow K(\os_v)\longrightarrow K(F_v)
\]
together with $K(k_v)/p\simeq H\F_p$ and $LK(k_v)=0$, shows that the integer-ring and field comparisons agree in degrees at least two.
For complex components, use Suslin's finite-coefficient comparison.
For real components at $p=2$, use the comparison of real algebraic $K$-theory with $KO$; at odd primes, restriction and transfer for $\mathbb R\subset\mathbb C$, together with invertibility of $2$, reduce to the complex comparison.
At ramified real components, $p$ is odd and $N_{\mathbb H/\mathbb R}/2$ reduces the problem to the same calculation over $\mathbb R$ \cite{QuillenFinite,ClausenMathewMorrow,SuslinClosed,KaroubiWeibelReal}.
Finite Moore spectra are dualizable and hence commute with products; the same range is preserved in the comparisons involving filtered colimits and positive covers.
The cofiber modulo $p$ is therefore bounded above.
The vanishing proved above and conservativity of reduction modulo $p$ on $K(1)$-local spectra give
\[
 LB_R\xrightarrow{\sim}L\Theta_{R,p}.
\]
For the conservativity assertion, if $T/p=0$, then $p$ is invertible on $T$ but acts by zero on $K(1)$, so $K(1)\otimes T=0$; if $T$ is also $K(1)$-local, then $T=0$.

Tensoring with perfect modules over the common model is compatible with all local factors and transition maps, so $\Theta_{R,p}$ can be constructed in the category of $LK(U_0)$-modules.
Evaluation against $y$ gives $\Theta_{D,p}\to\Theta_{F,p}$ with local components $m_{x|_v}^{-1}$.
Taking products and the colimit of the local homotopies used to construct $h_0$ gives a commutative square with the maps $B_R\to\Theta_{R,p}$.
Hence $Lh_0$ is the height-one localization of this evaluation.

Tensoring with a perfect left $\mathcal A$-module preserves the exact structures on discrete modules, adelic blocks, and locally compact modules, and commutes naturally with localization.
Hence tensoring with perfect modules over the common model also acts on the localization sequence of Theorem~\ref{thm:paper-lca-source}, giving a single diagram of the rows and their cofibers.
In the gluing above, we have prescribed the homotopies on the completed components of this diagram; if $p\nmid\deg D$, the relevant component is the one with $\deg D$ inverted, where $y=[\mathcal A^{\mathrm{op}}]/\deg D$ and we use the projection-formula homotopy for $N_{D/F}/\deg D$.
Lifting at the finitely many exceptional places preserves them, and the rational correction returning to the fixed $g_D$ vanishes under $L$.
Applying $L$ kills the overlap; for $p\mid\deg D$ the $p$-completed component survives and for $p\nmid\deg D$ the component with $\deg D$ inverted, so the whole square agrees with the prescribed evaluation.
Thus $Lc^{\mathrm{lc}}_D$ is evaluation against the same $y$.
\end{proof}

\subsection{The suspension of the comparison fiber and splitting in degrees at least two}
\label{subsec:paper-lca-fiber}
The global splitting in Section~\ref{sec:global-spectrum} passes to the cofibers of the comparison square.

\begin{theorem}
\label{thm:paper-lca-fiber-splitting}
Let $Z=\fib(g_D)$.
The comparison of Theorem~\ref{thm:paper-lca-comparison} satisfies
\begin{equation}\label{eq:paper-lca-fiber-suspension}
 \fib(c^{\mathrm{lc}}_D)\simeq\Sigma Z.
\end{equation}
There is also a noncanonical equivalence whose first projection is $\tau_{\ge2}c^{\mathrm{lc}}_D$,
\begin{equation}\label{eq:paper-lca-splitting}
 \tau_{\ge2}K(\operatorname{LCA}_D)
 \simeq\tau_{\ge2}K(\operatorname{LCA}_F)\oplus\Sigma Z.
\end{equation}
On the second summand, the localization boundary is the natural inclusion $\Sigma Z\to\Sigma\tau_{\ge2}K(D)$.
\end{theorem}
\begin{proof}
By the $3\times3$ lemma, the vertical fibers in~\eqref{eq:paper-lca-comparison-square} form a fiber sequence
\[
 Z\longrightarrow\fib(h_D)\longrightarrow\fib(c^{\mathrm{lc}}_D).
\]
Since $h_D$ is an equivalence, the middle term vanishes, giving \eqref{eq:paper-lca-fiber-suspension}.

We have checked after~\eqref{eq:paper-lca-positive-row} that applying $\tau_{\ge2}$ to all three terms still gives a fiber sequence.
By Theorem~\ref{thm:paper-global-split}, we may write
\[
 \tau_{\ge2}K(D)\simeq\tau_{\ge2}K(F)\oplus Z,
 \qquad\tau_{\ge2}g_D=\operatorname{pr}_1.
\]
After identifying the adelic terms by $\tau_{\ge2}h_D$, we check, including the homotopies, that $\lambda_D\simeq(\lambda_F,0)$.
Let $s$ be a section of $\tau_{\ge2}g_D$ and $i\colon Z\to\tau_{\ge2}K(D)$ the fiber inclusion, and use the decomposition given by $s\oplus i$.
Restricting the specified homotopy of the square along $s$ gives $\lambda_Ds\simeq\lambda_Fg_Ds\simeq\lambda_F$.
Restricting the same homotopy along $i$ and composing with the fiber nullhomotopy $g_Di\simeq0$ gives $\lambda_Di\simeq\lambda_Fg_Di\simeq0$.
Thus $\lambda_D\simeq(\lambda_F,0)$ in the arrow category, and this equivalence also identifies $c^{\mathrm{lc}}_D$ on cofibers with the first projection.
Taking cofibers gives
\[
 \cofib(\lambda_F,0)\simeq\cofib(\lambda_F)\oplus\Sigma Z,
\]
which proves~\eqref{eq:paper-lca-splitting}.
The final assertion follows from the definition of the cofiber boundary.
\end{proof}

\begin{corollary}
\label{cor:paper-lca-gem}
Let $e_\ell=\ell^{v_\ell(\deg D)}$ and $G_{\ell,j}=\pi_jZ_{(\ell)}$, as in Theorem~\ref{thm:paper-global-fiber-gem}.
Then
\[
 \fib(c^{\mathrm{lc}}_D)\simeq
 \bigoplus_{\ell\mid\deg D}\ \bigoplus_{j=4}^{e_\ell-3}
       \Sigma^{j+1}H(G_{\ell,j}),
\]
and, for $n\ge2$,
\[
 K_n(\operatorname{LCA}_D)\cong
 K_n(\operatorname{LCA}_F)\oplus
       \bigoplus_{\ell\mid\deg D}G_{\ell,n-1}.
\]
Here we set $G_{\ell,j}=0$ outside the stated range, and we omit $\ell=2$ from the sums if there is real ramification.
In particular, $c^{\mathrm{lc}}_D$ is an isomorphism in degrees $1$ through $4$, and the first possible difference is in degree $5$.
If $e_3=9$, then at $\ell=3$ the only possible differences are three finite groups, in degrees $5$, $6$, and $7$.
\end{corollary}
\begin{proof}
We combine Theorem~\ref{thm:paper-global-fiber-gem} with~\eqref{eq:paper-lca-fiber-suspension} and~\eqref{eq:paper-lca-splitting}.
The degree-one isomorphism follows from the long exact sequence of the fiber sequence $\fib(c^{\mathrm{lc}}_D)\to Y_D\to Y_F$.
No assertion is made that any of these finite groups is nonzero.
\end{proof}

\begin{theorem}\label{thm:paper-lca-k2}
For every central division algebra $D$ over $F$, there is an isomorphism
\[
 K_2(\operatorname{LCA}_D)\cong K_2(\operatorname{LCA}_F)
\]
with integral coefficients, also when $D$ is ramified at a real place.
In particular, \cite[Conjecture 10.2]{BHRAdeles} holds.
\end{theorem}
\begin{proof}
By \cite[Theorem 10.5]{BHRAdeles}, it suffices that, for $D$ and for every completion $D_v$, the reduced norm on $K_2$ is injective and its image is the kernel of the Hilbert symbols at the real places at which the algebra is ramified.
For $D$ this is Theorem~\ref{thm:paper-global-norm}, and for $D_v$ at a finite place $v$ it is Theorem~\ref{thm:paper-local-norm}.
At an infinite place at which $D$ splits, it holds by Morita equivalence.
At a real place at which $D$ is ramified, $D_v$ is Morita equivalent to the Hamilton quaternions $\mathbb H$, of index $2$.
By \cite[Corollary 2]{KahnLevine2010}, whose hypothesis is supplied by the norm residue isomorphism theorem \cite{Voevodsky}, $\Nrd_2$ is injective on $K_2(\mathbb H)$ with image the kernel of $[\mathbb H]\cup-\colon K_2(\mathbb R)\to H^4(\mathbb R,\Z/2)$.
Since $[\mathbb H]\cup(a)\cup(b)\ne0$ exactly when $a$ and $b$ are both negative, this kernel is the kernel of the Hilbert symbol.
\end{proof}

\subsection{Selmer \texorpdfstring{$K$}{K}-homology and id\`ele classes}\label{subsec:paper-lca-selmer}
We use Selmer $K$-homology to compare $c^{\mathrm{lc}}_D$ with Clausen's pro-$p$ Artin map and compute the map on id\`ele classes (Theorem~\ref{thm:paper-lca-artin}).

Let $p$ be a prime, assumed odd if $D$ is ramified at a real place.
Write $L=L_{K(1),p}$.
We write $\mathrm{lc}_\Z$ for the bounded derived category of second-countable locally compact abelian groups $\operatorname{LCA}_{\Z,\aleph}$, as in \cite[Definitions 0.2 and 3.2]{ClausenArtin}, $\omega=\omega_{K(1),p}$ for the $K(1)$-local dualizing object, $\mathbb D_p(T)=\operatorname{map}_{\mathrm{Sp}}(T,\omega)$ for the function spectrum.
For a small idempotent-complete $\Z$-linear stable $\infty$-category $\mathcal P$, Clausen's $p$-adic Selmer $K$-homology $dK^{\mathrm{Sel}}_p(\mathcal P)$ is the homotopy pushout of $\mathbb D_pK(\mathcal P)\leftarrow\mathbb D_p\mathrm{TC}(\mathcal P)\to d_{\mathrm{TC}}\mathrm{TC}(\mathcal P)$, where the left map is dual to the cyclotomic trace and $d_{\mathrm{TC}}$ is Clausen's duality for $\mathrm{TC}$ \cite[Definition 2.1]{ClausenArtin}; if $p$ is invertible in $\mathcal P$, it is equivalent to $\mathbb D_pK(\mathcal P)$ \cite[Proposition 2.11(1)]{ClausenArtin}.
We write $\mathbf j\in\pi_0dK^{\mathrm{Sel}}_p(\mathrm{lc}_\Z)$ for Clausen's fundamental class \cite[Definitions 1.4 and 4.2]{ClausenArtin}.
Under our convention of right $D$-modules, the perfect complexes used for evaluation are left $D$-modules, that is, objects of $\operatorname{Perf}(D^{\mathrm{op}})$.
Accordingly, define the Selmer objects by
\[
 \mathsf S_{D,p}=dK^{\mathrm{Sel}}_p(D^{\mathrm{op}}),
 \qquad
 \mathsf S_{F,p}=dK^{\mathrm{Sel}}_p(F).
\]
The evaluations below are performed on second-countable locally compact modules and transferred to the whole category using the $K$-equivalence of the inclusion in Lemma~\ref{lem:paper-lca-countability}.
Indeed, a second-countable right $D$-module $M_0$ defines an exact functor $\operatorname{Perf}(D^{\mathrm{op}})\to\mathrm{lc}_{\Z}$ by $P\mapsto M_0\otimes_D P$.
This correspondence induces an exact functor
\[
 D^b(\operatorname{LCA}_{D,\aleph})\longrightarrow
 \operatorname{Fun}_{\Z}
 \bigl(\operatorname{Perf}(D^{\mathrm{op}}),\mathrm{lc}_{\Z}\bigr).
\]
Pulling back Clausen's fundamental class along this evaluation gives
\[
 \operatorname{Art}_{D,p}\colon Y_D\longrightarrow\mathsf S_{D,p}.
\]
This is the composite of the Artin map for $\mathrm{lc}_{\operatorname{Perf}(D^{\mathrm{op}})}=\operatorname{Fun}_{\Z}(\operatorname{Perf}(D^{\mathrm{op}}),\mathrm{lc}_{\Z})$ with the $K$-theory of the exact functor above.
No equivalence between the stable categories themselves is required \cite[Definitions 3.2 and 4.2, Remarks 4.3--4.4]{ClausenArtin}.

\begin{theorem}\label{thm:paper-lca-artin}
Let $x$ be the class at $p$ on the common model and $y$ its Morita dual, fixed after Lemma~\ref{lem:paper-lca-common-model}.
Then $\rho_{D,p}=\mathbb D_p(m_y)\colon\mathsf S_{D,p}\to\mathsf S_{F,p}$ is an equivalence, and the comparison $c^{\mathrm{lc}}_D\colon Y_D\to Y_F$ of Theorem~\ref{thm:paper-lca-comparison} fits into a homotopy-commutative square
\begin{equation}\label{eq:paper-lca-artin-square}
 \begin{array}{ccc}
 Y_D&\xrightarrow{\operatorname{Art}_{D,p}}&\mathsf S_{D,p}\\
 {\scriptstyle c^{\mathrm{lc}}_D}\downarrow&&\downarrow{\scriptstyle\rho_{D,p}}\\
 Y_F&\xrightarrow{\operatorname{Art}_{F,p}}&\mathsf S_{F,p}.
 \end{array}
\end{equation}
Moreover, the restriction of $\operatorname{Art}_{D,p}$ to $\fib(c^{\mathrm{lc}}_D)\simeq\Sigma Z$ is zero.
In terms of Clausen's isomorphism $e_{F,p}\colon G_F^{\mathrm{ab},\wedge p} \xrightarrow{\sim}\pi_1\mathsf S_{F,p}$, the map on id\`ele classes is given, for $[(a_v)_v]\in\pi_1Y_D=\prod'_vK_1(D_v)/K_1(D)$, by
\begin{equation}\label{eq:paper-lca-reduced-artin}
 e_{F,p}^{-1}\,\pi_1(\rho_{D,p}\operatorname{Art}_{D,p})
 \bigl([(a_v)_v]\bigr)
 =
 \operatorname{rec}_{F,p}
 \bigl([(\Nrd_{D_v/F_v}(a_v))_v]\bigr).
\end{equation}
Here $a_v\in K_1(D_v)$ is assumed to come from the integral model at all but finitely many finite places.
The map $\operatorname{rec}_{F,p}$ is the global Artin map normalized so that a uniformizer at an unramified place maps to the geometric Frobenius, followed by the projection to the maximal pro-$p$ quotient, and $e_{F,p}$ is normalized as in \cite[Definition 2.37]{ClausenArtin}.
\end{theorem}

We prove Theorem~\ref{thm:paper-lca-artin} after two lemmas, which share the following choices.
Let $x\in\pi_0LK(D)$ and $y\in\pi_0LK(D^{\mathrm{op}})$ be the base changes to $F$ of the class at $p$ on the common model and its Morita dual, fixed after Lemma~\ref{lem:paper-lca-common-model}.
We also fix the homotopy making~\eqref{eq:paper-lca-comparison-square} commute, chosen as in Remark~\ref{rem:paper-lca-choice}, which is compatible with the localization sequences.

\begin{lemma}\label{lem:paper-lca-selmer-evaluation}
With these choices, $m_y\colon LK(F)\to LK(D^{\mathrm{op}})$ is an equivalence, and so is
\[
 \rho_{D,p}=\mathbb{D}_p(m_y)\colon\mathsf S_{D,p}\longrightarrow\mathsf S_{F,p}.
\]
The Selmer $K$-theory of the bilinear functor $(M_0,P)\mapsto M_0\otimes_DP$ defines a map $\Psi_y\colon LY_D\to LY_F$, evaluation against $y$, and $Lc^{\mathrm{lc}}_D\simeq\Psi_y$ as maps of localization sequences.
\end{lemma}
\begin{proof}
Put $M=LK(D)$ and $M'=LK(D^{\mathrm{op}})$.
Since the categories over a number field are $\Z[1/p]$-linear, \cite[Proposition 2.11]{ClausenArtin} gives
\[
 \mathsf S_{D,p}\simeq \mathbb{D}_p(M'),\qquad
 \mathsf S_{F,p}\simeq \mathbb{D}_p(LK(F)).
\]
By the choice of $y$, the map $m_y\colon LK(F)\to M'$ is an equivalence and $\mu(-,y)\simeq m_x^{-1}\colon M\to LK(F)$.
If $p\nmid \deg D$, then $y=[D^{\mathrm{op}}]/\deg D$ and $\mu(-,y)\simeq N_{D/F}/\deg D$, where $N_{D/F}$ is the norm used in Section~\ref{subsec:paper-lca-comparison}.
Thus $\rho_{D,p}=\mathbb D_p(m_y)$ is an equivalence.

In general, $y$ is not represented by an actual perfect module, so its action is defined on localized $K$-theory.
Put $\mathcal L_D=D^b(\operatorname{LCA}_{D,\aleph})$ and $\mathcal L_F=D^b(\operatorname{LCA}_{F,\aleph})$.
Tensoring with a finitely generated projective left $D$-module is defined on locally compact modules as a direct summand of a finite direct sum, and preserves second countability and exact sequences.
It therefore gives an exact bilinear functor
\[
 \mathcal L_D\otimes_F\operatorname{Perf}(D^{\mathrm{op}})
       \longrightarrow\mathcal L_F,
 \qquad (M_0,P)\longmapsto M_0\otimes_D P.
\]
The categories above and $\mathrm{lc}_{\operatorname{Perf}(D^{\mathrm{op}})}$ are $\Z[1/p]$-linear, so their $p$-completed $\operatorname{TC}$ vanishes, as for Selmer $K$-homology in \cite[Proposition 2.11]{ClausenArtin}.
For each of these categories $\mathcal P$, the definition of Selmer $K$-theory therefore gives
\[
 K^{\mathrm{Sel}}_p(\mathcal P)\simeq LK(\mathcal P).
\]
Thus the Selmer $K$-theory of this functor defines a pairing of spectra over $LK(F)$,
\[
 LK(\mathcal L_D)\otimes_{LK(F)}M'\longrightarrow LK(\mathcal L_F).
\]
Evaluation against $y$ gives $\Psi_y$ by substituting $m_y\colon LK(F)\to M'$ in the second variable.

Local maximal orders are hereditary and have no negative $K$-groups, so every factor in~\eqref{eq:paper-lca-restricted-spectrum}, and hence the filtered colimit, is connective.
Since $K(D)$ and $K(F)$ are connective, the localization sequences show that $K(\operatorname{LCA}_D)$ and $K(\operatorname{LCA}_F)$ are connective as well.
Taking positive covers thus removes only $H(K_0)$, which vanishes under height-one localization.
Using Lemma~\ref{lem:paper-lca-countability}, we therefore have $LK(\mathcal L_D)\simeq LY_D$ and $LK(\mathcal L_F)\simeq LY_F$.
With the homotopy making~\eqref{eq:paper-lca-comparison-square} commute fixed before this lemma, Lemma~\ref{lem:paper-lca-periodic-evaluation} gives, as maps of entire localization sequences,
\[
                         Lc^{\mathrm{lc}}_D\simeq \Psi_y.
\]
\end{proof}

\begin{lemma}\label{lem:paper-lca-artin-compatibility}
With the choices of Lemma~\ref{lem:paper-lca-selmer-evaluation},
\[
 \rho_{D,p}\operatorname{Art}_{D,p}^{L}
 \simeq\operatorname{Art}_{F,p}^{L}\Psi_y
 \simeq\operatorname{Art}_{F,p}^{L}Lc^{\mathrm{lc}}_D,
\]
where the superscript $L$ denotes factorization through height-one localization.
\end{lemma}
\begin{proof}
The Artin map of \cite[Definition 4.2]{ClausenArtin} is the composite of the Selmer $K$-theory evaluation pairing with the fundamental class $\mathbf j$.
We first define $\mathbf j_F\colon LK(D^b(\operatorname{LCA}_{F,\aleph}))\to\omega$ as the pullback of $\mathbf j$ along the forgetful functor to $\mathrm{lc}_\Z$, and we compare both sides after composing with the same $\mathbf j_F$.
Associativity of evaluation,
\[
 (M_0\otimes_D P)\otimes_F P'
 \simeq M_0\otimes_D(P\otimes_F P'),
\]
is a natural isomorphism of bilinear functors for this functor and the action of $\operatorname{Perf}(F)$, where $P'\in\operatorname{Perf}(F)$.
Passing to Selmer $K$-theory and pairing with $\mathbf j$, the two adjoints in question are
\[
 LY_D\otimes M'\longrightarrow LY_F\xrightarrow{\mathbf j_F}\omega
 \qquad\text{and}\qquad
 LY_F\otimes LK(F)\longrightarrow LY_F\xrightarrow{\mathbf j_F}\omega,
\]
where $M'=LK(D^{\mathrm{op}})$ and the first arrows are the evaluation pairing and the action of $LK(F)$.
Precomposed with $1\otimes m_y$ and with $\Psi_y\otimes1$, respectively, they agree as \emph{maps of spectra} $LY_D\otimes LK(F)\to\omega$, by balancedness of the evaluation pairing over $LK(F)$:
\[
 \operatorname{ev}\circ
   (\operatorname{Art}_{D,p}^{L}\otimes m_y)
 \simeq
 \operatorname{ev}\circ
   ((\operatorname{Art}_{F,p}^{L}\Psi_y)\otimes\mathrm{id}_{LK(F)}).
\]
By the internal-Hom adjunction and $\rho_{D,p}=\mathbb{D}_p(m_y)$,
\[
 \rho_{D,p}\operatorname{Art}_{D,p}^{L}
 \simeq\operatorname{Art}_{F,p}^{L}\Psi_y
 \simeq\operatorname{Art}_{F,p}^{L}Lc^{\mathrm{lc}}_D.
\]
The second equivalence is Lemma~\ref{lem:paper-lca-selmer-evaluation}.
\end{proof}

\begin{proof}[Proof of Theorem~\ref{thm:paper-lca-artin}]
Lemma~\ref{lem:paper-lca-artin-compatibility}, applied with the choices and the equivalence $\rho_{D,p}$ of Lemma~\ref{lem:paper-lca-selmer-evaluation}, gives $\rho_{D,p}\operatorname{Art}_{D,p}^{L}\simeq\operatorname{Art}_{F,p}^{L}Lc^{\mathrm{lc}}_D$.
Both Selmer objects are height-one local, so precomposing with $Y_D\to LY_D$ gives~\eqref{eq:paper-lca-artin-square}.
The choices of $x$, $y$, and the commutativity homotopy for the localization sequence are fixed simultaneously; independence of the choice of $\rho_{D,p}$ or $c^{\mathrm{lc}}_D$ is not asserted \cite[Definition 4.2 and Remarks 4.3--4.4]{ClausenArtin}.

We can check even more directly that the restriction of $\operatorname{Art}_{D,p}$ to $\fib(c^{\mathrm{lc}}_D)$ vanishes.
By Theorem~\ref{thm:paper-global-fiber-gem}, $Z^\wedge_p$ is a finite GEM with finite $p$-group coefficients.
Height-one localization is unchanged by $p$-completion and annihilates each Eilenberg--Mac Lane factor, so $L(\Sigma Z)=0$.
Consequently,
\[
 \operatorname{Map}_{\mathrm{Sp}}(\Sigma Z,\mathsf S_{D,p})
 \simeq
 \operatorname{Map}_{\mathrm{Sp}}(L\Sigma Z,\mathsf S_{D,p})
 \simeq *.
\]
The restriction of $\operatorname{Art}_{D,p}$ to $\fib(c^{\mathrm{lc}}_D)$ therefore vanishes, even independently of the choice of splitting.

Finally, since $\pi_0X_D=0$, the fiber sequence~\eqref{eq:paper-lca-positive-row} gives $\pi_1Y_D=\coker\bigl(K_1(D)\to\prod'_vK_1(D_v)\bigr)$, and similarly for $F$, so $\pi_1(c^{\mathrm{lc}}_D)$ is the map induced by $(h_D)_1$.
If there is real ramification, we first form the integral cokernel; since the target has $p$ odd, $2$ is invertible there and the formula extends uniquely to $\Z[1/2]$ coefficients.
Here $\nu_{1,D_v}=\Nrd_{D_v/F_v}$ for every $v$: at finite places by Theorem~\ref{thm:paper-local-norm}, at split infinite places because the Morita comparison on $K_1$ is the determinant, and at a ramified real place because $\pi_1N_{\mathbb H/\mathbb R}=2\Nrd_{\mathbb H/\mathbb R}$ on $K_1(\mathbb H)$.
Hence the case $j=1$ of~\eqref{eq:paper-lca-all-positive-norms} in Section~\ref{subsec:paper-lca-comparison} gives $\pi_1(c^{\mathrm{lc}}_D)([(a_v)_v])=[(\Nrd_{D_v/F_v}(a_v))_v]$.
On the field side, the class of an id\`ele in $\pi_1Y_F$ is represented by multiplication on $\mathbb A_F$, which is the description used in Clausen's comparison map \cite[Proposition 3.28(1) and Corollary 3.29(1)]{ClausenArtin}.
With $e_{F,p}$ defined as the negative of the Pontryagin-dual edge map, \cite[Theorem 0.9]{ClausenArtin} identifies $e_{F,p}^{-1}\circ\pi_1(\operatorname{Art}_{F,p})$ with $\operatorname{rec}_{F,p}$, as normalized in the theorem, on id\`ele classes \cite[Definition 2.37 and Lemma 4.5]{ClausenArtin}.
Combining this with~\eqref{eq:paper-lca-artin-square}, we obtain~\eqref{eq:paper-lca-reduced-artin}.
\end{proof}

\appendix
\section{Frobenius and comparison with finite coefficients}\label{app:second-frobenius}
This appendix treats the characteristic-primary part over fields $F$ of characteristic $p$ with $[F:F^p]=p$, which include local and global function fields; such a field $F$ is fixed throughout.
Its main result, Proposition~\ref{prop:app-frobenius}, compares twisted and ordinary cycle complexes in weights $w\ge1$ when $D$ has $p$-power period, and shows that $K_i(D;\Z/p^a)=0$ for $i\ge2$ for every central division algebra $D$; this is the input used at the characteristic prime in Sections~\ref{sec:paper-local-norm}, \ref{sec:paper-global-weight-one}, \ref{sec:global-chow}, and~\ref{sec:global-spectrum}.

The comparison of cycle complexes uses the coordinate Frobenius on algebraic simplices, whose pullback produces multiplicities divisible by the index.
The $K$-theoretic consequence also requires the support tower over an imperfect field, which we establish in Lemma~\ref{lem:app-twisted-tower}.
For a central division algebra $D$ over $F$ and $w\ge0$, by the reduced rank (Definition~\ref{def:paper-twisted-chow}), the reduced cycle map $\rho$ identifies $z^w(F,\bullet;D)$ with the subcomplex of $z^w(F,\bullet)$ consisting of the cycles whose coefficient at each generator $[W]$ is divisible by $\ind(D_{\kappa(W)})$.

\begin{lemma}\label{lem:app-frobenius-homotopy}
Let $D$ be a central division algebra over $F$ whose period $q=p^c$ is a power of $p$, and for $n\ge0$ let $f_n:\Delta^n\to\Delta^n$ be the finite flat map that raises every barycentric coordinate to the $q$th power.
In each weight $w\ge1$, the pullbacks $f_n^*$ send the ordinary cycle complex $z^w(F,\bullet)$ into the subcomplex $\rho(z^w(F,\bullet;D))$ of cycles whose coefficient at each generator $[W]$ is divisible by $\ind(D_{\kappa(W)})$.
The identity of the ordinary cycle complex is simplicially homotopic to $f_\bullet^*$ through maps preserving this subcomplex.
\end{lemma}
\begin{proof}
Let $W$ be a generator of codimension $w$, let $L$ be the reduced function field of the inverse image of $W$, and let $e$ be the multiplicity of this inverse image.
The map $f_n$ is finite, flat, and radicial of degree $q^n$, and the function field of $W$ is identified with $FL^q$; hence $e[L:FL^q]=q^n$.
Since $L$ has transcendence degree $n-w$ over $F$ and $[F:F^p]=p$, invariance of the $p$-rank under finite extensions gives $[L:L^q]=q^{1+n-w}$.
The $q$th-power isomorphism carries $LF^{1/q}/L$ to $FL^q/L^q$, so $[L:FL^q]=[L:L^q]/[LF^{1/q}:L]$.
Combining these equalities, we obtain
\[
 e=q^{w-1}[LF^{1/q}:L].
\]
The last degree is a multiple of $\ind(D_L)$: under the $q$th-power isomorphism $F^{1/q}\to F$, restriction acts on Brauer classes by multiplication by $q$, since it raises the defining $2$-cocycles to the $q$th power, so $LF^{1/q}$ splits $D_L$.
Since $w\ge1$, we obtain $\ind(D_L)\mid e$, so $f_n^*$ sends the ordinary cycle complex into the subcomplex with index-divisible coefficients.

For a monotone map $\alpha:[n]\to[1]$, put $s=\sum_{\alpha(i)=1}t_i$ and $h_\alpha(t)_j=(1-s)t_j+st_j^q$.
The sum of the image coordinates is $(1-s)+s\sum_jt_j^q=(1-s)+s(\sum_jt_j)^q=1$.
On a geometric fiber, the equation $\sum_{\alpha(i)=1}h_\alpha(t)_i=s-s^2+s^{q+1}$ leaves only finitely many possibilities for $s$, and for each fixed $s$ only finitely many for each coordinate.
Thus $h_\alpha$ and, by the same calculation, its restrictions to every face are quasi-finite maps between smooth simplices of the same dimension; since the sources are Cohen--Macaulay and the targets are regular, they are flat by miracle flatness \cite[Tag 00R4]{StacksProject}.
Flatness on each face and the composition formula for Tor lengths of properly intersecting cycles along Cartier faces give the compatibility of $h_\alpha^*$ with the face maps \cite[\S\S5.6, 5.9]{KahnLevine2010}.
The maps commute with the sums of coordinates occurring in degeneracies because the $q$th-power map is additive.
Thus $h_\alpha^*$ induces a simplicial homotopy between the identity and $f_\bullet^*$.
Since the index over every residue field after flat pullback divides the original index, $h_\alpha^*$ preserves the subcomplex with index-divisible coefficients.
\end{proof}

\begin{lemma}\label{lem:app-twisted-tower}
Let $k$ be an infinite field, and let $A$ be a central simple algebra over $k$.
For every $a\ge0$, the $a$th layer of the homotopy coniveau tower of $K(-;A)$ at $\Spec k$ is identified, through the cycle map of \cite[\S6.2]{KahnLevine2010}, with the spectrum associated with the twisted cycle complex of $A$ in weight $a$.
In particular, the strongly convergent spectral sequence of \cite[Remark 6.1.2, Corollary 6.1.4]{KahnLevine2010} is available over $k$ without assuming that $k$ is perfect.
\end{lemma}
\begin{proof}
Perfectness enters the cited identification only through the geometric purity of \cite[Corollaries 4.2.3--4.2.4]{LevineHomotopyConiveau} and through \cite[\S5.2]{LevineHomotopyConiveau}; we replace the first by localization in twisted $G$-theory and use only the arguments of the second.
The other inputs, namely \cite[Corollary 3.2.2, Theorem 4.1.1]{LevineHomotopyConiveau}, the moving lemma \cite[Theorem 3.3.5]{LevineMovingLemma}, and \cite[Corollary 6.3.3, Proposition B.5]{KahnLevine2010} together with the computation in the proof of \cite[Lemma 6.4.3]{KahnLevine2010}, require only an infinite base field.
Every residue field of a $k$-scheme contains $k$ and is therefore infinite, and the towers below over a residue field $\kappa(z)$ are taken with $\kappa(z)$ as base field.

For a $k$-scheme $X$ and $d\in\Z$, let $G_{(d)}(X;A)$ be the realization of the simplicial spectrum whose $n$-simplices are the colimit, over closed subsets $W\subset X\times\Delta^n$ with $\dim(W\cap(X\times\Delta^f))\le d+\dim\Delta^f$ for every face $\Delta^f$, of the $G$-theory of $A$-modules supported on $W$.
For $X$ of pure dimension $m$, write $G_a(X)=G_{(m-a)}(X;A)$ and $\mathrm{gr}_a=\cofib(G_{a+1}\to G_a)$.
For $b\ge0$, the \emph{$b$th outer stage} $G_a^{(b)}$ of $G_a$ is the realization of the simplicial spectrum whose $n$-simplices are the colimit, over closed subsets $W\subset X\times\Delta^n$ with $\operatorname{codim}_{X\times\Delta^f}(W\cap(X\times\Delta^f))\ge b$ for every face $\Delta^f$, of $\fib\bigl(G_a(X\times\Delta^n)\to G_a((X\times\Delta^n)\setminus W)\bigr)$; the support conditions inside $G_a$ and the outer support $W$ are imposed separately.
More generally, for a functor $E$ to which this outer construction applies, we write $E^{(b)}$ for its $b$th outer stage and $E^{(b/b+1)}=\cofib(E^{(b+1)}\to E^{(b)})$ for its \emph{$b$th outer layer}; in particular $E^{(0/1)}=\cofib(E^{(1)}\to E^{(0)})$.
We put $Q_{a,b}=G_a^{(b/b+1)}$.
For a closed subscheme $Z\subset X$ of pure codimension $c$ in $X$, \cite[Corollary 6.3.3]{KahnLevine2010} gives $(G_a)_Z(X)\simeq G_{a-c}(Z)$.
Put $Z_n=Z\times\Delta^n$; taking the fiber of restriction twice, the outer support becomes $W\cap Z_n$, which has codimension at least $b-c$ on each face of $Z_n$, and conversely every such support on $Z_n$ arises from some $W$.
From now on, let $X$ be smooth of pure dimension over $k$; the closed subscheme $Z$ need not be smooth.
Admissible supports on $U=X\setminus Z$ need not have admissible closures in $X$; we first compute the double fiber with the supports inherited from $X$ and then, after realization, compare with all admissible supports on $U$ by \cite[Corollary 3.2.2]{LevineHomotopyConiveau} applied to $E=G_a$.
Then $G_a$ is a coniveau stage of $K(-;A)$ and the excision of \cite[Theorem 4.1.1]{LevineHomotopyConiveau} applies.
Applying the same localization to each outer support $W\subset Z\times\Delta^n$ shifts the codimension conditions on every face from $(a,b)$ to $(a-c,b-c)$; closed pushforward and pullback to faces are compatible because the product squares obtained by pulling back $Z\times\Delta^n$ to the faces of the simplex are Tor-independent.
Filtering by codimension, measured in $X$, and passing to generic points, the supported term of the $b$th outer layer at a point $z$ of codimension $c$ is the term with indices $(a-c,b-c)$ over $\kappa(z)$.
To pass to generic points, we remove closed subsets of higher codimension from the ambient smooth $X$, take differences of closed supports, and pass to the limit over open subsets; this does not require the closed strata to be smooth.
This replaces the geometric purity used in \cite[Corollaries 4.2.3--4.2.4]{LevineHomotopyConiveau}.

We claim that the outer layer $Q_{a,b}$ vanishes for $0\le b<a$, by induction on $a$.
At a point $z$ of codimension $c>0$, the supported term of $Q_{a,b}$ is $Q_{a-c,b-c}$ over $\kappa(z)$.
It vanishes if $b-c<0$, since the support conditions of the two adjacent outer stages then coincide, and by induction if $b-c\ge0$, since then $0\le b-c<a-c<a$.
Shrinking the dimension conditions on finitely many faces to a dense open subset justifies the passage to generic points.
For $b=0$, the semilocal-simplex computation in the proof of \cite[Lemma 5.2.1]{LevineHomotopyConiveau}, which uses supports avoiding the vertices and the moving lemma, applied to the outer $0$th layer $G_a^{(0/1)}$ with both the inner and outer face conditions imposed, gives $Q_{a,0}=0$ directly.
For $b>0$, the induction shows that $Q=Q_{a,b}$ is birational, so its supported fibers in positive codimension vanish, and $\mathbb A^1$-invariance gives $Q(K)\simeq Q^{(0/1)}(K)$.
By exactness of taking layers, the right-hand side is the cofiber of $(G_a^{(b+1)})^{(0/1)}\to(G_a^{(b)})^{(0/1)}$, and both terms vanish by the same semilocal-simplex computation, applied to $G_a^{(b)}$ and $G_a^{(b+1)}$ with both the inner and outer face conditions imposed.
This is the argument of \cite[Proposition 5.2.2]{LevineHomotopyConiveau}, with the geometric Gersten terms replaced by the supported terms above.
Separately, the supported term of $\mathrm{gr}_a$ at a point of codimension $c>a$ vanishes, since $G_{a-c}$ and $G_{a+1-c}$ are then both $G$-theory and the dimension conditions agree.

For smooth $X$, applying the outer $a$th layer to $G_a\to K(-;A)$ gives an equivalence on supports of codimension at least $a$.
Applying it to $G_a\to\mathrm{gr}_a$ also gives an equivalence, because the contribution of $G_{a+1}$ vanishes by the claim; the outer layers with $b<a$ vanish as well, since $Q_{a,b}=Q_{a+1,b}=0$.
Since $\mathrm{gr}_a$ vanishes in outer codimension at least $a+1$, the outer $a$th layer recovers $\mathrm{gr}_a$.
To compute it, note that in each simplicial degree $n$, the vanishing in higher codimension removes intersections of supports, and additivity for finitely many disjoint supports gives the direct sum of $\mathrm{gr}_0(\kappa(w);A_{\kappa(w)})$ over the points $w$ of codimension $a$ in $X\times\Delta^n$ whose closures satisfy all the face conditions.
Each term is identified with $K(\widehat\Delta^\bullet;A_{\kappa(w)})$, the $K$-theory of the semilocal simplices over $\kappa(w)$, localized for supports avoiding the vertices.
The semilocal-simplex computation in the proof of \cite[Lemma 6.4.3]{KahnLevine2010} gives $\mathrm{gr}_0(\kappa(w);A_{\kappa(w)})\simeq H\bigl(K_0(A_{\kappa(w)})\bigr)$: the other degrees vanish by the computation with the full semilocal simplex and its boundary, and $\pi_0=K_0(A_{\kappa(w)})$ because the two vertex evaluations of a semilocal $1$-simplex agree in reduced rank.
This computation uses only the $K_1$-regularity of the boundary from Vorst's theorem and \cite[Proposition B.5]{KahnLevine2010}, surjectivity of $\operatorname{GL}$ from the semilocal ring to the boundary, and injectivity on $K_0$, all of which hold over an arbitrary residue field.
We use only this computation of $HK_0$, and not the sheaf identification in that lemma, which relies on \cite[Theorem 3.2.4]{KahnLevine2010}.
The map on $\pi_0$ is restriction to the generic point, which is the cycle map of \cite[\S6.2]{KahnLevine2010}.
In summary, for smooth $X$ we have
\[
 K(-;A)^{(a/a+1)}(X)\xleftarrow{\ \sim\ }G_a^{(a/a+1)}(X)\xrightarrow{\ \sim\ }\mathrm{gr}_a^{(a/a+1)}(X)\xleftarrow{\ \sim\ }\mathrm{gr}_a(X).
\]
\end{proof}

\begin{proposition}\label{prop:app-frobenius}
Let $D$ be a central division algebra over $F$.
\begin{enumerate}
\item If the period of $D$ is a power of $p$, then in each weight $w\ge1$ the reduced cycle map from the twisted cycle complex of $D$ to the ordinary cycle complex of $F$ is a quasi-isomorphism with integral coefficients.
\item For $i\ge2$ and every $a\ge1$, one has $K_i(D;\Z/p^a)=0$, and multiplication by $p$ is an isomorphism on $K_i(D)$.
\end{enumerate}
\end{proposition}
\begin{proof}
(1) As recalled before Lemma~\ref{lem:app-frobenius-homotopy}, $\rho$ identifies the twisted cycle complex with the subcomplex $\rho(z^w(F,\bullet;D))$ of $z^w(F,\bullet)$; the identification is compatible with the face maps by \cite[Lemma 5.9.1]{KahnLevine2010}.
By Lemma~\ref{lem:app-frobenius-homotopy}, the identity of the ordinary cycle complex is simplicially homotopic to $f_\bullet^*$ through maps preserving this subcomplex, so the homotopy descends to the quotient, on which $f_\bullet^*$ induces zero.
On the alternating face complex, it gives a chain homotopy from the identity of the quotient to zero, so the quotient is acyclic.

(2) Since $[F:F^p]=p>1$, the field $F$ is imperfect and hence infinite.
The prime-to-$p$ Brauer components can be removed by restriction to a finite separable extension of degree prime to $p$ that splits them; this preserves $[F:F^p]=p$, and for comparison modulo $p$ the restriction is a split injection, since the composite of restriction and corestriction is multiplication by the degree.
By Morita equivalence we may then assume that $D$ has $p$-power period.
By Lemma~\ref{lem:app-twisted-tower}, the strongly convergent spectral sequence for the twisted cycle complexes of $D$ is available over $F$, and the connectivity estimates for the tower remain valid with Moore coefficients.
In weight $w\ge1$, part~(1) identifies the mod-$p$ layers with the commutative ones, which are concentrated in degree $w$, where they are $K_w^M(F)/p$ \cite{GeisserLevine}.
For $w\ge2$, they vanish: $K_w^M(F)/p$ embeds into $\Omega^w_F$ by the $d\log$ map \cite[Corollary (2.8)]{BlochKato1986}, and $\Omega^w_F=0$ since $[F:F^p]=p$.
The weight-zero layer is concentrated in degree zero.
Hence $K_i(D;\Z/p)=0$ for $i\ge2$.
The coefficient sequences for $\Z/p\to\Z/p^{a+1}\to\Z/p^a$ give $K_i(D;\Z/p^a)=0$ for $i\ge2$ and all $a\ge1$ by induction on $a$.
In the Bockstein sequence $0\to K_i(D)/p\to K_i(D;\Z/p)\to K_{i-1}(D)[p]\to0$, the vanishing for $i\ge2$ gives $K_i(D)/p=0$ for $i\ge2$ and $K_j(D)[p]=0$ for $j\ge1$.
Thus multiplication by $p$ is an isomorphism on $K_i(D)$ for $i\ge2$.
\end{proof}

\section{Relative \texorpdfstring{$\mathrm{TC}$}{TC} and integral lifts}\label{app:second-tc}
The proof of Theorem~\ref{thm:second-two-wild} compares a class in $K_2$ of a maximal order with $\mathrm{TC}$ at a finite truncation and then represents it by an integral class whose image in $K_2(D)$ is controlled; the proof of Theorem~\ref{thm:second-three-coker} uses the results on deep principal units and relative $K_1$.
Sections~\ref{app:second-tc-supports} and~\ref{app:second-tc-separation} provide the trace comparison at finite truncations and show that the image in $K_2(D)$ of a deep congruence kernel consists of transfers from a fixed splitting field; the remaining subsections show that, in a fixed truncation, images coming from sufficiently deep truncations are represented by integral classes.

\subsection{Supports and congruence towers}\label{app:second-tc-supports}
At finite truncations of a maximal order, we use the relative cyclotomic trace to compare integral relative $K$-classes.

\begin{lemma}\label{lem:app-nilpotent-trace}
Let $\Lambda$ be a unital associative ring in which $p$ is nilpotent, and let $I\subset\Lambda$ be a nilpotent two-sided ideal.
Then each $K_i(\Lambda,I)$ is a $p$-primary torsion group of bounded exponent, with a bound depending on $i$, and the homotopy groups of $K(\Lambda,I)$ are unchanged by $p$-completion.
The relative trace
\[
 K(\Lambda,I)^\wedge_p \xrightarrow{\ \sim\ } \mathrm{TC}(\Lambda,I)^\wedge_p
\]
is an equivalence; in particular $K_i(\Lambda,I)\simeq\pi_i\mathrm{TC}(\Lambda,I)^\wedge_p$ for all $i$.
\end{lemma}
\begin{proof}
The first assertion is Theorem A of Geisser--Hesselholt \cite{GeisserHesselholtNilpotent}, and the second follows from it.
The equivalence is the relative trace comparison for nilpotent ideals \cite{DundasGoodwillieMcCarthy}.
\end{proof}

Let $F$ be a complete discretely valued field, let $D$ be a central division algebra over $F$, and let $P=\mathcal O_D$ be its maximal order; let $J$ be the Jacobson radical of $P$, and let $\Delta=P/J$ be the residue division algebra of $D$.
Localization and d\'evissage give the fiber sequence
\[
 K(\Delta)\longrightarrow K(P)\longrightarrow K(D).
\]
Let $\pi$ be a uniformizer of the center, and put $P_s=P/\pi^sP$ and $J_s=J/\pi^sP$.
Then $J_s$ is nilpotent and $P_s/J_s=\Delta$.
If $p$ is the residue characteristic, then $p$ is nilpotent in $P_s$, so Lemma~\ref{lem:app-nilpotent-trace} applies to $(P_s,J_s)$, as well as to Iwahori and commutative truncations, and identifies relative classes at a fixed truncation integrally.
This comparison cannot be applied directly to the kernel $\pi^sP$ of $P\to P_s$.
To compare relative terms in the congruence tower with terms with supports, we use their respective continuity results \cite{DundasGoodwillieMcCarthy,GeisserHesselholtNilpotent,DundasMorrow,ClausenMathewMorrow}.

\begin{lemma}\label{lem:app-iwahori-derived-cut}
Let $R$ be a discrete valuation ring with uniformizer $\pi$, let $C$ be a maximal order finite free over $R$, and put $\Delta=C/J(C)$.
Let $H$ be the $r\times r$ Iwahori order with $C$ on and above the diagonal and $J(C)$ below it.
Let $J_H=J(H)$, and let $f$ be the last diagonal idempotent.
Let $H_s=H/\pi^sH$ and $C_s=C/\pi^sC$.
The transition map $Z_s\to Z_j$ between the relative corner cofibers
\[
 Z_s=\cofib\bigl(K(C_s,\Delta)
       \longrightarrow K(H_s,H/J_H)\bigr)
\]
is nullhomotopic whenever $s>j$.
Before truncation, the relative corner map $K(C,\Delta)\to K(H,H/J_H)$ is an equivalence.
\end{lemma}
\begin{proof}
The proof uses only that $C$ is finite free over $R$ and that $\pi C\subset J(C)$.
In $Z_s$, the map $K(C_s)\to K(H_s)$ is $-\otimes_{C_s}fH_s$ for right modules, and $K(\Delta)\to K(\Delta^r)$ is the inclusion of the last factor, where $H/J_H=\Delta^r$.
First consider the orders before truncation.
The right $C$-module $Hf$ over the corner $fHf=C$ is free, so $Hf\otimes^{\mathbb L}_CfH=Hf\otimes_CfH$, and multiplication $Hf\otimes_C fH\to H$ is injective.
The cokernel $\mathcal E$ of this multiplication map is the upper triangular algebra over $\Delta$ on the remaining $r-1$ vertices, so the derived quotient $H/^{\mathbb L}(HfH)$ by $f$ is equivalent to $\mathcal E$.
By its universal property of making $f$ a zero object, the derived quotient commutes with derived base change along the center.
Hence
\[
 H_s/^{\mathbb L}(H_sfH_s)
 \simeq\mathcal E_s,\qquad\mathcal E_s=\mathcal E\otimes_R^{\mathbb L}R/\pi^sR.
\]
Since $H,C$ are flat over $R$, the corner on the left is the actual ring $C_s$.
For a differential graded algebra $\mathcal E'$, such as $\mathcal E$ and $\mathcal E_s$, we write $K(\mathcal E')$ for the $K$-theory of the category $\operatorname{Perf}(\mathcal E')$ of perfect $\mathcal E'$-modules.
This identification is compatible with the transition maps \cite{BraunChuangLazarev}.

The two-term Koszul model of $R/\pi^sR$ gives $\mathcal E\oplus\mathcal E\epsilon_s$ as an $\mathcal E$-algebra model of $\mathcal E_s$, with $|\epsilon_s|=1$, $\epsilon_s^2=0$, and $d\epsilon_s=0$.
For $s\ge j$, the transition is given by $\epsilon_s\mapsto\pi^{s-j}\epsilon_j$.
Since $\pi\mathcal E=0$, for $s>j$ the transition factors as
\[
 \mathcal E_s\longrightarrow\mathcal E\longrightarrow\mathcal E_j,
\]
where the first map is the augmentation $\epsilon_s\mapsto0$ and the second is the inclusion in degree zero, a section of the augmentation of $\mathcal E_j$.
Thus $K(\mathcal E_s,\mathcal E)\to K(\mathcal E_j,\mathcal E)$ is nullhomotopic.

We compare the localization sequence for the derived quotient with the corner localization sequence for the residue ring $H/J_H=\Delta^r$.
The $K$-theory of an upper triangular algebra is the $K$-theory of its diagonal, so $K(\mathcal E)\simeq\bigoplus_{i=1}^{r-1}K(\Delta)$.
Placing $K(C_s)\to K(H_s)\to K(\mathcal E_s)$ over $K(\Delta)\to K(\Delta^r)\to K(\Delta^{r-1})$, with the right vertical map the augmentation followed by the diagonal projection, the vertical fibers give $Z_s\simeq K(\mathcal E_s,\mathcal E)$ compatibly with the transition maps.
This gives the required nullhomotopy.
Before truncation, the same comparison has zero cofiber.
\end{proof}

\begin{lemma}
\label{lem:app-nilpotent-etale-descent}
Let $R$ be a commutative Henselian local ring with residue field $k$ of characteristic $p$, let $\Lambda$ be a unital associative $R$-algebra, and let $I\subset\Lambda$ be a nilpotent ideal.
Let $a\ge1$.
For a finite \'etale $k$-algebra $U$, let $R_U$ be its Henselian lift, and write $\Lambda_U=\Lambda\otimes_RR_U$ and $I_U=I\otimes_RR_U$.
Let $\mathcal X$ be the sheaf of spectra defined by finite \'etale base change from
\[
 U\longmapsto K(\Lambda_U,I_U)/p^a.
\]
In particular, $\mathcal X(k)=K(\Lambda,I)/p^a$.
Let $\mathcal X_{\bar k}$ be the geometric stalk of $\mathcal X$.
Then there is a natural short exact sequence
\[
 0\longrightarrow H^1(k,\pi_3\mathcal X_{\bar k})
 \longrightarrow\pi_2\mathcal X(k)
 \longrightarrow(\pi_2\mathcal X_{\bar k})^{G_k}
 \longrightarrow0.
\]
\end{lemma}
\begin{proof}
Hesselholt--Larsen--Lindenstrauss prove \'etale base change for an arbitrary associative $R$-algebra \cite[Proposition 3.1, Addendum 3.2]{HesselholtLarsenLindenstrauss2020}.
Their Corollary~3.3 gives, for a finite Galois \'etale extension $R\to R'$ with group $G$,
\[
 \mathrm{TR}^n_i(\Lambda;p)\simeq
 H^0(G,\mathrm{TR}^n_i(\Lambda_{R'};p)),\qquad
 H^{q}(G,\mathrm{TR}^n_i(\Lambda_{R'};p))=0\quad(q>0).
\]
These results of Hesselholt--Larsen--Lindenstrauss assume neither that the residue field is finite nor that $p\ne2$.
We compare finite Postnikov truncations and then take the inverse limit to obtain $\mathrm{TR}^n(\Lambda;p)\simeq\mathrm{TR}^n(\Lambda_{R'};p)^{hG}$.
Taking the limits along restriction and Frobenius, we obtain $\mathrm{TC}(\Lambda;p)\simeq\mathrm{TC}(\Lambda_{R'};p)^{hG}$.
Since homotopy fixed points preserve finite limits and finite colimits of spectra, this descent also holds for relative terms and Moore coefficients.

Since $I$ is nilpotent, the relative trace comparison gives $K(\Lambda,I)/p^a\simeq\mathrm{TC}(\Lambda,I)^\wedge_p/p^a$, so the finite \'etale descent above also holds for $\mathcal X$.
Commutation of $\mathrm{TC}$ modulo $p^a$ with filtered colimits implies that $\mathcal X_{\bar k}\simeq\operatorname*{colim}_{k'}\mathcal X(k')$ and that each homotopy group of $\mathcal X_{\bar k}$ is a discrete $G_k$-module, where $k'$ runs over the finite separable extensions of $k$ in $\bar k$ \cite{DundasGoodwillieMcCarthy}, \cite[Theorem 2.7]{ClausenMathewMorrow}.
The homotopy groups of $\mathcal X$ and of $\mathcal X_{\bar k}$ are killed by $p^{2a}$, being extensions of $\pi_jK/p^a$ by $\pi_{j-1}K[p^a]$; in particular they are $p$-primary.
By \cite[Corollary 3.3]{HesselholtLarsenLindenstrauss2020} the higher group cohomology vanishes; the Eilenberg--Mac Lane sheaves of each $\mathrm{TR}^n$ are sheaves, so finite Postnikov truncations are hypersheaves, and connectivity, limits defining $\mathrm{TC}$, relative fibers and Moore coefficients preserve this.
Since $k$ has characteristic $p$, $\cd_p(k)\le1$ \cite[II, \S2.2, Proposition 3]{SerreGalois}, so the descent spectral sequence for each finite Postnikov truncation has only cohomological degrees zero and one.
The truncated tails lose at most one degree of connectivity under descent, so their connectivity still tends to infinity, and the continuous descent spectral sequence converges to give the displayed short exact sequence.
\end{proof}

\subsection{Congruence separation for a fixed splitting field}\label{app:second-tc-separation}
Next we show that the image in $K_2(D)$ of a sufficiently deep congruence kernel of $K_2(P)$ consists of transfers from a fixed splitting field.
This allows us to work at a finite truncation in the main text.

\begin{lemma}\label{lem:app-congruence-separation}
Let $F$ be a complete discretely valued field, let $D/F$ be a central division algebra, and let $M/F$ be a finite separable splitting field.
For all sufficiently large $s$,
\[
 \im\bigl(\ker[K_2(P)\to K_2(P_s)]\longrightarrow K_2(D)\bigr)
 \ \subseteq\ \theta_M K_2(M).
\]
\end{lemma}

\begin{proof}
The local rings $P$ and $P_s$ have Bass stable rank $1$.
Indeed, at least one entry of a left-unimodular pair is a unit; if the first entry is not a unit, adding the second makes it a unit.
The noncommutative stability theorem of Suslin--Tulenbaev therefore allows us to fix a single $n\ge5$, independent of $s$, with $K_2(P)\simeq\ker(\operatorname{St}_n(P)\to E_n(P))$ and with the analogous isomorphisms for $P_s$ and $D$ \cite{SuslinTulenbaev}.

Enlarging $n$ if necessary, which preserves these stability isomorphisms, we may assume that Corollary~\ref{cor:paper-congruence-splitting} applies with $E=M$.
Since stabilization identifies the kernels of both $\operatorname{St}_n(D)\to E_n(D)$ and the stable universal central extension with $K_2(D)$, the group $\operatorname{St}_n(D)$ is the pullback of the stable extension to $E_n(D)$.
Put
\[
 \Gamma_m=E_n(D)\cap(1+\pi^m\operatorname{Mat}_n(P)).
\]
By Corollary~\ref{cor:paper-congruence-splitting}, for sufficiently large $m$ there is a group-homomorphic section $\sigma$ over $\Gamma_m$ of $\operatorname{St}_n(D)/\theta_MK_2(M)\to E_n(D)$.

We take $s\ge 2m$.
By the presentation of the Steinberg group, the kernel of $\operatorname{St}_n(P)\to\operatorname{St}_n(P_s)$ is the normal closure of the root elements $x_{ij}(t)$ with $t\in\pi^sP$.
Indeed, in the quotient by this normal closure the coefficients can be written as elements of $P_s$, and the three types of Steinberg relations over $P_s$ can be checked using arbitrary lifts.
For distinct $i,j,k$,
\[
 x_{ij}(t)=[x_{ik}(\pi^m),x_{kj}(\pi^{-m}t)].
\]
Let $\gamma\in\operatorname{St}_n(P)$, and let $u_1$ and $u_2$ be the matrix images of $\gamma x_{ik}(\pi^m)\gamma^{-1}$ and $\gamma x_{kj}(\pi^{-m}t)\gamma^{-1}$.
Since $s\ge2m$, the matrix images of $x_{ik}(\pi^m)$ and $x_{kj}(\pi^{-m}t)$ lie in $\Gamma_m$, and conjugation by the image of $\gamma$ preserves $E_n(D)\cap(1+\pi^m\operatorname{Mat}_n(P))$, so $u_1,u_2\in\Gamma_m$.
In a central extension, changing either lift by a central element leaves the commutator unchanged.
Therefore the image of $\gamma x_{ij}(t)\gamma^{-1}$ in $\operatorname{St}_n(D)/\theta_MK_2(M)$ is $[\sigma(u_1),\sigma(u_2)]\in\sigma(\Gamma_m)$.

Every element of $\ker[K_2(P)\to K_2(P_s)]$ belongs to this normal closure and has trivial matrix image.
Its image in $\operatorname{St}_n(D)/\theta_MK_2(M)$ therefore lies in $\sigma(\Gamma_m)$ and maps to the identity in $E_n(D)$.
By the section property this image is trivial, so the generic-fiber image of the element lies in $\theta_MK_2(M)$.
\end{proof}

\subsection{Deep principal units and relative classes at finite truncations}\label{app:second-deep-units}
Lemma~\ref{lem:app-deep-unit-norm} controls the reduced norm on deep principal units.

For relative classes $y,z$, the relative trace comparison for nilpotent ideals (Lemma~\ref{lem:app-nilpotent-trace}) shows that equality of their images in relative $\mathrm{TC}^\wedge_p$ at a fixed truncation implies equality of their integral $K_2$ images at that truncation.
For corrections from commutative subfields and degree-zero cyclic orbits, we use Lemma~\ref{lem:app-orbit-stable-image} and Corollary~\ref{cor:app-commutative-cut-representative} below.

\begin{lemma}\label{lem:app-deep-unit-norm}
Let $T$ be a complete discretely valued field, let $B/T$ be a central division algebra, put $C=\mathcal O_B$, and let $\pi$ be a central element of positive valuation in $T$.
Let $\mathcal U_s(C)=\im(1+\pi^sC\to K_1(C))$.
There are integers $s_0$ and $t$ such that for all $s\ge s_0$ the reduced norm is injective on $\mathcal U_s(C)$ and
\[
 1+\pi^{s+t}\mathcal O_T\subseteq\Nrd_1(\mathcal U_s(C))
       \subseteq1+\pi^s\mathcal O_T.
\]
Moreover, for each $j$, every sufficiently deep unit of reduced norm one can be written as $\prod_i[b_i,c_i]$, with fixed $b_i\in C^\times$ and $c_i\in1+\pi^jC$.
These depths can be chosen uniformly over finite unramified extensions for which the base change of $C$ is the maximal order of a division algebra of the same degree.
When $B/T$ is a central simple algebra and $C\subset B$ is a finite $\mathcal O_T$-subalgebra satisfying $C\otimes_{\mathcal O_T}T=B$, we still have injectivity, the two inclusions for the norm image, and the commutator representation.
\end{lemma}
\begin{proof}
Fix integral units $b_i$ spanning $B$ over $T$.
The commutator product has differential $\sum_i(\operatorname{Ad}b_i-1)c_i$ at the identity.
By nondegeneracy of the reduced trace pairing, the orthogonal complement of the image of the differential is the common centralizer $T$ of the $b_i$, so the differential surjects onto $\ker\operatorname{Trd}_{B/T}$.
After clearing denominators, we apply the inverse function theorem with a fixed right inverse and realize every sufficiently deep element of reduced norm one using $c_i$ of any prescribed depth.
Since commutators vanish in $K_1(C)$, every sufficiently deep element of reduced norm one has trivial class in $K_1(C)$; this proves the injectivity of the reduced norm on $\mathcal U_s(C)$.
We obtain the left inclusion from Hensel's lemma in a direction with nonzero reduced trace, and the right inclusion from the reduced characteristic polynomial with integral coefficients.
Since unramified extension does not change the valuations of the denominators of the fixed Jacobian and of its right inverse, one depth works for all such base changes.
For a general finite subalgebra $C$, its units still span $B$ over $T$, since $1+\pi^ma\in C^\times$ for every lattice basis element $a$ and all sufficiently large $m$.
The same reduced trace pairing and fixed-Jacobian argument therefore applies.
\end{proof}

\subsection{Stable images of relative finite-group orbits}
At a fixed truncation, we show that the image on $\pi_2$ of the relative homotopy orbits of a sufficiently deep truncation agrees with that of the integral relative homotopy orbits.

Throughout this subsection, let $F$ be a complete discretely valued field, let $T/F$ be a finite Galois extension, put $G=\Gal(T/F)$, and let $\pi\in F$ have positive valuation.
Let $B/T$ be a central division algebra and put $C=\mathcal O_B$.
Let $H$ be either $C$ or the Iwahori order used in the main text with corner $C$, and write $J_H=J(H)$, $\bar H=H/J_H$, and $H_s=H/\pi^sH$.
We assume that $G$ acts semilinearly, preserving the order, its residue, and its truncations, and inducing the specified Galois action on the generic center $T$.
Through the corresponding action on the module category, we also use weak actions represented by ring automorphisms and coherence isomorphisms given by units.
For an Iwahori cover, we likewise use the actual action on the algebra formed from the direct sum of all vertices, on which vertex translation acts; in that case $H$ denotes this algebra.
In the homogeneous-unit type of the main text, identifying the vertices by powers of the unit $x$ gives $H\simeq M_r(C)$, and the weak action transported to $C$ by Morita equivalence falls under the case $H=C$.
This Morita equivalence is compatible with the residue and with each truncation, and it carries vertex translation to the weak action on the corner; at the generic point the cover side is $M_r(B)$ and the corner side is $B$.
For both types, $H_T=H\otimes_{\mathcal O_T}T\simeq M_r(B)$ at the generic point, and $N_{T/F}\Nrd_{H_T/T}$ is $G$-invariant on $K_1(H_T)$, so it sends the boundary used below, which vanishes in the coinvariants, to $1$.

\begin{lemma}
\label{lem:app-orbit-stable-image}
Let $Y=K(H,\bar H)$ and $Y_s=K(H_s,\bar H)$.
For every $j$ and all sufficiently large $s$,
\[
 \im[\pi_2((Y_s)_{hG})\to\pi_2((Y_j)_{hG})]
 =\im[\pi_2(Y_{hG})\to\pi_2((Y_j)_{hG})]
.
\]
When $G=1$, the same conclusion holds for the absolute spectra $Y=K(H)$ and $Y_s=K(H_s)$.

If the residue characteristic is $p$, then for every $j$ and every $x\in\pi_2((\mathrm{TC}(H,\bar H)^\wedge_p)_{hG})$ there is $y\in\pi_2(Y_{hG})$ such that $x$ and $\operatorname{tr}(y)$ have the same image in $\pi_2((\mathrm{TC}(H_j,\bar H)^\wedge_p)_{hG})\simeq\pi_2((Y_j)_{hG})$.
\end{lemma}

\begin{proof}
The order $H$ has global dimension one, and the rings $H_s$ and $\bar H$ are Artinian, so the negative $K$-groups of $H$, $H_s$, and $\bar H$ vanish \cite[III, \S4 and Exercise 4.3]{WeibelKbook}.
Since $H$ is complete, idempotents lift modulo $\pi^sH$ and modulo $J_H$, which gives the isomorphisms on $K_0$; together with surjectivity on $K_1$, they imply that $K(H,\pi^sH)$, $Y$, and $Y_s$ are concentrated in degrees at least one.
Thus, taking homotopy orbits of the $G$-equivariant fiber sequence $K(H,\pi^sH)\to Y\to Y_s$, we see that the boundary of each element of $\pi_2((Y_s)_{hG})$ belongs to $K_1(H,\pi^sH)_G$ and can be represented by $u\in1+\pi^sH$ \cite[III, Lemma 2.4(1)]{WeibelKbook}.
This boundary maps to zero in $(\pi_1Y)_G$, so the $G$-invariant norm $\lambda=N_{T/F}\circ\Nrd_{H_T/T}$ gives $\lambda(u)=1$.

We apply the argument of Lemma~\ref{lem:app-deep-unit-norm} over $F$.
We choose fixed $b_i\in H^\times$ spanning $H_T$ over $T$, and consider the map to $\ker\lambda$
\[
 \Phi((c_i),(a_\sigma))=\prod_i[b_i,c_i]\,
                   \prod_{\sigma\in G}\sigma(a_\sigma)a_\sigma^{-1}.
\]
Its differential is $\sum_i(\operatorname{Ad}b_i-1)x_i+\sum_\sigma(\sigma-1)y_\sigma$.
For the nondegenerate pairing $(x,y)\mapsto\Tr_{T/F}\operatorname{Trd}_{H_T/T}(xy)$, the orthogonal complement of the image of this differential is $T$ after the commutator terms, and then $T^G=F$ after the $G$-terms.
The image is therefore $\ker(\Tr_{T/F}\operatorname{Trd}_{H_T/T})$.
Since $d\lambda=\Tr_{T/F}\circ\operatorname{Trd}_{H_T/T}$ and both traces are surjective, $d\lambda$ is nonzero even when the characteristic divides the degree, and $\ker\lambda$ is smooth at the identity.
Using the same inverse function theorem with a fixed Jacobian, we see that, for every $j$ and all sufficiently large $s$, every $u\in1+\pi^sH$ with $\lambda(u)=1$ can be written as $u=\Phi((c_i),(a_\sigma))$ with $c_i,a_\sigma\in1+\pi^jH$.

Taking $b=b_i$ and $x=(c_i-1)b_i^{-1}$ in the relative elementary relation $[1+bx]=[1+xb]$ \cite[III, Exercise 2.4]{WeibelKbook} shows that $[b_i,c_i]$ vanishes in $K_1(H,\pi^jH)$.
The class of $\sigma(a_\sigma)a_\sigma^{-1}$ also vanishes in $K_1(H,\pi^jH)_G$, so the original boundary is zero at truncation $j$.
By exactness, the image at truncation $j$ of the given element of $\pi_2((Y_s)_{hG})$ lifts from the orbits of $Y$.
The depth depends only on $j$ and the fixed data.
The reverse inclusion holds because the map from $Y$ to $Y_j$ factors through $Y_s$.

For $G=1$, the same argument applies with $Y=K(H)$ and $Y_s=K(H_s)$, starting from $\Nrd_1(u)=1$, since the boundary maps to zero in $K_1(H)$.
The element $\pi$ may be any central element of positive valuation.

In residue characteristic $p$, the relative $K$-groups of a finite truncation are $p$-primary torsion groups of bounded exponent in each degree by Lemma~\ref{lem:app-nilpotent-trace}.
Only finitely many terms contribute in each orbit degree, so $p$-completion does not change the homotopy groups of $(Y_s)_{hG}$.
We apply the trace equivalence of the same lemma at a deep truncation, and then use the equality of images just proved to obtain an integral representative at the fixed truncation.
\end{proof}

\begin{corollary}
\label{cor:app-integral-cut-stable}
Let $T$ be a complete discretely valued field, let $B/T$ be a central division algebra, put $C=\mathcal O_B$ and $\Delta=C/J(C)$, and let $\pi\in T$ be any central element of positive valuation.
Let $C_s=C/\pi^sC$.
For every $j$ and all sufficiently large $s$,
\[
 \im[K_2(C_s)\to K_2(C_j)]
 =\im[K_2(C)\to K_2(C_j)].
\]
The same equality holds for the relative groups $K_2(C_s,\Delta)$ and $K_2(C,\Delta)$.
\end{corollary}
\begin{proof}
We apply the absolute and relative conclusions of Lemma~\ref{lem:app-orbit-stable-image} with $G=1$ and $H=C$.
\end{proof}

\subsection{Integral representatives over commutative valuation rings}
Finally, we record the commutative case of Lemma~\ref{lem:app-orbit-stable-image}, which is used for corrections by commutative subfields.

\begin{corollary}
\label{cor:app-commutative-cut-representative}
Let $R$ be a complete discrete valuation ring with residue field $\kappa$ of characteristic $p$.
For $j\ge1$, let $R_j$ be the quotient of $R$ by the $j$th power of a fixed uniformizer.
For every $j\ge1$ and every $x\in\pi_2\mathrm{TC}(R,\kappa)^\wedge_p$ there is $y\in K_2(R,\kappa)$ such that $x$ and $\operatorname{tr}(y)$ have the same image in $\pi_2\mathrm{TC}(R_j,\kappa)^\wedge_p\simeq K_2(R_j,\kappa)$.
\end{corollary}
\begin{proof}
We apply the last conclusion of Lemma~\ref{lem:app-orbit-stable-image} with $G=1$ and $H=R$.
\end{proof}

\section{Common input for zero-cycles with coefficients and splitting transfers}\label{app:review-cycle-transfer}
The object of this appendix is to prove Corollary~\ref{cor:review-small-cd-splitting-transfer} and Lemma~\ref{lem:review-p-rank-two-frobenius-norm}, which together give Proposition~\ref{prop:review-splitting-transfer-norm-iso}: over a two-dimensional local field, the reduced norm maps the subgroup generated by splitting transfers isomorphically onto $K_2(F)$.
Away from the characteristic, we compute zero-cycles with Milnor $K$-coefficients on products of Severi--Brauer varieties; at the characteristic prime, we use the Frobenius.

Throughout this appendix, for a smooth variety $X$ over a field $k$, $M(X)$ denotes its effective motive with integral coefficients, and $\operatorname{Hom}$ denotes morphisms in the category of effective motives.
For a central simple algebra $B/k$, the sheaf with transfers $\Z_B$ takes the value $\ind(B_{k(U)})\Z\subset\Z$ on a smooth integral variety $U/k$, and $\Z_B(q)=\Z_B\otimes\Z(q)$ denotes its Tate twist.
We make motivic calculations with $\Z_{(\ell)}$ coefficients for a fixed prime $\ell$ different from the characteristic; if $\operatorname{char}k=p>0$, then $p$ is invertible in these coefficients.
The Frobenius calculation in Lemma~\ref{lem:review-p-rank-two-frobenius-norm} uses integral coefficients.
Recall (Definition~\ref{def:paper-transfer}) that $I_F(D)=\sum_E\theta_EK_2(E)$, the sum over all finite splitting fields $E/F$, separable or not, and that $\Nrd_2\theta_E=N_{E/F}$.

\begin{definition}\label{def:review-rost-zero-cycles}
For a finite-type $k$-scheme $Z$ and $j\ge0$, let $Z_{(a)}$ denote the set of points whose closures have dimension $a$.
The low-degree terms of the Rost complex are
\[
 C_a(Z,K_j^M)=\bigoplus_{z\in Z_{(a)}}K_{j+a}^M(\kappa(z))\quad(a=0,1),
 \qquad A_0(Z,K_j^M)=\coker(C_1\xrightarrow{\partial}C_0).
\]
The differential is the sum of the residue maps on the normalizations of the closures of one-dimensional points, followed by norms on residue fields.
Proper pushforward is defined by the residue-field norm when the image has the same dimension and by zero when the dimension drops; it commutes with the differential.
In particular, when $Z$ is proper over $k$, the structure map gives the degree map $\deg_Z:A_0(Z,K_j^M)\to K_j^M(k)$ \cite{Rost}.
\end{definition}

\begin{lemma}
\label{lem:review-proper-cycle-pushforward}
Let $f:X\to Y$ be a proper morphism of finite-type $k$-schemes, and let $j\ge0$.
Suppose that, for every closed point $y\in Y_{(0)}$, the map
\[
 A_0(X_y,K_j^M)\xrightarrow{\deg}K_j^M(\kappa(y))
\]
is an isomorphism, and that $A_0(X_\eta,K_{j+1}^M)\to K_{j+1}^M(\kappa(\eta))$ is surjective for every $\eta\in Y_{(1)}$ (the fibers $X_y$ and $X_\eta$ being regarded as schemes over $\kappa(y)$, resp.\ $\kappa(\eta)$).
Then $f_*:A_0(X,K_j^M)\to A_0(Y,K_j^M)$ is an isomorphism.
The same holds when both the hypotheses on the fibers and the conclusion are taken with $\Z_{(\ell)}$ coefficients.
\end{lemma}
\begin{proof}
Surjectivity follows by lifting the coefficient at each closed point of a zero-cycle from the corresponding fiber.
For injectivity, take $z\in C_0(X,K_j^M)$ and write $f_*z=\partial b$.
By assumption, the coefficient of $b$ at each one-dimensional point $\eta$ lifts from $K_{j+1}^M$ coefficients supported on finitely many closed points of $X_\eta$.
We take their closures in $X$, with the same coefficients over their function fields, to obtain a chain $\widetilde b\in C_1(X,K_j^M)$.
These closures are curves that are proper and generically finite over $\overline{\{\eta\}}$.
We have $f_*\partial\widetilde b=\partial b$ by the compatibility of residues on the normalizations of these closures with field norms, that is, by the chain-map property of proper pushforward.

Thus $f_*(z-\partial\widetilde b)=0$ in $C_0(Y,K_j^M)$.
Each part of $z-\partial\widetilde b$ on a closed fiber has degree zero, so injectivity of the degree map on the closed fibers makes it the boundary of a chain in $C_1(X_y,K_j^M)$.
Pushing these chains forward along the closed immersions $X_y\to X$ and adding them shows that $z-\partial\widetilde b$, and hence $z$, is a boundary.
\end{proof}

\begin{lemma}
\label{lem:review-prime-degree-cycle}
Suppose that $\ell\ne\operatorname{char}k$, $j\ge1$, and $\cd_\ell k\le j+1$, and let $B/k$ be a central simple algebra of degree $\ell$.
Then the degree map
\[
 A_0(\SB(B),K_j^M)_{(\ell)}\xrightarrow{\sim}K_j^M(k)_{(\ell)}
\]
is an isomorphism.
\end{lemma}
\begin{proof}
All isomorphisms in this proof are understood after $\ell$-localization.
First assume that $k$ is perfect, and put $X=\SB(B)$.
Consider the coniveau spectral sequence computing $\CH^{\ell-1+j}(X,j)$; with $n=\dim X+j$, its $E_1$-terms are $E_1^{a,b}=\bigoplus_{x\in X^{(a)}}H^{b-a}(\kappa(x),\Z(n-a))$.
On the target diagonal, degree minus weight equals $\dim X-a$, and at the source of $d_r$ it equals $r-1$, so the terms in codimension $a<\ell-1$ and the differentials $d_r$ with $r\ge2$ entering the target term vanish because the cohomological degree over a field exceeds the weight.
Only $K_j^M$ at closed points and residues from $K_{j+1}^M$ at one-dimensional points remain, so duality \cite[Corollary C.7.3]{KahnLevine2010} and the cancellation theorem give
\[
 A_0(X,K_j^M)=\CH^{\ell-1+j}(X,j)
       =\operatorname{Hom}(\Z,M(X)(j)[j]).
\]

By the finite tower \cite[(8.1)]{KahnLevine2010}, for which $f_\ell M(X)=0$, the $s$th slice of $M(X)$ is $\Z_{B^{\otimes(s+1)}}(s)[2s]$ for $0\le s\le\ell-1$.
For $s>0$, the contributing group
\[
 H^{j+2s}(k,\Z_{B^{\otimes(s+1)}}(j+s))
\]
and the group one degree higher vanish because degree exceeds weight.
Hence $A_0(X,K_j^M)\simeq H^j(k,\Z_B(j))$; on the bottom layer the structure map is $\Z_B\hookrightarrow\Z$, so the degree map corresponds to the reduced norm induced by this inclusion \cite[Lemmas 5.4.2, 8.1.1(3), Theorems 6.5.3, 7.4.2]{KahnLevine2010}.
The norm residue isomorphism theorem supplies the Beilinson--Lichtenbaum hypothesis of the Kahn--Levine exact sequence in every required weight \cite{Voevodsky}.
The sequence makes the degree map injective with cokernel contained in $H^{j+2}_{\et}(k,\Z/\ell(j+1))$, which vanishes by the cohomological dimension hypothesis, so the degree map is an isomorphism \cite[Theorem 8.2.2(2)]{KahnLevine2010}.

If $k$ is imperfect of positive characteristic, pass to its perfect closure; the finite \'etale sites are equivalent, so $\cd_\ell$ is preserved.
Zero-cycles with coefficients and their relations are defined by finitely many data and therefore descend to a finite purely inseparable stage.
Both injectivity and surjectivity descend to $k$ because restriction followed by transfer is multiplication by the degree of this stage, a power of $p$ that is invertible after $\ell$-localization.
\end{proof}

\begin{proposition}
\label{prop:review-small-cd-cycle-degree}
Suppose that $\ell\ne\operatorname{char}k$, $j\ge1$, and $\cd_\ell k\le j+1$.
For every finite product $Z/k$ of Severi--Brauer varieties of central simple algebras, the degree map
\[
 \deg_Z:A_0(Z,K_j^M)_{(\ell)}\xrightarrow{\sim}K_j^M(k)_{(\ell)}
\]
is an isomorphism.
\end{proposition}
\begin{proof}
We prove the assertion simultaneously for all pairs $(k,j)$ satisfying the hypotheses.
For an $a$-dimensional point of a finite-type $k$-scheme $W$, the residue field $\kappa$ satisfies $\cd_\ell\kappa\le\cd_\ell k+a\le(j+a)+1$ \cite[II, \S4.2, Proposition 11]{SerreGalois}.
Thus, when we apply Lemma~\ref{lem:review-proper-cycle-pushforward} to a projection, the base fields of both types of fibers satisfy the bound, with coefficients $K_j^M$ at closed points and $K_{j+1}^M$ at one-dimensional points.
For a product of factors $\SB(B_i)$ with $\deg B_i=\ell$, apply Lemma~\ref{lem:review-prime-degree-cycle} fiberwise and induct on the number of factors to see that the degree map is an isomorphism.

Now consider $X=\SB(B)$ for an arbitrary central simple algebra $B/k$.
We may pass to a finite separable extension of degree prime to $\ell$ that splits the other primary Brauer components, and adjoin $\mu_\ell$, which requires an extension of degree dividing $\ell-1$.
Restriction and transfer descend the degree isomorphism to the original field.
Assume henceforth that $[B]$ has period $\ell^r$, and induct on $r$.
For $r=0$, use the projective bundle formula.
We express $\ell^{r-1}[B]$ as a sum of classes of symbol algebras of degree $\ell$, and let $Y$ be the product of their Severi--Brauer varieties \cite{MerkurjevSuslin1982,Voevodsky}.

Every fiber of the projection $f:X\times Y\to X$ is a product of factors $\SB(B_i)$ with $\deg B_i=\ell$, so the preceding calculation and Lemma~\ref{lem:review-proper-cycle-pushforward} show that $f_*$ is an isomorphism.
Over every residue field $\kappa(y)$ of $Y$, the point $y$ splits each symbol algebra factor, so $\ell^{r-1}[B]=0$ over $\kappa(y)$.
The induction on the period therefore applies to the fibers of $g:X\times Y\to Y$ over closed and one-dimensional points, with coefficients $j$ and $j+1$, respectively, so $g_*$ is an isomorphism by the same lemma.
Transitivity of proper pushforward gives $\deg_X=\deg_Y\circ g_*\circ f_*^{-1}$, which is an isomorphism because $\deg_Y$ already is.
Finally, we adjoin one arbitrary algebra factor at a time and use the same fiberwise argument to prove the assertion for general finite products.
\end{proof}

\begin{corollary}
\label{cor:review-small-cd-splitting-transfer}
Suppose that $\ell\ne\operatorname{char}F$ and $\cd_\ell F\le3$.
For every central division algebra $D$ over $F$,
\[
 \Nrd_2:I_F(D)_{(\ell)}\xrightarrow{\sim}K_2(F)_{(\ell)}
\]
is an isomorphism.
\end{corollary}
\begin{proof}
Put $X=\SB(D)$.
Compose the Morita equivalence given by the universal right ideal at each closed point with proper pushforward to obtain $\bigoplus_{t\in X_{(0)}}K_2^M(\kappa(t))\to K_2(D)$.
On the normalization of the closure of each one-dimensional point, the pullback of the universal ideal splits the Azumaya algebra.
There the composite of the map from Milnor $K_3$ to Quillen $K_3$ with the localization boundary and proper pushforward is zero, since consecutive maps in the localization sequence compose to zero.
The Rost boundary relations therefore vanish, giving a map $\Theta:A_0(X,K_2^M)\to K_2(D)$.
Every closed point $t$ splits $D$, and the component of $\Theta$ at $t$ is the splitting transfer $\theta_{\kappa(t)}$.
Conversely, a rational point of $X$ over a finite splitting field $E/F$ maps to a closed point $t$, and we have $\theta_E=\theta_{\kappa(t)}N_{E/\kappa(t)}$.
Thus $\im\Theta=I_F(D)$.
Since the reduced norm is compatible with the field norm at each closed point, we have $\Nrd_2\Theta=\deg_X$ \cite{Rost,QuillenHigher,KahnLevine2010}.
By Proposition~\ref{prop:review-small-cd-cycle-degree} with $j=2$, the degree map is an isomorphism after $\ell$-localization, so the reduced norm on the image $I_F(D)_{(\ell)}$ of $\Theta_{(\ell)}$ is an isomorphism.
\end{proof}

\begin{lemma}
\label{lem:review-p-rank-two-frobenius-norm}
Suppose that $\operatorname{char}F=p$ and $[F:F^p]=p^2$.
Let $q=p^a$ and $L=F^{1/q}$, and let $\varphi:L\xrightarrow{\sim}F$ be $x\mapsto x^q$.
Then
\[
 N_{L/F}=\varphi_*:K_2(L)\xrightarrow{\sim}K_2(F).
\]
Moreover, for every central division algebra $D$ over $F$, the reduced norm gives an isomorphism
\[
 \Nrd_2:I_F(D)_{(p)}\xrightarrow{\sim}K_2(F)_{(p)}.
\]
\end{lemma}
\begin{proof}
We use Matsumoto's identification $K_2=K_2^M$.
Since $[L:F]=q^2$, one has $N_{L/F}\res_{L/F}=q^2$, and on symbols one also has $\varphi_*\res_{L/F}=q^2$.
Since symbols generate $K_2(L)$ and each $q^2\{u,v\}=\{u^q,v^q\}$ lies in the image of restriction, it follows that $q^2(N_{L/F}-\varphi_*)=0$.
Since $K_2(F)$ has no $p$-torsion, we have $N_{L/F}=\varphi_*$ \cite{Matsumoto}, \cite[\S2.4, Theorem 6]{Izhboldin2000}.

First suppose that $D$ is $p$-primary, and choose $q$ divisible by its period.
Pullback along the $q$th-power Frobenius acts on $\Br(F)$ as multiplication by $q$, since it raises the defining $2$-cocycles to the $q$th power.
Since the composite $F\subset L\xrightarrow{\varphi}F$ is $x\mapsto x^q$, we have $\varphi_*[D_L]=q[D]=0$ in $\Br(F)$, and $D_L$ splits because $\varphi$ is an isomorphism.
Since the definition of $I_F(D)$ includes finite purely inseparable splitting fields, the splitting transfer $\theta_L$ is defined.

For any finite extension $F'/F$, the tower formula and $[{F'}^p:F^p]=[F':F]$ give $[F':{F'}^p]=p^2$, so the first assertion, applied to $F'$, shows that $N_{{F'}^{1/q}/F'}$ is an isomorphism.
Since $F'L\subset{F'}^{1/q}$, this norm is the composite $N_{F'L/F'}N_{{F'}^{1/q}/F'L}$, so $N_{F'L/F'}$ is surjective.
In particular, for every finite splitting field $E/F$ and $z\in K_2(E)$, write $z=N_{EL/E}w$.
Transitivity of splitting transfers gives
\[
 \theta_Ez=\theta_{EL}w=\theta_LN_{EL/L}w.
\]
Consequently $I_F(D)=\theta_LK_2(L)$, and $\Nrd_2$ maps it isomorphically onto $K_2(F)$ since $\Nrd_2\theta_L=N_{L/F}$ is an isomorphism.

For general $D$, split the prime-to-$p$ Brauer components over a finite extension $F'/F$ of degree prime to $p$.
Restriction and transfer preserve the corresponding splitting transfer subgroups, and the composite of these two maps is multiplication by this degree.
We use $\Nrd_2\circ\cor_{F'/F}=N_{F'/F}\circ\Nrd_2$ \cite[Proposition 6.6.1(1)]{KahnLevine2010}.
For surjectivity, given $z\in K_2(F)_{(p)}$, the transfer of a preimage of $\res_{F'/F}z$ over $F'$ is a preimage of $[F':F]z$, and $[F':F]$ is a unit in $\Z_{(p)}$.
For injectivity, restriction is injective on $I_F(D)_{(p)}$ for the same reason, so injectivity over $F'$ suffices.
\end{proof}

\end{document}